\documentclass[12pt]{article}
\usepackage[utf8]{inputenc}
\usepackage[english]{babel}
\usepackage[T1]{fontenc}
\usepackage{hyperref}
\usepackage{amsmath}
\usepackage{amsthm}
\usepackage{amsfonts}
\usepackage{amssymb}
\numberwithin{equation}{section}
\usepackage{graphicx}
\usepackage{tikz}
\usetikzlibrary{calc}
\usepackage{framed}
\usepackage{caption}
\usetikzlibrary{arrows}
\tikzset{>=latex}
\usetikzlibrary{patterns}
\usepackage{authblk}
\usepackage[top=2.5cm, bottom=2.5cm, left=2cm, right=2cm]{geometry}
\usepackage{fancyhdr}
\usepackage{cuted}
\usepackage{appendix}
\usepackage{chngcntr}
\usepackage[all,cmtip]{xy}

\newtheorem{theorem}{Theorem}[section]
\newtheorem{proposition}[theorem]{Proposition}
\newtheorem{corollary}[theorem]{Corollary}
\newtheorem{lemma}[theorem]{Lemma}
\newtheorem*{notation}{Notation}
\newtheorem{definition}[theorem]{Definition}

\newtheorem{remark}[theorem]{Remark}
\theoremstyle{definition}

\newcommand{\N}{\mathbb{N}}

\newcommand{\R}{\mathbb{R}}
\newcommand{\C}{\mathcal{C}}
\newcommand{\Z}{\mathbb{Z}}

\newcommand{\Q}{\mathbb{Q}}

\newcommand{\Prob}{\mathbb{P}}
\newcommand{\E}{\mathbb{E}}
\newcommand{\cB}{\mathcal{B}}
\newcommand{\e}{\mathcal{E}}

\newcommand{\T}{\mathbb{T}}

\newcommand{\s}{\mathcal{S}}

\renewcommand{\epsilon}{\varepsilon}

\newcommand\xqed[1]{%
  \leavevmode\unskip\penalty9999 \hbox{}\nobreak\hfill
  \quad\hbox{#1}}
\newcommand\demo{\xqed{$\blacksquare$}}

\newcommand\yqed[1]{%
  \leavevmode\unskip\penalty9999 \hbox{}\nobreak\hfill
  \quad\hbox{#1}}
\newcommand\emptydemo{\yqed{$\square$}}

\newlength{\squareheight}

\DeclareRobustCommand{\squaredash}{%
  \mathbin{\text{\settoheight{\squareheight}{\mathstrut}\raisebox{0.0\squareheight}{%
    \tikz{\draw (0,0) rectangle (0.85\squareheight,0.85\squareheight);\draw(0,0) -- (0.85\squareheight,0.85\squareheight)}}}}}

\title{Large deviations for out of equilibrium correlations in the symmetric simple exclusion process}

\author{Thierry Bodineau\footnote{I.H.E.S., Universit\'e Paris-Saclay, CNRS, Laboratoire Alexandre Grothendieck. 
35 Route de Chartres, 91440 Bures-sur-Yvette (France). E-mail: {\tt bodineau@ihes.fr}.}\hspace{2cm} Benoit Dagallier\footnote{Courant Institute of Mathematical Sciences, NYU and DPMMS, University of Cambridge. E-mail: {\tt bd2543@nyu.edu}.}}

\date{}

\begin{document}

\maketitle

\begin{abstract}
For finite size Markov chains, the Donsker-Varadhan theory fully describes the large deviations of the time averaged empirical measure. 
We are interested in the extension of the Donsker-Varadhan theory for a large size non-equilibrium system:  the one-dimensional symmetric simple exclusion process connected with reservoirs at different densities. The Donsker-Varadhan functional encodes a variety of scales depending on the observable of interest. 
In this paper, we focus on the time-averaged two point correlations
and investigate the large deviations from the steady state behaviour.
To control two point correlations out of equilibrium, 
the key input is the construction of a simple approximation to the invariant measure. 
This approximation is quantitative in time and space as estimated through relative entropy bounds building on the work of Jara and Menezes~\cite{Jara2018}.   
\end{abstract}

\section{Introduction}
For a fluid in thermal equilibrium, 
spatial correlations are expected to have fast decay, 
in such a way that, roughly speaking, 
each macroscopic portion of the fluid is basically independent from the rest.
For fluids driven out of equilibrium, 
e.g. by contact with reservoirs at two different temperatures, 
the picture is quite different: 
the fluid settles in a steady state where heat and/or matter are transported at a macroscopic level.  
The transport induces long-range correlations, 
which can be modelled by a variety of approaches and that have been observed experimentally, 
see~\cite{Spohn1983,Garrido1990} and references therein. These general predictions are part of the results of the Macroscopic Fluctuation Theory (see the review~\cite{Bertini2015}), 
which proposes a framework to study out of equilibrium fluids at a macroscopic level. \\

The derivation, from a microscopic model, 
of the steady state correlations, 
which are of a genuinely dynamical nature, is usually a difficult problem. 
Rigorous results are mostly obtained for certain simple interacting particle systems on a lattice. 
The Symmetric Simple Exclusion Process connected with reservoirs (henceforth open SSEP) is a paradigmatic example for which this correlation structure can be analysed~\cite{Derrida2005,Landim2006}. 
In the open SSEP, defined in Section~\ref{sec_notations_and_results}, 
  particles follow symmetric random walks interacting by an exclusion rule on a finite subdomain of $\Z$. 
Reservoirs pump particles in and out of the system, 
 fixing a certain density of particles in their vicinity. 
When reservoirs are at the same density, the open SSEP dynamics is reversible. 
However, when connected with reservoirs which enforce a different density of particles, 
this dynamics settles in long time in a non-equilibrium steady state, characterised by a macroscopic current of particles. 
The strength of this current is proportional to the density difference between the reservoirs~\cite{Spohn1983,Eyink1990,Farfan2011}. 
Two-point correlations in the steady state are known exactly~\cite{Spohn1983}, 
as well as all higher cumulants~\cite{Derrida2007,Derrida2007a,Derrida2013} in dimension one. 
The correlation structure of the steady state of the open SSEP is conjectured to be representative of a large class of out of equilibrium systems~\cite{Spohn1983,Bertini2015}. 
However, much less is known rigorously about steady state correlations for general lattice gases.\\

Our goal is to estimate the asymptotic probability of observing a correlation structure that is different from the one of the steady state, 
thereby also gaining information on this invariant measure. 
When the value $N$ of the scaling parameter is fixed, 
this question has already received a comprehensive answer by Donsker and Varadhan~\cite{Donsker1975}. 
For a general, irreducible Markovian dynamics on a finite state space $\Omega_N$, 
they study the time empirical measure $\tilde\pi_T$, 
defined for each $T>0$ as a probability measure on the configuration space $\Omega_N$ by:
\begin{equation}
\forall \eta\in\Omega_N,\qquad 
\tilde \pi^T 
= 
\frac{1}{T}\int_0^T\delta_{\eta_t}dt
.
\end{equation}
The quantity $\tilde \pi^T(\eta)$ then corresponds to the proportion of time spent at a configuration $\eta\in\Omega_N$. 
A full large deviation principle with speed $T$ and rate function $I^N_{DV}$ is then provided in~\cite{Donsker1975} for the time empirical measure $\tilde \pi_T$, 
in the sense that, if $\mu^N$ is a probability measure on $\Omega_N$ and $\Prob$ denotes the probability associated with the dynamics:
\begin{equation}
\lim_{T\rightarrow\infty}\frac{1}{T}\log \Prob\big(\tilde\pi^T\approx \mu^N\big) = - I^N_{DV}(\mu^N),\label{eq_large_devs_DV}
\end{equation}
where $\approx$ means proximity in the weak topology of probability measures on $\Omega_N$.  
The rate function $I^N_{DV}$ vanishes only at the invariant measure $\pi^N_{inv}$ of the dynamics, 
and is defined through a complicated variational problem involving the generator $L$ of the Markov chain:
\begin{equation}
I^N_{DV}(\mu^N) 
:=
\sup_{h:\Omega_N\rightarrow\R} \mu^N\Big(e^{-h}(-L)\, e^h\Big)
.
\label{eq_DV_variational_principle}
\end{equation}

\medskip

We are interested in the macroscopic behaviour of the system, i.e. the large $N$ limit of the probability in~\eqref{eq_large_devs_DV}. 
To obtain these asymptotics, 
one possibility is to study the limit of~\eqref{eq_DV_variational_principle} when $N \to \infty$.
When the underlying dynamics is reversible, 
this can be carried out:
the variational problem~\eqref{eq_DV_variational_principle} can be solved, 
and $I^N_{DV}$ is expressed in terms of the Dirichlet form of the dynamics 
(see Section~\ref{sec_heuristics_correlation} below). 
Such computations are carried out in Section~\ref{sec_heuristics_correlation}, 
in the case of the open SSEP with reversible dynamics (where the scaling parameter $N$ is roughly the number of sites in the model). 
These computations highlight an important fact: 
one cannot naively take the large $N$ limit in~\eqref{eq_large_devs_DV} without losing information, 
because not all the information contained in a measure $\mu^N$ is stored at the same scale in $N$. 
By this we mean e.g. that observing a macroscopic density different from that of the steady state, 
or observing different two point-correlation but with the same density, are not events that have the same scaling in terms of $N$. 
Informally, it is shown in Section~\ref{sec_heuristics_correlation} that observing a macroscopic density profile different from the one of the invariant measure in the reversible open SSEP has a probability that scales like $e^{-TN^{-2}\cdot N}$ in the large $T$, then large $N$ limit, 
up to sub-exponential corrections. 
In contrast, changing the two-point correlation structure only requires a cost of order $e^{-TN^{-2}}$. 
To study the scaling limit of~\eqref{eq_large_devs_DV}, 
one therefore has to choose a scale. 
Out of equilibrium, the rate function $I^N_{DV}$ is not known explicitly, 
and the equilibrium heuristics cannot be used,  
but we prove in Theorem~\ref{theo_large_devs} below that scales are still separated in the same way. \\

In this article, 
we focus on the scale corresponding to two-point correlations, 
and quantify the probability of observing anomalous two-point correlations in the one dimensional, 
out of equilibrium open SSEP in the large $T,N$ limits. 
We establish a large deviation principle for the time-averaged two-point correlation field, 
in Theorem~\ref{theo_large_devs} below. 
We do not start from the Donsker-Varadhan asymptotics, 
but instead provide quantitative estimates on the dynamics, 
as a function of time and the system size. 
Note that density large deviations are well understood since the seminal paper~\cite{Kipnis1989} (see also Chapter 10 of~\cite{Kipnis1999} for a review and~\cite{Bertini2003} with reservoirs). 
The main difficulty of the article is to generalise these ideas to the estimate of two-point correlations, 
that are objects  living on a much finer scale than the density.
To illustrate this, recall that the
two-point correlations in the steady state of the open SSEP are long range, 
and scale like $O(N^{-1})$, 
compared to $O_N(1)$ for the density of particle at a given site. 
For this reason, 
and while model-dependent estimates on correlations have been obtained e.g. in~\cite{Spohn1983,Derrida2005,Landim2006,Goncalves2020}, 
to our knowledge there is no general method to study the out of equilibrium behaviour of the two-point correlation field in the long-time, large $N$ limits.

The proof of our result on two-point correlations, 
Theorem~\ref{theo_large_devs}, 
builds upon a refinement of the relative entropy method obtained by Jara and Menezes~\cite{Jara2018,Jara2020}.  
This method, originally introduced by Yau~\cite{Yau1991}, 
consists in quantifying, at each time and in terms of the relative entropy, 
the proximity of the law of the Markovian dynamics in an interacting particle system with a known reference measure. 
The idea behind the method is that, locally, 
the dynamics in large microscopic boxes equilibrates much faster than the typical time-scale at which the system evolves macroscopically. 
If one has an ansatz for the evolution of macroscopic variables of interest, say, henceforth, 
the density in a lattice gas; 
one then expects that the corresponding microscopic variables, 
when averaged over a sufficiently large microscopic box, are close to their macroscopic counterpart.  
This property, 
known as local equilibrium, 
has recently been shown quantitatively even for mesoscopic boxes, see~\cite{Goncalves2022}.

From the local equilibrium heuristics, 
one can build a reference measure in terms of the evolution of the macroscopic density only. 
If one is interested in the evolution of the macroscopic density or its fluctuations, 
it can be shown that local equilibrium holds and this reference measure is indeed a good enough approximation of the law of the dynamics, 
see Chapter 6 in~\cite{Kipnis1999},~\cite{Jara2018} and references therein. 
In particular, the reference measure does not need to contain any information on correlations.

To study two point correlations, however, 
the reference measure has to also contain information on the dynamical correlations. 
Adding such a correlation term in the reference measure is our key input.    
In the case of the open SSEP, 
since density fluctuations around the typical density profile at each time (and in the steady state~\cite{Landim2006}) are known to be Gaussian~\cite{Jara2018}, 
our candidates for reference measures are discrete Gaussian measures, 
see~\eqref{eq_def_bar_nu_g_intro_bis}. 
One expects that a good choice of discrete Gaussian measure will contain all leading order information about two point correlations. 
A similar observation was already present in~\cite{DemasiSmallDeviationsLocal1982}.  
It is made precise in Theorem~\ref{theo_entropic_problem} where we obtain,  
for a family of exclusion dynamics that occur in the proof of the large deviation result of Theorem~\ref{theo_large_devs}, 
a characterisation of long time, large $N$ correlations as the solution of a certain partial differential equation.\\

The approach used in this paper is not restricted to the symmetric simple exclusion process. 
Much like the usual relative entropy method, it can be used for a large class of one-dimensional diffusive interacting particle systems satisfying the so-called gradient condition (see Section 8 in~\cite{Jara2018}). 
In particular, 
very special features of symmetric simple exclusion such as the fact that correlations in the steady state are known are not used in the proof, see Section~\ref{sec_entropy_invariant_mes} for more details. 
There are however some technical difficulties to be expected when generalising the present approach. 
All these points are discussed further in Section~\ref{sec_conclusion_perspectives}. \\ 
Let us however mention that, at equilibrium, 
the behaviour of various $n$-point correlation fields 
has come under much scrutiny in the past few years. 
In~\cite{Assing2007}, 
a two-point correlation field is studied in the SSEP on $\Z$. It is not the same object as in~\cite{Goncalves2019}, 
where two point correlations are studied on the one-dimensional torus as a means to defining squares of distributions arising in certain ill-posed stochastic partial differential equations. 
In~\cite{Ayala2021} and~\cite{Chen2021}, 
interacting particle systems enjoying a self-duality property are investigated in all dimensions. 
In that context, equilibrium fluctuation fields involving $n$-point functions are investigated for any $n$.

In the same direction but using different techniques, 
long time large deviations for the density and the current have recently been considered in~\cite{Bertini2021}. 
Both the long diffusive time limit starting from the dynamical large deviation functional, 
and the long time, then large $N$ limits using Donsker and Varadhan's formula are investigated. 
This last limit is the same as the one studied in the present article at the level of two point correlations. 
The approach is however different: 
here, we provide a quantitative (i.e. non asymptotic) control on the dynamics at the microscopic level. 
On the other hand, 
in the case of the density and the current, 
the microscopic model considered in~\cite{Bertini2021} may have dynamical phase transitions. 
To capture this very subtle phenomenon, 
the Donsker-Varadhan variational principle~\eqref{eq_DV_variational_principle} in~\cite{Bertini2021} is solved indirectly in the large $N$ limit, 
by looking at large deviations at process level (so-called level 3 large deviations), 
then using a contraction argument. 
A related paper by the same authors~\cite{Bertini_FW2022} uses a similar approach in the settings of diffusions with small noise. 
Note that different kind of asymptotics for the Donsker-Varadhan functional have been considered, 
e.g. to study metastability~\cite{Bertini_metastab2022,Landim_metastab2022}.\\

The rest of the article is structured as follows. In Section~\ref{sec_notations_and_results}, 
we present the model and results. 
Section~\ref{sec_computations} gives the main microscopic tool for the study of two-point correlations: 
the relative entropy estimate when the reference measure is a certain discrete Gaussian measure. 
Properties of these measures are established in Appendix~\ref{app_discrete_gaussian_measures}. 
The relative entropy bounds allow for the computation of the Radon-Nikodym derivative between the open SSEP and the tilted processes introduced to estimate rare events. 
This requires sharp estimates collected in the appendices. 
The large deviations are then established, in Section~\ref{sec_large_devs} for the upper bound, 
and~\ref{sec_lower_bound} for the lower bound. 
For the lower bound, control of the open SSEP dynamics in long-time is obtained via the study of certain Poisson equations. 
Well-posedness of these equations is investigated in Appendix~\ref{app_Poisson}, 
while Appendix~\ref{app_sobolev_spaces} gathers useful topological facts.

\paragraph{Acknowledgements.}
This work has been motivated by many discussions with Bernard Derrida 
on the structure of correlations in non-equilibrium particle systems.
We are extremely grateful to him for sharing his insights. 
We would also like to thank Stefano Olla for  very useful suggestions and discussions at various stages of the writing process. Part of this work was carried out while B.D. was supported by the
European Research Council under the European Union’s Horizon 2020 research
and innovation programme (grant agreement No. 851682 SPINRG).
\section{Notations and results}\label{sec_notations_and_results}
\subsection{Notations and definition of the microscopic model}
\subsubsection{The microscopic model} 
For $N\in\N^*$, let $\Lambda_N :=\{-N+1,...,N-1\}$ and $\Omega_N = \{0,1\}^{2N-1}$. 
Elements of $\Omega_N$, denoted by the letter $\eta$, 
will be called configurations. 
We say that there is a particle at site $i\in\Lambda_N$ if $\eta_i=1$, and no particle if $\eta_i=0$. 
The variable $\eta_i$ is called the occupation number (of site $i$). 
On $\Omega_N$, we consider the dynamics given by the Symmetric Simple Exclusion Process connected to reservoirs at position $\pm N$ (henceforth open SSEP), 
which we now describe. 
For a survey of particle systems in contact with reservoirs, we refer to~\cite{Derrida2007a, Bertini2015}.

Let $\rho_-<\rho_+\in(0,1)$ be the densities of the left (for $\rho_-$) and the right ($\rho_+$) reservoirs. 
The open SSEP is defined through its generator $N^2L := N^2(L_0+L_-+L_+)$. 
It is made up of two parts, the bulk and boundary dynamics, corresponding to $L_0$ and $L_\pm$ respectively. 
The operators $L_0$ and $L_\pm$ act on $f:\Omega_N\rightarrow\R$ as follows:
\begin{align}
\forall\eta\in\Omega_N,\qquad N^2L_0 f(\eta) 
&= 
\frac{N^2}{2}\sum_{i<N-1}c(\eta,i,i+1)\big[f(\eta^{i,i+1}) - f(\eta)\big],\label{eq_def_L_0}\\
N^2L_\pm f(\eta) 
&= 
N^2L_+ f(\eta) +N^2L_- f(\eta)
,\nonumber\\
\quad 
N^2L_\epsilon f(\eta)
&= \frac{N^2}{2}c(\eta,\epsilon (N-1))\big[f (\eta^{\epsilon(N-1)})-f(\eta)\big]
,\qquad 
\epsilon\in\{-,+\}
\label{eq_def_L_pm}
.
\end{align}
Above, the jump rates $c$ are defined, for each $\eta\in\Omega_N$, by:
\begin{align}
c(\eta,i,i+1) 
&= 
\eta_{i+1}(1-\eta_i) + \eta_i(1-\eta_{i+1}),\qquad i<N-1,\nonumber\\
c(\eta,\pm(N-1)) 
&= 
(1-\rho_\pm)\eta_{\pm (N-1)} + \rho_{\pm}(1-\eta_{\pm (N-1)}),\label{eq_def_jump_rates}
\end{align}
and for $i,j\in\Lambda_N$ and $\eta\in\Omega_N$, 
the configurations $\eta^i$ and $\eta^{i,j}$ read:
\begin{equation}
\eta^i_\ell 
= 
\begin{cases}
\eta_\ell\quad &\text{if }\ell\neq i,\\
1-\eta_i\quad &\text{if }\ell=i,
\end{cases}
\qquad 
\eta^{i,j}_\ell 
= 
\begin{cases}
\eta_\ell\quad &\text{if }\ell\notin\{i,j\},\\
\eta_j\quad &\text{if }\ell=i,\\
\eta_i\quad &\text{if }\ell=j.
\end{cases}
\end{equation}
We write $\Prob^\eta,\E^\eta$ for the probability/expectation under this dynamics starting from $\eta\in\Omega_N$. \\

\subsubsection{The invariant measure and the correlation field} 

For each $N\in\N^*$, let $\pi^N_{inv}$ denote the unique invariant measure of the open SSEP. 
If $\rho_+=\rho_- = \rho\in(0,1)$, 
$\pi^N_{inv}$ is simply the Bernoulli product measure on $\Lambda_N$ with parameter $\rho$. 
If $\rho_-<\rho_+$, however, the measure $\pi^N_{inv}$ is not product. 
The average occupation number at each site was computed in~\cite{Spohn1983}: 
it is given in terms of an affine function $\bar\rho$, with:
\begin{equation}
\forall i\in\Lambda_N,\qquad 
\bar\rho(i/N)
:=  \E^{\pi^N_{inv}}\big[ \eta_i \big]
= \Big(1-\frac{i}{N}\Big)\frac{\rho_-}{2} + \Big(1+\frac{i}{N}\Big)\frac{\rho_+}{2}.\label{eq_def_profil_invariant}
\end{equation}
Note that, as $\bar\rho$ is affine, 
$\bar\rho' = (\rho_+-\rho_-)/2$ is a constant. 
As $\pi^N_{inv}$ is invariant, one has, for $F:\Omega_N\rightarrow\R$ and each $t\geq 0$:
\begin{equation}
\forall t\geq 0,\qquad 
\pi^N_{inv}(F) 
:= \E^{\pi^N_{inv}}\big[F(\eta_t)\big] = 
\sum_{\eta\in\Omega_N}\pi^N_{inv}(\eta)F(\eta) .
\end{equation}
In the following, the expectation $\E_\mu[F]$ with respect to  a measure $\mu$ is written $\mu(F)$.

The measure $\pi^N_{inv}$ exhibits long-range correlations. 
To make this statement precise, 
let us first define the main object of interest in this article, 
the \emph{correlation field} $\Pi^N$. 
It is a distribution, acting on test functions $\phi:(-1,1)^2\rightarrow\R$ according to:
\begin{equation}
\forall\eta\in\Omega_N,\qquad \Pi^N(\phi) = \Pi^N(\phi)(\eta) 
= \frac{1}{4N}\sum_{i\neq j\in\Lambda_N}\bar\eta_i\bar\eta_j \phi(i/N,j/N),
\quad \bar\eta_\cdot:= \eta_\cdot-\bar\rho(\cdot/N).
\label{eq_def_Pi}
\end{equation}
Note the scaling of~\eqref{eq_def_Pi} with $N$: 
the sum contains order $N^2$ terms, and the normalisation is only proportional to $N^{-1}$ as $\Pi^N(\phi)$ measures fluctuations of a central limit order.
Indeed the correlation field $\Pi^N$ is strongly related to the fluctuation field $Y^N$ defined for any  bounded test function $\psi:[-1,1]\rightarrow\R$ by:
\begin{equation}
\forall \eta\in\Omega_N,\qquad 
Y^N(\psi) 
= 
\frac{1}{\sqrt{N}}\sum_{i\in\Lambda_N}\bar\eta_i \psi(i/N)
.
\label{eq_def_fluctuations_sec2}
\end{equation}
Restricting to product test functions $\phi (x,y) = \psi_1 (x) \psi_2 (y)$, 
we note that:
\begin{equation}
\Pi^N(\phi)(\eta) 
= \frac{1}{4} Y^N(\psi_1) Y^N(\psi_2) - \frac{1}{4 N}
\sum_{i \in\Lambda_N} (\bar\eta_i)^2   \psi_1(i/N)  \psi_2(i/N).
\end{equation}
The fluctuation field has been extensively studied in the hydrodynamic scaling and 
at each time $t$, the fluctuation field $Y^N_{t}(\psi)$ can be proven to converge to a Gaussian random variable
(see Chapter 11 in~\cite{Kipnis1999} in the equilibrium case, or~\cite{Derrida2005,Landim2006,Jara2018} when $\rho_-\neq\rho_+$). 
The two-point correlations $Y_t(\psi_1)Y_t(\psi_2)$ are thus also of order $1$ in $N$ 
(see e.g.~\cite{Goncalves2019} for the reversible SSEP on the torus).

As mentioned in the introduction, the correlation field $\Pi^N$ is a natural quantity to study non equilibrium  properties. For SSEP, the average value of the correlation field under the invariant measure was obtained in~\cite{Spohn1983,Derrida2002} in the large $N$ limit :
\begin{equation}
\lim_{N\rightarrow\infty}\pi^N_{inv}\big(\Pi^N(\phi)\big) 
= 
\frac{1}{4}\int_{(-1,1)^2} \phi(x,y)k_0(x,y)\, dx\, dy,
\label{eq_def_correl_steady_state}
\end{equation}
where the kernel $k_0$ is a symmetric continuous  function on $[-1,1]^2$ given by :
\begin{equation}
\forall (x,y)\in(-1,1)^2,\qquad k_0(x,y) = -\frac{(\bar\rho')^2}{2}\Big[(1+x)(1-y){\bf 1}_{x\leq y}+(1+y)(1-x){\bf 1}_{x\geq y}\Big].
\label{eq_def_k_0}
\end{equation}
Thus when $\bar\rho' = (\rho_+-\rho_-)/2 \not = 0$,  
the system is driven out of equilibrium by the reservoirs 
 and long range correlations arise.
These correlations come from the diffusive transport of particles and, as a consequence, 
$k_0$ is obtained as the solution of the Laplace equation restricted to the triangle $\rhd = \{(x,y)\in [-1,1]^2 : x< y\}$ with fixed normal derivative on the diagonal $D = \{(x,x):x\in(-1,1)\}$: 
\begin{equation}
\forall (x,y) \in \rhd, \quad  \Delta k_0 (x,y) =0, \quad \text{with} \quad 
 \quad (\partial_1-\partial_2)k_0(x_\pm,x) = \pm( \bar\rho')^2,
 \quad x\in (-1,1)
 ,
\label{eq_Laplacien_k_0}
\end{equation}
with $0$ Dirichlet boundary conditions on $\{ x=-1 \}$ and $\{ y = 1 \}$.
 Above, we used the notation $x_{\pm}= \lim_{h\downarrow 0}(x\pm h)$. 
 
In the next section, we provide further insight on the correlation field and explain how $\Pi^N$ is related to the Donsker-Varadhan functional.

\begin{notation}
Throughout this paper, we write $\bar\rho_i := \bar\rho(i/N)$ for $i\in\Lambda_N$. 
More generally, for a function $\phi:[-1,1]^p\rightarrow\R$, $p\in\N^*$, 
we write $\phi_{i_1,...,i_p}$ for $\phi(i_1/N,...,i_p/N)$, $(i_1,...,i_p)\in\Lambda_N^p$. 
The letters $i,j,\ell...$, when used as indices, index elements of $\Lambda_N$; while $x,y,z$ are used for continuous variables.\\
More generally, when we speak of $n$-point correlations ($n\in\N^*$), 
it will always mean products of centred variables, 
the $\bar\eta$'s, of the form $\bar\eta_{i_1}...\bar\eta_{i_n}$ for some $(i_1,...,i_n)\in\Lambda_N^n$. 
When considering a trajectory $\eta_t\in\Omega_N,t\geq 0$, 
we write $\Pi^N_t$ for $\Pi^N(\cdot)(\eta_t)$.
\end{notation}

\subsection{Heuristic scaling of the Donsker-Varadhan functional}
\label{sec_heuristics_correlation}
In this section, we consider the SSEP dynamics  at equilibrium as given in \eqref{eq_def_L_0}--\eqref{eq_def_L_pm}, 
i.e. with two reservoirs at equal densities $\rho_-=\rho_+ = \rho\in(0,1)$.
The dynamics is reversible with respect to the Bernoulli product measure $\nu^N_\rho$ with parameter $\rho$. 
In this setting, the Donsker-Varadhan rate function is explicit and is given by the Dirichlet form of the dynamics~\cite{Donsker1975}: 
if $\mu^N$ is a probability measure on the state space $\Omega_N$, 
then one has for each fixed $N$:
\begin{align}
\lim_{T\rightarrow\infty}\frac{1}{T}\log \Prob^N_{\rho}
\Big(\frac{1}{T}\int_0^T\delta_{\eta_t}\, dt \approx \mu^N\Big)
&=
-
I^N_{DV}(\mu^N) 
\nonumber\\
&= 
-\nu^N_{\rho}\big(\sqrt{f} (-N^2L)\sqrt{f}\big),
\quad
f:= \frac{d\mu^N}{d\nu^N_{\rho}}
,
\label{eq_rate_function_DV_reversible}
\end{align}
where the $\approx$ sign denotes proximity in the weak topology of measures on $\Omega_N$. 

In the next two sections, 
we use the explicit form \eqref{eq_rate_function_DV_reversible} of the Donsker-Varadhan functional to derive its asymptotics as $N \to \infty$.
We will show that the scalings in $N$ are different  for the cost of observing an anomalous macroscopic density 
(in Section \ref{sec_shifting_density_intro}) or the cost of observing anomalous macroscopic two-point correlations 
(in Section \ref{sec_shifting_correl_intro}).
Thus different types of information are intertwined in the Donsker Varadhan functional and 
one has to zoom at the correct scaling to extract the relevant physical
information on a given observable.
The exact computation of Section \ref{sec_shifting_correl_intro} justifies our choice of focusing on the non-equilibrium large deviations of the   two-point correlation field $\Pi^N$ 
(see Section \ref{sec: Large deviations for time-averaged correlations}).

\subsubsection{Changing the macroscopic density}
\label{sec_shifting_density_intro}
Consider a smooth density profile $\hat\rho:[-1,1]\rightarrow(0,1)$. 
For simplicity, assume that the density close to the boundaries is unchanged, i.e. that $\hat\rho(x) = \rho$ for $x$ in an open neighbourhood of $\pm 1$. Define then:
\begin{align}
\mu^N = \bigotimes_{i\in\Lambda_N}\text{Ber}\big(\hat\rho(i/N)\big).
\end{align}
For a function $q:(-1,1)\rightarrow\R$, write for short $q_i := q(i/N)$ for $i\in\Lambda_N$. Introduce the chemical potential $\lambda$ and its discrete derivative $\partial^N\lambda_i$, defined by:
\begin{equation}
\forall x\in(-1,1),\qquad 
\lambda(x) 
:= 
\log\Big(\frac{\hat\rho(x)}{1-\hat\rho(x)}\Big),\quad \partial^N\lambda_i = N\big[\lambda_{i+1}-\lambda_i\big],\qquad i<N-1.\label{eq_def_lambda_intro}
\end{equation}
Recall that $\mu^N[\cdot]$ (or $\mu^N(\cdot)$) denotes expectation under the measure $\mu^N$. 
Elementary computations using \eqref{eq_rate_function_DV_reversible} then give, for each $N\in\N^*$:
\begin{align}
I^N_{DV}(\mu^N)
= 
\frac{N^2}{4}\mu^N\bigg[\sum_{i<N-1}c(\eta,i,i+1) \Big(\exp\Big[-\frac{(\eta_{i+1}-\eta_i)}{2N}\partial^N\lambda_i\Big]-1\Big)^2\bigg].
\label{eq_DV_densite}
\end{align}
Note that there is no contribution from the boundary dynamics, since $\hat\rho$ is constant and equals to $\rho$ close to $\pm 1$ by assumption.
 Expanding the right-hand side of \eqref{eq_DV_densite}, one finds, with the notation $\sigma(r) = r(1-r)$ for $r\in[0,1]$:
\begin{align}
\lim_{T\rightarrow\infty}\frac{1}{T}\log \Prob^N_\rho\Big(\frac{1}{T}\int_0^T\delta_{\eta_t}\, dt\approx \mu^N\Big) 
&= 
-\frac{1}{16}\mu^N\Big[\sum_{|i|<N-1}c(\eta,i,i+1)\big[\partial^N\lambda_i\big]^2\Big] +O_N(1)
\nonumber\\
&=
- \frac{N}{8}\int_{(-1,1)}\sigma(\hat\rho(x))|\nabla\lambda(x)|^2dx + O_N(1)
,
\label{eq_deviations_densite_DV}
\end{align}
where we used the smoothness of $\hat\rho$, and $\mu^N(c(\eta,i,i+1)) = 2\sigma(\hat\rho_i)+O(N^{-1})$ for each $i<N-1$, with the $O(N^{-1})$ uniform on $i$. It follows that a macroscopic change of density is observed with probability of order $e^{-TN}$ in the large $T$, then large $N$ limit.
\\

As a remark, notice that, up to factors of $N,T$, the right-hand side of \eqref{eq_deviations_densite_DV} is the same as the one given by the dynamical rate functional obtained in diffusive time in~\cite{Bertini2003}. 
To see it, recall that the rate functional $I_{SSEP}$ evaluated at the constant profile $\hat\rho$ on the time interval $[0,T]$ for $T>0$ is given by:
\begin{equation}
I_{SSEP}\big((\hat\rho(x)dx)_{t\leq T}) \big) 
= 
\frac{1}{2}\int_0^T \int_{(-1,1)}|\nabla h(t,x)|^2\sigma(\hat\rho(x))\, dx\, dt,
\label{eq_I_SSEP_with_h}
\end{equation}
where $h$ is the bias such that $h(t,\pm 1) = 0$ for each $t\in[0,T]$, and:
\begin{equation}
\partial_t\hat\rho = 0 = \frac{1}{2}\Delta\hat\rho - \nabla\cdot\big(\sigma(\hat\rho)\nabla h\big).
\end{equation}
In particular, integrating the divergence operator, there is a divergence-free function $j$ on $(-1,1)$ (the current), i.e. a constant in our one-dimensional setting, such that:
\begin{equation}
\nabla h = \frac{(1/2)\nabla\hat\rho + j}{\sigma(\hat\rho)}\quad \Rightarrow\quad I_{SSEP}\big((\hat\rho(x)dx)_{t\leq T}) \big)=\frac{T}{2}\int_{(-1,1)}\frac{\big((1/2)\nabla\hat\rho + j\big)^2}{\sigma(\hat\rho(x))}\, dx
.
\label{eq: 2.21}
\end{equation}
In the present case, $h(\cdot,\pm 1) = 0$ and $\nabla\lambda = \nabla \hat\rho/\sigma(\hat\rho)$ implies that $j=0$, with $\lambda$ defined in \eqref{eq_def_lambda_intro}. As a result:
\begin{equation}
\forall t\geq 0,x\in(-1,1),\qquad \nabla h(t,x) = \nabla\lambda(x)/2.
\end{equation}
Thus replacing $h$ by $\nabla\lambda/2$ in \eqref{eq_I_SSEP_with_h}, the functional 
\eqref{eq_deviations_densite_DV} is recovered.
In other words, the long time, large $N$ limit and the long diffusive time limits coincide :
\begin{align}
\lim_{N\rightarrow\infty}\lim_{T\rightarrow\infty}
\frac{1}{NT}\log &\, \Prob^N_\rho\Big(\frac{1}{T}\int_0^T\delta_{\eta_t}\, dt\approx \mu^N\Big) 
\nonumber\\
&\qquad = 
\lim_{T\rightarrow\infty}\lim_{N\rightarrow\infty}\frac{1}{NT}\log \Prob^N_\rho \Big((\pi_t^N)_{t\leq T} \approx (\hat\rho(x)dx)_{t\leq T}\Big)
. 
\label{eq_two_limits_coincide_density}
\end{align}

It is not difficult to adapt the proof of~\cite{Bertini2003}, derived in the hydrodynamic regime, to recover  \eqref{eq_two_limits_coincide_density}, i.e. to take the long time limit first before taking $N$ large. 
Indeed  the large deviation functional of the SSEP is convex, so that  the optimal way to observe an averaged density profile in the long diffusive time limit is obtained by a time-independent tilt : there is no dynamical phase transitions. 
The exchange of limits~\eqref{eq_two_limits_coincide_density} then remains valid in non equilibrium in the absence of a dynamical phase transition. 
The much more delicate proof that the two limits coincide even in the presence of a dynamical phase transition is carried out in~\cite{Bertini2021},  
where the joint deviations of the current and density are investigated.

\subsubsection{Changing the macroscopic correlations}
\label{sec_shifting_correl_intro}
Consider again the open SSEP at equilibrium at density $\rho\in(0,1)$. 
Our aim is to consider the large size asymptotics for a measure with the equilibrium density, but  different correlations.
Recall the definition~\eqref{eq_def_Pi} of the (off diagonal) correlation field $\Pi^N$, acting on a bounded test function $\phi:[-1,1]^2\rightarrow\R$ according to:
\begin{equation}
\forall \eta\in\Omega_N,\qquad 
\Pi^N(\phi) 
= 
\frac{1}{4N}\sum_{i\neq j \in\Lambda_N}\bar\eta_i \bar\eta_j \phi\Big(\frac{i}{N},\frac{j}{N}\Big).
\end{equation}
In view of the above discussion, to find a measure that is close to $\nu^N_\rho$ but with a different correlation structure, it is reasonable to look at:
\begin{equation}
\mu^N 
= 
\nu^N_{\rho,\phi} 
:= 
\frac{1}{\mathcal Z^{N}_{\rho,\phi}}\, e^{2\Pi^N(\phi)}\nu^N_{\rho}
,\qquad 
\mathcal Z^N_{\rho,\phi}\, \text{ a normalisation factor}
.
\label{eq_choice_mu_N_heuristic_correlations}
\end{equation}
Assume $\|\phi\|_\infty$ is sufficiently small and $\phi$ is smooth, symmetric, i.e.:
\begin{equation}
\forall (x,y)\in[-1,1]^2,\qquad \phi(x,y) = \phi(y,x).
\end{equation}
Assume also that $\phi(x,\cdot) = 0$ for $x$ in an open neighbourhood of $\pm 1$. 
Then the macroscopic density is still given by $\rho$, 
but the measure $\mu^N$ now features long-range correlations. 
Indeed, using the same kind of arguments as in Appendix~\ref{sec_correl_bar_nu_G} 
(see also~\cite{dagallierFluctuationsCorrelationsWeakly2023a}), 
one can show that there is a limiting kernel $k:\squaredash\rightarrow\R$ such that:
\begin{equation}
\sup_{i \not = j\in\Lambda_N}\big|\mu^N(\bar\eta_i\bar\eta_j) - k(i/N,j/N)\big| = o(N^{-1}),
\qquad 
\sup_{i \in\Lambda_N}\big|\mu^N(\bar\eta_i\bar\eta_i) - \sigma(\rho) \big| = o_N(1).
\label{eq_two_points_mu_N}
\end{equation}
The limiting covariance can therefore be described by an operator $C$ acting on $\psi_1,\psi_2\in\mathbb L^2((-1,1))$ according to 
(recall the definition~\eqref{eq_def_fluctuations_sec2} of the fluctuation field $Y^N$):
\begin{equation}
\lim_{N\rightarrow\infty}\mu^N\big(Y^N(\psi_1)Y^N(\psi_2)\big)
=
\int_{(-1,1)} \psi_1(x)(C\psi_2)(x)\, dx,
\end{equation}
with:
\begin{equation}
\forall x\in(-1,1),
\qquad 
C \psi_2(x) 
:= 
\sigma(\rho) \psi_2(x) + \int_{(-1,1)} k (x,y)\psi_2(y)\, dy
.
\label{eq: operateur C equilibre}
\end{equation}
The operator $C$ is obtained from $\phi$ as the inverse of the operator $U_\phi$, 
defined for $\psi\in\mathbb L^2((-1,1))$ by:
\begin{equation}
\forall x\in(-1,1),\qquad 
U_\phi\psi(x) :=  \sigma(\rho)^{-1}\psi(x) -\int_{(-1,1)} \phi(x,y)\psi(y)\, dy.
\end{equation}
Intuitively, this relation means that for large $N$, the density under the measure $\mu^N$ in \eqref{eq_choice_mu_N_heuristic_correlations} behaves as a Gaussian field of covariance $C$.
A similar structure will be derived out of equilibrium and the proofs will be given then, 
 see for instance Theorem~\ref{theo_entropic_problem} where it is stated that the out of equilibrium SSEP dynamics stays very close to a measure of the form $\mu^N$ at each time, 
in the sense that it has the same two-point correlation structure.
To summarise, the macroscopic density is still given by $\rho$, 
but the measure $\mu^N$ has now long-range correlations parametrised by $\phi$. 
Let us again compute the Donsker-Varadhan rate functional for the measure $\mu^N$ 
with large $N$. 
\begin{lemma}
\label{lemm_heur_correl}
For $\phi$ and $\mu^N$ as in~\eqref{eq_choice_mu_N_heuristic_correlations}, one has with the notation~\eqref{eq: operateur C equilibre}:
\begin{align}
&\lim_{T\rightarrow\infty}
\frac{1}{T}\log \Prob\Big(\frac{1}{T}\int_0^T\delta_{\eta_t}\, dt\approx \mu^N\Big) 
\nonumber\\
&\hspace{3cm}= 
-\frac{1}{8}\int_{(-1,1)} \sigma(\rho)\big<\partial_1 \phi(z,\cdot), C \partial_1\phi(z,\cdot)\big>\, dz +o_N(1)
,
\label{eq_formule_heuristique_avec_matrice_correl}
\end{align}
where $\big<\cdot,\cdot\big>$ is the scalar product in $\mathbb L^2((-1,1))$.
\end{lemma}
\begin{proof}
Plugging the expression~\eqref{eq_choice_mu_N_heuristic_correlations} of the measure $\mu^N$ 
in the Donsker-Varadhan functional \eqref{eq_rate_function_DV_reversible}, 
we obtain a formula 
similar to~\eqref{eq_DV_densite} where the variation of the  chemical potential is now replaced by 
the non local expression
\begin{equation}
\Pi^N(\phi)(\eta^{i,i+1})-\Pi^N(\phi)(\eta)
=
-\frac{(\eta_{i+1}-\eta_i)}{2N^2}\sum_{j\notin\{i,i+1\}}\bar\eta_j \partial^N_1 \phi_{i,j}
,
\end{equation}
with $\partial^N_1 \phi_{i,j} = N[\phi_{i+1,j}-\phi_{i,j}]$ bounded uniformly in $i,j,N$.  
Expanding the exponential, we get as in \eqref{eq_deviations_densite_DV}
\begin{align}
\lim_{T\rightarrow\infty}\frac{1}{T}\log\, &\Prob^N_\rho\Big(\frac{1}{T}\int_0^T\delta_{\eta_t}\, dt\approx \mu^N\Big) 
\nonumber\\
&\qquad
= 
-\frac{1}{16 N^2}\mu^N\Big[\sum_{\substack{|i|<N-1 \\ j, \ell \notin\{i,i+1\} }}c(\eta,i,i+1)
\bar\eta_j \bar\eta_{\ell} \; \partial^N_1 \phi_{i,j} \partial^N_1 \phi_{i,\ell} \Big] 
+o_N(1)
.
\label{eq_deviations_correl_DV}
\end{align}
The limit \eqref{eq_formule_heuristique_avec_matrice_correl} can be easily guessed from 
the above formula, 
as $2\sigma (\rho)$ arises from $c(\eta,i,i+1)$ and the correlations 
$\bar\eta_j \bar\eta_{\ell}$ are approximated by the limiting covariance \eqref{eq_two_points_mu_N} 
of $\mu^N$.
Similar computations will be carried out numerous times in Section~\ref{sec_computation_L_star}, 
so we give no details here.
\end{proof}

Compared with the asymptotics of the density large deviations \eqref{eq_deviations_densite_DV},
the cost \eqref{eq_formule_heuristique_avec_matrice_correl} of modifying only the correlations has a different scaling in $N$. 
Our aim is to derive similar results for systems driven out of equilibrium by reservoirs 
when the Donsker-Varadhan rate function is no longer given by the explicit formula \eqref{eq_rate_function_DV_reversible} but by the variational principle \eqref{eq_DV_variational_principle}. 
The computations in the reversible case highlight the fact that:
\begin{itemize}
	\item the whole correlation structure of the invariant measure is not contained at the same scale in $N$, and therefore:
	\item one has to focus on observables at a specific scale to get a non trivial limit when $N$ is large.
\end{itemize}
In the following, we focus on the large deviations of the correlation field 
$\Pi^N$ introduced in~\eqref{eq_def_Pi},
and generalise formula~\eqref{eq_formule_heuristique_avec_matrice_correl} to a non-equilibrium situation in Theorem~\ref{theo_large_devs}.

\subsection{The topology for correlations}\label{sec_topology}
Motivated by the heuristics of Section~\ref{sec_heuristics_correlation},
we focus in this article on the next scale after the density and consider the two point correlation field $\Pi^N$, 
defined in~\eqref{eq_def_Pi}. 
Zooming at the level of correlations amounts to rewriting the asymptotics of Lemma~\ref{lemm_heur_correl} as follows: compute the asymptotics of observing 
 a given correlation field $\Pi$ 
\begin{equation}
\Prob \Big(\frac{1}{T}\int_0^T \Pi^N_t\, dt \overset{weak^*}{\approx} \Pi\Big)
\quad \text{when } T \text{ and then } N \text{ are large}
.
\label{eq_proba_observer_k_general}
\end{equation}
In~\eqref{eq_proba_observer_k_general}, 
$\overset{weak^*}{\approx}$ means proximity in the weak$^*$ topology. 
In this section, we start by defining the functional space to which $\Pi^N$ belongs and the associated weak$^*$ topology.

Let us start with a few observations. By definition~\eqref{eq_def_Pi}, $\Pi^N$ can be seen as a linear form on several function spaces. 
Let $\square = (-1,1)^2$ and notice that $\Pi^N$ is symmetric in the following sense:
\begin{equation}
\forall \phi:\square\rightarrow\R,\qquad 
\Pi^N(\phi) 
= 
\Pi^N(\phi_s),\quad 
\phi_s(x,y) 
= 
\phi(x,y)/2 + \phi(y,x)/2,
\quad (x,y)\in\square
.
\label{eq_def_symmetric_of_a_function}
\end{equation}
In other words, $\Pi^N$ could really be defined on the triangle  $\{x,y\in(-1,1): x<y\}$, 
but we work on the square for symmetry reasons. 
Note also that any symmetric function $\phi$, 
i.e. $\phi=\phi_s$, 
that is $C^1$ on the whole of $\square$, 
satisfies:
\begin{align}
&\forall (x,y)\in \square,\qquad 
\partial_1\phi(x,y) 
= 
\partial_2\phi(y,x)
\nonumber\\
&\hspace{3cm}\quad \Rightarrow\quad 
\forall x\in(-1,1),\qquad 
\big(\partial_1-\partial_2\big)\phi(x,x) 
= 
0
.
\end{align}
\begin{figure}[h]
\centering
\includegraphics[width=5cm]{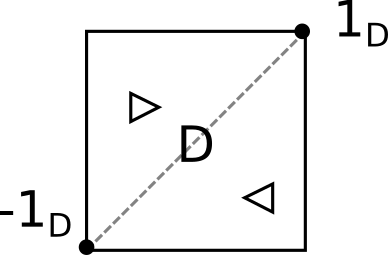}
\caption{\small 
The domain $\square = (-1,1)^2$ and the diagonal $D$ (in dashed lines) are depicted, 
with $\squaredash:= \square\setminus D$. 
According to notation \eqref{eq_notation_domain}, the lower part is $\lhd$ and the upper part $\rhd$.
The extremities of the diagonal $1_D$ and $-1_D$ are defined in \eqref{eq: extremities of the diagonal}.
}
\label{fig: domain}
\end{figure}
In view of~\eqref{eq_def_k_0}, 
the two-point correlations are symmetric functions with singularities on the diagonal $D$ of $\square$, 
defined by:
\begin{equation}
D := \big\{(x,x) : x\in(-1,1)\big\}.
\label{eq_def_diagonale_D}
\end{equation}
 We therefore cannot only consider symmetric $\phi$ that are smooth on the whole of $\square$. 
Finally, to account for the reservoirs, 
we require $\phi$ to be continuous on $\partial\square$ as well, 
and set $\phi_{|\partial\square} = 0$.

Let us now define the test functions $\Pi^N$ will act on. 
Split $\square$ as follows (see Figure \ref{fig: domain}) :
\begin{equation}
\square 
= 
\rhd \cup D\cup \lhd 
,\quad 
\rhd 
:= 
\{(x,y)\in\square : x<y\},\ 
\lhd 
:= 
\{(x,y)\in\square:x>y\},\quad 
\squaredash 
:= 
\rhd\cup\lhd
.
\label{eq_notation_domain}
\end{equation}
For $n\in\N$ and $p\geq 1$, 
let $\mathbb W^{n,p}(\squaredash) := \mathbb W^{n,p}(\rhd)\cap \mathbb W^{n,p}(\lhd)$ be the Sobolev space of functions with distributional derivatives up to order $n$ in $\mathbb L^2(\square)$. 
Properties of these spaces are recalled in Appendix~\ref{app_sobolev_spaces}. 
Note that $\mathbb L^p(\square) = \mathbb L^p(\squaredash)$ since the diagonal has vanishing two-dimensional Lebesgue measure. 
The difference between functions on $\square$ and $\squaredash$ arises in the integration by parts formula defining their weak derivatives. 
If $p=2$, we simply write $\mathbb H^n(\rhd):=\mathbb W^{n,2}(\rhd)$. 
Define then the set $\mathcal T$ of test functions:
\begin{equation}
\mathcal T
:=
\mathbb H^2(\squaredash)
= 
\mathbb H^2(\rhd)\cap \mathbb H^2(\lhd)
.
\label{eq_def_T}
\end{equation}
The set $\mathcal T$ is a separable Hilbert space, 
and $\mathcal T\subset C^0(\bar\rhd)\cap C^0(\bar\lhd)$ by Sobolev embedding, 
see Appendix~\ref{app_sobolev_spaces}, 
where $\bar\rhd,\bar\lhd$ respectively denote the closure of $\rhd,\lhd$. 
Denote then by $\mathcal T'$ the set of bounded linear forms on $\mathcal T$, 
and by $\mathcal T'_s\subset \mathcal T'$ the subset of those forms that are symmetric (recall~\eqref{eq_def_symmetric_of_a_function}):
\begin{equation}
\mathcal T'_s = \big\{ \Pi \in \mathcal T' : \forall \phi\in \mathcal T,\quad \Pi(\phi_{s}) = \Pi(\phi) \big\}.\label{eq_def_T_prime}
\end{equation}
To keep topology-related issues as simple as possible, 
we equip $\mathcal T'_s$ with the weak$^*$ topology, 
i.e. $(\Pi_n)_n\subset \mathcal T'_s$ converges to $\Pi\in\mathcal T'_s$ if and only if $\lim_{n\rightarrow\infty}\Pi_n(\phi) = \Pi(\phi)$ for each $\phi \in \mathcal T$ 
(or, equivalently, each $\phi \in \mathcal T_{\rhd}$). 
To avoid ambiguities, 
we write $(\mathcal T'_s,*)$ when we explicitly refer to the weak$^*$ topology.\\
As a bounded linear form on (a closed subset of) the Hilbert space $\mathbb H^2(\squaredash)$, 
the Riesz representation theorem allows each $\Pi \in \mathcal T'_s$ to be written as:
\begin{equation}
\Pi(\phi) := \frac{1}{4}\big<k_\Pi, \phi \big>,
\qquad \phi \in\mathcal T.
\label{eq_def_k_Pi}
\end{equation}
Above, $\big<\cdot,\cdot\big>$ denotes the standard scalar product on $\mathbb L^2(\squaredash)$ 
and duality pairing between elements of $\mathbb H^n(\squaredash)$ and $(\mathbb H^n(\squaredash))'$, $n\in\N$. 
The norm on $\mathbb L^2(\square)=\mathbb L^2(\squaredash)$ is denoted by $\|\cdot\|_2$. 
We use both $\Pi$ and $k_\Pi$ indifferently in the following.

\subsection{Large deviations for time-averaged correlations}
\label{sec: Large deviations for time-averaged correlations}
Our main result concerns the large deviation behaviour of the probability~\eqref{eq_proba_observer_k_general} in the large $T$, 
then large $N$ limits.
To state it, we need more notations.

Compared with the techniques used for the large deviations of the density  
(see e.g. Chapter 10 in~\cite{Kipnis1999}), producing atypical correlations requires  
modifying the jump rates of the dynamics by adding a long range interaction.
Thus we consider
the generator $L_h$ parametrised by a function $h:[-1,1]^2\rightarrow\R$ which is a \emph{non-local bias}, 
with the corresponding modified jump rates given for $\eta\in\Omega_N$, $i\in\{\pm(N-1)\}$ and $j<N-1$ by:
\begin{align}
c_{h}(\eta,i) 
&= 
c(\eta,i)\exp\big[\Pi^N(h)(\eta^i)-\Pi^N(h)(\eta)\big],\nonumber\\
c_h(\eta,j,j+1) 
&= 
c(\eta,j,j+1)\exp\big[\Pi^N(h)(\eta^{j,j+1})-\Pi^N(h)(\eta)\big]
.
\label{eq_def_jump_rates_H}
\end{align}
We write $\Prob_h,\E_h$ for the probability/expectation under this dynamics, 
and $\Prob^{\mu^N}_h,\E^{\mu^N}_h$ when starting from the measure $\mu^N$ on $\Omega_N$. 
The strategy is to find  the correct bias $h$ so that the rare event of observing the correlation $\Pi$ in \eqref{eq_proba_observer_k_general} becomes typical :
\begin{equation}
\Prob^{\mu^N}_h\Big(\frac{1}{T}\int_0^T \Pi^N_tdt \overset{weak^*}{\approx} \Pi\Big) 
=1 
\quad \text{when } T \text{ and then  }N \text{ are large}
.
\label{eq_proba_h_that_pi_n_close_to_k}
\end{equation}
Before stating our main Theorem \ref{theo_large_devs}, we have to introduce some restrictions, in particular on the size of the bias $h$.\\

\noindent\textbf{Main assumption and characterisation of the biases.} 
In theory, one could define the open SSEP dynamics with any value of the reservoir densities $\rho_-,\rho_+$, 
and consider any sufficiently regular $h\in\mathcal T$ and any $\Pi\in\mathcal T'_s$ in~\eqref{eq_proba_h_that_pi_n_close_to_k}. 
In practice, to focus on the key ideas and avoid many technical issues, 
we will restrict the range of $\rho_-,\rho_+$ as well as the size of the biases $h$, 
see Theorems~\ref{theo_large_devs}--\ref{theo_entropic_problem} below. 
This restriction is discussed in Section~\ref{sec_conclusion_perspectives}. 
Let us now define the set of biases.\\
Let $\pm 1_{D}$ denote the two corners of the triangle $\rhd$ 
(see Figure \ref{fig: domain}) corresponding to extremities of the diagonal $D$,  
defined in~\eqref{eq_def_diagonale_D}:
\begin{equation}
1_D 
= (1,1)
,\qquad 
-1_D 
= (-1,-1)
.
\label{eq: extremities of the diagonal}
\end{equation}
For $\epsilon>0$, define:
\begin{align}
\mathcal S(\epsilon) &= \Big\{ h\in \mathcal T:h\in \mathbb W^{4,p}(\rhd)\text{ for some }p>2, \ h \text{ is symmetric}, 
\nonumber\\
&\hspace{3cm}\|h\|_\infty,\|\partial_1 h\|_{\infty} \leq \epsilon,\ \lim_{(x,y)\rightarrow\pm 1_D}(\partial_1-\partial_2) h(x,y) = 0\Big\}
.
\label{eq_def_test_functions_s_B}
\end{align}
The condition on $(\partial_1-\partial_2)h$ at $\pm 1_D$ is purely technical. 
Introduce also the set $\mathcal S(\infty)$ of biases without size constraints:
\begin{equation}
\s(\infty)
:=
\bigcup_{\epsilon>0}\s(\epsilon)
.
\label{eq_def_s_infty}
\end{equation}

\noindent\textbf{The rate function.} 
Introduce the bilinear mapping $\mathcal M$ from $\mathbb L^2(\squaredash)^2$ to $\mathbb L^2(\squaredash)$, 
defined for $(u,v)\in\mathbb L^2(\squaredash)^2$ by:
\begin{equation}
\forall (x,y)\in\squaredash,\qquad 
\mathcal M(u,v)(x,y) = \int_{(-1,1)}u(z,x)\bar\sigma(z)v(z,y)dz,\label{eq_def_I}
\end{equation}
where $\bar\sigma$ is defined as 
\begin{equation}
\forall x\in(-1,1), \qquad
\bar\sigma(x) 
:= 
\sigma(\bar\rho(x))
\quad \text{with} \quad 
\sigma(r) = r(1-r), \quad r\in[0,1]
.
\label{eq_def_bar_sigma}
\end{equation}

For $\epsilon>0$, let $h$ belong to the set $\s(\epsilon)$ defined in~\eqref{eq_def_test_functions_s_B}. 
Introduce the functional $J_h$, 
defined for $\Pi\in \mathcal T'_s\cap \mathbb H^1(\squaredash)$ by:
\begin{align}
J_h(\Pi) 
&= 
-\frac{1}{2}\Pi\big(\Delta h + \mathcal M(\partial_1h,\partial_1h)\big) +\frac{1}{4}\int_{(-1,1)} \text{tr}_D( k_\Pi)(x) (\partial_2-\partial_1)h(x_+,x)\, dx \nonumber\\
&\quad + 
\frac{(\bar\rho')^2}{4}\int_{(-1,1)} h(x,x)\, dx - \frac{1}{8}\int_{\squaredash}\bar\sigma(x)\bar\sigma(y)\big[\partial_1 h(x,y)\big]^2\, dx\, dy
,
\label{eq_def_J_h}
\end{align}
and: 
\begin{equation}
J_h(\Pi) 
= 
+\infty \qquad\text{if }\Pi\notin \mathcal T'_s\cap \mathbb H^1(\squaredash)
.
\label{eq_J_h_infty_outside_H1}
\end{equation}
In~\eqref{eq_def_J_h}, tr$_D(k_\Pi)$ is the trace of $k_\Pi$ on the diagonal $D$ (defined in~\eqref{eq_def_diagonale_D}), 
with $k_\Pi$ related to $\Pi$ via~\eqref{eq_def_k_Pi}. 
It is well defined for $\Pi$ in  $\mathbb H^1(\squaredash)$, 
see Theorem 1.5.1.3 in~\cite{Grisvard2011}. 
Moreover, $\Pi\in\mathcal T'_s\cap \mathbb H^1(\squaredash)$ implies that $k_\Pi$ is symmetric, 
thus its trace on either side of the diagonal is the same and the notation tr$_{D}(k_\Pi)$ is not ambiguous.\\ 
Define then
the functionals $\mathcal I_\epsilon,\mathcal I_\infty:(\mathcal T'_s,*)\rightarrow\R_+$ $(\epsilon>0)$ as follows:
\begin{equation}
\mathcal I_\epsilon 
= 
\sup_{h\in\s(\epsilon)} J_h
,\qquad
\mathcal I_\infty 
=
\sup_{h\in \s(\infty)}J_h
,
\label{eq_def_rate_function}
\end{equation}
where $\epsilon>0$ stands for the restriction on the size of the biases.  
To demystify the expression~\eqref{eq_def_J_h} of $J_h$, 
let $\Pi = \frac{1}{4}\big<k,\cdot\big>\in\mathcal T'_s$ with $k\in C^3(\bar\rhd)$ a regular kernel. 
The associated correlation operator $C_k =\bar\sigma + k$ is defined for any test function  $\phi \in\mathbb L^2((-1,1)^2)$ by: 
\begin{equation}
C_k \phi (x) = \bar\sigma (x) \phi (x) +  \int_{(-1,1)}k (x,y) \phi (y)dy,
\qquad x\in[-1,1]
,
\label{eq_def_operateur_k}
\end{equation}
where the diagonal part $\bar\sigma (x)$ corresponds to the variance at a single site, 
and $k$ encodes the long range correlations.
In comparison, in the steady state (corresponding to $h=0$), the correlation operator in the large size limit is $C_{k_0} = \bar\sigma + k_0$, 
with the kernel $k_0$ introduced in~\eqref{eq_def_correl_steady_state}.

The following proposition connects the bias $h$ and the correlation kernel $k$. 
It is a classical result, proven in Appendix~\ref{app_EL}.
\begin{proposition}[Euler-Lagrange equation]
\label{prop_Euler-Lagrange equation}
Assume that the supremum in the definition~\eqref{eq_def_rate_function} is reached at some $h\in\s(\infty)$. 
Then $k$ satisfies the Euler-Lagrange equation: for each $\phi\in\mathcal T$, 
\begin{align}
&\frac{1}{2}\int_{\squaredash}\nabla (k-k_0)(x,y)\cdot \nabla\phi(x,y)\,dx\,dy 
+ \frac{1}{2}\int_{(-1,1)}\bar\sigma(z)\big<\partial_1 h(z,\cdot), C_k \partial_1\phi(z,\cdot)\big>\, dz 
=
0
.
\label{eq_EL}
\end{align}
The expression~\eqref{eq_def_rate_function} of the rate function $\mathcal I_\infty$ then simplifies for $\Pi=\frac{1}{4}\big<k,\cdot\big>$:
\begin{equation}
\mathcal I_\infty \big((1/4)\big<k,\cdot\big>\big)
=
\frac{1}{8}\int_{(-1,1)}\bar\sigma(z)\big<\partial_1 h(z,\cdot), C_k \partial_1h(z,\cdot)\big>\, dz 
.
\label{eq_rate_function_explicit}
\end{equation}
\end{proposition}
This generalises the asymptotics of the  Donsker-Varadhan functional of Lemma~\ref{lemm_heur_correl} in the reversible case.
In addition, the expression~\eqref{eq_rate_function_explicit} has the familiar form of an $\mathbb L^2$ norm of the gradient of the bias in a suitable weighted space in terms of the target distribution (see~\cite[Chapter 10]{Kipnis1999}). \\

The next theorem gives a large deviation result for the law of $\frac{1}{T}\int_0^T\Pi^N_t\, dt$. 
A more general claim is discussed in Section~\ref{sec_conclusion_perspectives}.
\begin{theorem}
\label{theo_large_devs}
Let $\rho_-\in(0,1)$. 
There is then $\epsilon_B=\epsilon_B(\rho_-)$, 
defined in Theorem~\ref{theo_entropic_problem} below, 
such that, 
if $\bar\rho'\leq \epsilon_B$, 
then the following holds. 
Let $\mathcal O,\mathcal C\subset (\mathcal T'_s,*)$ respectively be an open, closed set. Then: 
\begin{equation}
\limsup_{N\rightarrow\infty}\limsup_{T\rightarrow\infty}\frac{1}{T}\log \Prob^{\pi^N_{inv}}\Big(\frac{1}{T}\int_0^T\Pi^N_t\, dt\in \mathcal C\Big)
\leq 
-\inf_{\mathcal C} \mathcal I_{\epsilon_B}
,
\end{equation}
The rate function correctly captures the behaviour of the correlation field for kernels close to the kernel $k_0$ of the steady state, 
in the following sense. 
Let $\mathcal C_B\subset\mathcal T'_s$ be the set of correlation kernels associated with a regular bias $h\in\s(\epsilon_B)$:
\begin{align}
\mathcal C_B 
&=
\Big\{ k\in\mathcal T'_s : k\text{ solves the Euler-Lagrange}
\nonumber\\ 
&\hspace{3cm}\text{equation~\eqref{eq_EL} for some }h\in\s(\epsilon_B)\Big\}
.
\label{eq_ensemble_k_h}
\end{align}
Then, for $k\in\mathcal C_B$, 
one has $\mathcal I_\infty\big(\frac{1}{4}\big<k,\cdot\big>\big) = \mathcal I_{\epsilon_B}\big(\frac{1}{4}\big<k,\cdot\big>\big)$, 
and:
\begin{equation}
\liminf_{N\rightarrow\infty}\liminf_{T\rightarrow\infty}
\frac{1}{T}\log \Prob^{\pi^N_{inv}}\Big(\frac{1}{T}\int_0^T\Pi^N_t\, dt\in\mathcal O\Big)
\geq 
-\inf_{\mathcal O\cap \mathcal C_B} \mathcal I_{\epsilon_B}
=
-\inf_{\mathcal O\cap \mathcal C_B} \mathcal I_{\infty}.
\end{equation}
\end{theorem}
The statement of Theorem~\ref{theo_large_devs} and extensions 
(in particular a lower bound for non-regular kernels) are discussed in Section~\ref{sec_conclusion_perspectives}. 
The fact that, for $k\in \C_B$, 
the associated $h$ is a global optimiser and thus $\mathcal I_\infty,\mathcal I_{\epsilon_B}$ agree at $\frac{1}{4}\big<k,\cdot\big>$ is proven in Section~\ref{sec_dynamical_part_lower_bound}.
\begin{remark}
Recall that the kernel $k_0$ of the steady state $\pi^N_{inv}$ of the open SSEP in the large $N$ limit, defined in~\eqref{eq_def_k_0}, was observed to be smooth away from the diagonal. 
As Theorem~\ref{theo_large_devs} shows, 
by definition of the rate function (see~\eqref{eq_J_h_infty_outside_H1}),  
it is in fact a general property that the time-average of $\Pi^N_\cdot$ is much more regular than an element of $\mathcal T'_s$ when $T,N$ are large: it belongs to $\mathbb H^1(\squaredash)$.
\demo
\end{remark}
\subsection{The relative entropy method}\label{sec_relative_entropy_into}
In this section, we explain the method used to establish Theorem~\ref{theo_large_devs}, i.e. to study, for some $\Pi\in\mathcal T'_s$, the probability:
\begin{equation}
\Prob^{\pi^N_{inv}}\Big(\frac{1}{T}\int_0^T \Pi^N_tdt \approx \Pi\Big)
\quad 
\text{when }T\text{, then }N\text{ are large}
.
\label{eq_proba_observer_k_general_bis}
\end{equation}
Understanding correlations out of equilibrium ($\rho_-\neq \rho_+$) is notoriously difficult. 
Existing results in the literature deal either with the equilibrium case $\rho_-=\rho_+$ without bias (i.e $h=0$), 
see e.g.~\cite{Goncalves2019}; 
or use methods particular to the $h=0$ case, 
which cannot easily be generalised to $h\neq 0$~\cite{Goncalves2020}. 
The methods rely on explicit knowledge of the invariant measure of the dynamics. 
However, for $h\neq 0$, the invariant measure $\pi^N_{inv,h}$ of the tilted dynamics $\Prob_h$ defined in~\eqref{eq_def_jump_rates_H} is not known explicitly. 
Even for $h=0$, 
where the invariant measure $\pi^N_{inv}$ is well understood~\cite{Derrida2002}, 
its complexity makes the study of the probability in~\eqref{eq_proba_observer_k_general_bis} difficult. 
To study~\eqref{eq_proba_observer_k_general_bis} in full generality, we therefore need a different approach. \\

The key idea is to find an approximation of the invariant measure, 
which is sufficiently close to control the large time behaviour of correlations, yet simple enough to make explicit computations possible. 
The proximity of the law of the dynamics to this approximate invariant measure is quantified through the relative entropy method, 
relying on the beautiful generalisation by Jara and Menezes~\cite{Jara2018} of the ideas of Yau~\cite{Yau1991}. 
The relative entropy method is presented in more mathematical terms in Section~\ref{sec_entropy_and_FK}, 
and here we only describe it informally. 
In our context, the relative entropy method consists in finding a measure $\mu^N$ on the state space $\Omega_N$, 
that is both sufficiently simple to perform explicit computations, 
and as close as possible to the invariant measure $\pi^N_{inv,h}$. 
This closeness to the invariant measure aims at ensuring that, 
if the dynamics $\Prob_h$ starts from $\mu^N$ then, at time $t\geq 0$, 
the law $f_t\mu^N$ of the dynamics is still close to $\mu^N$. The proximity to the invariant measure is quantified by the relative entropy $H(f_t\mu^N|\mu^N)$. 
The level of precision needed on this relative entropy depends both on the quantity to study - e.g. the density, 
the density fluctuations, or in our case the correlations; 
and on the time range one wishes to probe. \\

Here, we improve on the estimates of Jara and Menezes~\cite{Jara2018}, 
obtaining sufficiently precise relative entropy estimates to study the probability~\eqref{eq_proba_observer_k_general_bis}.
To do so, for a bias $h$ in the set $\s(\epsilon_B)$, 
defined in~\eqref{eq_def_test_functions_s_B}, 
we compare the law of the dynamics at each time with a discrete Gaussian measure $\nu^N_{g_h}$ of the following form: for each $\eta$ in $\Omega_N$,
\begin{equation}
\nu^N_{g_h}(\eta) 
= 
\big(\mathcal Z^N_{g_h}\big)^{-1}\exp\big[2\Pi^N(g_h)\big]\bar\nu^N(\eta),
\quad \text{with}\quad 
\bar\nu^N(\eta) := \bigotimes_{i\in\Lambda_N}\text{Ber}(\bar\rho_i).
\label{eq_def_bar_nu_g_intro_bis}
\end{equation}
Above, for $\rho\in[0,1]$, 
$\text{Ber}(\rho)$ is the Bernoulli measure on $\{0,1\}$ with parameter $\rho$. 
The partition function $\mathcal Z^N_{g_h}$ is a normalisation factor.
The function $g_h:\squaredash\rightarrow\R$, 
which solves a partial differential equation depending on the dynamical bias $h$, 
is the function that allows us to minimise the entropy production $\partial_t H(f_t\nu^N_{g_h}|\nu^N_{g_h})$, as stated in Theorem~\ref{theo_entropic_problem} below.

Let us first consider $h=0$ and give a heuristic reason why  the invariant measure $\pi^N_{inv}$ can be approximated 
by a measure of the form~\eqref{eq_def_bar_nu_g_intro_bis} for a well chosen function $g_0 := g_{h=0}$. 
Under the product Bernoulli measure $\bar\nu^N$ 
(which has the same density profile as the steady state, but no correlations), 
for each test function $\phi\in C^0([-1,1])$, 
the fluctuation field $Y^N(\phi)=N^{-1/2}\sum_i \bar\eta_i\phi_i$ converges to a Gaussian field with covariance $\bar\sigma$. 
Assuming for a moment that each $\bar\eta_i$ is a continuous variable, 
the law of fluctuations under $\nu^N_g$ is thus close to a Gaussian field with covariance $(\bar\sigma^{-1} - g)^{-1}$. 
This claim is made rigorous in the forthcoming paper~\cite[Appendix A]{dagallierFluctuationsCorrelationsWeakly2023a}. 
Thus, in order for $\nu^N_g$ to have the same correlations $k_0$ 
(defined in~\eqref{eq_def_correl_steady_state}) as the steady state $\pi^N_{inv}$ when $N$ is large, 
one can take $g=g_0$ to be the inverse correlation kernel of the steady state:
\begin{equation}
C_{k_0} 
=
\bar\sigma + k_0 =: (\bar\sigma^{-1} - g_0)^{-1}
\quad\Leftrightarrow\quad g_0 := \bar\sigma^{-1} - (\bar\sigma + k_0)^{-1}
.
\label{eq_def_g_0}
\end{equation}
Note that $C_{k_0}$ is indeed invertible as, by construction, 
$k_0$ is a negative definite operator on $\mathbb L^2(\square)$:
\begin{align}
\forall \phi \in\mathbb L^2((-1,1))\setminus \{0\},
\qquad 
&\int_{(-1,1)^2} \phi(x)k_0(x,y) \phi(y)\, dx\, dy<0
.
\label{eq_k_0_g_0_negative_kernels}
\end{align}
The fact that $g_0$ is indeed a kernel operator 
follows from $g_0 C_{k_0} = \bar\sigma^{-1}k_0$, 
see the first theorem in section 4.6.1 in~\cite{Bukhvalov1996}. 
This identity also implies that $g_0$ is a negative definite kernel.
Moreover, as a function on $\squaredash$, 
$g_0$ inherits the regularity of $k_0$ 
(smooth all the way up to the diagonal $D$, 
but with normal derivative having a jump across $D$).  
The following proposition and theorem provide a systematic way of approximating the invariant measure for a 
non-local bias $h\in\s(\epsilon)$ for small enough $\epsilon>0$.
\begin{proposition}[Main equation]\label{prop_main_equation}
Let $0<\rho_-,\rho_+<1$, 
take a bias $h\in\s(\infty)$ and consider the following problem with unknown $g$, 
referred to as \emph{the main equation:}
\begin{equation}
\begin{cases}
&\displaystyle{\Delta (g-h)(x,y) + \frac{\bar\sigma'(x)}{\bar\sigma(x)}\partial_1(2g-h)(x,y) +\frac{\bar\sigma'(y)}{\bar\sigma(y)}\partial_2(2g-h)(x,y) }\hspace{0.75cm} \text{for }(x,y)\in \squaredash\\
&\hspace{0.2cm}\displaystyle{  + \int_{(-1,1)}\bar\sigma(z)\big[\partial_1 g(z,x)\partial_1 (g-h)(z,y)+ \partial_1 g(z,y)\partial_1 (g-h)(z,x)\big]dz = 0,}\\
&g  = 0 \quad \text{ on }\partial\square,\\
&(\partial_2-\partial_1) (h-g)(x_+,x)= (\partial_1-\partial_2)(h-g)(x_-,x) = \displaystyle{\frac{(\bar\rho')^2}{\bar\sigma(x)^2}} \quad \text{for }x\in (-1,1).
\end{cases}
\label{eq_main_equation}
\end{equation}
There is $\epsilon_0(\rho_-)>0$ and a function $\delta(\cdot)$ with $\lim_{x\downarrow 0} \delta(x) =0$ such that, 
for any $\epsilon\in(0, \epsilon_0(\rho_-)]$, 
$\bar\rho'\leq \epsilon$ and $h\in\s(\epsilon)$ implies that~\eqref{eq_main_equation} has a unique solution $g_h\in g_0+\s(\delta(\epsilon))$. 
\end{proposition}
Proposition~\ref{prop_main_equation} is proven in Appendix~\ref{app_Poisson}. 
The next theorem provides the key control on the dynamics and determines the parameter $\epsilon_B$ mentioned in Theorem~\ref{theo_large_devs}.
\begin{theorem}\label{theo_entropic_problem}
Let $0<\rho_-\leq \rho_+<1$. 
There is $\epsilon_B=\epsilon_B(\rho_-)$ such that, 
if $\rho'\leq \epsilon_B$ and $h\in\s(\epsilon_B)$, 
then the measure 
$\nu^N_{g_h}$ defined in~\eqref{eq_def_bar_nu_g_intro_bis} is a good approximation of the invariant measure of the dynamics $\Prob_h$ with bias $h$ in the following sense. 
Let $f_t\nu^N_{g_h}$ denote the law of $\Prob_h$ at time $t\geq 0$. 
There is then $C,K>0$ depending on $h,\rho_\pm$ such that:
\begin{equation}
\forall t\geq 0,\qquad 
H(f_t\nu^N_{g_h}|\nu^N_{g_h}) 
\leq 
e^{-Kt}H(f_0\nu^N_{g_h} |\nu^N_{g_h}) + \frac{C}{N^{1/2}}.
\end{equation}
\end{theorem}
Theorem~\ref{theo_entropic_problem} is proven in Section~\ref{sec_computations}. 
The conditions that the parameter $\epsilon_B$ must satisfy are summarised in Definition~\ref{def_epsilon_B_and_G}. 
\begin{remark}
	For each bias $h\in\s(\epsilon_B)$, 
	we prove in Appendix~\ref{app_equivalence_EL_main_equation} that obtaining a solution $g_h$ to the main equation~\eqref{eq_main_equation} with the desired regularity is equivalent to obtaining a classical solution $k_h$ of the Euler-Lagrange~\eqref{eq_EL} with the same regularity, 
	with $g_h$ and $k_h$ related through the following identity of the associated operators:
	\begin{equation}
	\bar\sigma + k_h 
	= 
	(\bar\sigma^{-1}-g_h)^{-1}
	.
	\label{eq_relation_k_h_g_h}
	\end{equation}
	In particular this validates the heuristics behind the choice~\eqref{eq_def_g_0} of $g_0$. 
	The following diagram summarises the relationships between $h$, $k$ and $g$.
	$$
\xymatrix{ k \ar[rr] \ar[rd]_{\eqref{eq_relation_k_h_g_h}} &&  h \ar[ll]_{\eqref{eq_EL}}  \ar[ld]\\ & g \ar[ru]_{\eqref{eq_main_equation}}  \ar[lu]  }
$$
\demo
\end{remark}
\begin{remark}
The entropy control of Theorem 2.7 has many interesting consequences. 
For instance, it allows one to compute two point correlations under the invariant measure. 
Indeed,  
if $\pi^N_{inv,h}$ is the invariant measure for $\Prob^N_h$, 
then for any bounded $f_1,f_2:(-1,1)\to\R$, 
\begin{align}
\pi^N_{inv,h}\big[Y^N(f_1)Y^N(f_2)\big] 
&=
\nu^N_{g_h}\big[Y^N(f_1)Y^N(f_2)\big] 
+o_N(1)
\nonumber\\
&=
\int_{(-1,1)} f_1(x) (\bar\sigma^{-1}-g_h\big)^{-1}f_2(x)
\, dx
+o_N(1)
.
\end{align}
By Pinsker's inequality, the relative entropy controls the total variation distance, which together with the uniform integrability of the correlation field, Lemma~\ref{lemm_moment_3_halves} implies the first equality. 
The second equality follows from the explicit knowledge of correlations under $\nu^N_{g_h}$, 
see Proposition A.2 in~\cite{dagallierFluctuationsCorrelationsWeakly2023a}. 
\demo
\end{remark}
\subsection{Conclusion and perspectives}
\label{sec_conclusion_perspectives}
\subsubsection{Extensions of the large deviation principle} 
The large deviation result of Theorem~\ref{theo_large_devs} is stated starting from the invariant measure $\pi^N_{inv}$. 
In fact any choice of initial condition is possible, with no change to the proof.

An advantage of our very precise, 
quantitative microscopic estimates is that Theorem~\ref{theo_large_devs} also holds if one takes a diverging sequence $T_N$ of times. 
Moreover, it is also possible to take the limits in the opposite order (large $N$, then large time). 
In this case the proof of Theorem~\ref{theo_large_devs} is slightly simpler as the relative entropy does not need to be controlled uniformly in time. 
The proof is otherwise nearly the same.

Theorem~\ref{theo_large_devs} gives large deviation bounds in the weak-$*$ topology. 
In fact, the theorem also holds in the strong dual topology.  This can be seen to hold with no change to the proof for the upper bound. 
For the lower bound, one needs to be more careful, 
so we chose to work in a weaker topology to avoid technicalities. 
An extension of the lower bound to non-regular correlation fields is also possible, 
and sketched in Section~\ref{sec_extension_non_regular}.
\subsubsection{Restriction on the biases and reservoir densities} 
In Theorem~\ref{theo_large_devs}, 
restrictions are imposed on both the biases $h$ and the difference $\rho_+-\rho_-$ of the reservoir densities.
\begin{itemize} 
\item 
Some restriction has to be imposed on the bias size $h$, otherwise one expects that not only the correlation structure, but also the typical density of the corresponding biased dynamics changes. 
Nevertheless, the conditions in~\eqref{eq_def_test_functions_s_B} are far from sharp.
\item 
The restriction on the slope $\bar \rho' = (\rho_+-\rho_-)/2$ has two technical purposes. 
First, this is convenient to prove existence of $g_h$ solving the main equation \eqref{eq_main_equation} 
(even though it should be possible to remove this assumption). 
Secondly, it is used in the derivation of the relative entropy bound of Theorem~\ref{theo_entropic_problem} to control the relative entropy uniformly in time: 
one has to make sure that error terms are bounded by $cN$ uniformly on the state space for sufficiently small $c$ in order to estimate their exponential moments, 
and $c$ depends on $\bar\rho'$.  
Again, one should be able to relax this assumption, 
for instance through a priori estimate on the size of the error terms using large deviation results for the density. 
All other uses of the fact that $\bar\rho'$ is small can be relaxed without additional work, but at the cost of more technicalities. 
\end{itemize}
Let us however stress again that the large deviation result of Theorem~\ref{theo_large_devs} is sharp for correlation kernels close to the kernel $k_0$ of the steady state~\eqref{eq_def_k_0}, 
as in that case the rate functions $\mathcal I_{\epsilon_B},\mathcal I_\infty$ (recall~\eqref{eq_def_rate_function}) coincide.
\subsubsection{The relative entropy method}\label{sec_2.6.3}
The relative entropy approach used to obtain large deviations can be applied to many other settings, 
still in dimension one. 
Let us list a few.\\

\indent A first direction is the study of more general diffusive gradient systems with the following restrictions. 
For systems more complicated than the SSEP,	
	the behaviour of correlations at the boundary may be problematic, 
	and it is unclear whether one could still get a relative entropy estimate. 
	However, results of the present paper should carry over to any diffusive gradient system on the torus 
	(see the discussion in Section 8.1. in~\cite{Jara2018}). 
	Additional work would be needed, 
	in particular some large deviation estimates for the density to ensure that higher order correlations on the diagonal (typically of the form $\sum_{i}\bar\eta_i\bar\eta_{i+1}\bar\eta_{i+2}$) cannot be very large.  
	The relative entropy bound would then only be of size $o_N(1)$, 
	instead of $O(N^{-1/2})$ as in Theorem~\ref{theo_entropic_problem}. \\
\indent In the paper~\cite{dagallierFluctuationsCorrelationsWeakly2023a}, 
the present refinement of the relative entropy method is used to study fluctuations in the WASEP on the one-dimensional discrete torus, 
constrained to produce a macroscopic current on the time interval $[0,T]$. 
The presence of this current is known to create a rich correlation structure. 
An expression of these correlations in the long time limit has been conjectured in~\cite{Bodineau2008}. 
In~\cite{dagallierFluctuationsCorrelationsWeakly2023a}, 
this conjecture is proved and the fluctuations of this process are described. 
In a related manner, 
Derrida and Sadhu~\cite{Derrida2019} study the density large deviations for a non equilibrium SSEP conditioned to have an atypical macroscopic current.\\

Another direction of inquiry concerns the extension of the relative entropy results to higher dimensions. 
This seems to be a complicated problem.  
One major hurdle is the fact that, 
in higher dimensions, correlations are not smooth on the diagonal (see e.g.~\cite{Spohn1983}), 
while the relative entropy method requires smoothness. 
\subsubsection{Entropy of the invariant measure}
\label{sec_entropy_invariant_mes}
The relative entropy estimate of Theorem~\ref{theo_entropic_problem} can be used to recover the asymptotic behavior
of the entropy of the steady state found in~\cite{Derrida2007}. 
The entropy of $\pi^N_{inv}$ is by definition:
\begin{equation}
S(\pi^N_{inv}) 
=
-\sum_{\eta\in\Omega_N}\pi^N_{inv}(\eta)\log\pi^N_{inv}(\eta)
\geq 
0
.
\end{equation}
In~\cite{Derrida2007}, 
the authors compute the first and second order contributions to this entropy as a function of the size of the system. 
The leading order term is given by the entropy of a Bernoulli measure with the correct density profile, 
and heuristics showing that the next order is the relative entropy of a Gaussian with appropriate covariance are provided. 
The computations rely on the precise knowledge of $n$-point correlation functions between each individual lattice sites and for each integer $n$. 
These are obtained by recursion equations specific to the open exclusion process. 

In the range of $\bar\rho'$ in which it applies, 
Theorem~\ref{theo_entropic_problem} gives the same result as we explain next.  
Contrary to~\cite{Derrida2007}, however, 
no knowledge of $n$-point correlations for $n$ larger than $2$ are required.

Let us now see why. 
Recall that $\bar\nu^N$ is the Bernoulli product measure~\eqref{eq_def_bar_nu_g_intro_bis}. 
Then:
\begin{align}
S(\bar\nu^N)
-S(\pi^N_{inv})  
&=
\sum_{\eta\in\Omega_N}(\pi^N_{inv}-\bar\nu^N)(\eta)\log\bar\nu^N(\eta) + \sum_{\eta\in\Omega_N}\pi^N_{inv}(\eta)\log\Big(\frac{\pi^N_{inv}(\eta)}{\bar\nu^N(\eta)}\Big)
\nonumber\\
&=
H(\pi^N_{inv}|\bar\nu^N)
,
\end{align}
where the last equality is obtained by noticing that the first term vanishes, since for each $\eta\in\Omega_N$:
\begin{equation}
\log\bar\nu^N(\eta)
=
\sum_{i\in\Lambda_N}\eta_i\log\Big(\frac{\bar\rho_i}{1-\bar\rho_i}\Big)+ \sum_{i\in\Lambda_N} \log(1-\bar\rho_i)
,
\end{equation}
and each $\eta_i$ ($i\in\Lambda_N$) has the same average under both $\pi^N_{inv}$ and $\bar\nu^N$. 
Theorem~\ref{theo_entropic_problem} gives $H(\pi^N_{inv}|\nu^N_{g_0})=O(N^{-1/2})$ 
(with $\nu^N_{g_0}$ defined in~\eqref{eq_def_bar_nu_g_intro_bis}). 
To use that result, let us turn $\bar\nu^N$ into $\nu^N_{g_0}$:
\begin{align}
H(\pi^N_{inv}|\bar\nu^N)
&=
\sum_{\eta\in\Omega_N}\pi^N_{inv}(\eta)\log\Big(\frac{\pi^N_{inv}(\eta)}{\nu^N_{g_0}(\eta)}\Big)
-\sum_{\eta\in\Omega_N}\pi^N_{inv}(\eta)\log\Big(\mathcal Z^N_{g_0}e^{-2\Pi^N(g_0)}\Big)
\nonumber\\
&=
H(\pi^N_{inv}|\nu^N_{g_0})
-\log\mathcal Z^N_{g_0} + 2\pi^N_{inv}(\Pi^N(g_0))
.
\end{align}
The first term is bounded by $O(N^{-1/2})$. 
The other two terms can be computed explicitly. 
Indeed, the measure $\pi^N_{inv}$ in the last term can be replaced by $\nu^N_{g_0}$, 
again using Theorem~\ref{theo_entropic_problem} to argue that $\|\pi^N_{inv}-\nu^N_{g_0}\|_{TV}=o_N(1)$ and that $(\Pi^N(g_0))_{N}$ is uniformly integrable under $\pi^N_{inv}$ (see Lemma~\ref{lemm_moment_3_halves}). 

Building on the fact that $\nu^N_{g_0}$ is very close to a Gaussian measure (it approximately satisfies Wick theorem, see Appendix A in~\cite{dagallierFluctuationsCorrelationsWeakly2023a}), 
which one can show through elementary but tedious computations in the spirit of Lemma~\ref{lemm_bound_Z_g}, 
one obtains in particular:
\begin{equation}
\sup_{i\in \Lambda_N}\big|\nu^N_{g_0}((\bar\eta_i)^2)
- C^N_{g_0}(i,i)\big|
=
o_N(1)
,\qquad 
\sup_{i\neq j\in \Lambda_N} N\big|\nu^N_{g_0}(\bar\eta_i\bar\eta_j)
- C^N_{g_0}(i,j)\big| 
=
o_N(1)
,
\end{equation}
with, interpreting $g_0$ below as the matrix $((g_0)_{i,j})$:
\begin{equation}
C^N_{g_0}(i,j) 
:= 
\Big(\bar\sigma^{-1} {\bf 1}_{i =j} -\frac{g_0}{N}{\bf 1}_{i\neq j}\Big)^{-1}(i,j)
.
\end{equation}
Adding and subtracting a diagonal term, this implies:
\begin{align}
\nu^N_{g_0}(2\Pi^N(g_0)) 
&=
-\frac{1}{2}\text{Tr}\big[C^N_{g_0}\big(C^N_{g_0}\big)^{-1}\big]
+\frac{1}{2}\sum_i\nu^N_{g_0}(\bar\eta_i^2)\bar\sigma_i^{-1}
+o_N(1)
=
o_N(1)
,
\end{align}
where we used $\sup_i|\nu^N_{g_0}(\bar\eta_i^2) - \bar\sigma_i|=o_N(1)$ as proven in Lemma~\ref{lemm_bound_correlations_bar_nu_G}. 
Similarly, write:
\begin{align}
\log \mathcal Z^N_{g_0} 
&=
\int_0^1\partial_t \log \mathcal Z^N_{tg_0}\, dt
=
\int_0^1\nu^N_{tg_0}\big(2\Pi^N(g_0)\big)\, dt
\nonumber\\
&=
\int_0^1 \Big[\frac{1}{2N}\text{Tr}\Big( C^N_{tg_0} (g_0{\bf 1}_{i\neq j})\Big)+o_N(1)\Big]\, dt
.
\end{align}
A careful control of the error terms shows that the integral of $o_N(1)$ is still $o_N(1)$. 
Taking out $\bar\sigma^{1/2}$ on each side and defining:
\begin{equation}
M^\sigma(i,j)
:=
\frac{1}{N}\bar\sigma^{1/2}_i (g_0)_{i,j}{\bf 1}_{i\neq j}\bar\sigma^{1/2}_j
,\qquad
i,j\in\Lambda_N^2
,
\end{equation}
we end up with, writing also $\text{sp}$ for the spectrum of $M^\sigma$:
\begin{align}
\log \mathcal Z^N_{g_0}
&=
\frac{1}{2}\int_0^1\text{Tr}\Big( (\text{id} - tM^\sigma)^{-1}M^\sigma\Big)\, dt +o_N(1)
=
\frac{1}{2}\int_0^1\sum_{\lambda\in\text{sp}}\frac{\lambda}{1-t\lambda}\, dt +o_N(1)
\nonumber\\
&=
\frac{1}{2}\sum_{\lambda\in\text{sp}}\log\Big(\frac{1}{1-\lambda}\Big)+o_N(1)
.
\end{align}
Thus:
\begin{equation}
\log\mathcal Z^N_{g_0}
=
\frac{1}{2}\log\det\Big(\bar\sigma^{-1/2}C^N_{g_0}\bar\sigma^{-1/2}\Big) 
+o_N(1)
.
\end{equation}
This implies the result of~\cite{Derrida2007}:
\begin{equation}
S(\pi^N_{inv})
=
S(\bar\nu^N) 
-
\frac{1}{2}\log\det\Big(\bar\sigma^{-1/2}C^N_{g_0}\bar\sigma^{-1/2}\Big) 
+o_{N}(1)
.
\end{equation}
In fact all error terms can be shown to be $O(N^{-1/2})$, 
consistent with the fact that the next order correction to $S(\pi^N_{inv})$ should come from three point correlations. 
The above computation does not use any special feature of $\pi^N_{inv}$ and in particular generalises to the invariant measures of the tilted dynamics $\Prob^N_h$ 
as well as, more interestingly, 
all models for which an equivalent of Theorem~\ref{theo_entropic_problem} holds (having $o_N(1)$ relative entropy bound rather than $O(N^{-1/2})$ is enough). 
This is the case for one-dimensional Glauber + Kawasaki dynamics on a torus as will be shown in future work. 
More generally it should be the case for all diffusive one-dimensional gradient models on the torus as discussed in Section~\ref{sec_2.6.3}.
\section{Main ingredient: the entropic estimate}\label{sec_computations}
In this section, 
we provide the key microscopic estimates to study the long-time behaviour of the process $(\Pi^N_t)_{t\geq 0}$, 
i.e. we prove Theorem~\ref{theo_entropic_problem}. 
We follow the strategy of Jara and Menezes~\cite{Jara2018}--\cite{Jara2020} and improve the controls by tuning correlations of the reference measure. 
The same kind of computations give the expression of the Radon-Nikodym derivative $D_h = d\Prob_h/d\Prob$ for $h\in \s(\epsilon_B)$, 
stated in Proposition~\ref{prop_expression_der_radon_nykodym} at the end of the section.
\subsection{The relative entropy method and Feynman-Kac inequality}\label{sec_entropy_and_FK}
Let us now recall the main features of the relative entropy method.  
Let $(\omega_t)_{t\geq 0}$ be a Markov chain on a state space $\Omega$, assumed to be finite for simplicity. 
Let $\Prob, \E$ denote the associated probability and expectation. 
Let $V:\Omega\rightarrow\R$. 
One would like to estimate quantities of the form:
\begin{equation}
\E\big[V(\omega_t)\big],\qquad t>0. \label{eq_motivation_relative_entropy}
\end{equation}
The entropy inequality provides a tool to estimate~\eqref{eq_motivation_relative_entropy}. 
Let $\mu$ be any probability measure on $\Omega$ and $f_t\mu$ be the law of $\omega_t$, $t\geq 0$. 
Then, for any $\gamma>0$:
\begin{equation}\label{eq_entropy_inequality}
\forall t\geq 0,\qquad 
\E\big[V(\omega_t)\big] 
\leq 
\gamma^{-1}H(f_t\mu|\mu)+ \gamma^{-1}\log \mu\big(\exp[\gamma V]\big)
.
\end{equation}
\begin{remark}
The symbol $\E$ always denotes dynamical expectations. 
In contrast, static expectations with respect to a measure $\mu$ are denoted by $\mu[\cdot]$ (or $\mu(\cdot)$).
\demo
\end{remark}
Through~\eqref{eq_entropy_inequality}, 
the dynamical estimate of $V$ in~\eqref{eq_motivation_relative_entropy} is reduced to a static problem: 
a relative entropy estimate, and a concentration-of-measure result under $\mu$. 
The relative entropy method aims at finding a measure $\mu$ which satisfies the following two criteria: 
exponential moments must be under control, 
and the relative entropy $H(f_t\mu|\mu)$, $t\geq 0$ must be sufficiently small, 
the size depending on the kind of observables $V$ one is interested in.\\
The analysis in~\cite{Jara2018}-\cite{Jara2020} greatly improves the existing method to control $H(f_t\mu|\mu)$, $t\geq 0$. As a starting point, Jara and Menezes revisit Yau's entropy bounds in the following form.
\begin{lemma}[Lemma A.1 in~\cite{Jara2018}]\label{lemm_entropy_FK}
Let $(\omega_t)_{t\geq 0}$ be a Markov chain on a finite state space $\Omega$, 
with jump rates $(c(\omega,\omega'))_{(\omega,\omega')\in\Omega^2}$. 
Denote by $L$ its generator and by $\Gamma$ the corresponding carré du champ operator:
\begin{equation}
\forall \omega\in\Omega,\forall f : \Omega\rightarrow\R,\qquad  \Gamma f(\omega) = \sum_{\omega'\in\Omega}c(\omega,\omega')\big[f(\omega')-f(\omega)\big]^2.\label{eq_carre_du_champ_general}
\end{equation}
Let $\mu$ be a probability measure on $\Omega$ satisfying $\inf_{\omega\in \Omega}\mu(\omega)>0$. 
Let $f_t\mu$ be the law of the process $(\omega_s)_{s\geq 0}$ at time $t\geq 0$. 
Then:
\begin{align}
\forall t\geq 0,\qquad \partial_tH(f_t\mu|\mu) \leq -\mu\big(\Gamma (\sqrt{f_t})\big) + \mu\big(f_tL^*{\bf 1}\big),\label{eq_entropy_dissipation}
\end{align}
where $L^*$ is the adjoint of $L$ in $\mathbb L^2(\mu) = \{f:\Omega\rightarrow\R : \mu(f^2)<\infty\}$. It acts on $f:\Omega\rightarrow\R$ according to:
\begin{equation}
\forall \omega\in\Omega,\qquad 
L^*f(\omega) 
= 
\sum_{\omega'\in\Omega}\Big[c(\omega',\omega)f(\omega')\frac{\mu(\omega')}{\mu(\omega)} - c(\omega,\omega')f(\omega)\Big]
.
\label{eq_def_L_star_H_general_mesure}
\end{equation}
\end{lemma}
Since the adjoint $L^*$ is known explicitly in terms of $\mu$,~\eqref{eq_entropy_dissipation} provides a way to estimate $\partial_tH(f_t\mu|\mu)$, $t\geq 0$. An estimate of $H(f_t\mu|\mu)$ follows by applying the entropy and Gronwall inequalities.

The same estimates used to bound the relative entropy will allow us, 
together with a log-Sobolev inequality (see Lemma~\ref{lemm_LSI_sec3}), 
to get estimates on exponential moments. 
This is a consequence of a bound involving the Feynman-Kac formula, 
stated now for future reference.
\begin{lemma}[Feynman-Kac inequality, Lemma A.2. in~\cite{Jara2018}]\label{lemm_FK}
For $V : \Omega\rightarrow\R$ and $T\geq 0$, 
\begin{equation}\label{eq_FK_sous_mu_general}
\log \E^{\mu}\Big[\exp\Big(\int_0^T V(\omega_t)\, dt\Big)\Big] 
\leq 
T\sup_{f\geq 0 :\mu(f)=1}\Big\{\mu(fV) - \frac{1}{2}\mu\big(\Gamma(\sqrt{f})\big) + \frac{1}{2}\mu\big(fL^*{\bf 1}\big)\Big\}
.
\end{equation}
\end{lemma}
The main difficulty to use~\eqref{eq_entropy_dissipation}--\eqref{eq_FK_sous_mu_general} is to control the term $L^*{\bf 1}$. 
Notice that $L^*{\bf 1}=0$ if and only if $\mu=\pi$ is the invariant measure. 
The quantity $L^*{\bf 1}$ thus appears as a way to quantify the proximity of $\mu$ to the invariant measure $\pi$. This serves as an informal guiding principle for the choice of $\mu$:
\begin{equation}
\text{for each }f\geq 0\text{ with }\mu(f)=1,\quad \mu\big(fL^*{\bf 1}\big)\text{ must be  small}.\label{eq_L_star_1_must_be_small}
\end{equation}
We are going to apply Lemma~\ref{lemm_entropy_FK} to the dynamics $\Prob_h$, defined in~\eqref{eq_def_jump_rates_H}, for $h\in \s(\epsilon_B)$. 
Let us emphasize again that "small" must always be understood in comparison with the size of the $V$'s one wishes to estimate as in~\eqref{eq_motivation_relative_entropy}. 
Typically, if one wants to study the hydrodynamic limit of the density of particles in a $d$-dimensional open SSEP on a lattice of side-length $N$, 
related observables are of the form $\sum_{i}\eta_i\phi_i\approx N^d$ for a test function $\phi$. 
The rule of thumb is then that one needs $o(N^d)$ bounds on the relative entropy. 
These can be achieved if $\mu$ is a product measure with the same densities as those of the invariant measure at each site, i.e:
\begin{equation}
\mu= \bar\nu^N = \bigotimes_{i\in\Lambda_N}\text{Ber}(\bar\rho_i),\label{eq_def_mu_product}
\end{equation}
with Ber$(\rho)$ the Bernoulli measure on $\{0,1\}$ with parameter $\rho\in(0,1)$, 
and $\bar\rho$ the steady state density profile in the large $N$ limit, 
see~\eqref{eq_def_profil_invariant}. 
In contrast, consider the fluctuation field:
\begin{equation}
\forall \phi:(-1,1)\rightarrow\R,\qquad 
Y^N_t(\phi) 
= 
\frac{1}{N^{1/2}}\sum_{i\in \Lambda_N}\bar\eta_i(t)\phi(i/N)
.
\end{equation}
The typical observables to study $Y^N_\cdot$ should be of the form $N^{1/2}Y^N(\phi)$ which is typically of order $N^{1/2}$, 
so one needs $o(N^{1/2})$ bounds on the relative entropy (in fact $o(N^{d/2})$ bounds in dimension $d$). 
Remarkably, while one could expect that some information on the correlation structure of the invariant state should be necessary to study $Y^N$, 
Jara and Menezes~\cite{Jara2018} managed to obtain such bounds on the relative entropy by still taking $\mu$ product as in~\eqref{eq_def_mu_product}. 
To do so, they set up a general renormalisation scheme to bound $\mu(fL^*{\bf 1})$, for a $\mu$-density $f$, in terms of the carré du champ, and objects that can be estimated by the entropy inequality. 
Precisely, they manage to prove bounds of the form:
\begin{equation}
H(f_t\mu|\mu)
\leq 
C(T) a_d(N)N^{d-2}
\qquad
\text{with}\qquad
a_d(N) 
=
\begin{cases}
N\quad &\text{if }d=1\\
\log N\quad &\text{if }d=2\\
1\quad &\text{if }d\geq 3.
\end{cases}
\label{eq_entropy_bounds_JM}
\end{equation}
This is enough to study fluctuations in dimension $d<4$. 
They also argue that these bounds are the best possible when $\mu$ is product.\\

Let us come back to the study of the correlation process $\Pi^N_\cdot$ in the (one-dimensional) open SSEP. 
Observables, of the form $\Pi^N(\phi)$, should be bounded with $N$, so we need $o_N(1)$ bounds on the relative entropy at each time. 
The measure $\mu$ therefore cannot be taken product: one needs to include information on the correlations under the invariant measure in $\mu$. 
With~\eqref{eq_L_star_1_must_be_small} in mind, 
we look for $\mu$ that has both the same density at each site, 
and the same two-point correlations as the invariant measure - which are in general not known - when $N$ is large. 
We tune these unknown correlations in an indirect way, 
taking a smooth function $g:\squaredash\rightarrow\R$, 
and looking for the optimal choice of $g$ such that the measure $\mu = \nu^N_g$ satisfies~\eqref{eq_L_star_1_must_be_small}, 
with:
\begin{equation}
\nu^N_{g} := \frac{1}{\mathcal Z^N_g}\exp\big[2\Pi^N(g)\big]\bar\nu^N,\qquad \mathcal Z^N_g\text{ a normalisation factor, and }\bar\nu^N\text{ as in}~\eqref{eq_def_mu_product},\label{eq_def_nu_g_dans_methode_entropique}
\end{equation}
For each bias $h\in\s(\epsilon)$ for sufficiently small $\epsilon>0$, the optimal $g=g_h$ arises as the solution of a certain partial differential equation, 
that we call the main equation~\eqref{eq_main_equation}. 
For this $g_h$, the method of~\cite{Jara2018}, 
adapted to this context, 
and a logarithmic-Sobolev inequality yield the bound of Theorem~\ref{theo_entropic_problem}:
\begin{equation}
\exists C,K>0,\forall t\geq 0,\qquad H(f_t\nu^N_{g_h}|\nu^N_{g_h})
\leq 
e^{-Kt}H(f_0\nu^N_{g_h}|\nu^N_{g_h}) + CN^{-1/2}
.
\label{eq_estimee_entropique_generale}
\end{equation}
The exponent $-1/2$ improves on~\eqref{eq_entropy_bounds_JM}. 
It is related to the size of exponential moments of three-point and four-point correlation functions for product measures 
(see Section~\ref{sec_Ising_measure}).  
It thus cannot be improved without adding a correction to $\nu^N_{g_h}$. 
The proof of~\eqref{eq_estimee_entropique_generale} is the main technical result of the article. 
The precise statement of the result is the content of the next two lemmas, a more comprehensive reformulation of Theorem~\ref{theo_entropic_problem}. 
For $h\in\s(\infty)$ 
(defined in~\eqref{eq_def_s_infty}) 
and $f:\Omega_N\rightarrow\R_+$, 
let $\Gamma_h(\sqrt{f})$ be the carré du champ operator associated with the generator $L_h$ biased by $h$, 
with jump rates $c_h$ defined in~\eqref{eq_def_jump_rates_H}:
\begin{align}
\forall\eta\in\Omega_N,\qquad 
\Gamma_h(\sqrt f)(\eta) 
&= 
\frac{1}{4}\sum_{i<N-1}c_h(\eta,i,i+1)\big[\sqrt f(\eta^{i,i+1})-\sqrt f(\eta)\big]^2
\nonumber \\
&\quad
+ \frac{1}{4} \sum_{i\in\{\pm (N-1)\}}c_h(\eta,i)\big[\sqrt f(\eta^{i})-\sqrt f(\eta)\big]^2
.
\label{eq_def_carre_du_champ}
\end{align}
\begin{lemma}[Log-Sobolev inequality, adapted from~\cite{Goncalves2021}]\label{lemm_LSI_sec3}
Let $0<\rho_-\leq \rho_+<1$. 
There is a constant $C_{LS} = C_{LS}(\rho_\pm)$ such that, 
for each $\epsilon,\epsilon'\in(0,1/4)$, 
each $h\in\s(\epsilon)$ and each $g\in g_0+\s(\epsilon')$, 
the following inequality holds. 
For any density $f$ for $\nu^N_g$:
\begin{equation}
\forall N\in\N^*,\qquad
H(f\nu^N_g|\nu^N_g)
\leq 
C_{LS}N^2 \nu^N_g\big(\Gamma_h(\sqrt{f})\big)
.
\end{equation}
\end{lemma}
The $\epsilon_B$ appearing in the next lemma is the same as the one of Theorems~\ref{theo_large_devs}--\ref{theo_entropic_problem}.
\begin{lemma}[Approximation of the invariant measure]\label{lemm_L_star_as_e_plus_carre_du_champ}
Let $0<\rho_-\leq \rho_+<1$. 
Let $h\in\s(\infty)$ and assume that 
the main equation~\eqref{eq_main_equation} has a solution $g_h\in g_0+\s(\infty)$
(the set $\s(\infty)$ is defined in~\eqref{eq_def_test_functions_s_B}). 
For $N\in\N^*$, 
the reference measure $\nu^N_{g_h}$ is defined by~\eqref{eq_def_nu_g_dans_methode_entropique}. 

There is then $\epsilon_B= \epsilon_B(\rho_-)>0$ such that, 
if $\bar\rho':= \frac{\rho_+-\rho_-}{2}\leq \epsilon_B$ and $h\in\s(\epsilon_B)$, 
then $g_h$ is a negative kernel 
(as defined in~\eqref{eq_k_0_g_0_negative_kernels}) 
and the following is true. \\
\noindent$\bullet$ For any $N\in\N^*$, 
there is a function $\e:\Omega_N\rightarrow\R$ such that, 
for any $\nu^N_{g_h}$-density $f$, 
the adjoint $L^*_h$ of $L_h$ in $\mathbb L^2(\nu^N_{g_h})$ satisfies:
\begin{align}
\nu^N_{g_h}\big(fN^2L^*_h{\bf 1}\big) \leq \nu^N_{g_h}(f\e) +\frac{N^2}{2}\nu^N_{g_h}\big(\Gamma_h(\sqrt{f})\big).\label{eq_bound_adjoint_dansl_lemma_L_star}
\end{align}
\noindent$\bullet$ There are constants $\gamma>8C_{LS},C>0$ depending on $h$ and the reservoir densities $\rho_\pm$, 
such that:
\begin{equation}
\forall N\in\N^*,\qquad 
\gamma^{-1}\log\nu^N_{g_h}\Big(\exp\big[\gamma|\e|\big]\Big) 
\leq 
\frac{C}{N^{1/2}}
.
\label{eq_e_controllable_lemm_L_star}
\end{equation}
\noindent$\bullet$ As a consequence, 
for any density $f$ for $\nu^N_{g_h}$,
\begin{align}
\nu^N_{g_h}\big(fL^*_h{\bf 1}\big) - N^2\nu^N_{g_h}\big(\Gamma_h(\sqrt{f})\big)
&\leq 
\frac{H(f\nu^N_{g_h}|\nu^N_{g_h})}{8C_{LS}} + \frac{C}{N^{1/2}}- \frac{N^2}{2}\nu^N_{g_h}\big(\Gamma_h(\sqrt{f})\big)
\nonumber\\
&\leq 
-\frac{H(f\nu^N_{g_h}|\nu^N_{g_h})}{8C_{LS}} + \frac{C}{N^{1/2}}- \frac{N^2}{4}\nu^N_{g_h}\big(\Gamma_h(\sqrt{f})\big)
,
\label{eq_bound_adjoint_lemm_relative_entropy}
\end{align}
and if $f_t\nu^N_{g_h}$ denotes the law of the dynamics at time $t$, then:
\begin{align}
\forall t\geq 0,\forall N\in\N^*, \qquad 
\partial_t H(f_t\nu^N_{g_h}|\nu^N_{g_h}) 
&\leq 
-\frac{H(f_t\nu^N_{g_h}|\nu^N_{g_h})}{8C_{LS}} + \frac{ C}{ N^{1/2}}
.
\label{eq_formule_entropie_lemm_L_star}
\end{align}
\end{lemma}
Let us now prove Theorem~\ref{theo_entropic_problem} using Lemmas~\ref{lemm_LSI_sec3}--\ref{lemm_L_star_as_e_plus_carre_du_champ}. 
Let $g_h$ solve the main equation~\eqref{eq_main_equation}. 
Recalling~\eqref{eq_formule_entropie_lemm_L_star} and applying Gronwall inequality to $t\mapsto H(f_t\nu^N_{g_h}|\nu^N_{g_h})$ yields Theorem~\ref{theo_entropic_problem}:
\begin{equation}\label{eq_application_LSI_sec3_2}
\forall t\geq0,\qquad 
H(f_t\nu^N_{g_h}|\nu^N_{g_h})
\leq  
H(f_0\nu^N_{g_h}|\nu^N_{g_h})e^{-(8C_{LS})^{-1} t} + \frac{8C_{LS} C}{N^{1/2}}\Big(1-e^{-(8C_{LS})^{-1} t}\Big)
.
\end{equation}
Lemma~\ref{lemm_LSI_sec3} is proven in Appendix~\ref{app_LSI}. 
It is a direct adaptation of the proof of~\cite{Goncalves2021}. 
The proof of the key ingredient, 
Lemma~\ref{lemm_L_star_as_e_plus_carre_du_champ}, 
takes up the next four subsections. 
The fact that $g_h$ is a negative kernel if $h\in\s(\epsilon)$ and $\bar\rho'\leq \epsilon$ for small enough $\epsilon>0$ is a consequence of Proposition~\ref{prop_solving_P_triangle}, 
where it is shown that $\|g_h-g_0\|_2$ vanishes with $\epsilon$.

A function $h\in\s(\infty)$ and a negative kernel $g\in g_0+\s(\infty)$ are fixed throughout the rest of the section, 
and we highlight where we need to restrict the size of $h$ and $\bar\rho'$. 
Recall that $h,g$ are symmetric functions, and that their restrictions to $\rhd$ are in $C^3(\bar\rhd)$. 
We do not a priori assume that $g$ solves the main equation~\eqref{eq_main_equation}, 
and explain along the proof where this comes into play.

\subsection{Estimates on $L^*_h{\bf 1}$}
\label{sec_computation_L_star}
To prove Lemma~\ref{lemm_L_star_as_e_plus_carre_du_champ}, 
we need to compute $L^*_h{\bf 1}$. 
The computation of $L^*_h{\bf 1}$ will give rise to many different objects which, 
roughly speaking, 
will either contribute to leading order in $N$, 
or be sub-leading order error terms.
In this section, 
we define precisely how to estimate the size of a function in terms of $N$, 
and formulate criteria to identify which terms are error terms.%
\subsubsection{Size of error terms}
Consider a density $f$ for $\nu^N_g$, 
and a function $X_N:\Omega_N\rightarrow\R$. 
Our main tool to estimate the scaling of $X_N$ with $N$ is the entropy inequality:
\begin{equation}
\forall \gamma>0,\qquad 
\nu^N_{g}(f|X_N|) 
\leq 
\frac{H(f\nu^N_{g}|\nu^N_{g})}{\gamma} + \frac{1}{\gamma}\log \nu^N_{g}\Big(\exp\big[\gamma|X_N|\big]\Big)
.
\label{eq_entropy_ineq_321}
\end{equation}
Informally, we will say that $X_N$ is small if its moment generating function under $\nu^N_g$ vanishes with $N$ for $\gamma$ in a neighbourhood of $0$. 
This is the kind of characterisation of smallness that is used to estimate the size of the function $\e$ in Lemma~\ref{lemm_L_star_as_e_plus_carre_du_champ}. 
In some cases, typically when dealing with the effect of the reservoirs, 
we will encounter an $X_N$ that is not small, 
but can be transformed into some $\tilde X_N$ that is indeed small, 
up to a cost estimated by the carré du champ operator. 
The next definition formalises these considerations 
and examples are given in Lemma~\ref{lemm_size_controllable} below.
\begin{definition}\label{def_controllability}
Let $a_N\in \R^*_+$, $N\in\N^*$. 
A family $X_N:\Omega_N\rightarrow\R,N\in\N^*$ of functions is said to be:
\begin{itemize}
	\item \emph{Controllable with size $a_N$} if there are $\gamma,K$ independent of $N$ such that:
	\begin{equation}
\forall N\in\N^*,\qquad 
\frac{1}{\gamma}\log \nu^N_{g}\Big(\exp\big[\gamma|X_N|\big]\Big) 
\leq 
K a_N
.
\label{eq_controllable_rv}
	\end{equation}
	By convention, if $a_N(p)$ depends on an additional parameter $p$, 
	then the constant $K$ in~\eqref{eq_controllable_rv} will not depend on $p$. 
	By the entropy inequality~\eqref{eq_entropy_ineq_321},~\eqref{eq_controllable_rv} implies, 
	for each density $f$ for $\nu^N_g$:
\begin{equation}
\forall N\in\N^*,\qquad 
\nu^N_g\big(f|X_N|\big)
\leq 
\frac{H(f\nu^N_g|\nu^N_g)}{\gamma} + Ka_N
.
\label{eq_consequence_controllable}
\end{equation}	
	\item \emph{$\Gamma$-controllable with size $a_N$} if one can transform $X_N$, 
	using the carré du champ $\Gamma_h$, 
	into a controllable function with size $a_N$. 
	More precisely: 
	$X_N$ is $\Gamma$-controllable with size $a_N$ if there are controllable functions $\tilde X^{N}_{\pm}$ with size $a_N$ such that, 
	for each $\delta>0$, 
 	each $N\in\N^*$ and each density $f$ for $\nu^N_{g}$:
	\begin{equation}\label{eq_def_Gamma_controllability}
\nu^N_{g}(f(\pm X_N)) 
\leq 
\delta N^2 \nu^N_{g}\big(\Gamma_h(f^{1/2})\big) + \frac{1}{\delta}\nu^N_{g}(f \tilde X^{N}_{\pm})
.
	\end{equation}
	The entropy inequality~\eqref{eq_entropy_ineq_321} then again implies that there are $\gamma>0,K>0$ independent of $\delta,N$ such that:
	\begin{equation}
\nu^N_{g}(f(\pm X_N))
\leq 
\delta N^2 
\nu^N_{g}\big(\Gamma_h(f^{1/2})\big) +\frac{H(f\nu^N_g|\nu^N_g)}{\gamma\delta}+ \frac{K a_N}{\gamma \delta}
.
\label{eq_consequence_Gamma_controllable}
	\end{equation}
	\item An \emph{error term with size $a_N$}, or \emph{error term} for short, 
	if it is either controllable or $\Gamma$-controllable with size $a_N$, and $a_N = o_N(1)$.
\end{itemize}
\end{definition}
\begin{remark}
$(\Gamma$-$)$controllability behaves well with respect to multiplication by a small constant in the following sense. 
Assume that $X_N$ is $(\Gamma$-$)$controllable with size $a_N$ 
and let $b_N\in[0,1]$ ($N\in\N^*$) satisfy $b_N=o_N(1)$. 
Then $b_N X_N$ is $(\Gamma$-$)$controllable with size $b_N a_N$ by Jensen inequality.
\demo
\end{remark}
To illustrate the notion of controllability, the following proposition, 
proven in Appendix~\ref{prop_Boltzmann_Gibbs}, 
states its consequence on the dynamical behaviour of observables. 
\begin{proposition}\label{prop_Boltzmann_gibbs_sec3}
Let $h\in\s(\epsilon_B)$ with $\epsilon_B$ given by Lemma~\ref{lemm_L_star_as_e_plus_carre_du_champ}, 
and let $g_h$ be the associated solution of the main equation~\eqref{eq_main_equation} as in Lemma~\ref{lemm_L_star_as_e_plus_carre_du_champ}. 
Let $E^N:\Omega_N\rightarrow\R$ be an error term with size $a_N =o_N(1)$, 
and let $F^N$ be ($\Gamma$)-controllable with size $1$. 
There are then $\gamma,C$ and $\gamma',C'>0$ independent of $N,T$ such that:
\begin{align}
\forall T>0,\qquad 
\frac{1}{T}\log\E^{\nu^{N}_{g_h}}\Big[\exp\Big|\gamma\int_0^T  E^N(\eta_t)dt\Big|\Big] 
&\leq 
C a_N,\qquad\nonumber\\ 
\sup_N\frac{1}{T}\log \E^{\nu^{N}_{g_h}}\Big[\exp\Big|\gamma'\int_0^T  F^N(\eta_t)dt\Big|\Big] 
&\leq 
C'
.
\end{align}
\end{proposition}
Let $\phi:\Lambda_N\rightarrow\R$. 
To determine whether a field is an error term or not, one must keep in mind the following heuristics: the measures $\nu^N_g$ are discrete Gaussian measures, in the sense that the fluctuation field $Y^N(\phi)$, which reads:
\begin{equation}
Y^N(\phi) 
:= 
\frac{1}{N^{1/2}}\sum_{i\in\Lambda_N}\bar\eta_i \phi(i)
,
\label{eq_def_fluct}
\end{equation}
is close to a Gaussian random variable when $N$ is large, provided $\|\phi\|_{\infty}<\infty$. 
In particular, $\lambda \mapsto \nu^N_g(\exp[\lambda Y^N(\phi)])$ is bounded uniformly in $N$ in a neighbourhood $[0,\gamma(\phi))$ of $0$ for some $\gamma(\phi)>0$. 
By~\eqref{eq_controllable_rv}, 
this means that $Y^N(\phi)$ is controllable with size $1$ 
(or size $C(\phi)$ for some constant $C(\phi)>0$ if we want to keep track of the dependence on $\phi$).
In analogy with Gaussian random variables, one can prove that $Y^N(\phi)^2$ is controllable with size $1$, but $Y^N(\phi)^n$ for $n\geq 3$ is not. 
Similarly, the quantity: 
\begin{equation}
Z^N(\phi) 
:= 
\frac{1}{N^{1/2}}\sum_{i<N-1}\bar\eta_i\bar\eta_{i+1}\phi(i)\label{eq_def_Z_n_fluct}
\end{equation}
should not have worse concentration properties than $Y^N(\phi)$. 
Contrary to genuine Gaussian random variables, however, 
$Y^N(\phi)$ and $Z^N(\phi)$ are bounded, 
by $C\|\phi\|_\infty N^{1/2}$ for some $C>0$. 
As a result, it is always possible to find $a_N$ small enough such that $a_N(Y^N(\phi))^n$ is controllable with size $1$. This discussion is summarised in the next lemma.
\begin{lemma}\label{lemm_size_controllable}
For $n\in\N^*$, 
let $\phi_n :\Lambda_N^n\rightarrow\R$ satisfy $\sup_{N}\|\phi_n\|_{\infty}<\infty$. 
Define, for $N\in\N^*$ and either $J = \{0\}$ or $J = \{0,1\}$, the functions (an empty product is by convention equal to $1$):
\begin{align}\label{eq_def_rv_lemm_size_controllable}
\forall \eta\in\Omega_N,\qquad 
X^{\phi_{n}}_{n,J}(\eta) 
&= 
\sum_{i_0<N-|J|}\sum_{i_1,...,i_{n-1}\in\Lambda_N}\phi_{n}(i_0,...,i_{n-1})\Big(\prod_{j\in J}\bar\eta_{i_0+j}\Big)\Big(\prod_{a=1}^{n-1}\bar\eta_{i_a}\Big) ,\nonumber\\
U_0^\epsilon(\eta) 
&= 
\bar\eta_{\epsilon(N-1)},\quad 
U_1^{\epsilon}(\eta) 
= 
\bar\eta_{\epsilon(N-1)}\sum_{i\neq \epsilon(N-1)}\bar\eta_{i}\phi_1(i),\qquad\epsilon\in\{-,+\}
.
\end{align}
Then:
\begin{itemize}
	\item the function $U_0^\pm$ is $\Gamma$-controllable with size $N^{-1}$, 
	and $U_1^{\pm}$ is $\Gamma$-controllable with size $\|\phi_1\|_{\infty}N^{-1}$. 
	Moreover, $N^{-1/2}X^{\phi_1}_{1,J}$ is controllable with size $\|\phi_1\|_{\infty}^2$.  
	\item For $n\geq 2$, 
	there are constants $\gamma_n>0$ that depend only on $n$, 
	but not on $\phi_n$, 
	such that:
\begin{equation}\label{eq_true_controllability}
\log\nu^N_g\Big(\exp\Big[\frac{\gamma_n}{\|\phi_n\|_\infty} N^{-(n-1)}X^{\phi_n}_{n,J}\Big]\Big) 
\leq 
\frac{3}{N^{\frac{n-2}{2}}}
.
	\end{equation}	
	This implies (but is stronger) that 
	the function $N^{-(n-1)}X^{\phi_n}_{n,J}$ is controllable with size $\|\phi_n\|_{\infty} N^{-(n-2)/2}$. 
	In addition, 
	if $\phi_n(i_0,...,i_{n-1})$ vanishes whenever an index appears twice in the collection $i_0+J, i_1,...,i_{n-1}$, 
	then one can replace $\|\phi_n\|_\infty$ by the weaker norm $\|\phi_n\|_{2,N}$ in~\eqref{eq_true_controllability}, 
	where:
	\begin{equation}\label{eq_def_norm_2_N}
	\|\phi_n\|_{2,N} := \Big(\frac{1}{N^{n}}\sum_{i_0,...,i_{n-1}}\phi_n(i_0,...,i_{n-1})^2\Big)^{1/2}
	.
	\end{equation}
	\item For $n=2$ and $|J|>1$, say $J=\{0,1\}$, 
	the previous estimate can be improved: $N^{-1}X^{\phi_2}_{2,\{0,1\}}$ is $\Gamma$-controllable with size $\|\phi_n\|_{\infty}N^{-1/2}$ 
	(or $\|\phi_n\|_{2,N}N^{-1/2}$ if $\phi_{2}(i,i) = 0 = \phi_2(i,i+1)$ for each $i$).
\end{itemize}
\end{lemma}
\begin{remark}
To help clarify the definition of $X^{\phi_n}_{n,J}$, 
take $J=\{0,1\}$ and, for $n\in\N^*$, 
let $\phi_n(i_0,...,i_{n-1}) = \prod_{\ell=0}^{n-1}\phi_1(i_\ell)$. 
Then:
\begin{equation}
X^{\phi_{n}}_{n,J} = N^{n/2}(Y^N(\phi_1))^{n-1}Z^N(\phi_1),
\end{equation}
with $Y^N(\phi)$, $Z^N(\phi)$ defined in~\eqref{eq_def_fluct}-\eqref{eq_def_Z_n_fluct}. 
Moreover, 
for any $\phi_2:\Lambda^2_N\rightarrow\R$ with $\phi_2(i,i) = 0$ for $i\in\Lambda_N$, 
then the correlation field \eqref{eq_def_Pi} can be recovered:
$X^{\phi_2}_{2,\{0\}} = 4N\Pi^N(\phi_2)$.
\demo
\end{remark}
Lemma~\ref{lemm_size_controllable} is proven in Appendix~\ref{sec_Ising_measure} for the variables $X^{\phi_n}_{n,J}$. The statement for $U_0^\pm,U_1^\pm$ is proven in Appendix~\ref{sec_boundary_correl_dirichlet}, 
and the last item corresponds to Proposition~\ref{prop_error_terms}.
\subsubsection{Bounding the entropy}
After classifying each term arising in the computation of $L^*_h{\bf 1}$ according to the categories in Definition~\ref{def_controllability}, 
we will obtain the estimate~\eqref{eq_bound_adjoint_dansl_lemma_L_star}: for a good choice $g_h$ depending on $h$ and any density $f$ for $\nu^N_{g_h}$,
\begin{equation}
\nu^N_{g_h}\big(fL^*_h{\bf 1}\big) 
\leq 
\nu^N_{g_h}\big(f\e\big) + \frac{N^2}{2}\nu^N_{g_h}\big(\Gamma_h(\sqrt{f})\big)
.
\end{equation}
Above, 
$\e$ will be a controllable error term with size $N^{-1/2}$: 
there is $\gamma,C>0$ depending on $\rho_\pm,h$ such that
\begin{equation}\label{eq_bound_adjoint_to_prove}
\forall N\in\N^*,\qquad 
\log\nu^N_{g_h}\Big(\exp\big[\gamma|\e|\big]\Big) 
\leq 
\frac{C}{N^{1/2}}
.
\end{equation}
To turn this estimate on the adjoint into the entropy estimate of Theorem~\ref{theo_entropic_problem}, 
there is one more step, 
which involves the log-Sobolev inequality of Lemma~\ref{lemm_LSI_sec3}, 
see~\eqref{eq_application_LSI_sec3_2}.
To apply the log-Sobolev inequality, 
we need to ensure that the constant $\gamma$ in~\eqref{eq_bound_adjoint_to_prove} is sufficiently large compared to the log-Sobolev constant $C_{LS}$, 
in the sense that we want to have:
\begin{align}
&\nu^N_{g_h}\big(fL^*_h{\bf 1}\big) - N^2\nu^N_{g_h}\big(\Gamma_h(f^{1/2})\big)
\nonumber\\
&\qquad\leq 
\frac{H(f\nu^N_{g_h}|\nu^N_{g_h})}{\gamma} + \frac{1}{\gamma}\log\nu^N_{g_h}\Big(\exp\big[\gamma|\e|\big]\Big) 
- \frac{N^2}{2}\nu^N_{g_h}\big(\Gamma_h(f^{1/2})\big)
\nonumber\\
&\qquad\leq
- \frac{H(f\nu^N_{g_h}|\nu^N_{g_h})}{8C_{LS}} + \frac{C}{\gamma N^{1/2}} 
-\frac{N^2}{4}\nu^N_{g_h}\big(\Gamma_h(f^{1/2})\big)
,
\end{align}
where $(8C_{LS})^{-1}$ in the last line could be replaced by any smaller positive number. 
The last bound requires $\gamma>8C_{LS}$, 
and is the reason for the introduction of the following terminology.
\begin{definition}\label{def_LS_type}
\begin{itemize}
	\item A $(\Gamma$-$)$controllable error term $X_N$ of size $a_N$ is said to be of \emph{vanishing type} if 
there is a sequence $(\gamma_N)_N$ of positive numbers such that $\lim_N\gamma_N=\infty$, 
and $\gamma_N X_N$ is $(\Gamma$-$)$controllable with size $1$.
	\item A function $X_N$ (error term or not) is said to be of \emph{large type} if, for some $\gamma>1$ independent of $N$, 
	$\gamma X_N$ has unbounded exponential moment as $N\rightarrow\infty$.
	\item A controllable function is said to be of \emph{LS type} if one can take $\gamma>2^{10}C_{LS}$ in the definition~\eqref{eq_controllable_rv} of controllability 
	(this is in particular always true for $N$ large in the case of error terms of vanishing type). 
	A $\Gamma$-controllable function is said to be of LS type if the controllable functions $Y^N_\pm$ associated via~\eqref{eq_def_Gamma_controllability} are of LS type.
\end{itemize}
\end{definition}
\begin{remark}
The constant $2^{10}$ is not at all optimal 
(any large enough constant would also work), 
but is sometimes convenient later on.
\demo
\end{remark}
Let us clarify this notion by classifying the functions appearing in Lemma~\ref{lemm_size_controllable}.
\begin{lemma}\label{lemm_type}
For $n\in\N$, $J=\{0\}$ or $\{0,1\}$ and 
$\phi_n:\Lambda_N^n\rightarrow\R$ with $\sup_N\|\phi_n\|_\infty<\infty$, 
recall the definitions of $X^{\phi_n}_{n,J},U^\pm_n$ in Lemma~\ref{lemm_size_controllable}.
\begin{itemize}
	\item $N^{1/2}U_0^\pm, N^{1/2}U_1^{\pm}$ are $\Gamma$-controllable with size $1$, 
	so both $U_0^\pm$ and $U^1_\pm$ are error terms of vanishing type. 
	Similarly, 
	for any sequence $\epsilon_N>0$ $(N\in\N^*)$ with $\epsilon_N=o_N(1)$, 
	$\epsilon_NN^{-1/2}X^{\phi_1}_{1,J}$ and $\epsilon_N N^{-(n-1)}X^{\phi_n}_{n,J}$ for $n\geq 2$ are error terms of vanishing type.
	\item For $n\geq 2$,  
	$\gamma N^{-(n-1)}X^{\phi_n}_{n,J}$ has exponential moment bounded with $N$ only if $\gamma>0$ is small enough. 
Thus $N^{-(n-1)}X^{\phi_n}_{n,J}$ is of large type.
	\item For $n\geq 2$, 
	there is a numerical constant $\zeta_n>0$ such that $\frac{\zeta_n}{\|\phi_n\|_\infty} N^{-(n-1)}X^{\phi_n}_{n,J}$ is of LS type. \\
	If in addition $\phi_n(i_0,...,i_{n-1})$ vanishes whenever the same index appears twice in the collection $i_0+J,i_1,...,i_{n-1}$, 
	then the same statement is true replacing $\|\phi_n\|_\infty$ with $\|\phi_n\|_{2,N}$ (defined in~\eqref{eq_def_norm_2_N}).
	\item If a function $X_N$ is of LS type and $\epsilon\in(0,1)$, 
	then $\epsilon X_N$ is of LS type. 
	Moreover, $(1+\alpha)X_N$ is of LS type for sufficiently small $\alpha\in(0,1)$ independent of $N$.
\end{itemize}
\end{lemma}
Lemma~\ref{lemm_type} is obtained as a consequence of the proof of Lemma~\ref{lemm_size_controllable}, 
carried out in Corollary~\ref{coro_estimate_W_section_Ising}.

To see why the claim of Lemma~\ref{lemm_type} is reasonable,
consider again $Y^N(\phi_1) = N^{-1/2}\sum_{i}\bar\eta_i\phi_1(i)$   
($\phi_1$ is bounded) as an approximately Gaussian random variable when $N$ is large. 
Similarly, one should see $N^{-(n-1)}X^{\phi_n}_{n,J}$ as approximately $N^{-(n-2)/2}$ times a product of $n$ Gaussian random variable ($n\geq 2$) for bounded $\phi_n$. 
For $n\geq 3$, 
the prefactor $N^{-(n-2)/2}$ vanishes, 
thus $N^{-(n-1)}X^{\phi_n}_{n,J}$ should become an error term. 
Error terms, however, 
are defined in terms of smallness of exponential moments, 
and the moment generating function $M_n(\gamma)$ of a product of $n$ Gaussian random variables is unbounded for $\gamma$ large enough ($n=2$) or for any $\gamma>0$ if $n\geq 3$. 
For $n\geq 3$, 
the vanishing prefactor $N^{-(n-2)/2}$ when considering $N^{-(n-1)}X^{\phi_n}_{n,J}$ as a product of $n$ Gaussians is precisely right to ensure that its exponential moments are bounded with $N$ in a neighbourhood of $0$. \\

For an error term of LS type, 
a straightforward adaptation of the proof of Proposition~\ref{prop_Boltzmann_gibbs_sec3} 
(see Section~\ref{sec_proof_boltzmann_gibbs})
yields the following result.
\begin{corollary}\label{coro_Boltzmann_gibbs_sec3}
Let $E^N$ be a $(\Gamma$-$)$controllable function with size $a_N=O_N(1)$, 
and assume that $E^N$ is of LS type. 
Then one can take $\gamma = 1$ in Proposition~\ref{prop_Boltzmann_gibbs_sec3}, 
i.e. there is $C=C(\rho_\pm)>0$ such that:
\begin{equation}
\forall N\in\N^*,\forall T>0\qquad 
\frac{1}{T}\log\E^{\nu^{N}_{g_h}}\Big[\exp\Big|\int_0^T  E^N(\eta_t)dt\Big|\Big] 
\leq 
C a_N
.
\end{equation}
In particular, 
this applies to $E^N = \zeta_2\|\phi\|_{2,N}^{-1}\Pi^N(\phi)$ with $a_N =1$ for any bounded $\phi:\squaredash\rightarrow\R$, 
with $\zeta_2$ the constant in item 3 of Lemma~\ref{lemm_type}.
\end{corollary}
Equipped with these notations, 
we now turn to computing $L^*_h{\bf 1}$, 
in the next two subsections. 
We split $L^*_{h}$ into $L^*_{h} = L^*_{h,0} + L^*_{h,\pm}$, 
the adjoint dynamics respectively in the bulk and at the boundaries, 
and study each contribution separately. 
\subsection{Adjoint at the boundary}\label{sec_adjoint_boundary}
In this section, we compute $L^*_{h,\pm}{\bf 1}$, the part of the adjoint $L^*_h{\bf 1}$ of $L_h$ 
with respect to $\nu^{N}_g$
corresponding to the dynamics at the boundary. 
By~\eqref{eq_def_L_star_H_general_mesure}, 
it reads:
\begin{align}
&N^2L^*_{h,\pm}{\bf 1}(\eta)
=
\frac{N^2}{2}\sum_{i\in\{\pm (N-1)\}}\Big[ c_h(\eta^i,i)\frac{\nu^N_{g}(\eta^{i})}{\nu^N_{g}(\eta)}-c_h(\eta,i)\Big]
\nonumber\\
&\quad = 
\frac{N^2}{2}\sum_{i\in\{\pm (N-1)\}}\bigg[c(\eta^i,i)\Big(\frac{\bar\rho_i}{1-\bar\rho_i}\Big)^{1-2\eta_i} \exp\Big[\frac{(1-2\eta_i)}{2N}\sum_{j\neq i}\bar\eta_j (2g-h)_{i,j}\Big]
\nonumber\\
&\hspace{5cm}
-c(\eta,i)\exp\Big[\frac{(1-2\eta_i)}{2N}\sum_{j\neq i}\bar\eta_j h_{i,j}\Big]\bigg]
.
\label{eq_starting_point_L_star_boundary}
\end{align}
The jump rates $c(\eta,i)$, $i\in\{\pm(N-1)\}$ are defined in~\eqref{eq_def_jump_rates}. To compute~\eqref{eq_starting_point_L_star_boundary}, recall that both $h,g$ satisfy $h(\pm 1,\cdot) = 0 = g(\pm 1,\cdot)$ by hypothesis. It follows that the arguments of the exponentials in~\eqref{eq_starting_point_L_star_boundary} are bounded by $O(N^{-1})$. Moreover, introduce for $i\in\Lambda_N = \{ -N +1, \dots, N-1\}$ and at the boundaries:
\begin{equation}
\lambda_i :=\log\Big(\frac{\bar\rho_i}{1-\bar\rho_i}\Big)
\quad \text{as well as} 
\quad \lambda_{\pm N} = \log\Big(\frac{\rho_\pm}{1-\rho_\pm}\Big).\label{eq_def_lambda}
\end{equation}
With this notation and by reversibility, one has for $i\in\{\pm(N-1)\}$:
\begin{equation}
c(\eta^i,i) 
= 
c(\eta,i) \Big(\frac{\rho_{\text{sign(i)}}}{1-\rho_{\text{sign(i)}}}\Big)^{2\eta_i -1}
= 
c(\eta,i)\exp \Big[(2\eta_i -1) \lambda_{\text{sign(i)} N} \Big]
.
\end{equation}
As a result, with the notation $\partial^N\lambda_\ell = N ( \lambda_{\ell+1} - \lambda_\ell)$, we get :
\begin{align}
&c(\eta^i,i)\Big(\frac{\bar\rho_i}{1-\bar\rho_i}\Big)^{1-2\eta_i} 
\nonumber\\
&\qquad = 
c(\eta,i)\exp\Big[\frac{(1-2\eta_i)}{N}\big({\bf 1}_{i= -(N-1)}\partial^N\lambda_{-N} - {\bf 1}_{i =N-1}\partial^N\lambda_{N-1}\big)\Big]
.
\label{eq_c_eta_i_i_fonction_of_c_eta_i}
\end{align}
Write for short the boundary term as:
\begin{equation}
\cB_i := {\bf 1}_{i= -(N-1)}\partial^N\lambda_{-N} - {\bf 1}_{i =N-1}\partial^N\lambda_{N-1}.
\end{equation}
Using~\eqref{eq_c_eta_i_i_fonction_of_c_eta_i} and the existence of $C_{h,g,\bar\rho} = C(\|h\|_\infty,\|g\|_\infty,\|\bar\rho\|_\infty)>0$ such that 
\begin{equation}
|e^x-1-x^2/2| \leq C_{h,g,\bar\rho}|x|^3
\quad \text{for} \quad |x|\leq 2(\|h\|_{\infty}+2\|g\|_\infty + \|\lambda\|_\infty),
\end{equation}
Equation~\eqref{eq_starting_point_L_star_boundary} can be expanded with an error term:
\begin{align}
&\delta^{N,1}_{\pm}(\eta) := 
N^2L^*_{h,\pm}{\bf 1}(\eta)
 - \frac{N}{2}\sum_{i\in\{\pm (N-1)\}}c(\eta,i)(1-2\eta_i)\Big[\frac{1}{N}\sum_{j\neq i}\bar\eta_j N(g-h)_{i,j} +  \cB _i \Big]
\label{eq_expansion_L_boundary}\\
&\hspace{1.5cm}
- \frac{1}{4}\sum_{i\in\{\pm(N-1)\}}c(\eta,i) \bigg[
\Big(\frac{1}{2 N}\sum_{j\neq i}\bar\eta_j N(2g-h)_{i,j} +  \cB_i \Big)^2
- 
\Big(\frac{1}{ 2 N}\sum_{j\neq i}\bar\eta_j Nh_{i,j}\Big)^2
\bigg] 
,\nonumber
\end{align}
and $\delta^{N,1}_{\pm}$ satisfies $|\delta^{N,1}_{\pm}|\leq C_{h,g,\bar\rho}/N$. 
To compute the two terms in the left-hand side of~\eqref{eq_expansion_L_boundary}, 
let $i\in\{\pm(N-1)\}$ and let us first rewrite the jump rate in terms of $\bar\eta_i$ in two different ways. One has:
\begin{align}
c(\eta,i) := (1-\eta_i)\rho_{\text{sign}(i)} + (1-\rho_{\text{sign}(i)})\eta_i = a_{\text{sign}(i)} + \bar\eta_i(1-2\rho_{\text{sign}(i)}),\label{eq_computation_c_eta_i_from_def}
\end{align}
with: 
\begin{equation}
a_{\text{sign}(i)} := (1-\bar\rho_i)\rho_{\text{sign}(i)} + (1-\rho_{\text{sign}(i)})\bar\rho_i.\label{eq_def_a_pm}
\end{equation}
Moreover, it also holds that:
\begin{align}
(1-2\eta_i)c(\eta,i) &= -(\eta_i-\rho_{\text{sign}(i)}) = -\bar\eta_i+ (\rho_{\text{sign}(i)}-\bar\rho_i)\nonumber \\
&=- \bar\eta_i -\frac{1}{N} \big[{\bf 1}_{i= -(N-1)}\partial^N\bar\rho_{-N} - {\bf 1}_{i =N-1}\partial^N\bar\rho_{N-1}\big].\label{eq_computation_c_eta_i}
\end{align}
Using~\eqref{eq_computation_c_eta_i} in the first sum in~\eqref{eq_expansion_L_boundary}, one finds:
\begin{align}
\frac{N}{2}&\sum_{i\in\{\pm (N-1)\}}c(\eta,i)(1-2\eta_i)\Big[\frac{1}{N}\sum_{j\neq i}\bar\eta_j N(g-h)_{i,j} 
\nonumber\\
&\hspace{6cm}+ \big({\bf 1}_{i= -(N-1)}\partial^N\lambda_{-N} - {\bf 1}_{i =N-1}\partial^N\lambda_{N-1}\big)\Big] 
\nonumber\\
& = 
-\frac{N\partial^N\lambda_{-N}}{2}\Big[\bar\eta_{-(N-1)} + \frac{\partial^N\bar\rho_{-N}}{N}\Big] + \frac{N\partial^N\lambda_{N-1}}{2}\Big[\bar\eta_{N-1} -\frac{\partial^N\bar\rho_{N-1}}{N}\Big]+ \delta^{N,2}_\pm(\eta)
,
\label{eq_line_1_L_star_boundary}
\end{align}
where $\delta^{N,2}_\pm$ is an error term that reads:
\begin{align}
&\delta^{N,2}_\pm(\eta) 
\\
&\quad= 
-\frac{1}{2}\sum_{i\in\{\pm (N-1)\}}\Big(\bar\eta_i + \frac{1}{N}\big[{\bf 1}_{i= -(N-1)}\partial^N\bar\rho_{-N} - {\bf 1}_{i =N-1}\partial^N\bar\rho_{N-1}\big]\Big)\sum_{j\neq i}\bar\eta_j N(g-h)_{i,j}
.
\nonumber
\end{align}
Since $N(g-h)_{i,j}$ is of order 1 for $i$ close to the boundary,
the term involving $\bar\eta_i$ above is of the same form as $U_1^\pm$ in Lemma~\ref{lemm_size_controllable}, 
and recall that $U_1^\pm$ is $\Gamma$-controllable with size $N^{-1}$ and of vanishing type (as defined in Lemma~\ref{lemm_type}). 
The other term is of the form $N^{-1/2}Y^N(\phi)$ for a bounded $\phi$ (recall~\eqref{eq_def_fluct}), 
and $Y^N(\phi)$ is controllable with size $1$, 
thus $N^{-1/2}Y^N(\phi)$ is also of vanishing type. 
It follows that $\delta^{N,2}_\pm$ is $\Gamma$-controllable with size $N^{-1}$ and of vanishing type.\\
Consider now the second sum in~\eqref{eq_expansion_L_boundary}. Using~\eqref{eq_computation_c_eta_i_from_def} and recalling the definition~\eqref{eq_def_a_pm} of $a_\pm$, it reads:
\begin{align}
\frac{1}{4}&\sum_{i\in\{\pm(N-1)\}}c(\eta,i)\bigg[ 
\Big(\frac{1}{2 N}\sum_{j\neq i}\bar\eta_j N(2g-h)_{i,j} 
\nonumber\\
&\hspace{2.8cm}+ \big({\bf 1}_{i= -(N-1)}\partial^N\lambda_{-N} - {\bf 1}_{i =N-1}\partial^N\lambda_{N-1}\big)\Big)^2 
- \Big(\frac{1}{2N}\sum_{j\neq i}\bar\eta_j Nh_{i,j}\Big)^2\bigg] 
\nonumber\\
&\hspace{2.5cm}= 
\frac{a_-(\partial^N\lambda_{-N})^2}{4} + \frac{a_+(\partial^N\lambda_{N-1})^2}{4} + \delta^{N,3}_{\pm}(\eta)
,
\label{eq_line_2_L_star_boundary}
\end{align}
where $\delta^{N,3}_\pm$ is an error term that contains all other contributions:
\begin{align}
&\delta^{N,3}_\pm(\eta) 
= 
\frac{1}{4}\sum_{i\in\{\pm(N-1)\}}c(\eta,i)\bigg[
\Big(\frac{1}{2N}\sum_{j\neq i}\bar\eta_j N(2g-h)_{i,j}\Big)^2
- \Big(\frac{1}{2N}\sum_{j\neq i}\bar\eta_j Nh_{i,j}\Big)^2\bigg]
\nonumber\\
&\quad + 
\frac{1}{4  N}\sum_{i\in\{\pm(N-1)\}}c(\eta,i)\Big[\big({\bf 1}_{i= -(N-1)}\partial^N\lambda_{-N} - {\bf 1}_{i =N-1}\partial^N\lambda_{N-1}\big)\sum_{j\neq i}\bar\eta_j N(2g-h)_{i,j}\Big]
\nonumber\\
&\quad + 
\bar\eta_{-(N-1)}(1-2\rho_-)\frac{(\partial^N\lambda_{-N})^2}{4}
+\bar\eta_{N-1}(1-2\rho_+)\frac{(\partial^N\lambda_{N-1})^2}{4}
.
\end{align}
With the notations of Lemma~\ref{lemm_size_controllable} and bounding $c(\eta,\cdot)$ by $1$, the first line of $\delta^{N,3}_\pm$ is bounded by a term of the form $N^{-2}X^{\phi_2}_{2,\{0\}}$, 
thus $N^{-1}$ times a quantity controllable with size $1$. 
Similarly, the second line is bounded by a term of the form $N^{-1/2}|Y^N(\phi)|$ (recall that $Y^N$ is defined in~\eqref{eq_def_fluct}) for bounded $\phi$, 
and $|Y^N(\phi)|$ is controllable with size $1$. 
Finally, the third line is of the form $U_0^\pm$, 
and $U_0^{\pm}$ is $\Gamma$-controllable with size $N^{-1}$ and of vanishing type. 
The quantity $\delta^{N,3}_\pm$ is therefore $\Gamma$-controllable with size $N^{-1/2}$ and of vanishing type.\\

Putting together~\eqref{eq_line_1_L_star_boundary} and~\eqref{eq_line_2_L_star_boundary}, we have obtained the following expression of the adjoint at the boundary:
\begin{align}
N^2L^*_{h,\pm}{\bf 1}(\eta) &= -\frac{N\partial^N\lambda_{-N}}{2}\Big[\bar\eta_{-(N-1)} + \frac{\partial^N\bar\rho_{-N}}{N}\Big] + \frac{N\partial^N\lambda_{N-1}}{2}\Big[\bar\eta_{N-1} -\frac{\partial^N\bar\rho_{N-1}}{N}\Big] \nonumber\\
&\quad + \frac{a_-\big(\partial^N\lambda_{-N}\big)^2}{4} + \frac{a_+\big(\partial^N\lambda_{N-1}\big)^2}{4} + \sum_{q=1}^3\delta^{N,q}_\pm(\eta).
\end{align}
It remains to notice that the constant terms in the last equation compensate each other to obtain the final expression for $N^2L^*_{h,\pm}{\bf 1}(\eta)$. 
Indeed, for each $i\in\{\pm(N-1)\}$, a Taylor expansion yields:
\begin{equation}
a_{\text{sign}(i)} = 2\bar\sigma_i +O(N^{-1}),\qquad \partial^N\lambda_i = \frac{\partial^N\bar\rho_i}{\bar\sigma_i} +O(N^{-1}).
\end{equation}
It follows that there is a configuration-independent error term $\delta^{N,4}$, 
with $\delta^{N,4} = O(N^{-1})$, 
such that:
\begin{align}
N^2L^*_{h,\pm}{\bf 1}(\eta) = -\frac{N\partial^N\lambda_{-N}}{2}\bar\eta_{-(N-1)}+ \frac{N\partial^N\lambda_{N-1}}{2}\bar\eta_{N-1} + \delta^{N}_{\pm}(\eta),\qquad \delta^N_{\pm} := \sum_{q=1}^4\delta^{N,q}_\pm(\eta).\label{eq_final_adjoint_boundary}
\end{align}
The quantities $\delta^N_-,\delta^N_+$ are, by definition, 
$\Gamma$-controllable with size $N^{-1}$ and of vanishing type.
\subsection{Adjoint in the bulk}\label{sec_adjoint_bulk}
We now compute $L^*_{h,0}{\bf 1}$. For each $i<N-1$, define $B^h_i,C^h_i,D^h_i$ as follows: 
\begin{equation}
B^h_i(\eta) = \frac{1}{2N}\sum_{j\notin\{i,i+1\}}\bar\eta_j\partial^N_1 h_{i,j},
\qquad D^h_i(\eta) = \frac{\partial^N\bar\rho_{i} \, h_{i,i+1}}{2N},
\qquad C^h_i = B^h_i + D^h_i,\label{eq_def_B_x_D_x}
\end{equation}
where, for $u : \Lambda_N^2\rightarrow\R$:
\begin{equation}
\partial^N_1 u(i,j) = N[u(i+1,j) - u(i,j)],
\qquad i<N-1,j\in\Lambda_N
.
\end{equation}
With these definitions,
\begin{equation}
\forall i<N-1,\qquad \Pi^N(h)(\eta^{i,i+1})-\Pi^N(h)(\eta) = -\frac{(\eta_{i+1}-\eta_i)}{N}C^h_i(\eta).
\end{equation}
Define similarly $C^g_\cdot$, and notice that, since $h,g$ are regular:
\begin{equation}
\sup_{N\in\N^*}\sup_{\eta\in\Omega_N}\sup_{i<N-1} \big( |C^h_i(\eta)| + |C^g_i(\eta)|\big)<\infty.\label{eq_C_g_C_h_bounded}
\end{equation}
By~\eqref{eq_def_L_star_H_general_mesure}, the adjoint $L^*_{h,0}$ in the bulk reads, by definition:
\begin{equation}
N^2L^*_{h,0}{\bf 1}(\eta) = \frac{N^2}{2}\sum_{i<N-1}\Big[ c_h(\eta^{i,i+1},i,i+1)\frac{\nu^N_{g}(\eta^{i,i+1})}{\nu^N_{g}(\eta)}-c_h(\eta,i,i+1)\Big].
\end{equation}
With the above notations, this becomes:
\begin{align}
N^2L^*_{h,0}{\bf 1}(\eta) &= \frac{N^2}{2}\sum_{i<N-1} c(\eta,i,i+1)\bigg[\exp\Big[\frac{(\eta_{i+1}-\eta_i)}{N}(C^h_i-2C_i^{g})\Big] \Big[\frac{\bar\rho_i(1-\bar\rho_{i+1})}{\bar\rho_{i+1}(1-\bar\rho_i)}\Big]^{\eta_{i+1}-\eta_i} \nonumber\\
&\hspace{5cm}-\exp\Big[-\frac{(\eta_{i+1}-\eta_i)}{N}C^h_i\Big]\bigg].\label{eq_adjoint_bulk_init}
\end{align}
To compute~\eqref{eq_adjoint_bulk_init}, recall the definition of $\lambda$:
\begin{equation}
\forall i\in\Lambda_N,\qquad \lambda_i := \log\Big(\frac{\bar\rho_i}{1-\bar\rho_i}\Big).
\end{equation}
Notice that $C^h_\cdot - 2C^g_\cdot = C^{h-2g}_\cdot$. Moreover, $c(\eta,i,i+1) = (\eta_{i+1}-\eta_i)^2$ for each $i<N-1$. With these notations,~\eqref{eq_adjoint_bulk_init} reads:
\begin{align}
N^2L^*_{h,0}{\bf 1}(\eta) &= \frac{N^2}{2}\sum_{i<N-1}(\eta_{i+1}-\eta_i)^2\bigg[\exp\Big[\frac{(\eta_{i+1}-\eta_i)}{N}(C^{h-2g}_i-\partial^N\lambda_i)\Big] \nonumber\\
&\hspace{5cm}-\exp\Big[-\frac{(\eta_{i+1}-\eta_i)}{N}C^h_i\Big]\bigg].\label{eq_adjoint_bulk_init_1}
\end{align}
To compute~\eqref{eq_adjoint_bulk_init_1}, we expand the above exponentials. Write $(N^2L^*_{h,0}{\bf 1})_{\text{order }p}$ for the term of order $p\in\N$. From the existence of $C_{h,g} = C(\|h\|_\infty,\|g\|_\infty)>0$ such that $|e^{x}-1-x-x^2/2-x^3/6|\leq C_{h,g}x^4$ when $|x|\leq 2(\|h\|_\infty + 2\|g\|_\infty)$, one has $|\delta^{N}_{0,\,\text{order}\geq 4}(\eta)|\leq 2C_{h,g}/N$, with $\delta^{N}_{0,\,\text{order}\geq 4}(\eta)$ given by:
\begin{align}
N^2L^*_{h,0}{\bf 1}(\eta)&=  \frac{N}{2}\sum_{i<N-1}(\eta_{i+1}-\eta_i)\Big[2C^{h-g}_i - \partial^N\lambda_i\Big]\label{eq_first_line_L_star_bulk}\\
&\quad +  \frac{1}{4}\sum_{i<N-1}(\eta_{i+1}-\eta_i)^2\Big[\big(C^{h-2g}_i-\partial^N\lambda_i\big)^2-\big(C^h_i\big)^2\Big] \label{eq_second_line_L_star_bulk}\\
&\quad + \frac{1}{12N}\sum_{i<N-1}(\eta_{i+1}-\eta_i)\Big[\big[C^{h-2g}_i-\partial^N\lambda_i\big]^3+\big(C^h_i\big)^3\Big] + \delta^N_{0,\,\text{order}\geq 4}(\eta). \label{eq_L_star_bulk_new_0} 
\end{align}
The sum in the last line~\eqref{eq_L_star_bulk_new_0} will later be found to be an error term, in Section~\ref{sec_three_point}. The important terms are therefore the sums in~\eqref{eq_first_line_L_star_bulk}-\eqref{eq_second_line_L_star_bulk}, which we will see impose conditions on the choice of $g$. \\
To highlight the structure of $L^*_{h,0}{\bf 1}$, let us rewrite the sums in~\eqref{eq_first_line_L_star_bulk}-\eqref{eq_second_line_L_star_bulk} by grouping together terms involving $n$-point correlations, $n\in\N^*$. By~\eqref{eq_def_B_x_D_x}, $C$ is the sum of $B$, which involves one-point correlations (i.e. one $\bar\eta$); and of $D$, which is configuration-independent, like $\lambda$. Moreover, the sum in~\eqref{eq_first_line_L_star_bulk} will have to be integrated by parts to remove the $N$ factor. To do so, write:
\begin{equation}
\forall i<N-1,\qquad \eta_{i+1}-\eta_i = \bar\eta_{i+1}-\bar\eta_i + \bar\rho_{i+1}-\bar\rho_i = \bar\eta_{i+1}-\bar\eta_i + N^{-1}\partial^N\bar\rho_i.\label{eq_eta_as_bar_eta_plus_bar_rho_for_IPP}
\end{equation}
The sum in~\eqref{eq_first_line_L_star_bulk} therefore contains constant terms, fluctuations and two-point correlations. Let us similarly analyse the second line~\eqref{eq_second_line_L_star_bulk}. The jump rate $(\eta_{i+1}-\eta_i)^2$, $i<N-1$ can be expressed in terms of $\bar\eta_i$ and $\bar\eta_{i+1}$ as follows:
\begin{equation}
\forall i<N-1,\qquad (\eta_{i+1}-\eta_i)^2 = a_i + \sigma'(\bar\rho_i)\bar\eta_{i+1}+\sigma'(\bar\rho_{i+1})\bar\eta_i-2\bar\eta_i\bar\eta_{i+1},\label{eq_dvplt_jump_rate}
\end{equation}
where $a_i = \bar\rho_{i+1}(1-\bar\rho_i) + \bar\rho_i(1-\bar\rho_{i+1})$, $i<N-1$. The sum in~\eqref{eq_second_line_L_star_bulk} therefore involves constant terms and $n$-point correlations for each $1\leq n\leq 4$. In Section~\ref{sec_three_point}, we prove that three-point and four-point correlations lead to an error term $\delta^N_{0,3-4}$, while the sum in the third line~\eqref{eq_L_star_bulk_new_0} is an error term $\delta^N_{0,\,\text{order }3}$. The adjoint in the bulk thus reads:
\begin{equation}
N^2L^*_{h,0}{\bf 1} = \text{ Const } + \text{ Fluct } + \text{ Corr }+\delta^N_{0,3-4}+  \delta^{N}_{0,\,\text{order }3} + \delta^N_{0,\,\text{order}\geq 4}, \label{eq_L_star_as_F_plus_Corr_plus_Con}
\end{equation}
where Const, Fluct, Corr respectively denote the constant terms, the fluctuations and the correlations. The expression of these terms is given in the next three sections. Informally, these are small only when $\nu^N_g$ satisfies specific conditions. Namely, the fluctuations term Fluct is small because, by definition, $\nu^N_g$ has the same average occupation number as the invariant measure in the large $N$ limit, see Section~\ref{sec_fluctuations_L_star}. On the other hand, the correlations Corr are small provided $g$ solves the partial differential equation~\eqref{eq_main_equation}, as shown in Section~\ref{sec_correlations_L_star}. Finally, the constant Const is small provided all other terms are, as established in Section~\ref{subsec_constant_terms}. We will repeatedly use the following estimates (recall the definition~\eqref{eq_def_B_x_D_x} of $D$):
\begin{align}
&\sup_{i<N-1}|D_i| = O(N^{-1}),
\nonumber\\
&\sup_{i<N-1}\Big|\frac{a_i\partial^N\lambda_i}{2}-\partial^N\bar\rho_i\Big| 
= 
O(N^{-1}) 
= \sup_{i<N-1}\Big|\partial^N\lambda_i - \frac{\partial^N\bar\rho_i}{\bar\sigma_i}\Big|
.
\label{eq_D_and_lambda_moins_rho_prime_petits}
\end{align}
\subsubsection{The fluctuations}\label{sec_fluctuations_L_star}
Here, we estimate the fluctuations term Fluct in~\eqref{eq_L_star_as_F_plus_Corr_plus_Con}, which we recall accounts for all terms with a single $\bar\eta$ in the two sums~\eqref{eq_first_line_L_star_bulk}--\eqref{eq_second_line_L_star_bulk}. 
Recalling~\eqref{eq_dvplt_jump_rate}, it reads:
\begin{align}
\text{Fluct }&=\frac{1}{2}\sum_{i<N-1}\bigg[N(\bar\eta_{i+1}-\bar\eta_i)\big(2D^{h-g}_i-\partial^N\lambda_i\big)+2\partial^N\bar\rho_iB^{h-g}_i 
\nonumber\\
&\hspace{1.5cm}+\frac{a_i}{2}\big[2B_i^{h-2g}\big(D^{h-2g}_i-\partial^N\lambda_i\big)-2B_i^hD_i^h\big]
\nonumber \\
&\hspace{1.5cm}+\frac{1}{2}\big[\sigma'(\bar\rho_i)\bar\eta_{i+1}+\sigma'(\bar\rho_{i+1})\bar\eta_i\big]\big[(D_i^{h-2g}-\partial^N\lambda_i)^2 - (D_i^h)^2\big]\bigg]
.
\label{eq_line_2_L_star_new}
\end{align}
To estimate the size of each term above, 
recall from Lemma~\ref{lemm_size_controllable} that a term of the form $Y^N(\phi)$ (defined in~\eqref{eq_def_fluct}), 
with $\phi:(-1,1)\rightarrow\R$ bounded, is controllable with size $1$. 
Using~\eqref{eq_D_and_lambda_moins_rho_prime_petits},~\eqref{eq_line_2_L_star_new} thus turns into:
\begin{align}
\text{Fluct }&= \frac{1}{2}\sum_{i<N-1}\Big[N(\bar\eta_{i+1}-\bar\eta_i)\big(2D_i^{h-g}-\partial^N \lambda_i\big)+2\partial^N\bar\rho_i B_i^{h-g}-a_i\partial^N\lambda_i B_i^{h-2g}\nonumber\\
&\hspace{3cm} + \frac{(\partial^N\lambda_i)^2}{2}\big[\sigma'(\bar\rho_i)\bar\eta_{i+1}+\sigma'(\bar\rho_{i+1})\bar\eta_i\big]\Big] + \delta^{N,1}_{0,1}(\eta),\label{eq_line_2_L_star_new_1}
\end{align}
where $\delta^{N,1}_{0,1}(\eta)$ reads:
\begin{align}
\delta^{N,1}_{0,1}(\eta) 
&:= 
\frac{1}{2}\sum_{i<N-1}\bigg[a_i\big[B^{h-2g}_i D_i^{h-2g}-B_i^gD_i^g\big] 
\\
&\hspace{2cm}+
\big[\bar\sigma'(\bar\rho_i)\bar\eta_{i+1} + \sigma'(\bar\rho_{i+1})\bar\eta_i\big]\big[(D_i^{h-2g})^2-2D_i^{h-2g}\partial^N\lambda_i - (D_i^g)^2\big]\bigg]
\nonumber
.
\end{align}
Recall from~\eqref{eq_def_B_x_D_x} that each $B^\cdot_i$ for $i<N-1$ is of the form $N^{-1/2}Y^N(\phi)$ for a bounded function $\phi$, and that $\sup_i|D_i| = O(N^{-1})$. As a result, all terms composing $\delta^{N,1}_{0,1}$ are of the form $N^{-1/2}Y^N(\psi)$ or $N^{-3/2}Y^N(\psi)$ for a bounded $\psi$, 
and $Y^N(\psi)$ is controllable with size $1$ by Lemma~\ref{lemm_size_controllable}. 
It follows that $\delta^{N,1}_{0,1}$ is controllable with size $N^{-1}$ and of vanishing type (recall Definitions~\ref{def_controllability}--\ref{def_LS_type}).\\ 

Let us compute~\eqref{eq_line_2_L_star_new_1}. 
We start by integrating its first term by parts.
From~\eqref{eq_def_B_x_D_x} and the regularity of $h,g$, one draws, for each $i<N-1$:
\begin{equation}
D^{h-g}_{i-1}-D^{h-g}_i = -\frac{1}{2N^2}\Big[(h-g)_{i,i+1}\Delta^N\bar\rho_i + N\big[(h-g)_{i,i+1}-(h-g)_{i-1,i}\big]\partial^N\bar\rho_{i-1}\Big],\label{eq_gradient_D^h_x}
\end{equation}
where:
\begin{equation}
\Delta^N\bar\rho_i = \partial^N (\partial^N\bar\rho_{i-1}) 
= 
N^2[\bar\rho_{i+1}+\bar\rho_{i-1}-2\bar\rho_i]
.
\end{equation}
\begin{remark}
In the present case, $\Delta^N\bar\rho_i = 0$ and $\partial^N\bar\rho_i = \bar\rho'$ for each $i$, 
which could be used to simplify~\eqref{eq_gradient_D^h_x} and several other expressions below. 
We chose not to use the properties of $\bar\rho$ until the end of this section to highlight the structure of the Fluct term: it will contain a discrete PDE involving $\bar\rho$. 
Had we defined $\nu^N_g$ in terms of a density function $\rho:[-1,1]\rightarrow[0,1]$, 
this PDE would determine the choice of density $\rho$ in the measure $\nu^N_g$ in order to obtain an optimal bound on the adjoint.
\demo
\end{remark}
As a result, 
$\sup_{i}\big|D^{h-g}_{i-1}-D^{h-g}_i\big| = O(N^{-2})$. 
Moreover, $g(\pm 1,\cdot) = 0 = h(\pm 1,\cdot)$, which implies that $D^{h-g}_{-(N-1)} = O(N^{-2}) = D^{h-g}_{N-2}$. 
An integration by parts therefore turns~\eqref{eq_line_2_L_star_new_1} into:
\begin{align}
\text{Fluct }
&= 
\frac{1}{2}\sum_{|i|<N-1}\bar\eta_i\Delta^N\lambda_i
-\frac{N}{2}\bar\eta_{N-1}\partial^N\lambda_{N-2}+\frac{N}{2}\bar\eta_{-(N-1)}\partial^N\lambda_{-(N-1)}  \label{eq_line_2_L_star_new_2}\\
&\quad+
\sum_{i<N-1}\Big[\partial^N\bar\rho_i B_i^{h-g}-\frac{a_i\partial^N\lambda_i}{2} B_i^{h-2g}+ \frac{(\partial^N\lambda_i)^2}{4}\big[\sigma'(\bar\rho_i)\bar\eta_{i+1}+\sigma'(\bar\rho_{i+1})\bar\eta_i\big]\Big] 
\nonumber\\
&\quad 
+ \sum_{q=1}^2\delta^{N,q}_{0,1}(\eta)
,
\nonumber
\end{align}
where $\delta^{N,2}_{0,1}$ is a controllable error term with size $N^{-1}$ of vanishing type, 
as it reads:
\begin{equation}
\delta^{N,2}_{0,1}(\eta) = \frac{1}{N}\sum_{|i|<N-1}\bar\eta_i N^2\big[D^{h-g}_{i-1}-D^{h-g}_i\big] + \bar\eta_{N-1}D^{h-g}_{N-2} - \bar\eta_{-(N-1)}D^{h-g}_{-(N-1)}.
\end{equation}
Using~\eqref{eq_D_and_lambda_moins_rho_prime_petits} to express $a_\cdot\partial^N\lambda_\cdot$ in terms of $\partial^N\bar\rho_\cdot$ in the terms involving $B$ in the second line,~\eqref{eq_line_2_L_star_new_2} becomes:
\begin{align}
\text{Fluct }&= \frac{1}{2}\sum_{|i|<N-1}\bar\eta_i\Delta^N\lambda_i -\frac{N}{2}\bar\eta_{N-1}\partial^N\lambda_{N-2}+\frac{N}{2}\bar\eta_{-(N-1)}\partial^N\lambda_{-(N-1)} \nonumber\\
&+\sum_{i<N-1}\Big[\partial^N\bar\rho_i B_i^{g}+\frac{(\partial^N\lambda_i)^2}{4}\big[\sigma'(\bar\rho_i)\bar\eta_{i+1}+\sigma'(\bar\rho_{i+1})\bar\eta_i\big]\Big] + \sum_{q=1}^3\delta^{N,q}_{0,1}(\eta),\label{eq_line_2_L_star_new_3}
\end{align}
with $\delta^{N,3}_{0,1}(\eta)$ of the form $N^{-1/2}Y^N(\phi)$ (recall~\eqref{eq_def_fluct}) with $\phi$ bounded, 
thus controllable with size $N^{-1}$ and of vanishing type. 
Indeed, recalling the estimates~\eqref{eq_D_and_lambda_moins_rho_prime_petits} and the definition~\eqref{eq_def_B_x_D_x} of $B$, 
$\delta^{N,3}_{0,1}$ reads:
\begin{equation}
\delta^{N,3}_{0,1}(\eta) 
= 
-\frac{1}{2N}\sum_{j\in\Lambda_N}\bar\eta_j \bigg(\frac{1}{N}\sum_{i\notin\{j-1,j,N-1\}}N\Big[\frac{a_i\partial^N\lambda_i}{2}-\partial^N\bar\rho_i\Big]\partial^N_1(h-2g)_{i,j}\bigg),
\end{equation}
and the term between parenthesis is bounded uniformly in $j\in\Lambda_N$ and $N\in\N^*$. 
Let us now compute the term involving $B$ in~\eqref{eq_line_2_L_star_new_3}. 
Recall that $g\in C^3(\bar\rhd)$ and is symmetric, 
thus $\sup_{|j|<N-1} N|g_{j-1,j}-g_{j+1,j}|$ is bounded uniformly in $N$. 
Integrating by parts and using $g(\pm 1,\cdot) =0$, 
and the symmetry of $g$ in the last line below; 
we find:
\begin{align}
\sum_{i<N-1}&\partial^N\bar\rho_i B_i^{g} = \frac{1}{2N}\sum_{i<N-1}\partial^N\bar\rho_i\sum_{j\notin\{i,i+1\}}\bar\eta_j\partial^N_1 g_{i,j} = \frac{1}{2N}\sum_{j\in\Lambda_N} \bar\eta_j\sum_{i\notin\{j-1,j,N-1\}}\partial^N\bar\rho_i\partial^N_1 g_{i,j}\nonumber\\
&= - \frac{1}{2N}\sum_{j\in\Lambda_N} \bar\eta_j \bigg[\sum_{\substack{|i|<N-1 \\|j-i|>1}}\Delta^N\bar\rho_i g_{i,j} + \partial^N\bar\rho_{j-2}(Ng_{j-1,j}) - \partial^N\bar\rho_{j+1}(Ng_{j+1,j}) \nonumber\\
&\hspace{5cm}+ \partial^N\bar\rho_{N-2} (Ng_{N-1,j}) - \partial^N\bar\rho_{-N+1} (Ng_{-N+1,j})\bigg]\nonumber\\
&=: - \frac{1}{2N}\sum_{j\in\Lambda_N} \bar\eta_j \sum_{\substack{|i|<N-1 \\|j-i|>1}}\Delta^N\bar\rho_i g_{i,j} + \delta^{N,4}_{0,1}(\eta) \nonumber\\
&=: - \frac{1}{2}\sum_{j\in\Lambda_N} \bar\eta_j \big(N^{-1}M_g\big)(\Delta\bar\rho_\cdot)(j) + \delta^{N,5}_{0,1}(\eta),\label{eq_IPP_sur_B}
\end{align}
where $\delta^{N,4}_{0,1},\delta^{N,5}_{0,1}$ are of the form $N^{-1/2}Y^N(\phi)$ ($Y^N(\phi)$ is defined in~\eqref{eq_def_fluct}) for bounded $\phi$, 
and therefore $\delta^{N,4}_{0,1},\delta^{N,5}_{0,1}$ are controllable with size $N^{-1}$ and of vanishing type. 
$M_g$ is the matrix $(g_{i,j})_{(i,j)\in\Lambda_N^2}$ and $(\Delta\bar\rho_\cdot)$ the vector $(\Delta\bar\rho_i)_{i\in\Lambda_N}$, 
so that $\delta^{N,5}_{0,1}$ accounts for the replacement of $\Delta^N\bar\rho$ by $\Delta \bar\rho$ (this cost vanishes in our case since $\bar\rho$ is linear, 
but we do not use this fact at this point), 
as well as the addition of missing terms in the sum on $i$:
\begin{align}
\delta^N_{0,5}(\eta)
&= 
- \frac{1}{2N}\sum_{j\in\Lambda_N} \bar\eta_j \sum_{\substack{-N+1<i<N-1 \\|j-i|>1}}\big[\Delta^N\bar\rho_i-\Delta\bar\rho_i\big] g_{i,j}
\nonumber\\
&\quad - \frac{1}{2N}\sum_{j\in\Lambda_N} \bar\eta_j \sum_{i\in\{\pm(N-1),j,j\pm 1\}}\Delta\bar\rho_i g_{i,j}
.
\end{align}
Consider now the sums involving $\lambda_\cdot$ in~\eqref{eq_line_2_L_star_new_3}. Elementary computations give, for each $i<N-1$:
\begin{align}
\Delta^N\lambda_i = N[\partial^N\lambda_i-\partial^N\lambda_{i-1}] = \frac{\Delta^N\bar\rho_i}{\bar\sigma_i} - \frac{\big(\partial^N\bar\rho_i\big)^2\sigma'(\bar\rho_i) }{(\bar\sigma_i)^2} +  \epsilon^N_i,\qquad \sup_{|i|<N-1}|\epsilon^N_i| = O(N^{-1}).
\end{align}
By~\eqref{eq_D_and_lambda_moins_rho_prime_petits}, we also know $\sup_i|\partial^N\lambda_i-\partial^N\bar\rho_i/\bar\sigma_i| = O(N^{-1})$. As a result:
\begin{align}
\frac{1}{2}\sum_{|i|<N-1}\bar\eta_i\Delta^N\lambda_i 
&+\sum_{i<N-1}\frac{(\partial^N\lambda_i)^2}{4}\big[\sigma'(\bar\rho_i)\bar\eta_{i+1}+\sigma'(\bar\rho_{i+1})\bar\eta_i\big] 
\nonumber\\
&-\frac{N}{2}\bar\eta_{N-1}\partial^N\lambda_{N-2}+\frac{N}{2}\bar\eta_{-(N-1)}\partial^N\lambda_{-(N-1)} 
\nonumber\\
&\hspace{-1cm}=  
\frac{1}{2}\sum_{i\in\Lambda_N}\bar\eta_i\frac{\Delta\bar\rho_i}{\bar\sigma_i} -\frac{N}{2}\bar\eta_{N-1}\partial^N\lambda_{N-2}+\frac{N}{2}\bar\eta_{-(N-1)}\partial^N\lambda_{-(N-1)} + \delta^{N,6}_{0,1}(\eta)
,
\label{eq_recast_lambda_fluctuations}
\end{align}
where $\delta^{N,6}_{0,1}$ reads:
\begin{align}
\delta^{N,6}_{0,1}(\eta) 
&:= 
\frac{1}{2}\sum_{|i|<N-1}\bar\eta_i\Big[\Delta^N\lambda_i + \frac{1}{4}\big[(\partial^N\lambda_{i-1})^2\sigma'(\bar\rho_{i-1})+(\partial^N\lambda_i)^2\sigma'(\bar\rho_{i+1})\big]-\frac{\Delta\bar\rho_i}{\bar\sigma_i}\Big]\nonumber\\
&\quad
+ \frac{\bar\eta_{N-1}(\partial^N\lambda_{N-2})^2\sigma'(\bar\rho_{N-2})}{4} + \frac{\bar\eta_{-(N-1)}(\partial^N\lambda_{-(N-1)})^2\sigma'(\bar\rho_{-(N-2)})}{4}\nonumber\\
&\quad
 - \bar\eta_{N-1}\frac{\Delta\bar\rho_{N-1}}{\bar\sigma_{N-1}} - \bar\eta_{-(N-1)}\frac{\Delta\bar\rho_{-(N-1)}}{\bar\sigma_{-(N-1)}}.
\end{align}
The function $\delta^{N,6}_{0,1}$ involves terms of the form $\bar\eta_{\pm(N-1)}u(\pm(N-1))$ with $u$ bounded (the last two lines), 
and $N^{-1/2}Y^N(\phi)$ for bounded $\phi:(-1,1)\rightarrow\R$ (the first line). 
It follows that $\delta^{N,6}_{0,1}$ is $\Gamma$-controllable with size $N^{-1}$ and of vanishing type by Lemma~\ref{lemm_type}. 
Note that the last line actually vanishes since $\Delta\bar\rho = 0$. 
This last line is an error term anyway, 
so we do not need this fact.\\
Putting~\eqref{eq_line_2_L_star_new_3},~\eqref{eq_IPP_sur_B} and~\eqref{eq_recast_lambda_fluctuations} together, 
we have computed the fluctuations term~\eqref{eq_line_2_L_star_new} in $N^2L^*_{h,0}{\bf 1}$: 
\begin{align}
\text{Fluct }&= \frac{1}{2}\sum_{i\in\Lambda_N}\bar\eta_i\big(\bar\sigma^{-1}-N^{-1}M_g\big)\big(\Delta\bar\rho_\cdot)(i) -\frac{N}{2}\bar\eta_{N-1}\partial^N\lambda_{N-2}\nonumber\\
&\qquad +\frac{N}{2}\bar\eta_{-(N-1)}\partial^N\lambda_{-(N-1)} + \sum_{q=1}^6\delta^{N,q}_{0,1}(\eta).
\end{align}
Since $\bar\rho$ is the steady state profile satisfying $\Delta\bar\rho_\cdot = 0$, 
the first sum vanishes, and:
\begin{equation}
\text{Fluct }= -\frac{N}{2}\bar\eta_{N-1}\partial^N\lambda_{N-2}+\frac{N}{2}\bar\eta_{-(N-1)}\partial^N\lambda_{-(N-1)}  + \delta^N_{0,1}(\eta),\qquad \delta^N_{0,1} :=\sum_{q=1}^6\delta^{N,q}_{0,1}.\label{eq_final_expression_fluct_part_L_star}
\end{equation}
By definition of error terms, see Definition~\ref{def_controllability}, 
we have proven the following: 
for any $\theta,\gamma>0$, there is $C(\theta,\gamma)>0$ such that, for any density $f$ for $\nu^N_g$ and $N$ larger than some $N(\gamma)$:
\begin{align}
\nu^N_g\big(f\cdot\text{Fluct}\big) 
&\leq 
\frac{N}{2}\nu^N_g\big(f[-\bar\eta_{N-1}\partial^N\lambda_{N-2}+\bar\eta_{-(N-1)}\partial^N\lambda_{-(N-1)} ]\big)
\nonumber\\
&\quad + \theta N^2\nu^N_g\big(\Gamma_h(\sqrt{f})\big) + \frac{H(f\nu^N_g|\nu^N_g)}{\gamma} + \frac{C(\theta,\gamma)}{N}.\label{eq_Fluct_sous_f_nu_N_g_final}
\end{align}
The expectation in the first line of~\eqref{eq_Fluct_sous_f_nu_N_g_final} is not an error term, 
but it will cancel out with the boundary term obtained in~\eqref{eq_final_adjoint_boundary}. 
This cancellation will appear in Section~\ref{sec_conclusion_L_star}, 
where all contributions to $L^*_h{\bf 1}$ are summed.
\subsubsection{The correlations}\label{sec_correlations_L_star}
In this section, we compute the Corr term in~\eqref{eq_L_star_as_F_plus_Corr_plus_Con} and obtain the partial differential equation that an optimal $g$ must solve. Recall that Corr corresponds to all terms in~\eqref{eq_first_line_L_star_bulk}-\eqref{eq_second_line_L_star_bulk} that involve products of two $\bar\eta$'s. It reads:
\begin{align}
\text{Corr } &= \sum_{i<N-1}\bigg[N(\bar\eta_{i+1}-\bar\eta_i)B^{h-g}_i -a_iB^{g}_iB^{h-g}_i -\frac{1}{2}\bar\eta_i\bar\eta_{i+1}\big[\big(D_i^{h-2g}-\partial^N\lambda_i\big)^2 - (D_i^h)^2\big] \nonumber\\
&\hspace{1cm}+ \frac{1}{2}\big[\sigma'(\bar\rho_i)\bar\eta_{i+1}+\sigma'(\bar\rho_{i+1})\bar\eta_i\big]\big[B_i^{h-2g}\big(D_i^{h-2g}-\partial^N\lambda_i\big)-B_i^hD_i^h\big]\bigg]\label{eq_line_3_L_star_new}
\end{align}
Recall from Lemma~\ref{lemm_size_controllable} that terms of the form $\Pi^N(u)$, or of the form:
\begin{equation}
N^{-1/2}X^{v}_{1,\{0,1\}} 
= 
N^{-1/2}\sum_{i<N-1}\bar\eta_i\bar\eta_{i+1}v_i
,
\end{equation} 
are controllable with size $1$ as soon as the test functions $u,v$ are bounded. 
Multiplying them by $\epsilon_N$ with $\epsilon_N = o_N(1)$ therefore turns them into controllable error terms with size $\epsilon_N$, 
$\epsilon_N^2$ respectively by Lemma~\ref{lemm_size_controllable}. 
As in Section~\ref{sec_fluctuations_L_star}, 
we first use the estimate~\eqref{eq_D_and_lambda_moins_rho_prime_petits} 
on the size of $D$ to remove some terms from~\eqref{eq_line_3_L_star_new}:
\begin{align}
\text{Corr } 
&= 
\sum_{i<N-1}\Big[N(\bar\eta_{i+1}-\bar\eta_i)B^{h-g}_i -a_iB^{g}_iB^{h-g}_i -\frac{1}{2}\bar\eta_i\bar\eta_{i+1}\big(\partial^N\lambda_i\big)^2\nonumber\\
&\hspace{2.5cm}- 
\frac{1}{2}\big[\sigma'(\bar\rho_i)\bar\eta_{i+1}+\sigma'(\bar\rho_{i+1})\bar\eta_i\big]\partial^N\lambda_iB_i^{h-2g}\Big] + \delta^{N,1}_{0,2}(\eta),\label{eq_line_3_L_star_new_1}
\end{align}
with $\delta^{N,1}_{0,2}$ controllable with size $N^{-1}$ of vanishing type, defined by:
\begin{align}
\delta^{N,1}_{0,2}(\eta)&:= -\frac{1}{2N}\sum_{i<N-1}\bar\eta_i\bar\eta_{i+1}N\big[(D_i^{h-2g})^2-2D_i^{h-2g}\partial^N\lambda_i - (D_i^h)^2\big] \nonumber\\
&\qquad +\frac{1}{2N}\sum_{i<N-1}\big[\sigma'(\bar\rho_i)\bar\eta_{i+1}+\sigma'(\bar\rho_{i+1})\bar\eta_i\big]N\big[B_i^{h-2g}D_i^{h-2g}-B_i^hD_i^h\big].
\end{align}
Let us integrate by parts the term involving $\bar\eta_{i+1}-\bar\eta_i$ in~\eqref{eq_line_3_L_star_new_1}, $i<N-1$. To do so, notice first that, for each $i$ with $|i|<N-1$:
\begin{align}
N\big[B^{h-g}_{i-1} -B^{h-g}_{i}\big] 
&= 
N\big[B^{g-h}_{i} -B^{g-h}_{i-1}\big] 
\nonumber\\
&=
\frac{1}{2N}\sum_{j:|j-i|>1} \bar\eta_j \Delta^N_1 (g-h)_{i,j}
\nonumber\\
&\qquad
+\frac{1}{2}\big[\bar\eta_{i-1}\partial^N_1 (g-h)_{i,i-1}-\bar\eta_{i+1}\partial^N_1 (g-h)_{i-1,i+1}\big]
,
\end{align}
with $\Delta^N_1u(i,j)= \partial^N_1(\partial^N_1u(i-1,j))$ for $u:\Z^2\rightarrow\R$ and $(i,j)\in\Z^2$. 
As a result:
\begin{align}
&\sum_{i<N-1}N(\bar\eta_{i+1}-\bar\eta_i)B^{h-g}_i 
\nonumber\\
&\hspace{2cm}= 
N\bar\eta_{N-1}B^{h-g}_{N-2} -N\bar\eta_{-(N-1)}B^{h-g}_{-(N-1)} + \sum_{|i|<N-1}\bar\eta_i N\big[B_{i-1}^{h-g}-B_i^{h-g}\big]
\nonumber \\
&\hspace{2cm}=  
\delta^{N,2}_{0,2}(\eta)+ \frac{1}{2N}\sum_{|i|<N-1}\sum_{j:|j-i|>1} \bar\eta_i\bar\eta_j \Delta^N_1 (g-h)_{i,j} 
\label{eq_line_3_L_star_new_2}\\
&\hspace{2cm}\quad
+ \frac{1}{2}\sum_{|i|<N-1}\bar\eta_i\big[\bar\eta_{i-1}\partial^N_1 (g-h)_{i,i-1}-\bar\eta_{i+1}\partial^N_1 (g-h)_{i-1,i+1}\big]
,
\nonumber
\end{align}
where $\delta^{N,2}_{0,2}(\eta) := N\bar\eta_{N-1}B^{h-g}_{N-2}  -N\bar\eta_{-(N-1)}B^{h-g}_{-(N-1)}$. 
$\delta^{N,2}_{0,2}$ involves correlations between the reservoir and the bulk, 
of the same form as the function $U^\pm_1$ defined in Lemma~\ref{lemm_size_controllable}. 
$\delta^{N,2}_{0,2}$ is thus $\Gamma$-controllable with size $N^{-1}$ and of vanishing type. 
On the other hand, the last line of~\eqref{eq_line_3_L_star_new_2} becomes, changing indices:
\begin{align}
\sum_{|i|<N-1}&\bar\eta_i\big[\bar\eta_{i-1}\partial^N_1 (g-h)_{i,i-1}-\bar\eta_{i+1}\partial^N_1 (g-h)_{i-1,i+1}\big] \nonumber\\
&= \sum_{-(N-1)<i<N-2}\bar\eta_i\bar\eta_{i+1}\big[\partial^N_1 (g-h)_{i+1,i}-\partial^N_1 (g-h)_{i-1,i+1}\big] +\delta^{N,3}_{0,2}(\eta),
\end{align}
with:
\begin{align}
\delta^{N,3}_{0,2}(\eta) 
\, &:=
\bar\eta_{-(N-1)}\bar\eta_{-(N-2)}\partial^N_1(g-h)_{-(N-2),-(N-1)} 
\nonumber\\
&\qquad 
- \bar\eta_{N-2}\bar\eta_{N-1}\partial^N_1(g-h)_{N-3,N-1}
.
\end{align}
The quantity $\delta^{N,3}_{0,2}$ is again a $\Gamma$-controllable error term with size $N^{-1}$ of vanishing type, 
as it is of the same form as the function $U_1^\pm$ of Lemma~\ref{lemm_size_controllable}. 
Equation~\eqref{eq_line_3_L_star_new_2} thus reads:
\begin{align}
&\sum_{i<N-1}N(\bar\eta_{i+1}-\bar\eta_i)B^{h-g}_i = \frac{1}{2N}\sum_{|i|<N-1}\sum_{j:|j-i|>1} \bar\eta_i\bar\eta_j \Delta^N_1 (g-h)_{i,j}
\nonumber\\
&\quad 
+\frac{1}{2}\sum_{-(N-1)<i<N-2}\bar\eta_i\bar\eta_{i+1}\big[\partial^N_1 (g-h)_{i+1,i}-\partial^N_1 (g-h)_{i-1,i+1}\big]+\sum_{q=2}^3\delta^{N,q}_{0,2}(\eta)
.
\label{eq_Neumann_IPP_L_star}
\end{align}
The other terms in~\eqref{eq_line_3_L_star_new_1} are simpler. Indeed, recall that:
\begin{equation}
\sup_{i<N-1}|a_i-2\bar\sigma_i| = O(N^{-1}),\qquad \sup_{i<N-1}\Big|\partial^N\lambda_i -\frac{\partial^N\bar\rho_i}{\bar\sigma_i}\Big| = O(N^{-1}).
\end{equation}
Using these estimates in~\eqref{eq_line_3_L_star_new_1}, Corr becomes:
\begin{align}
&\text{Corr }
= 
\frac{1}{2N}\sum_{|i|<N-1}\sum_{j:|j-i|>1} \bar\eta_i\bar\eta_j \Delta^N_1 (g-h)_{i,j}\nonumber\\
&\quad+
\frac{1}{2}\sum_{-(N-1)<i<N-2}\bar\eta_i\bar\eta_{i+1}\Big[\partial^N_1 (g-h)_{i+1,i}-\partial^N_1 (g-h)_{i-1,i+1} - \frac{(\partial^N\bar\rho_i)^2}{(\bar\sigma_i)^2}\Big] \nonumber\\
&\quad- 
\sum_{i<N-1}\Big[\big[\sigma'(\bar\rho_i)\bar\eta_{i+1}+\sigma'(\bar\rho_{i+1})\bar\eta_i\big]\frac{\partial^N\bar\rho_i}{2\bar\sigma_i}B_i^{h-2g} +2\bar\sigma_i B^{g}_iB^{h-g}_i\Big]  +\sum_{q=1}^4\delta^{N,q}_{0,2}(\eta),\label{eq_line_3_L_star_new_3}
\end{align}
where $\delta^{N,4}_{0,2}$ reads:
\begin{align}
&\delta^{N,4}_{0,2}(\eta) 
= 
-\frac{1}{2N}\sum_{i<N-1}\big[\sigma'(\bar\rho_i)\bar\eta_{i+1}+\sigma'(\bar\rho_{i+1})\bar\eta_i\big]\Big[\partial^N\lambda_i - \frac{\partial^N\bar\rho_i}{\bar\sigma_i}\Big]B_i^{h-2g}  
\nonumber\\
&\qquad 
- \frac{1}{N}\sum_{i<N-1}N\big[a_i-2\bar\sigma_i\big]B_i^gB_i^{h-g}
-\frac{1}{2N}\sum_{-(N-1)<i<N-2}\bar\eta_i\bar\eta_{i+1}N\Big[(\partial^N\lambda_i)^2-\frac{(\partial^N\bar\rho_i)^2}{\bar\sigma_i^2}\Big]
\nonumber\\
&\qquad 
-\frac{1}{2}\Big[\bar\eta_{-(N-1)}\bar\eta_{-(N-2)}\big(\partial^N\lambda_{-(N-1)}\big)^2 + \bar\eta_{N-2}\bar\eta_{N-1}\big(\partial^N\lambda_{N-2}\big)^2\Big]
.
\end{align}
The last line comes from the fact that the sum involving $\bar\eta_i\bar\eta_{i+1}$ in~\eqref{eq_line_3_L_star_new_1} and in~\eqref{eq_Neumann_IPP_L_star} do not have the same range. 
It is of the same form as $U_1^+$ in Lemma~\ref{lemm_size_controllable}, 
thus is $\Gamma$-controllable term with size $N^{-1}$ and of vanishing type by Lemmas~\ref{lemm_size_controllable}--\ref{lemm_type}. 
The first line and the first term of the second line above are of the form $N^{-2}X^{\phi_1}_{2,\{0\}}$, 
while the second term of the second line reads $N^{-1/2}X^{\phi_2}_{1,\{0,1\}}$; 
for bounded tensors $\phi_1,\phi_2$. 
As a result, $\delta^{N,4}_{0,2}$ is $\Gamma$-controllable with size $N^{-1}$ and of vanishing type by Lemmas~\ref{lemm_size_controllable}--\ref{lemm_type}. \\

To conclude on the expression of the correlations, 
it remains to take care of the two terms involving $B$ in~\eqref{eq_line_3_L_star_new_3}. 
Recalling the definition~\eqref{eq_def_B_x_D_x} of $B$, 
using the regularity of $h,g$ and $\bar\rho$ and changing indices, 
one can write:
\begin{align}
\frac{1}{2}&\sum_{i<N-1}\big[\sigma'(\bar\rho_i)\bar\eta_{i+1}+\sigma'(\bar\rho_{i+1})\bar\eta_i\big]\frac{\partial^N\bar\rho_i}{\bar\sigma_i}B_i^{2g-h} \\
&\hspace{2.5cm}
= 
\frac{1}{2N}\sum_{|i|<N-1}\frac{\sigma'(\bar\rho_i)\partial^N\bar\rho_i}{\bar\sigma_i}\bar\eta_{i}\sum_{j\notin\{i,i+1\}}\bar\eta_j\partial^N_1(2g-h)_{i,j} + \delta^{N,5}_{0,2}(\eta)
,
\nonumber
\end{align}
where $\delta^{N,5}_{0,2}$ is an error term that reads:
\begin{align}
&\delta^{N,5}_{0,2}(\eta) 
= 
\frac{1}{2N}\sum_{|i|<N-1}\bar\eta_iN\big[\sigma'(\bar\rho_{i+1})-\sigma'(\bar\rho_i)\big]\frac{\partial^N\bar\rho_i}{\bar\sigma_i}B_i^{2g-h} \nonumber\\
&\qquad+ \frac{1}{2N}\sum_{|i|<N-1}\bar\eta_i N\Big[\sigma'(\bar\rho_{i-1})\frac{\partial^N\bar\rho_{i-1}}{\bar\sigma_{i-1}}B_{i-1}^{2g-h} - \sigma'(\bar\rho_i)\frac{\partial^N\bar\rho_i}{\bar\sigma_i}B_i^{2g-h}\Big]\nonumber\\
&\qquad + \bar\eta_{N-1}\frac{\sigma'(\bar\rho_{N-2})\partial^N\bar\rho_{N-2}}{2\bar\sigma_{N-2}}B_{N-2}^{2g-h} + \bar\eta_{-(N-1)}\frac{\sigma'(\bar\rho_{-(N-2)})\partial^N\bar\rho_{-(N-1)}}{2\bar\sigma_{-(N-1)}}B^{2g-h}_{-(N-1)}.
\end{align}
$\delta^{N,5}_{0,2}$ is of the form $N^{-1}X^{\phi_2}_{2,\{0\}}$ for a bounded $\phi_2$ for the first two lines, 
and $N^{-1}U_1^\pm$ for the third line. 
By Lemmas~\ref{lemm_size_controllable}--\ref{lemm_type}, 
$\delta^{N,5}_{0,2}$ is therefore $\Gamma$-controllable with size $N^{-1}$ and of vanishing type. \\

Finally, recall that $(\bar\eta_\cdot)^2 = \bar\sigma_\cdot + \sigma'(\bar\rho_\cdot)\bar\eta_\cdot$. 
Separating diagonal and off-diagonal contributions, 
the term involving $B^g_\cdot B_\cdot^{h-g}$ in~\eqref{eq_line_3_L_star_new_3} reads:
\begin{align}
-2\sum_{i<N-1}\bar\sigma_i B^g_iB_i^{h-g} &= \frac{1}{2N^2}\sum_{i<N-1}\bar\sigma_i \sum_{j,\ell\notin\{i,i+1\}}\bar\eta_j\bar\eta_\ell \partial^N_1 g_{i,j}\partial^N_1(g-h)_{i,\ell}\nonumber\\
&= \frac{1}{2N}\sum_{\substack{|j|<N-1 \\ \ell\neq j}}\bar\eta_j\bar\eta_\ell\Big(\frac{2}{2N}\sum_{i\notin\{j-1,j,\ell-1,\ell,N-1\}}\bar\sigma_i \partial^N_1 g_{i,j}\partial^N_1(g-h)_{i,\ell}\Big)\nonumber\\
&\quad + \frac{1}{2N^2}\sum_{i<N-1}\sum_{j\notin\{i,i+1\}}\bar\sigma_i\bar\sigma_j\partial^N_1 g_{i,j}\partial^N_1(g-h)_{i,j} + \delta^{N,6}_{0,2}(\eta),
\end{align}
where $\delta^{N,6}_{0,2}$ is the sum of error terms of the form $N^{-1/2}Y^N(\phi)$ (recall~\eqref{eq_def_fluct}) and $U^\pm_1$, 
thus $\delta^{N,6}_{0,2}$ is $\Gamma$-controllable with size $N^{-1}$ and of vanishing type. 
It is given by :
\begin{align}
\delta^{N,6}_{0,2}(\eta) 
&= 
\frac{1}{2N^2}\sum_{i<N-1}\sum_{j\notin\{i,i+1\}}\bar\sigma_i  \sigma'(\bar\rho_j)\bar\eta_j\partial^N_1 g_{i,j}\partial^N_1(g-h)_{i,j}\nonumber\\
&\quad +
\sum_{j\in\{\pm(N-1)\}}\sum_{\ell\neq j}\bar\eta_j\bar\eta_\ell\Big(\frac{2}{2N}\sum_{i\notin\{j-1,j,\ell-1,\ell,N-1\}}\bar\sigma_i \partial^N_1 g_{i,j}\partial^N_1(g-h)_{i,\ell}\Big).
\end{align}
The correlations~\eqref{eq_line_3_L_star_new} have so far been rewritten as follows:
\begin{align}
\hspace{-0.2cm}\text{Corr } 
&=
\sum_{\substack{|i|<N-1 \\ j\neq i}}\bar\eta_i\bar\eta_j\Bigg\{\frac{1}{2N}\bigg[{\bf 1}_{|i-j|>1}\Delta^N_1(g-h)_{i,j} 
+ {\bf 1}_{j\neq i+1}\partial^N\bar\rho_i\frac{\sigma'(\bar\rho_i)}{\bar\sigma_i}\partial^N_1 (2g-h)_{i,j} \label{eq_PDE_sur_g_minus_h_first_line}\\
&\hspace{3.5cm}+ 
\frac{1}{N}\sum_{\ell\notin\{i-1,i,j-1,j,N-1\}}\bar\sigma_\ell \partial^N_1 g_{\ell,i}\partial^N_1 (g-h)_{\ell,j}\bigg]\label{eq_PDE_sur_g_minus_h_second_line}\\
&\hspace{2cm}+
{\bf 1}_{j = i+1<N-1}\frac{1}{2}\bigg[\partial^N_1 (g-h)_{i+1,i}-\partial^N_1 (g-h)_{i-1,i+1} - \frac{\partial^N\bar\rho_i}{(\bar\sigma_i)^2}\bigg]\Bigg\}\label{eq_third_line_L_star_complete_with_omega_x_x+1}\\
&\quad 
+ \frac{1}{2N^2}\sum_{i<N-1}\sum_{j\notin\{i,i+1\}}\bar\sigma_i\bar\sigma_j\partial^N_1 g_{i,j}\partial^N_1(g-h)_{i,j} + \sum_{q=1}^6\delta^{N,q}_{0,2}(\eta).\label{eq_last_line_L_star_constants_and_errors}
\end{align}
We claim that the curly bracket is a discrete version of the partial differential equation~\eqref{eq_main_equation}. 
To see it, first use the symmetry of $g,h$ and exchange $i,j$. 
Recall then that, by assumption, $h,g\in W^{4,s}(\squaredash)$ for some $s>2$. 
By Sobolev embedding, $W^{4,s}(\squaredash)\subset C^3(\bar\rhd)\cap C^3(\bar\lhd)$, 
see Appendix~\ref{app_sobolev_spaces}. 
As a result, 
approximating discrete derivatives by continuous ones and the Riemann sum in~\eqref{eq_PDE_sur_g_minus_h_first_line} by an integral, 
there is an error term $\delta^{N,7}_{0,2}$, 
controllable with size $N^{-1}$ and of vanishing type, 
of the form $N^{-1}\Pi^N(u)+N^{-1}X^{v}_{1,\{0,1\}}$ for bounded $u,v$, 
such that the curly bracket spanning~\eqref{eq_PDE_sur_g_minus_h_first_line}--\eqref{eq_PDE_sur_g_minus_h_second_line}--\eqref{eq_third_line_L_star_complete_with_omega_x_x+1} equals:
\begin{align}
&\Pi^N\Big(\Delta(g-h) + \frac{(\bar\sigma)'}{\bar\sigma}\partial_1(2g-h) + \partial_2(2g-h)\frac{(\bar\sigma)'}{\bar\sigma}\Big)\nonumber\\
&\qquad
+ 
\Pi^N\Big( \mathcal M(\partial_1g,\partial_1(g-h)) + \mathcal M(\partial_1(g-h),\partial_1g)\Big) \nonumber\\
&\qquad
+ 
\sum_{i<N-1}\bar\eta_i\bar\eta_{i+1}\bigg[\partial_1(g-h)_{i_+,i}-\partial_1(g-h)_{i_-,i} -\frac{(\bar\rho')^2}{\bar\sigma_i^2}\bigg]+ \delta^{N,7}_{0,2}(\eta)
,\label{eq_terme_correspondant_a_EDP_dans_Corr}
\end{align}
with $\mathcal M$ the bilinear operator defined in~\eqref{eq_def_I} and the convention, for $(x,y)\in\squaredash$\,:
\begin{equation}
w\phi(x,y) = w(x)\phi(x,y),
\quad
\phi w(x,y) = \phi(x,y)w(y),
\qquad \phi:\squaredash\rightarrow\R, 
w:(-1,1)\rightarrow\R
. 
\end{equation}
If $g$ is chosen as the solution $g_h$ of the main equation~\eqref{eq_main_equation}, 
then the right-hand side of~\eqref{eq_terme_correspondant_a_EDP_dans_Corr} reduces to $\delta^{N,7}_{0,2}(\eta)$, 
whence:
\begin{equation}
\text{Corr } 
= 
\frac{1}{2N^2}\sum_{i<N-1}\sum_{j\notin\{i,i+1\}}\bar\sigma_i\bar\sigma_j\partial^N_1 (g_h)_{i,j}\partial^N_1(g_h-h)_{i,j} 
+ 
\delta^N_{0,2}(\eta),
\label{eq_final_expression_correl_part_L_star}
\end{equation}
with:
\begin{equation}
\delta^N_{0,2}
=
\sum_{q=1}^7\delta^{N,q}_{0,2}
.
\end{equation}
The first term is independent of the configuration, 
and will be estimated in Section~\ref{subsec_constant_terms}. 
By definition of error terms (see Definition~\ref{def_controllability}), 
for any $\theta,\gamma>0$, there is thus $C(\theta,\gamma)>0$ such that, 
for any density $f$ for $\nu^N_{g_h}$ and any $N$ large enough depending on $\gamma$:
\begin{align}
\nu^N_{g_h}\big(f\cdot \text{Corr}\big) 
&\leq 
\frac{1}{2N^2}\sum_{i<N-1}\sum_{j\notin\{i,i+1\}}\bar\sigma_i\bar\sigma_j\partial^N_1 (g_h)_{i,j}\partial^N_1(g_h-h)_{i,j} 
\nonumber\\
&\qquad 
+ \theta N^2\nu^N_{g_h}\big(\Gamma_h(\sqrt{f})\big)
+
\frac{H(f\nu^N_{g_h}|\nu^N_{g_h})}{\gamma} + \frac{C(\theta,\gamma)}{N}
.
\label{eq_Corr_sous_f_nu_N_g_final}
\end{align}
\begin{remark}
The choice of $g=g_h$ cancelling the curly bracket~\eqref{eq_PDE_sur_g_minus_h_first_line}--\eqref{eq_third_line_L_star_complete_with_omega_x_x+1} 
is optimal in the following sense: 
for another $g$, 
the first two terms in the right-hand side of~\eqref{eq_terme_correspondant_a_EDP_dans_Corr} are not even error terms, but only controllable with size $1$. 
As a result, for $g\neq g_h$,~\eqref{eq_Corr_sous_f_nu_N_g_final} can at best only be true with $\gamma_\theta N^{-1/2}, C_\theta$ respectively replacing $\gamma_\theta, C_\theta N^{-1/2}$ in the right-hand side, 
which breaks the $o_N(1)$ bound on the entropy obtained in Section~\ref{sec_conclusion_L_star}.
\demo
\end{remark}
\subsubsection{Higher-order correlations}\label{sec_three_point}
In this section, 
we again assume that $g=g_h$ solves the main equation~\eqref{eq_main_equation}, 
and we estimate the third order term $(N^2L^*_h{\bf 1})_{\text{order 3}}$ in the development~\eqref{eq_L_star_bulk_new_0} of the exponentials making up the adjoint in the bulk, 
as well as three-point and four-point correlations arising in~\eqref{eq_second_line_L_star_bulk}. 
The choice $g=g_h$ is not necessary here: higher-order correlations behave similarly under all measures $\nu^N_{g}$ with suitable $g$. 
We keep $g=g_h$ to avoid confusions.\\
 
Consider first $(N^2L^*_h{\bf 1})_{\text{order 3}}$, which by~\eqref{eq_L_star_bulk_new_0} reads:
\begin{align}
(N^2L^*_h{\bf 1}(\eta))_{\text{order 3}} =  \frac{1}{12N}\sum_{i<N-1}(\eta_{i+1}-\eta_i)\Big[\big[C^{h-2g_h}_i-\partial^N\lambda_i\big]^3+\big[C^h_i\big]^3\Big].
\end{align}
For $i<N-1$, write $\eta_{i+1}-\eta_i = \bar\eta_{i+1}-\bar\eta_i + N^{-1}\partial^N\bar\rho_i$ as before, 
and recall from~\eqref{eq_C_g_C_h_bounded} that $C^{h-2g_h}_\cdot$, $C_\cdot^h$ are bounded with $N$. As a result:
\begin{equation}
\Big|(N^2L^*_h{\bf 1}(\eta))_{\text{order 3}} -  \frac{1}{12N}\sum_{i<N-1}(\bar\eta_{i+1}-\bar\eta_i)\Big[\big[C^{h-2g_h}_i-\partial^N\lambda_i\big]^3+\big[C^h_i\big]^3\Big]\Big|\leq \frac{C(h)}{N}.\label{eq_bound_partielle_order_3_L_star_bulk}
\end{equation} 
One need not even integrate by parts to find that the sum in~\eqref{eq_bound_partielle_order_3_L_star_bulk} is an error term. 
Indeed, developing the cubes and recalling that $C = B +D$ (see~\eqref{eq_def_B_x_D_x}) with $ND$ bounded, one finds:
\begin{align}
\frac{1}{12N}\sum_{i<N-1}(\bar\eta_{i+1}-\bar\eta_i)\Big[\big[C^{h-2g_h}_i-\partial^N\lambda_i\big]^3+\big[C^h_i\big]^3\Big] 
= 
\sum_{n=1}^4N^{-n}X^{\phi_n}_{n,\{0\}}(\eta),
\end{align}
for bounded $\phi_n$, $1\leq n\leq 4$, depending on $\rho_\pm,h$. 
By Lemmas~\ref{lemm_size_controllable}--\ref{lemm_type}, 
$N^{-n}X^{\phi_n}_{n,\{0\}}$ is controllable with size $N^{-3/2}$ at most for $1\leq n\leq 4$, 
and of vanishing type. 
This observation and~\eqref{eq_bound_partielle_order_3_L_star_bulk} yield:
\begin{equation}
(N^2L^*_h{\bf 1}(\eta))_{\text{order 3}} :=\delta^{N}_{0,\,\text{order 3}}(\eta)
,\quad \big|\delta^{N}_{0,\,\text{order 3}}(\eta)\big| 
\leq 
\frac{C(h)}{N} + \Big|\sum_{n=1}^4N^{-n}X^{\phi_n}_{n,\{0\}}(\eta)\Big|
,
\label{eq_def_delta_order_3}
\end{equation}
and $\delta^{N}_{0,\,\text{order }3}$ is controllable with size $N^{-1}$ and of vanishing type. \\
 
Consider now three-and four-point correlations arising in the terms~\eqref{eq_L_star_bulk_new_0}--\eqref{eq_eta_as_bar_eta_plus_bar_rho_for_IPP}--\eqref{eq_dvplt_jump_rate}. 
They read:
\begin{align}
\delta^{N}_{0,3-4}(\eta) 
:= 
&\sum_{i<N-1}\bigg[-\big[\sigma'(\bar\rho_i)\bar\eta_{i+1}+\sigma'(\bar\rho_{i+1})\bar\eta_i\big]B^{g_h}_i B^{h-g_h}_i
\nonumber\\
&\qquad +
\bar\eta_i\bar\eta_{i+1}B_i^{g_h}\big(D_i^{2h-2g_h}-\partial^N\lambda_i\big) 
+ \bar\eta_i\bar\eta_{i+1}B_i^{h-g_h}\big(D_i^{2g_h}+\partial^N\lambda_i\big) \bigg]
\label{eq_line_4_L_star_new}
\\
&\, +
2\sum_{i<N-1}\bar\eta_i\bar\eta_{i+1}B^{g_h}_i B^{h-g_h}_i\label{eq_line_5_L_star_new}
.
\end{align}
Recall that $\sup_i|D^\phi_i|=O(N^{-1})$ for any bounded $\phi:[-1,1]^2\rightarrow\R$, 
and that $\sup_i|\partial^N\lambda_i -\bar\rho'/\bar\sigma_i| = O(N^{-1})$. 
There are then bounded functions $\phi_{n},\psi_n :(-1,1)^n\rightarrow\R$ ($n\in\{2,3\}$) such that, 
in the notations of Lemma~\ref{lemm_size_controllable}: 
\begin{align}
\sum_{i<N-1}\Big[-\big[\sigma'(\bar\rho_i)\bar\eta_{i+1}+\sigma'(\bar\rho_{i+1})\bar\eta_i\big]B^{g_h}_i B^{h-g_h}_i 
&=:
N^{-2}X^{\phi_3}_{3,\{0\}},
\\
\sum_{i<N-1} \bar\eta_i\bar\eta_{i+1}\Big[
B_i^{g_h}D_i^{2h-2g_h}+  B_i^{h-g_h}D_i^{2g_h}
\Big]
&=: 
N^{-2}X^{\phi_2}_{2,\{0,1\}},
\\
\sum_{i<N-1} \bar\eta_i\bar\eta_{i+1}B_i^{h-2g_h}\partial^N\lambda_i
&=:
N^{-1} X^{\psi_2}_{2,\{0,1\}},
\label{eq_def_psi_2_error_terms_3_point}
\\
2\sum_{i<N-1}\bar\eta_i\bar\eta_{i+1}B^{g_h}_i B^{h-g_h}_i
&=:
N^{-2}X^{\psi_3}_{3,\{0,1\}}
,
\end{align}
so that $\delta^{N}_{0,3-4}$ reads:
\begin{equation}\label{eq_splitting_higher_order_error_terms}
\delta^{N}_{0,3-4} 
= 
 N^{-2}X^{\phi_2}_{2,\{0,1\}} + N^{-2}X^{\phi_3}_{3,\{0\}} 
 + N^{-2}X^{\psi_3}_{3,\{0,1\}} + N^{-1}X^{\psi_2}_{2,\{0,1\}}.
\end{equation}
The first three terms are controllable error terms with size $N^{-1/2}$ at most, 
with the first term of vanishing type. 
The second and third term, involving $\phi_3,\psi_3$ are of large type, 
and $\phi_3,\psi_3$ are bounded as follows:
\begin{equation}
\max\Big\{\|\phi_3\|_\infty,\|\psi_3\|_\infty\Big\}
\leq 
\frac{1}{2}\|\partial_1 g_h\|_\infty\|\partial_1 (g_h-h)\|_\infty
.
\end{equation}
By taking $\bar\rho'\leq \epsilon$ and $h\in\s(\epsilon)$ for sufficiently small $\epsilon$, 
the right-hand side can be made as small as needed, see Appendix~\ref{app_choice_epsilonB}. 
Lemma~\ref{lemm_type} therefore ensures that $N^{-2}X^{\phi_3}_{3,\{0\}}$ and $N^{-2}X^{\psi_3}_{3,\{0,1\}}$ are of LS type.\\
The last term in~\eqref{eq_splitting_higher_order_error_terms} is $N^{-1}X^{\psi_2}_{2,\{0,1\}}$ given in~\eqref{eq_def_psi_2_error_terms_3_point}. 
$\psi_2$ satisfies:
\begin{equation}\label{eq_bound_psi_2_three_point_correl}
\|\psi_2\|_\infty
\leq 
\frac{C\bar\rho'\|\partial_1(h-2g_h)\|_\infty }{\min\{\rho_-(1-\rho_-),\rho_+(1-\rho_+)\}}
\underset{\bar\rho'\rightarrow 0}{\longrightarrow} 0
.
\end{equation}
At first glance, the second item of Lemma~\ref{lemm_size_controllable} only yields that $N^{-1}X^{\psi_2}_{2,\{0,1\}}$ is controllable with size $\|\psi_2\|_\infty$.  
The last item of Lemma~\ref{lemm_size_controllable} in fact yields that $N^{-1}X^{\psi_2}_{2,\{0,1\}}$ is $\Gamma$-controllable with size $\|\psi_2\|_\infty N^{-1/2}$. 
This improvement is obtained in Appendix~\ref{sec_estimate_W_3}, 
through the renormalisation scheme of Jara and Menezes~\cite{Jara2018}.
In view of the bound~\eqref{eq_bound_psi_2_three_point_correl}, 
taking $\bar\rho'\leq \epsilon$ and $h\in\s(\epsilon)$ for sufficiently small $\epsilon$ ensures that $N^{-1}X^{\phi_2}_{2,\{0,1\}}$ is of LS type. 
All terms in the decomposition~\eqref{eq_splitting_higher_order_error_terms} have now been treated, 
and we conclude that, 
for small enough $\epsilon>0$ such that $\bar\rho'\leq\epsilon,h\in\s(\epsilon)$, 
$\delta^{N}_{0,3-4}$ is a $\Gamma$-controllable error term of LS type.
Thus, for all large enough $N$ and each $\delta>0$:
\begin{align}
\nu^N_g\big(f\delta^N_{0,3-4}\big) 
&\leq 
\delta N^2 \nu^N_{g}\big(\Gamma_h(\sqrt{f})\big)
+ \frac{H(f\nu^N_g|\nu^N_g)}{2^{10}C_{LS}\delta} 
+ \frac{C'(g,h)}{N^{1/2}\delta }
.
\end{align}
\begin{remark}
In addition to the boundary terms, 
the term $N^{-1}X^{\psi_2}_{2,\{0,1\}}$ is the only one for which a renormalisation scheme using the carre du champ is needed.
\demo
\end{remark}
\subsubsection{The constant terms}\label{subsec_constant_terms}
Here, 
we prove that the configuration-independent terms appearing in the full adjoint $N^2L^*_{h}{\bf 1}$ are small when $g=g_h$, 
with $g_h$ the solution of the main equation~\eqref{eq_main_equation}. 
The terms in question correspond to various constant terms bounded by $O(N^{-1})$ encountered in the previous subsections and the computation of $N^2L^*_{h,\pm}{\bf 1}$, 
which already are error terms; 
and the sum of the constant term in~\eqref{eq_final_expression_correl_part_L_star}, 
as well as the Const term of~\eqref{eq_L_star_as_F_plus_Corr_plus_Con}, 
which reads:
\begin{equation}
\text{Const } := \sum_{i<N-1}\Big[ \partial^N\bar\rho_i\Big(D^{h-g}_i-\frac{\partial^N\lambda_i}{2}\Big)+\frac{a_i}{4}\big[\big(D^{h-2g}_i-\partial^N\lambda_i\big)^2 - (D_i^h)^2\big]\Big]. \label{eq_line_1_L_star_new}
\end{equation}
By definition of $N^2L^*_h{\bf 1}$, one has:
\begin{equation}
\nu^N_{g_h}\big(N^2L^*_h{\bf 1}\big) = \nu^N_{g_h}\big(N^2L_h{\bf 1}\big) = 0.
\end{equation}
We can also estimate $\nu^N_{g_h}(N^2L^*_h{\bf 1})$ through the expression~\eqref{eq_final_adjoint_boundary} of the adjoint at the boundary and the expansion~\eqref{eq_L_star_bulk_new_0} of the adjoint in the bulk. Indeed, Lemma~\ref{lemm_bound_correlations_bar_nu_G} yields the estimates:
\begin{align}
&\forall n\geq 1,\qquad 
\sup_{\substack{J\subset\Lambda_N \\ |J|=n}}\nu^N_{g_h}\Big(\prod_{j\in J}\bar\eta_j\Big) 
= 
O(N^{-n/2}),
\\
&\forall j\notin\{\pm(N-1)\},\qquad 
\nu^N_{g_h}(\bar\eta_{\pm (N-1)}\bar\eta_j) 
= O(N^{-2}),
\qquad 
\nu^N_{g_h}(\bar\eta_{\pm(N-1)}) = O(N^{-3/2})
.
\nonumber
\end{align}
These bounds can be used on the error term $\delta^N_{\pm}$ involving the adjoint at the boundary~\eqref{eq_final_adjoint_boundary}, 
the error term $\delta^N_{\text{order}\geq 3,0}$ defined in~\eqref{eq_def_delta_order_3}, 
the error term $\delta^N_{3-4,0}$ accounting for three and four point correlations defined in~\eqref{eq_line_4_L_star_new}--\eqref{eq_line_5_L_star_new}, 
the estimates~\eqref{eq_final_expression_fluct_part_L_star} of the fluctuations 
and~\eqref{eq_final_expression_correl_part_L_star} of the correlations. 
They yield:
\begin{align}
\nu^N_{g_h}\big(N^2L^*_h{\bf 1}\big) = 0 &= \text{ Const } +  \frac{1}{2N^2}\sum_{i<N-1}\sum_{j\notin\{i,i+1\}}\bar\sigma_i\bar\sigma_j\partial^N_1 g_{i,j}\partial^N_1(g-h)_{i,j} + O(N^{-1}) \nonumber\\
&=:\delta^N_{0,0} +O(N^{-1}).
\end{align}
The configuration-independent terms $\delta^{N}_{0,0}$ arising in $N^2L^*_{h,0}{\bf 1}$ are thus bounded by $O(N^{-1})$.
\subsection{Conclusion}\label{sec_conclusion_L_star}
Let us put together the estimates obtained so far to conclude the proof of Lemma~\ref{lemm_L_star_as_e_plus_carre_du_champ}. 
The expression of the adjoint at the boundary was obtained in~\eqref{eq_final_adjoint_boundary} for general $g$, 
while the adjoint in the bulk has been estimated in the last three sections, 
provided one takes $g = g_h$, with $g_h$ the solution of the main equation~\eqref{eq_main_equation}. 
One has therefore:
\begin{align}
N^2L^*_h{\bf 1}(\eta) 
&= 
N^2L^*_{h,\pm}{\bf 1}(\eta) + N^2L^*_{h,0}{\bf 1}(\eta) \nonumber\\
&= 
\frac{N}{2}\bar\eta_{-(N-1)}\Big(\partial^N\lambda_{-N}-\partial^N\lambda_{-(N-1)}\Big) -\frac{N}{2}\bar\eta_{N-1}\Big(\partial^N\lambda_{N-1}-\partial^N\lambda_{N-2} \Big) +\delta^N(\eta)\nonumber\\
&= 
\frac{\bar\eta_{-(N-1)}}{2}\Delta^N\lambda_{-(N-1)} -\frac{\bar\eta_{N-1}}{2}\Delta^N\lambda_{N-1} + \delta^N(\eta),
\end{align}
with $\delta^N$ a $\Gamma$-controllable error term with size $N^{-1/2}$, given by:
\begin{equation}
\delta^N(\eta) = \delta^N_{\pm}(\eta)+ \delta^N_{0,0} + \delta^N_{0,1}(\eta) + \delta^N_{0,2}(\eta) + \delta^N_{0,3-4}(\eta)+ \delta^N_{0,\,\text{order }3}(\eta) +\delta^N_{0,\,\text{order}\geq 4}(\eta). 
\end{equation}
Since $\sup_N|\Delta^N\lambda_{\pm(N-1)}|<\infty$, the quantity:
\begin{equation}
\frac{\bar\eta_{-(N-1)}}{2}\Delta^N\lambda_{-(N-1)} -\frac{\bar\eta_{N-1}}{2}\Delta^N\lambda_{N-1}
\end{equation}
is $\Gamma$-controllable with size $N^{-1}$ by Lemma~\ref{lemm_small_boundary_correl_sous_bar_nu_G}, 
and of vanishing type. 
It follows that $N^2L^{*}_h{\bf 1}$ is $\Gamma$-controllable with size $N^{-1/2}$, 
and equal to a sum of terms of vanishing type, 
plus $\delta^N_{0,3-4}$ which is of LS type.
Thus, by Definition~\ref{def_controllability} of $\Gamma$-controllability and taking $\delta = 1/2$ there, 
there is a controllable function $\e$ with size $N^{-1/2}$ and $C(\rho_\pm,h)$ satisfying, 
for some $\gamma>8C_{LS}$:
\begin{equation}
\frac{1}{\gamma}\log \nu^N_{g_h}\Big(\exp\big[\gamma|\e|\big]\Big)
\leq 
\frac{C(\rho_\pm,h)}{N^{1/2}}
,
\end{equation}
and for 
any density $f$ for $\nu^N_{g_h}$:
\begin{align}
\nu^N_{g_h}\big(fN^2L^*_h{\bf 1}\big)  
&\leq 
\nu^N_{g_h}(f\e) + \frac{N^2}{2}\nu^N_{g_h}\big(\Gamma_h(\sqrt{f})\big)\nonumber\\
&\leq \frac{H(f\nu^N_{g_h}|\nu^N_{g_h})}{8C_{LS}} + \frac{N^2}{2}\nu^N_{g_h}\big(\Gamma_h(\sqrt{f})\big) + \frac{C(\rho_\pm,h)}{N^{1/2}}.
\end{align}
This concludes the proof of Lemma~\ref{lemm_L_star_as_e_plus_carre_du_champ}.\emptydemo
\subsection{The Radon-Nikodym derivative}\label{sec_der_RD}
The computations in the previous subsections can be used to obtain an expression of the Radon-Nikodym derivative ${\rm D}_h = d\Prob_h/d\Prob$ ($h\in\s(\infty)$) on each fixed time interval. 
By definition, for trajectories up to time $T>0$, 
${\rm D}_h$ reads (see Appendix A.7 in~\cite{Kipnis1999}):
\begin{equation}
\log {\rm D}_h((\eta_t)_{t\leq T}) = \Pi^N_T(h)- \Pi^N_0(h) - N^2\int_0^Te^{-\Pi^N_t(h)}Le^{\Pi^N_t(h)}\, dt
.
\label{eq_def_der_radon_nikodym}
\end{equation}
The correlation field $\Pi^N$ is defined in~\eqref{eq_def_Pi}. 
\begin{proposition}\label{prop_expression_der_radon_nykodym}
Let $h\in\s(\epsilon_B)$, and recall that, 
for $u,v\in\mathbb L^2(\squaredash)$ and $(x,y)\in\squaredash$: 
\begin{equation}
\mathcal M(u,v)(x,y)
= 
\int_{(-1,1)}u(z,x)\bar\sigma(z)u(z,y)dz
.
\end{equation}
Then, for each $\eta\in\Omega_N$:
\begin{align}
N^2e^{-\Pi^N(h)}&Le^{\Pi^N(h)}(\eta) 
= 
\frac{1}{2}\Pi^N\Big(\Delta h+\mathcal M(\partial_1 h,\partial_1h)\Big) - \frac{(\bar\rho')^2}{4}\int_{(-1,1)}h(x,x)\, dx \label{eq_bulk_term_final}
\\
&+ \frac{1}{4}\sum_{i<N-1}\bar\eta_i\bar\eta_{i+1}\big(\partial_1 h_{i_+,i}-\partial_1h_{i_-,i}\big) + \frac{1}{8}\int_{\squaredash}\bar\sigma(x)\bar\sigma(y)\big[\partial_1 h(x,y)\big]^2dxdy + \epsilon^N(h),
\nonumber 
\end{align}
where $\epsilon^N(h)$ is a $\Gamma$-controllable error term with size $N^{-1/2}$ and of LS type. 
It thus satisfies by Corollary~\ref{coro_Boltzmann_gibbs_sec3}, 
for some $C(h,\rho_\pm)>0$, 
each large enough $N$ and each $T>0$:
\begin{equation}
\E^{\nu^N_{g_h}}_h\Big[\exp\Big[\int_0^T\epsilon^N_t(h)\, dt\Big]\Big] 
\leq 
\exp\Big[\frac{C(h,\rho_\pm)T}{N^{1/2}}\Big]
.
\end{equation}
\end{proposition}
\begin{remark}\label{remark_normal_derivative}
A bias $h\in\s(\epsilon_B)$ is a symmetric function by definition. 
As a result, for each $(x,y)\in\squaredash$:
\begin{equation}
\partial_1 h(x,y) = \partial_2 h(y,x)\quad \Rightarrow\quad \partial_1 h(x_+,x) - \partial_1 h(x_-,x) = (\partial_1-\partial_2)h(x_+,x).
\end{equation}
The first term in the second line of~\eqref{eq_bulk_term_final} thus corresponds to a contribution of the derivative of $h$ normal to the diagonal.
\demo
\end{remark}
The following corollary will be useful in the proof of lower bound large deviations.
\begin{corollary}\label{coro_comparaison_D_H_martingale}
Consider the $\Prob_h$-martingale $M^{N,\phi}$, defined for $T\geq 0$ and $\phi\in\mathcal T$ ($\mathcal T$ is defined in~\eqref{eq_def_T}) by:
\begin{equation}
M^{N,\phi}_T = \Pi^N_T(\phi) - \Pi^N_0(\phi) - \int_0^T N^2 L_h\Pi^N_t(\phi)dt.
\end{equation}
If additionally $\phi$ is a symmetric function in $C^3(\bar\rhd)$, 
there is a $\Gamma$-controllable error term $\tilde \epsilon^N(h,\phi)$ with size $N^{-1/2}$ such that, 
for any $T\geq 0$:
\begin{align}
M^{N,\phi}_T 
&= 
\Pi^N_T(\phi) - \Pi^N_0(\phi)-\frac{1}{2}\int_0^T\Pi^N_t\Big(\Delta \phi+ \mathcal M(\partial_1 \phi,\partial_1 h)+\mathcal M(\partial_1h,\partial_1\phi)\Big)dt
\nonumber\\
&\quad 
+ \frac{1}{4}\int_0^T\sum_{i<N-1}\bar\eta_i(t)\bar\eta_{i+1}(t)\big(\partial_1 \phi_{i_+,i}-\partial_1\phi_{i_-,i}\big)dt + \frac{(\bar\rho')^2T}{4}\int_{(-1,1)}\phi(x,x)\, dx
\nonumber\\
&\quad 
- \frac{T}{4}\int_{\squaredash}\bar\sigma(x)\bar\sigma(y)\partial_1 \phi(x,y)\partial_1 h(x,y)\, dx\, dy +\int_0^T\tilde{\epsilon}^N_t(h,\phi)\, dt
.
\end{align}
When $\phi=h$, one has in particular:
\begin{align}
M^{N,h}_T - \log {\rm D}_h((\eta_t)_{t\leq T}) 
&= 
-\frac{1}{2}\int_0^T\Pi^N_t\big(\mathcal M(\partial_1 h,\partial_1 h)\big)dt 
- \frac{T}{8}\int_{\squaredash}\bar\sigma(x)\bar\sigma(y)\big[\partial_1 h(x,y)\big]^2dxdy 
\nonumber\\
&\quad 
+ \int_0^T\hat\epsilon_t^N(h)dt
,
\end{align}
for a $\Gamma$-controllable error term $\hat \epsilon^N(h)$ with size $N^{-1/2}$ of LS type.
\end{corollary}

\section{Long-time behaviour: upper bound}\label{sec_large_devs}
In this section, we establish the upper bound in Theorem~\ref{theo_large_devs}: 
for the $\epsilon_B$ of Theorem~\ref{theo_entropic_problem} and any closed set $\mathcal F$ in $(\mathcal T'_s,*)$,
\begin{equation}
\limsup_{N\rightarrow\infty}\limsup_{T\rightarrow\infty}\frac{1}{T}\log \Prob^{\pi^N_{inv}}\Big(\frac{1}{T}\int_0^T\Pi^N_t\, dt\in \mathcal F\big)
\leq 
-\inf_{\Pi\in \mathcal F} \mathcal I_{\epsilon_B}(\Pi),\label{eq_upper_bound_in_sec_upper_bound}
\end{equation}
with $\mathcal I_{\epsilon_B}$ the functional defined in~\eqref{eq_def_rate_function}. 
A bound for compact sets is established in Section~\ref{sec_upper_bound_compact_sets}, 
relying on a regularity estimate in the space $\mathbb H^1(\squaredash)$, 
proven in Section~\ref{sec_regularity_estimate}. 
Exponential tightness, 
which yields the bound for closed sets, is obtained in Section~\ref{sec_exp_tightness}.\\

Before we start, 
let us make some remarks and fix notations. 
For $N\in\N^*,T>0$, a probability measure $\mu^N$ on the state space $\Omega_N$ and each measurable set $\mathcal B\subset \mathcal T'_s$, 
introduce the notation:
\begin{equation}
\Q^{\mu^N}_T\big(\mathcal B\big)
:=
\Prob^{\mu^N}\Big(\frac{1}{T}\int_0^T\Pi^N_t\, dt \in \mathcal B\Big)
.
\label{eq_def_Q_T}
\end{equation}
For short, we will also write:
\begin{equation}
\mathcal B^T 
= 
\mathcal B^T_N
:=
\Big\{ \frac{1}{T}\int_0^T\Pi^N_t\, dt\in\mathcal B\Big\}
.
\end{equation}
Recall that, 
for $N\in\N^*$, $\nu^N_{g_0}$ is the discrete Gaussian measure~\eqref{eq_def_bar_nu_g_intro_bis} built from the inverse correlation kernel $g_0$ of the steady state of the open SSEP.  Changing initial condition from $\pi^N_{inv}$ to the measure $\nu^N_{g_0}$ 
(defined in~\eqref{eq_def_nu_g_dans_methode_entropique}, 
with $g_0$ given by~\eqref{eq_def_g_0}) 
has a cost independent of $T$:
\begin{equation}
\forall T>0,\forall N\in\N^*,\qquad
\Q^{\pi^N_{inv}}_T(\mathcal B)
\leq 
\max_{\eta\in\Omega_N}\frac{\pi^N_{inv}(\eta)}{\nu^N_{g_0}(\eta)}\times \Q^{\nu^N_{g_0}}_T(\mathcal B)
.
\end{equation}
Upon taking $T^{-1}\log$ and the large $T$, 
then large $N$ limit, 
the initial condition of the dynamics does not change the value of the left-hand side in~\eqref{eq_upper_bound_in_sec_upper_bound}. 
The dynamics will therefore be started from the measure $\nu^N_{g_0}$.
\subsection{Upper bound for open and compact sets}\label{sec_upper_bound_compact_sets}
To estimate~\eqref{eq_def_Q_T}, 
we consider dynamics with tilted two-point correlations as in~\eqref{eq_def_jump_rates_H} and use the martingale method presented in Chapter 10 of~\cite{Kipnis1999}.  
It relies on the computation of the Radon-Nikodym derivatives $D_h = d\Prob_h/d\Prob$, 
$h\in\s(\epsilon_B)$. 
In Section~\ref{sec_a_first_upper_bound}, 
a first upper bound on compact sets with a rate function $\tilde{\mathcal I}_{\epsilon_B}\leq \mathcal I_{\epsilon_B}$ is established. 
The bound is then improved to $\mathcal I_{\epsilon_B}$ in Section~\ref{sec_refinement_upper_bound} with the help of the regularity estimate of Section~\ref{sec_regularity_estimate}.
\subsubsection{A first upper bound}\label{sec_a_first_upper_bound}
Here, we build a functional $\tilde{\mathcal I}_{\epsilon_B}:(\mathcal T'_s,*)\rightarrow\R_+$ such that, 
if $\mathcal{K}\subset (\mathcal T'_s,*)$ is a compact set,
\begin{equation}
\limsup_{N\rightarrow\infty}\limsup_{T\rightarrow\infty}\frac{1}{T}\log \Q_T^{\nu^N_{g_0}}(\mathcal K) \leq -\inf_{\mathcal K}\tilde {\mathcal I}_{\epsilon_B}.\label{eq_first_upper_bound}
\end{equation}
We first prove an upper bound for general Borel sets, 
then specify to compact sets. 
Let $h\in\s(\epsilon_B)$. In Proposition~\ref{prop_expression_der_radon_nykodym}, we proved that, for any $T>0$ and any trajectory $(\eta(t))_{t\leq T}$:
\begin{align}
\log {\rm D}_h\big((\eta(t)\big)_{t\leq T}) 
&= 
\Pi^N_T(h)-\Pi^N_0(h) - \frac{1}{2}\int_0^T\Pi^N_t\Big(\Delta h + \mathcal M(\partial_1h,\partial_1h)\Big)\, dt 
\label{eq_def_RD_sec_upper_bound}\\
&\hspace{-0.5cm}
- \frac{1}{4}\int_0^T\sum_{i<N-1}\bar\eta_i(t)\bar\eta_{i+1}(t)\big(\partial_1 h_{i_+,i}-\partial_1h_{i_-,i}\big) + \frac{(\bar\rho')^2T}{4}\int_{(-1,1)}h(x,x)\, dx
\nonumber\\
&\hspace{-0.5cm}
- \frac{T}{8}\int_{\squaredash}\bar\sigma(x)\bar\sigma(y)\big[\partial_1 h(x,y)\big]^2\, dx\, dy - \int_0^T\epsilon^N_t(h)\, dt
,
\nonumber
\end{align}
with $\epsilon^N(h)$ a $\Gamma$-controllable error term with size $N^{-1/2}$ of LS type 
(see Definitions~\ref{def_controllability}--\ref{def_LS_type}). 
For any Borel set $\mathcal B$ in $\mathcal T'_s$, 
one can thus write:
\begin{equation}\label{eq_change_measure_dynamical}
\Q^{\nu^N_{g_0}}_T(\mathcal B) 
=
\E^{\nu^N_{g_0}}_{h}\Big[{\bf 1}_{\mathcal{B}^T} (D_{h})^{-1}\Big],
\end{equation}
and the point is to build the functional $\tilde{\mathcal{I}}_{\epsilon_B}$ from $D_h^{-1}$ ($h\in\s(\epsilon_B)$). \\

\noindent\textbf{Closed expression.} 
Up to the error term $\epsilon^N(h)$, 
the expression~\eqref{eq_def_RD_sec_upper_bound} is nearly closed in terms of the distribution $\int_0^T \Pi^N_t dt$ applied to regular test functions. 
The only problematic term is the diagonal term involving the $\bar\eta_i\bar\eta_{i+1}$, $i<N-1$. 
Call it $W^{\text{Neu}}_h$:
\begin{equation}\label{eq_def_Neumann_term_sec4}
W^{\text{Neu}}_h(\eta) 
=
\frac{1}{4}\sum_{i<N-1}\bar\eta_i\bar\eta_{i+1}\big(\partial_1^N h_{i_+,i} - \partial^N_1h_{i_-,i}\big),
\qquad\eta\in\Omega_N
.
\end{equation}
In Appendix~\ref{sec_Neumann}, 
we estimate the cost of rewriting $W^{\text{Neu}}_h$ in terms of the correlation field $\Pi^N$ applied to a smooth test function supported in a small neighbourhood of the diagonal $D$. 
To state this estimate, 
consider a function $\chi^\epsilon\in C^\infty(\bar\square)$ with $\chi^\epsilon = 0$ on $\partial\square$, 
$0\leq \chi^\epsilon \leq 2/\epsilon$, 
such that $\chi^\epsilon(x,\cdot)$ is supported on $(x,x+\epsilon)\cap(-1,1)$ for each $x\in(-1,1)$, and:
\begin{equation}
\forall x<1-\epsilon,\qquad 
\int_{(x,x+\epsilon)}\chi^\epsilon(x,y)\, dy = 1,
\qquad\forall x\geq 1-\epsilon,\qquad 
\int_{(x,x+\epsilon)}\chi^\epsilon(x,y)\, dy \leq 1.\label{eq_mass_repartition_chi_epsilon}
\end{equation}
Define then $\mathcal N^\epsilon_h(x,y)$ for $(x,y)\in\squaredash$ as follows:
\begin{align}
\mathcal N_h^\epsilon(x,y) &= \frac{\bar\sigma(x)}{\bar\sigma(y)}\chi^\epsilon(x,y)\big(\partial_1-\partial_2)h(x_+,x).\label{eq_def_mathcal_N_h_epsilon}
\end{align}
Then $\mathcal N_h^\epsilon\in \mathcal T$ 
(defined in~\eqref{eq_def_T}), 
and we prove the following in Proposition~\ref{prop_neumann_condition_averaging}. 
Let $\theta>0$ and 
$\epsilon\in(0,1)$ be smaller than some $\epsilon_0(\rho_\pm,h,\theta)>0$. 
There are constants $C_1(\rho_\pm,h,\theta),C_2(\rho_\pm,h,\epsilon,\theta)>0$ such that, 
for each $N$ large enough depending on $\epsilon,\theta$ and each $T>0$:
\begin{align}
&\E^{\nu^N_{g_h}}_h\bigg[\exp\Big[\theta\int_0^T \big[W^{\text{Neu}}_h(\eta_t) - \Pi^N_t\big(\mathcal N^h_\epsilon\big)\big]dt\Big] \bigg] 
\nonumber\\
&\hspace{3cm}\leq
\exp\Big[C_1(h,\rho_{\pm},\theta)\epsilon^{1/2} T + \frac{C_2(h,\rho_{\pm},\epsilon,\theta)T}{N^{1/2}}\Big]
.
\label{eq_estimate_replacement_neumann_sec4}
\end{align}
The same bound is valid starting from $\nu^N_{g_0}$ up to an additional $e^{C(\rho_\pm,h) N}$ factor in the right-hand side (in fact only a factor $e^{C(\rho_\pm,h)}$, 
i.e. bounded with $N$, 
but this does not matter for our purposes). 
Since we take $T$ large first, this $e^{C(\rho_\pm,h) N}$ multiplicative constant is not a problem. 
 
Introduce then the following continuous functional $J_h^\epsilon$ on $(\mathcal T'_s,*)$:
\begin{align}
\forall \Pi\in\mathcal T'_s,\qquad J^\epsilon_h(\Pi) 
&= 
-\frac{1}{2}\Pi\Big(\Delta h + \mathcal M(\partial_1h,\partial_1 h)\Big)
-\Pi\big(\mathcal N^\epsilon_h\big) \nonumber\\
&\quad + 
\frac{(\bar\rho')^2}{4}\int_{(-1,1)} h(x,x)\, dx
- \frac{1}{8}\int_{\squaredash}\bar\sigma(x)\bar\sigma(y)\big[\partial_1 h(x,y)\big]^2\, dx\, dy.\label{eq_def_J_h_avec_epsilon}
\end{align}
With this definition, 
the Radon-Nikodym derivative~\eqref{eq_def_RD_sec_upper_bound} becomes:
\begin{align}
\log {\rm D}_h\big((\eta(t)\big)_{t\leq T}) 
&= 
\Pi^N_T(h)-\Pi^N_0(h) + TJ^\epsilon_h\Big(\frac{1}{T}\int_0^T\Pi^N_t\, dt\Big) \nonumber\\
&\quad 
- \int_0^T \big[W^{\text{Neu}}_h(\eta_t) - \Pi^N_t\big(\mathcal N^h_\epsilon\big)\big]\, dt - \int_0^T\epsilon^N_t(h)\, dt
.
\end{align}
\noindent\textbf{Estimates on error terms.} 
Let us estimate the contribution of $\epsilon^N(h) + W^{\text{Neu}}_h(\eta_t) - \Pi^N_t\big(\mathcal N^h_\epsilon\big)$.  
By Proposition~\ref{prop_expression_der_radon_nykodym}, 
$\epsilon^N(h)$ will not contribute to the large deviations as it has small exponential moment. 
Indeed, Proposition~\ref{prop_expression_der_radon_nykodym} 
together with the last item of Lemma~\ref{lemm_type} give that, for some small enough $\alpha>0$:
\begin{equation}\label{eq_size_exp_mom_epsilon_sec_4}
\limsup_{N\rightarrow\infty}\limsup_{T\rightarrow\infty}\frac{1}{T}\log \Prob^{\nu^N_{g_h}}\Big[\exp\Big[(1+\alpha)\int_0^T\epsilon^N_t(h)\, dt\Big]\Big] 
= 0
.
\end{equation}
Since a change of initial condition has no influence in the large $T$ limit, 
the above estimate is also valid starting from $\nu^N_{g_0}$. 
Thus, by Hölder inequality:
\begin{align}
\frac{1}{T}\log &\, \Prob_h^{\nu^N_{g_0}}\bigg[\exp\Big[\int_0^T \big[W^{\text{Neu}}_h(\eta_t) - \Pi^N_t\big(\mathcal N^h_\epsilon\big)\big]dt + \int_0^T\epsilon^N_t(h)dt\bigg]\nonumber\\
&\quad\leq
\frac{\alpha}{(1+\alpha)T}\log \Prob_h^{\nu^N_{g_0}}\bigg[\exp\Big[\int_0^T \frac{(1+\alpha)}{\alpha}\big[W^{\text{Neu}}_h(\eta_t)- \Pi^N_t\big(\mathcal N^h_\epsilon\big)\big]dt\Big]\bigg] \nonumber\\
&\qquad+ \frac{1}{(1+\alpha)T}\log \Prob^{\nu^N_{g_0}}\bigg[\exp\Big[ \int_0^T(1+\alpha)\epsilon^N_t(h)dt\bigg].
\end{align}
Combining~\eqref{eq_estimate_replacement_neumann_sec4}--\eqref{eq_size_exp_mom_epsilon_sec_4}, 
we find that there is a constant $C= C(h,\rho_\pm,\alpha)>0$ such that:
\begin{equation}\label{eq_estimate_epsilon_plus_Neumann}
\limsup_{N\rightarrow\infty}\limsup_{T\rightarrow\infty}
\frac{1}{T}\log \Prob^{\nu^N_{g_0}}\Big[\exp\Big[\int_0^T \big[W^{\text{Neu}}_h(\eta_t) - \Pi^N_t\big(\mathcal N^h_\epsilon\big)\big]dt + \int_0^T\epsilon^N_t(h)\, dt\Big]
\leq
C\epsilon^{1/2}
.
\end{equation}
We can now try to obtain the upper bound~\eqref{eq_first_upper_bound}. 
Take a Borel set $\mathcal B$. 
Starting from~\eqref{eq_change_measure_dynamical}, 
one has:
\begin{align}\label{eq_borne_sup_1}
&\frac{1}{T}\log\Q_T^{\nu^N_{g_0}}(\mathcal B) 
=
\frac{1}{T}\log\E^{\nu^N_{g_0}}_{h}\Big[{\bf 1}_{\mathcal{B}^T} (D_{h})^{-1}\Big]\\
&\quad\leq
\frac{1}{T}\log\E^{\nu^N_{g_0}}_h\bigg[
\exp\Big[-\Pi^N_T(h)+\Pi^N_0(h)+\int_0^T \big[W^{\text{Neu}}_h(\eta_t) - \Pi^N_t\big(\mathcal N^h_\epsilon\big)\big]dt + \int_0^T\epsilon^N_t(h)dt\Big]
\bigg]\nonumber
\\&\hspace{10cm}+
\sup_{\Pi\in \mathcal B}\big(-J^\epsilon_h(\Pi)\big)\nonumber\\
&\quad\leq 
\frac{C(h)N}{T} + \frac{1}{T}\log\E^{\nu^N_{g_0}}_h\bigg[
\exp\Big[\int_0^T \big[W^{\text{Neu}}_h(\eta_t) - \Pi^N_t\big(\mathcal N^h_\epsilon\big)\big]dt + \int_0^T\epsilon^N_t(h)dt\Big]
\bigg]
\nonumber\\
&\hspace{10cm}+ 
\sup_{\Pi\in \mathcal B}\big(-J^\epsilon_h(\Pi)\big)
.\nonumber
\end{align}
Taking the large $T$, then large $N$ limits and estimating the error terms via~\eqref{eq_estimate_epsilon_plus_Neumann}, 
one finds that,
for all $\epsilon>0$ smaller than some $\epsilon_0=\epsilon_0(\rho_\pm,h,\alpha)$:
\begin{equation}
\limsup_{N\rightarrow\infty}\limsup_{T\rightarrow\infty}\frac{1}{T}\log\Q_T^{\nu^N_{g_0}}(\mathcal B) 
\leq 
C(h,\rho_\pm)\epsilon^{1/2} + \sup_{\Pi\in \mathcal B}\big(-J^\epsilon_h(\Pi)\big)
.
\end{equation}
Taking the infimum on $\epsilon\in(0,\epsilon_0)$ and $h\in \s(\epsilon_B)$ yields a first bound:
\begin{align}\label{eq_borne_sup_2}
\limsup_{N\rightarrow\infty}\limsup_{T\rightarrow\infty}\frac{1}{T}\log\Q_T^{\nu^N_{g_0}}(\mathcal B)
\leq 
\inf_{h\in\s(\epsilon_B)}\inf_{\epsilon\in(0,\epsilon_0)}\sup_{\Pi\in \mathcal B}\big(-J^\epsilon_h(\Pi)\big)
.
\end{align}
We now work out a way to 
exchange supremum and infima in~\eqref{eq_borne_sup_2} when $\mathcal B$ is a compact set. 
The argument is standard and relies on Lemmas 3.2 and 3.3 in Appendix 2 of~\cite{Kipnis1999}. 
Let $\mathcal K\subset (\mathcal{T}',*)$ be compact. 
We wish to prove:
\begin{equation}
\limsup_{N\rightarrow\infty}\limsup_{T\rightarrow\infty}\frac{1}{T}\log\Q_T^{\nu^N_{g_0}}(\mathcal K)
\leq 
\sup_{\mathcal K}\inf_{h\in\s(\epsilon_B)}\liminf_{\epsilon\rightarrow0}\big(-J^\epsilon_h\big).\label{eq_borne_sup_compact_1}
\end{equation}
Since $(J^\epsilon_{h})_{h,\epsilon}$ is a family of continuous functionals on $(\mathcal T'_s,*)$, 
Lemmas 3.2 and 3.3 in Appendix 2 of~\cite{Kipnis1999} allow for the exchange of the infima on $h,\epsilon$ and the supremum on (open covers of) $\mathcal K$:
\begin{align}
\limsup_{T\rightarrow\infty}\frac{1}{T}\log\Q_T^{\nu^N_{g_0}}(\mathcal K) 
\leq
\sup_{\mathcal K}  \inf_{h\in\s(\epsilon_B)}\inf_{\epsilon\in(0,\epsilon_0)}\big(-J^\epsilon_{h}\big)
.
\label{eq_borne_sup_3}
\end{align}
Since $\inf_{\epsilon<\epsilon_0}$ increases when $\epsilon_0$ shrinks,
~\eqref{eq_borne_sup_3} yields:
\begin{align}
\limsup_{N\rightarrow\infty}\limsup_{T\rightarrow\infty}\frac{1}{T}\log\Q_T^{\nu^N_{g_0}}(\mathcal K) 
&\leq
\sup_{\mathcal K}  \inf_{h\in\s(\epsilon_B)}\liminf_{\epsilon\rightarrow 0}\big(-J_h^\epsilon\big)
.
\end{align}
This yields a first bound on compact sets and proves~\eqref{eq_first_upper_bound}:
\begin{equation}
\limsup_{N\rightarrow\infty}\limsup_{T\rightarrow\infty}\frac{1}{T}\log\Q_T^{\nu^N_{g_0}}(\mathcal K)
\leq 
-\inf_{\mathcal K}\tilde{\mathcal I}_{\epsilon_B},\qquad \tilde{\mathcal I}_{\epsilon_B} = \sup_{h\in\s(\epsilon_B)}\limsup_{\epsilon\rightarrow 0} J_h^\epsilon
.
\label{eq_first_bound_compact_set}
\end{equation}
\subsubsection{Refinement of the upper bound to restrict to more regular correlations}\label{sec_refinement_upper_bound}
To improve the expression for the functional $\tilde{\mathcal{I}}_{\epsilon_B}$, 
we would like to take the limit in $\epsilon$ in~\eqref{eq_first_bound_compact_set}. 
This is however not possible in general. 
Indeed, recall that $\Pi\in \mathcal T'_s$ is a distribution, 
and taking $\epsilon$ to $0$ amounts to asking for $\Pi$ to have a well-defined trace on the diagonal $D$ of the square $\squaredash$. 
This is possible only if $\Pi$ has some regularity.\\
In this section, 
we explain how to improve the upper bound~\eqref{eq_first_bound_compact_set} so that it is finite only on correlation fields $\Pi$ that have a well defined trace. 
More precisely, 
write $\Pi = \frac{1}{4}\big<k_\Pi,\cdot\big>\in\mathcal T'_s$. 
We prove that one can restrict to $\Pi$'s with $k_\Pi\in\mathbb H^1(\squaredash)$.
Note that, as $\Pi\in \mathcal T'_s\cap\mathbb H^1(\squaredash)$, 
$k_\Pi$ is a symmetric function, 
and the traces on both sides of the diagonal coincide. 
We thus write $\text{tr}_D(k_\Pi)$ for the trace of $k_\Pi$ on the diagonal without ambiguity. 
By definition of the trace, one has then:
\begin{equation}
\lim_{\epsilon\rightarrow 0}\Pi\big(\mathcal N^\epsilon_h\big) = \lim_{\epsilon\rightarrow 0}\frac{1}{4}\big<k_\Pi,\mathcal N^\epsilon_h\big> = \frac{1}{4}\int_{(-1,1)}\text{tr}(k_\Pi)(x,x) (\partial_2-\partial_1)h(x_+,x) dx.\label{eq_limite_mathcal_N_converges_to_trace}
\end{equation}
Thus, for $\Pi\in\mathcal T'_s\cap\mathbb H^1(\squaredash)$, 
$\limsup_{\epsilon\rightarrow 0}|J_h^\epsilon(\Pi)- J_h(\Pi)|=0$, 
with $J_h$ the functional equal to $+\infty$ outside of $\mathcal T'_s\cap \mathbb H^1(\squaredash)$, 
and:
\begin{align}
\forall \Pi\in\mathcal T'_s\cap \mathbb H^1(\squaredash),\qquad J_h(\Pi) 
&= 
-\frac{1}{2}\Pi\Big(\Delta h + \mathcal M(\partial_1h,\partial_1 h)\Big)
\nonumber\\
&\quad+
\frac{1}{4}\int_{(-1,1)}\text{tr}_{D}(k_\Pi)(\partial_2-\partial_1)h(x_+,x)dx
\label{eq_def_J_h_dans_upper_bound}\\
&\quad + \frac{(\bar\rho')^2}{4}\int_{(-1,1)} h(x,x)dx
- \frac{1}{8}\int_{\squaredash}\bar\sigma(x)\bar\sigma(y)\big[\partial_1 h(x,y)\big]^2dxdy
.
\nonumber
\end{align}
As a result, for $\Pi\in \mathcal T'_s\cap \mathbb H^1(\squaredash)$, 
$\tilde{\mathcal{I}}_{\epsilon_B}(\Pi) = \mathcal{I}_{\epsilon_B}(\Pi)$, 
with $\tilde {\mathcal{I}}_{\epsilon_B}$ defined in~\eqref{eq_first_bound_compact_set}, 
and $\mathcal I_{\epsilon_B}\geq \tilde{\mathcal{I}}_{\epsilon_B}$ the larger rate function defined in~\eqref{eq_def_rate_function}. 
The goal of the section is summarised in the following lemma.
\begin{lemma}\label{lemm_upper_bound_compact_sets_true}
Let $\mathcal K\subset (\mathcal T'_s,*)$ be a compact set. Then:
\begin{equation}
\limsup_{N\rightarrow\infty}\limsup_{T\rightarrow\infty}\frac{1}{T}\log\Q_T^{\nu^N_{g_0}}(\mathcal K)\leq -\inf_{\mathcal K}\mathcal I_{\epsilon_B}.\label{eq_better_upper_bound_I}
\end{equation}
\end{lemma}
To prove Lemma~\ref{lemm_upper_bound_compact_sets_true}, consider the functionals $\mathcal Q$ and $\mathcal Q_\phi$, defined for each $\phi\neq 0$ in $C^\infty_{c}(\squaredash)$, the set of compactly supported, $C^\infty$ functions on $\squaredash$, by:
\begin{equation}
\forall \Pi\in \mathcal T'_s,\qquad \mathcal Q_\phi(\Pi) = \frac{\Pi(\partial_1 \phi)}{\|\phi\|_2},\qquad \mathcal Q = \sup_{\phi\in C^\infty_{c}(\squaredash)\setminus\{0\}}\mathcal Q_\phi,\label{eq_def_regularity_functional}
\end{equation}
Note that $\mathcal Q$ is weak$^*$ lower semi-continuous on $(\mathcal T'_s,*)$. 
Indeed, it is a supremum over the $\mathcal Q_\phi$ for $\phi\in C^\infty_c(\squaredash)\setminus\{0\}$, 
and $\mathcal Q_\phi$ is weak$^*$ continuous since 
$\|\phi\|_2\mathcal Q_\phi(\Pi)$ is the evaluation of $\Pi$ at $\partial_1 \phi\in\mathcal T$ for $\Pi\in\mathcal T'_s$.\\
It is classical that $\mathcal Q$ controls the regularity of elements of $\mathcal T'_s$, 
as stated in the next lemma.
\begin{lemma}\label{lemm_properties_energie}
For $\Pi\in\mathcal T'_s$, $\Pi$ is in fact a function in $\mathbb H^1(\squaredash)$ if and only if $\mathcal Q(\Pi)<\infty$. 
\end{lemma}
\begin{proof}[Proof of Lemma~\ref{lemm_upper_bound_compact_sets_true}]
We now begin the proof of the large deviation bound~\eqref{eq_better_upper_bound_I}. 
Consider a sequence $\phi_j\in C_{c}^\infty(\squaredash)\setminus\{0\},j\in\N^*$, 
dense for the norm of $\mathbb H^3(\squaredash)$. 
Introduce then, for each $\ell\in\N^*$ and each $A>0$:
\begin{equation}
U(\ell,A) = \Big\{\max_{1\leq j \leq \ell}\mathcal Q_{\phi_j} \leq A\Big\}.\label{eq_def_U_ell_lambda_controle_energie}
\end{equation}
In Proposition~\ref{prop_energy_estimate}, 
we prove the existence of $C=C(\rho_\pm)>0$ such that, 
for $A$ larger than some $A_0>0$ and each $\ell\in\N^*$:
\begin{equation}
\limsup_{N\rightarrow\infty}\limsup_{T\rightarrow\infty}\frac{1}{T}\log\Q_T^{\nu^N_{g_0}}\big((U(\ell,A))^c\big)
\leq 
-CA
.\label{eq_energy_estimate_sec_upper_bound_compact_sets}
\end{equation}
Notice also that, for $\Pi\in\mathcal T'_s$, if $(\phi_{j_n})$ converges to $\phi\in C^\infty_c(\squaredash)$ in the norm of $\mathbb H^3(\squaredash)$, 
then $\partial_1 \phi_{j_n}$ converges to $\partial_1 \phi$ in $\mathcal T$, so that $\lim_n \mathcal Q_{\phi_{j_n}}(\Pi) = \mathcal Q_{\phi}(\Pi)$, and:
\begin{equation}
\mathcal Q(\Pi) = \sup_{j\in\N^*}\mathcal Q_{\phi_j}(\Pi).\label{eq_energie_is_sup_sur_dense_sequence}
\end{equation}
The regularity of the correlation fields can therefore be controlled through the $\mathcal Q_{\phi_j}$ ($j\in\N^*$). 
Let $\ell\in\N^*$ and $A>A_0$. 
Recall the notation: 
\begin{equation}
U(\ell,A)^T := \Big\{\frac{1}{T}\int_0^T\Pi_t\, dt \in U(\ell,A)\Big\},
\end{equation}
so that $\Prob(U(\ell,A)^T) = \Q_T(U(\ell,A))$. 
Let $\mathcal B\subset (\mathcal T'_s,*)$ be a Borel set. 
For ease of writing, 
let us abbreviate limits in $N,T$ as follows:
\begin{equation}
\limsup_{N,T\rightarrow\infty} 
:= 
\limsup_{N\rightarrow\infty}\limsup_{T\rightarrow\infty}
.
\end{equation}
We again estimate $\Q_T^{\nu^N_{g_0}}(\mathcal B)$, 
starting from the bound:
\begin{align}
\limsup_{N,T\rightarrow\infty}\frac{1}{T}\log \Q_T^{\nu^N_{g_0}}(\mathcal B) 
&\leq
\max \Big\{\limsup_{N,T\rightarrow\infty}\frac{1}{T}\log \Q_T^{\nu^N_{g_0}}\big(\mathcal B\cap U(\ell,A)\big), 
\nonumber\\
&\hspace{4cm}\limsup_{N,T\rightarrow\infty}\frac{1}{T}\log \Q_T^{\nu^N_{g_0}}\big(U(\ell,A)^c\big)\Big\}
\nonumber\\
&\leq 
\max \Big\{\limsup_{N,T\rightarrow\infty}\frac{1}{T}\log \Q_T^{\nu^N_{g_0}}\big(\mathcal B\cap U(\ell,A)\big), 
-CA\Big\}
,
\end{align}
where~\eqref{eq_energy_estimate_sec_upper_bound_compact_sets} was used to obtain the second inequality.
Proceeding as in~\eqref{eq_borne_sup_1}--\eqref{eq_borne_sup_2} to estimate the first probability, 
we find:
\begin{align}\label{eq_borne_sup_4}
\limsup_{N,T\rightarrow\infty}\frac{1}{T}\log\Q_T^{\nu^N_{g_0}}(\mathcal B)
\leq 
\max\Big\{C(h,\rho_\pm)\epsilon
+ \sup_{\Pi\in \mathcal B\cap U(\ell,A)}\big(-J^\epsilon_h(\Pi)\big), -CA\Big\}
.
\end{align}
For each admissible $h,\epsilon,\ell,A$, 
let $J_{h,\ell,A}^\epsilon$ be equal to $+\infty$ in $U(\ell,A)^c$, 
and on $U(\ell,A)$:
\begin{equation}
 J_{h,\ell,A}^\epsilon := 
\max\Big\{C(h,\rho_\pm)\epsilon
- J^\epsilon_h(\Pi), -CA\Big\}.
\end{equation}
Minimising~\eqref{eq_borne_sup_4} on $\epsilon\in(0,\epsilon_0)$, $h\in\s(\epsilon_B)$, $\ell\in\N^*$ and $A>A_0$, 
we have therefore obtained the upper bound:
\begin{equation}\label{eq_upper_bound_5}
\limsup_{N,T\rightarrow\infty}\frac{1}{T}\log\Q_T^{\nu^N_{g_0}}(\mathcal B)
\leq 
\inf_{A,\ell,h,\epsilon} \sup_{\mathcal B}\bar J_{h,\ell,A}^\epsilon
.
\end{equation}
Let us obtain a bound on compact sets from~\eqref{eq_upper_bound_5}. 
Let $\mathcal K\subset (\mathcal T'_s,*)$ be compact. 
Since $U(\ell,A)$ is weak$^*$ closed by continuity of each $\mathcal Q_{\phi_j}$ 
($1\leq j \leq \ell$), 
$(J_{h,\ell,A}^\epsilon)_{h,\epsilon,\ell,A}$ is a family of weak$^*$ upper semi-continuous functionals. 
Lemmas A.2.3.2 and A.2.3.3 in~\cite{Kipnis1999} thus give as before:
\begin{align}
\limsup_{N,T\rightarrow\infty}\frac{1}{T}\log\Q_T^{\nu^N_{g_0}}(\mathcal K) 
\leq 
\sup_{\mathcal K}\inf_{h\in\s(\epsilon_B)}\inf_{\ell,A,\epsilon}J^\epsilon_{h,\ell,A},\label{eq_borne_sup_en_cours_estim_energie_0}
\end{align}
with the infimum still taken on $\ell\in\N^*$, $A>A_0$ and $\epsilon\in(0,\epsilon_0)$.\\
One can again bound the infimum on $\epsilon$ by a liminf:
\begin{align}\label{eq_borne_sup_en_cours_estim_energie_0_bis}
\limsup_{N,T\rightarrow\infty}\frac{1}{T}\log\Q_T^{\nu^N_{g_0}}(\mathcal K) 
\leq 
\sup_{\mathcal K}\inf_{h\in\s(\epsilon_B)}\inf_{\ell,A}\max\big\{\liminf_{\epsilon\rightarrow 0}(-J^\epsilon_{h,\ell,A}),-CA\big\}
.
\end{align}
Let $A>A_0$. As $U(\ell,A)\subset U(\ell',A)$ if $\ell\leq \ell'$, 
the argument of the supremum on $\mathcal K$ in~\eqref{eq_borne_sup_en_cours_estim_energie_0_bis} is equal to $-\infty$ when evaluated at any $\Pi\notin \bigcap_{\ell\in\N^*}U(\ell,A)$. 
By definition of $U(\ell,A)$ in~\eqref{eq_def_U_ell_lambda_controle_energie} and using~\eqref{eq_energie_is_sup_sur_dense_sequence}, 
one has:
\begin{equation}
\bigcap_{\ell\in\N^*}U(\ell,A) = \{\mathcal Q\leq A\}.
\end{equation}
Equation~\eqref{eq_borne_sup_en_cours_estim_energie_0_bis} thus becomes:
\begin{align}
\limsup_{N,T\rightarrow\infty}\frac{1}{T}\log \Q_T^{\nu^N_{g_0}}(\mathcal K) \leq \sup_{\mathcal K}\inf_{h\in\s(\epsilon_B)}\inf_{A>A_0}\max\big\{\liminf_{\epsilon\rightarrow0}(-J^\epsilon_{h,A}),-CA\big\},\label{eq_borne_sup_en_cours_estim_energie_1}
\end{align}
with, for each $h,\epsilon,A$, 
$J^\epsilon_{h,A} = +\infty$ on $\{\mathcal Q>A\}$, 
and $J^\epsilon_{h,A} = J_h^\epsilon$ on $\{\mathcal Q\leq A\}$. 
Consider again $A>A_0$. 
For each $\Pi\in \{\mathcal Q\leq A\}$, 
the $k_\Pi$ associated with $\Pi$ via $\Pi(\cdot)=\frac{1}{4}\big<k_\Pi,\cdot\big>$ belongs to $\mathbb H^{1}(\squaredash)$ by Lemma~\ref{lemm_properties_energie}. 
In particular, by~\eqref{eq_limite_mathcal_N_converges_to_trace}, if $h\in\s(\epsilon_B)$ and $\Pi\in\{\mathcal Q\leq A\}$,
\begin{equation}
\liminf_{\epsilon\rightarrow 0}(-J^\epsilon_{h,A}(\Pi)) = -J_{h,A}(\Pi).
\end{equation}
Above, $J_{h,A} = J_h$ on $\{\mathcal Q\leq A\}$, $J_{h,A} = +\infty$ outside, and $J_h$ is defined in~\eqref{eq_def_J_h_dans_upper_bound}. Equation~\eqref{eq_borne_sup_en_cours_estim_energie_1} thus becomes:
\begin{align}
\limsup_{N,T\rightarrow\infty}\frac{1}{T}\log \Q_T^{\nu^N_{g_0}}(\mathcal K) \leq \sup_{\mathcal K}\inf_{h\in\s(\epsilon_B)}\inf_{A>A_0}\max\big\{-J_{h,A},-CA\big\}.
\end{align}
Finally, note that $J_{h,A}\geq J_h$ for $A>A_0$, 
since $J_h$ may be finite on $\{\mathcal Q>A\}$ while $J_{h,A}$ may not. 
Lemma~\ref{lemm_upper_bound_compact_sets_true} is proven:
\begin{equation}
\limsup_{N\rightarrow\infty}\limsup_{T\rightarrow\infty}\frac{1}{T}\log \Q_T^{\nu^N_{g_0}}(\mathcal K)
\leq 
-\inf_{\mathcal K} \mathcal I_{\epsilon_B},
\qquad \mathcal I_{\epsilon_B} := \sup_{h\in\s(\epsilon_B)}J_h
.\label{eq_upper_bound_compact sets}
\end{equation}
\end{proof}
\subsection{Regularity estimate}\label{sec_regularity_estimate}
In this section, we prove the energy estimate~\eqref{eq_energy_estimate_sec_upper_bound_compact_sets}. 
The key argument is the following proposition.
\begin{proposition}\label{prop_energy_estimate}
Let $0<\rho_-\leq \rho_+<1$ and assume $\bar\rho'\leq\epsilon_B$, 
with $\epsilon_B$ given by Theorem~\ref{theo_entropic_problem}.  Let $A>0$ be larger than some fixed $A_0>0$ and let $\phi \in C^2_{c}(\squaredash)$, 
where the subscript $c$ stands for compactly supported. 
There are constants $C_1(\rho_\pm),C_2(\rho_\pm)>0$ and $C_3(\rho_\pm,\phi)$ such that,
for each $T>0$, 
each $A>0$ and each $N$ larger than some $N_0(\phi)$:
\begin{align}
&\Prob^{\nu^N_{g_0}}\Big(\Big|\frac{1}{T}\int_0^T\Pi^N_t(\partial_1 \phi)dt\Big|> A\|\phi\|_2\Big) 
\nonumber\\
&\hspace{3cm}\leq 
2\exp\Big[-C_1(\rho_\pm)\big(A - C_2(\rho_\pm)\big) T  + \frac{C_3(\rho_\pm,\phi)T}{N^{1/2}}\Big]
.
\label{eq_controle_energie}
\end{align}
In particular,
\begin{equation}
\limsup_{N\rightarrow\infty}\limsup_{T\rightarrow\infty}\frac{1}{T}\log \Prob^{\nu^N_{g_0}}\Big(\Big|\frac{1}{T}\int_0^T\Pi^N_t(\partial_1 \phi)\, dt\Big|> A\|\phi\|_2\Big) 
\leq 
- C_1(\rho_\pm)\big(A - C_2(\rho_\pm)\big).
\end{equation}
\end{proposition}
Assuming the proposition,~\eqref{eq_energy_estimate_sec_upper_bound_compact_sets} immediately follows by a union bound and the appropriate limits.
\begin{proof}[Proof of Proposition~\ref{prop_energy_estimate}]
Up to considering $-\phi$, 
it is enough to prove the result without the absolute value and the factor $2$ in front of the probability. 
Up to taking $N$ large enough depending on $\phi$, 
we can assume that the support of $\phi$ is contained in $\{ z\in \squaredash:d(z,\partial\ \squaredash)>2/N\}$, 
so that $\phi_{i,i+1} = 0 = \partial_1 \phi_{i,i+1}$ for each $i<N-1$. 
We may also assume without loss of generality that $\phi$ is symmetric owing to the identity $\Pi^N(\phi) = \Pi^N(\check\phi)$, 
with $\check\phi(x,y) = \phi(y,x)$ for $x,y\in\squaredash$. 
The idea is to use Feynman-Kac inequality (Lemma~\ref{lemm_FK}) and a microscopic integration by parts to rewrite,  
for each density $f$ for $\nu^N_{g_0}$, 
the average $\nu^N_{g_0}\big(f\|\phi\|_2^{-1}\Pi^N(\partial_1 \phi)\big)$ as $\nu^N_{g_0}\big(f\Pi^N(F(\|\phi_2\|_2^{-1}\phi))\big)$ plus a term involving the carré du champ, 
for some function $F$.
The term $F(\|\phi\|_2^{-1}\phi)$ now involves only $\|\phi\|_2^{-1}\phi$,
and not its first partial derivative. 
The average $\nu^N_{g_0}\big(f\Pi^N(F(\|\phi_2\|_2^{-1}\phi))\big)$ is then estimated through the entropy inequality.\\

By Feynman-Kac inequality~\eqref{eq_FK_sous_mu_general} together with the bound~\eqref{eq_bound_adjoint_lemm_relative_entropy} on the adjoint, 
one has, 
for some $C=C(\rho_\pm)>0$ and some $\kappa>0$ to be chosen later:
\begin{align}
&\frac{1}{T}\log \Prob^{\nu^N_{g_0}}\Big(\frac{1}{T}\int_0^T\|\phi\|_2^{-1}\Pi^N_t(\partial_1 \phi)\, dt>  A\Big) 
\label{eq_FK_energy}\\
&\hspace{1cm}\leq 
-A\kappa
+ \sup_{f\geq 0 :\nu^N_{g_0}(f)=1}
\Big\{
\kappa\|\phi\|_2^{-1} \nu^N_{g_0}\big(f\Pi^N(\partial_1\phi)\big)  -\frac{N^2}{8}\nu^N_{g_0}\big(\Gamma(f^{1/2})\big)
+\frac{C(\rho_\pm)}{N^{1/2}}
\Big\}
.
\nonumber
\end{align}
To estimate the right-hand side,
let us write out $\|\phi\|_2^{-1}\Pi^N(\partial_1 \phi)$. 
Fix $\eta\in\Omega_N$.
\begin{align}
\|\phi\|_2^{-1}\Pi^N(\partial_1 \phi) 
&= 
\frac{1}{4\|\phi\|_2 N}\sum_{i<N-1}\bigg[\sum_{j\notin\{i,i+1\}}\bar\eta_i\bar\eta_j\partial^N_1 \phi_{i,j} +\bar\eta_i\bar\eta_{i+1}\partial^N_1\phi_{i,i+1}\bigg] 
\nonumber\\
&\quad
+\frac{1}{N\|\phi\|_2}\Pi^N(b)
\label{eq_discretisation_Pi_partial_1_phi}
,
\end{align}
where $b = N[\partial_1\phi - \partial^N_1\phi]$ is the discretisation error, which is bounded in $N$. 
By assumption on the compact support of $\phi$ in $\squaredash$, $\partial_1 \phi_{i,i+1} = 0$ for each $i<N-1$, 
and~\eqref{eq_discretisation_Pi_partial_1_phi} becomes:
\begin{align}
\|\phi\|_2^{-1}\Pi^N(\partial_1 \phi) 
=
\frac{1}{4\|\phi\|_2 N}\sum_{i<N-1}\sum_{j\notin\{i,i+1\}}\bar\eta_i\bar\eta_j\partial^N_1 \phi_{i,j} +\frac{1}{N\|\phi\|_2}\Pi^N(b)
.
\end{align}
By Lemma~\ref{lemm_size_controllable}, 
as $b$ is bounded, 
the second term in the right-hand side above is controllable with size $N^{-1}$, 
and of vanishing type. 
Let us rewrite the first term through an integration by parts:
\begin{align}
&\frac{1}{4\|\phi\|_2}\sum_{i<N-1}\sum_{j\notin\{i,i+1\}}\bar\eta_i\bar\eta_j\big(\phi_{i+1,j}-\phi_{i,j}\big) 
\nonumber\\
&\quad=  
\frac{1}{4\|\phi\|_2}\sum_{|i|<N-1}\sum_{j:|j-i|>1}(\bar\eta_{i-1}-\bar\eta_i)\bar\eta_j\phi_{i,j} 
+ \frac{1}{4\|\phi\|_2}\sum_{|i|<N-1}\bar\eta_{i-1}(\bar\eta_{i+1}\phi_{i,i+1}-\bar\eta_i\phi_{i,i-1})
\nonumber\\
&\quad= 
\frac{1}{4\|\phi\|_2}\sum_{|i|<N-1}\sum_{j:|j-i|>1}(\bar\eta_{i-1}-\bar\eta_i)\bar\eta_j\phi_{i,j} 
=:
S
,
\label{eq_iPP_continue_preuve_energie}
\end{align}
where the first equality makes use of $\phi_{\pm(N-1),\cdot} = 0$, while the second equality follows from $\phi_{i,i\pm1} = 0$ for each $|i|<N-1$.
To estimate the supremum in~\eqref{eq_FK_energy}, 
we see from~\eqref{eq_iPP_continue_preuve_energie} that we have to estimate $\nu^N_{g_0}(fS)$. 
This is done through the integration by parts Lemma~\ref{lemm_IPP_forme_dir}. 
This lemma is formulated with the variables $\omega_i = \bar\eta_i/\bar\sigma_i$ for $i\in\Lambda_N$, 
so we first rewrite~\eqref{eq_iPP_continue_preuve_energie} in terms of these variables. 
For $|i|<N-1$, using the identity:
\begin{align}
\bar\eta_{i-1} - \bar\eta_{i} 
&= 
\bar\sigma_{i-1}(\omega_{i-1}-\omega_i) - (\bar\sigma_i-\bar\sigma_{i-1})\omega_{i}
\nonumber\\
&= 
\bar\sigma_{i-1}(\omega_{i-1}-\omega_i) - \frac{\bar\rho'}{N}\Big(\sigma'(\bar\rho_i) + \frac{\bar\rho'}{N}\Big)\omega_i
\end{align}
and the convention:
\begin{equation}
q\phi(x,y) = q(x)\phi(x,y)\quad \text{for}\quad q:(-1,1)\rightarrow\R\ \text{and}\ (x,y)\in\squaredash\,
,
\end{equation}
the right-hand side $S$ of~\eqref{eq_iPP_continue_preuve_energie} can be rewritten as follows: 
\begin{align}
S 
&= 
\frac{1}{4\|\phi\|_2}\sum_{|i|<N-1}\sum_{j:|j-i|>1}\bar\sigma_{i-1}(\omega_{i-1}-\omega_i)\bar\eta_j\phi_{i,j} - \frac{1}{\|\phi\|_2}\Pi^N\Big(\frac{\bar\rho'}{\bar\sigma}\Big(\sigma'(\bar\rho)+\frac{\bar\rho'}{N}\Big)\phi\Big)
\nonumber\\
&=: 
S'+ \Pi^N(Y^{(0)}),\qquad Y^{(0)} = -\frac{1}{\|\phi\|_2}\frac{\bar\rho'}{\bar\sigma}\Big(\sigma'(\bar\rho)+\frac{\bar\rho'}{N}\Big)\phi
.
\end{align}
$\Pi^N(Y^{(0)})$ is of the form $\|\phi\|_2^{-1}\Pi^N(q\phi)$ with $q$ bounded independently of $\phi$, 
i.e. the kind of expression we were after. 
It remains to estimate $S'$. 
Define, for $|i|<N-1$, a function $v_i$ on $\Omega_N$ as follows:
\begin{align}
\forall\eta\in\Omega_N,\qquad v_i(\eta) = \frac{1}{4}\sum_{j:|j-i|>1}\bar\eta_j\bar\sigma_{i-1} \frac{\phi_{i,j}}{\|\phi\|_2}.\label{eq_def_v_energie}
\end{align}
Recall also the notation $C^{g_0}_\cdot$ defined in~\eqref{eq_def_B_x_D_x}:
\begin{equation}
\forall \eta\in\Omega_N,\forall i<N-1,\qquad \Pi^N(g_0)(\eta^{i,i+1})-\Pi^N(g_0)(\eta) =:-\frac{(\eta_{i+1}-\eta_i)}{N}C^{g_0}_i.
\end{equation}
With these notations, we can apply the integration by parts Lemma~\ref{lemm_IPP_forme_dir} to each $|i|<N-1$ with $u=v_i$, and obtain the existence of $C>0$ such that, for each $\delta>0$ and each density $f$ for $\nu^N_{g_0}$:
\begin{align}
\nu^N_{g_0}(fS') 
&\leq 
\delta N^2\nu^N_{g_0}\big(\Gamma(f^{1/2})\big) + \frac{C}{\delta N^2} \sum_{|i|<N-1}\nu^N_{g_0}\big(f(v_i)^2\big)\nonumber\\
&\quad 
-\sum_{|i|<N-1}(\bar\rho_{i}-\bar\rho_{i-1})\nu^N_{g_0}\big(\omega_{i-1}\omega_{i}fe^{-(\eta_{i}-\eta_{i-1})C^{g_0}_{i-1}/N}v_i\big)\label{eq_estimate_T_1_T_2_preuve_energie}\\
&\quad
+ \sum_{|i|<N-1}\nu^N_{g_0}\Big(\big(\omega_{i}-\omega_{i-1})\big(1-e^{-(\eta_{i}-\eta_{i-1})C^{g_0}_{i-1}/N}\big)fv_i\Big)
.\nonumber
\end{align}
We express each  term appearing in~\eqref{eq_estimate_T_1_T_2_preuve_energie} as $\Pi^N(G(\phi))$ plus error terms, 
for explicit $G$'s. 
Consider first the second term on the first line. 
By definition~\eqref{eq_def_v_energie} of $v$ and using $(\bar\eta_\cdot)^2 = \bar\sigma_\cdot + \sigma'(\bar\rho_\cdot)\bar\eta_\cdot$, 
it reads:
\begin{align}
\frac{C}{16\delta N}&\nu^N_{g_0}\Big(f\sum_{j,\ell\in\Lambda_N} \frac{\bar\eta_j\bar\eta_\ell}{N}\sum_{\substack{|i|<N-1:\\ |i-j|>1, |i-\ell|>1}}\bar\sigma_i^2\frac{\phi_{i,j}\phi_{i,\ell}}{\|\phi\|_2^2}\Big)\nonumber\\
&= 
\frac{1}{\delta}\nu^N_{g_0}\big(f\Pi^N(Y^{(1)})\big) 
+ \frac{C}{16\delta }\int_{\squaredash}\bar\sigma(x)\bar\sigma(y)^2\frac{\phi(x,y)^2}{\|\phi\|_2^2}dxdy 
+ \nu^N_{g_0}\big(f\theta^{N,1}_\delta(\phi/\|\phi\|_2)\big)
,
\label{eq_Y_1_plus_erreur_preuve_energie}
\end{align}
with $Y^{(1)}$ the function recording the off-diagonal, $\ell\neq j$ contribution: 
\begin{equation}\label{eq_def_Y_1_preuve_energie}
\forall (x,y)\in\squaredash,\qquad 
Y^{(1)}(x,y) 
= 
\frac{C}{4}\int_{(-1,1)}\bar\sigma(z)^2\frac{\phi(z,x)\phi(z,y)}{\|\phi\|_2^2}dz
.
\end{equation}
The error term $\theta^{N,1}_\delta(\|\phi\|_2^{-1}\phi)$ in~\eqref{eq_Y_1_plus_erreur_preuve_energie} involves discretisation errors and the diagonal, 
$\ell=j$ contributions. 
It is given for $\eta\in\Omega_N$ by:
\begin{align}
\theta^{N,1}(\|\phi\|_2^{-1}\phi)(\eta) 
&= 
\frac{C}{16\delta }\int_{\squaredash}\bar\sigma(x)\bar\sigma(y)^2\frac{\phi(x,y)^2}{\|\phi\|_2^2}dxdy - \frac{C}{16\delta N^2}\sum_{\substack{|i|<N-1 \\ |j-i|>1}}\bar\sigma_i \bar\sigma_j^2\frac{\phi_{i,j}^2}{\|\phi\|_2^2} \nonumber\\
&\quad 
+\frac{C}{16\delta N^2}\nu^N_{g_0}\Big(f\sum_{\substack{|i|<N-1 \\ |j-i|>1}}\bar\eta_i\sigma'(\bar\rho_i)\bar\sigma_j^2\frac{\phi_{i,j}^2}{\|\phi\|_2^2}\Big) + \frac{1}{N}\Pi^N(c),
\end{align}
where $c$ is a discretisation error arising in the replacement of~\eqref{eq_Y_1_plus_erreur_preuve_energie} by $Y^{(1)}$. 
We write $\theta^{N,1}(\|\phi\|_2^{-1}\phi)$ to emphasise the fact that $\theta$ only depends on $\|\phi\|_2^{-1}\phi$ rather than $\phi$ itself. 
The first line of $\theta^{N,1}(\|\phi\|_2^{-1}\phi)$ is configuration-independent, 
and bounded by $C(\phi)/N$. 
The first sum on the second line is of the form $N^{-1/2}Y^N(u)$ for $u:(-1,1)\rightarrow\R$ bounded, 
with $Y^N(u)$ the fluctuations defined in~\eqref{eq_def_fluct}. 
$\theta^{N,1}(\|\phi\|_2^{-1}\phi)$ is therefore controllable with size $N^{-1}$ by Lemma~\ref{lemm_size_controllable}, 
and of vanishing type by Lemma~\ref{lemm_type}. 
For later use, 
note that the middle term in~\eqref{eq_Y_1_plus_erreur_preuve_energie} is bounded by $C/(2^{10}\delta)$ for all large enough $N$, 
as $\bar\sigma\leq 1/4$. \\
Consider now line $2$ of~\eqref{eq_estimate_T_1_T_2_preuve_energie}. 
Using the identity $e^{x} = 1+\int_0^1xe^{tx}\, dt$ for $x\in\R$, there is $C(\rho_\pm)>0$ such that:
\begin{align}
\Big|\sum_{|i|<N-1}(\bar\rho_{i}-\bar\rho_{i-1})&\nu^N_{g_0}\big(\omega_{i-1}\omega_{i}fe^{-(\eta_{i}-\eta_{i-1})C^{g_0}_{i-1}/N}v_i\big) 
-\frac{\bar\rho'}{N}\sum_{|i|<N-1}\nu^N_{g_0}\big(\omega_{i-1}\omega_{i}fv_i\big) \Big|
\nonumber\\
&\quad \leq
\frac{C(\rho_\pm)}{N^2}\sum_{|i|<N-1}\nu^N_{g_0}\big(f\big|v_iC^{g_0}_{i-1}\big|\big)
.
\end{align}
The second sum in the left-hand side already involves three point correlations. 
It is shown to be $\Gamma$-controllable with size $N^{-1/2}$ in Proposition~\ref{prop_error_terms}, 
and is of large type. 
The second term is of the form $N^{-2}X^{v_2}_{2,\{0\}}$ in the notations of Lemma~\ref{lemm_size_controllable}, 
and therefore controllable with size $N^{-1}$ and of vanishing type. 
It follows that the first sum if a $\Gamma$-controllable error term with size $N^{-1/2}$ and of large type, 
depending only on $\|\phi\|_2^{-1}\phi$ rather than $\phi$ itself. 
We rewrite it as follows:
\begin{equation}
\theta^{N,2}(\|\phi\|_2^{-1}\phi)\big)
:=
-\sum_{|i|<N-1}(\bar\rho_{i}-\bar\rho_{i-1})\nu^N_{g_0}\big(\omega_{i-1}\omega_{i}fe^{-(\eta_{i}-\eta_{i-1})C^{g_0}_{i-1}/N}v_i\big) 
.
\end{equation}
By Lemma~\ref{lemm_type}, 
there is a numerical constant $\gamma_2>0$ such that $\gamma_2 \theta^{N,2}(\|\phi\|_2^{-1}\phi)$ is of LS type.\\
Consider finally line $3$ of~\eqref{eq_estimate_T_1_T_2_preuve_energie}. 
Let $C(\bar\rho)>0$ be such that:
\begin{equation}
\forall |i|<N-1,\qquad |\omega_i-\omega_{i-1}|\leq C(\bar\rho).
\end{equation}
Using this time the existence of $c(g_0)>0$ such that $|e^x-1-x|\leq c(g_0)x^2$ for all $|x|\leq 2\|g_0\|_\infty$, 
one can write, 
for each $\eta\in\Omega_N$:
\begin{align}
\sum_{|i|<N-1}\big(\omega_{i}-\omega_{i-1})\big(1-e^{-(\eta_{i}-\eta_{i-1})C^{g_0}_{i-1}/N}\big)v_i
&\leq 
\frac{1}{N}\sum_{|i|<N-1}(\omega_i - \omega_{i-1})(\eta_i-\eta_{i-1})C_{i-1}^{g_0}v_i \nonumber\\
&\quad + 
\frac{c(g_0)C(\bar\rho)}{N}\sum_{|i|<N-1}\frac{|v_i|}{N}\big(C^{g_0}_{i-1}\big)^2
.
\label{eq_decomp_third_line_energie}
\end{align}
The last term is an average over $i$ of terms of the form $N^{-3}|X^{w^i_3}_{3,\{0\}}|$ with the notations of Lemma~\ref{lemm_size_controllable}, 
where the $w^i_3$ satisfy $\sup_{N,i}\sup_{\Lambda_N^3}|w^i_3|<\infty$. 
It is therefore controllable with size $N^{-3/2}$ and of vanishing type. 
To estimate the first sum in the right-hand side of~\eqref{eq_decomp_third_line_energie}, 
we use the following elementary identity, 
valid for each $|i|<N-1$:
\begin{equation}\label{eq_omega_i_moins_omega_i_moins_un_fois_difference_eta_i}
(\omega_i-\omega_{i-1})(\eta_i-\eta_{i-1}) 
= 
\Big[2+(1-\bar\rho_i-\bar\rho_{i-1})[\omega_i + \omega_{i-1}]- (\bar\sigma_{i-1}+\bar\sigma_{i})\omega_{i-1}\omega_i\Big]
.
\end{equation}
This identity can be obtained by making the following observation:
\begin{equation}
\forall i\in\Lambda_N,\qquad \eta_i\omega_i := \eta_i\frac{(\eta_i-\bar\rho_i)}{\bar\sigma_i} = \frac{\eta_i}{\bar\sigma_i}(1-\bar\rho_i) = \frac{\bar\eta_i}{\bar\rho_i} +1.
\end{equation}
Looking at~\eqref{eq_decomp_third_line_energie}, 
we see that the term $C_{i-1}^{g_0}v_i$ already contains two-point correlations for each $|i|<N-1$. 
We therefore claim that only the constant term in the identity~\eqref{eq_omega_i_moins_omega_i_moins_un_fois_difference_eta_i} will give something that is not an error term in~\eqref{eq_decomp_third_line_energie}. 
More precisely, 
we claim that one can obtain the following bound for line $3$ of~\eqref{eq_decomp_third_line_energie}:
\begin{align}
&\sum_{|i|<N-1}\nu^N_{g_0}\Big(\big(\omega_{i}-\omega_{i-1})\big(1-e^{-(\eta_{i}-\eta_{i-1})C^{g_0}_{i-1}/N}\big)fv_i\Big) 
\label{eq_bound_Y_2_energy}\\
&\hspace{1cm}
\leq 
\nu^N_{g_0}\big(f\Pi^N(Y^{(2)})\big)+ \frac{1}{4}\int_{\squaredash}\bar\sigma(x)\bar\sigma(y)\frac{\phi(x,y)}{\|\phi\|_2}\partial_1 g_0\, (x,y)dxdy + \nu^N_{g_0}\big(f\theta^{N,3}(\|\phi\|_2^{-1}\phi)\big)
,
\nonumber
\end{align}
where $\theta^{N,3}(\|\phi\|_2^{-1}\phi)$ is controllable with size $N^{-1/2}$, 
of large type, and again depends on $\|\phi\|_2^{-1}\phi$ only rather than on $\phi$. 
By Lemma~\ref{lemm_type}, 
there is thus $\gamma_3>0$ such that $\theta^{N,3}(\|\phi\|_2^{-1}\phi)$ is of LS type.\\
There is $C(g_0) = C(\rho_\pm)>0$ independent of $\phi$ bounding the middle term in~\eqref{eq_bound_Y_2_energy}, 
which comes from the constant term in the identity $(\bar\eta_\cdot)^2=\bar\sigma_\cdot + \sigma'(\bar\rho_\cdot)\bar\eta_\cdot$, 
and $Y^{(2)}$ is defined as:
\begin{align}
\forall (x,y)\in\squaredash,\qquad 
Y^{(2)}(x,y) 
=  
\int_{(-1,1)} \bar\sigma(z)\partial_1 g_0(z,x)\frac{\phi(z,y)}{\|\phi\|_2}\, dz
.
\label{eq_def_Y_2_preuve_energie}
\end{align}
To summarise, 
we have established the following. 
If, for $\delta>0$, $V_\delta$ is the quantity: 
\begin{equation}
V_\delta = \Pi^N(\delta^{-1} Y^{(0)}+Y^{(1)}+Y^{(2)})
, 
\end{equation}
then, for each $N$ large enough depending on $\phi$:
\begin{equation}\label{eq_energy_estimate_after_IPP}
\nu^N_{g_0}\big(f\|\phi\|_2^{-1}\Pi^N(\partial_1 \phi)\big)
\leq 
\delta N^2\nu^N_{g_0}\big(\Gamma(\sqrt{f})\big) 
+ 
\nu^N_{g_0}(fV_\delta) + C(\rho_\pm) + \zeta_\delta^N(\|\phi\|_2^{-1}\phi)
.
\end{equation}
The quantity $\zeta_\delta^N(\|\phi\|_2^{-1}\phi)$ is a $\Gamma$-controllable error term of size $N^{-1/2}$ of large type. 
One can then choose the $\kappa$ appearing in~\eqref{eq_FK_energy} equal to
some $\kappa_0(\rho_\pm)>0$ independent of $\phi,\delta$ such that, for $N$ large enough depending on $\phi$: 
\begin{equation}
\kappa_0(\rho_\pm)\zeta_\delta^N(\|\phi\|_2^{-1}\phi)\text{ is of LS type}
.
\end{equation}
To conclude the estimate of the right-hand side of~\eqref{eq_FK_energy} using~\eqref{eq_energy_estimate_after_IPP}, 
it remains to bound the average of $V_\delta$.  
Fix $\delta=1/16$. 
Recall that $Y^{(0)}$ has $2$-norm bounded by $4$ and,  
by Cauchy-Schwarz inequality, 
that $Y^{(1)}$ has $2$-norm bounded by $2^{-6}\delta^{-1}C =C/4$ and $Y^{(2)}$ by $\sqrt{2}\|\nabla g_0\|_2/8$. 
The $Y^{(i)}$ ($0\leq i\leq 2$) thus have $2$-norm bounded independently of $\phi$.  
Lemma~\ref{lemm_size_controllable} then implies that there is $\kappa_0'(\rho_\pm)>0$ independent of $\phi$ and a numerical constant $C$ such that $\kappa_0'(\rho_\pm) V_{1/16}$ is of LS type. 
Taking $\kappa(\rho_\pm):= \min\{\kappa_0(\rho_\pm),\kappa'_0(\rho_\pm)\}$, 
there are then constants $C(\rho_\pm),C(\rho_\pm,\phi)>0$ such that, for each $N$ large enough depending on $\phi$: 
\begin{align}
&\kappa(\rho_\pm)\nu^N_{g_0}\Big(f\big[V_{1/16} + \zeta^N_{1/16}(\|\phi\|_2^{-1}\phi)\big]\Big)
\nonumber\\
&\hspace{3cm}\leq 
\frac{H(f\nu^N_{g_0}|\nu^N_{g_0})}{2^9 C_{LS}} + C(\rho_\pm) + \frac{C(\rho_\pm,\phi)}{N^{1/2}} 
+ \frac{N^2}{16}\nu^N_{g_0}\big(\Gamma(f^{1/2})\big)
.
\end{align}
Injecting this bound in~\eqref{eq_FK_energy} and using the log-Sobolev inequality of Lemma~\ref{lemm_LSI_sec3} concludes the proof:
\begin{align}
\frac{1}{T}\log \Prob^{\nu^N_{g_0}}\Big(\frac{1}{T}\int_0^T\|\phi\|_2^{-1}\Pi^N_t(\partial_1 \phi)\, dt>  A\Big) 
\leq
-A\kappa(\rho_\pm) + C(\rho_\pm) + \frac{C(\rho_\pm,\phi)}{N^{1/2}}
.
\end{align}
\end{proof}
\subsection{Exponential tightness}\label{sec_exp_tightness}
In this section, 
we prove that the upper bound~\eqref{eq_upper_bound_compact sets} is also valid for closed sets. 
We refer to Appendix~\ref{app_sobolev_spaces} for a characterisation of compact sets in $(\mathcal T'_s,*)$, 
and establish exponential tightness in the next proposition.
\begin{proposition}\label{prop_exp_tightness}
For each large enough $A>0$, 
there is a compact set $\mathcal K_A\subset (\mathcal T'_s,*)$ such that:
\begin{equation}
\sup_N\limsup_{T\rightarrow\infty}\frac{1}{T}\log \Prob^{\nu^N_{g_0}}\Big(\frac{1}{T}\int_0^T\Pi^N_t\, dt \notin \mathcal K_A\Big) 
\leq 
-A
.
\end{equation}
\end{proposition}
\begin{proof}
In Appendix~\ref{app_sobolev_spaces}, 
we prove that there is a norm $\|\cdot\|_{\mathbb{T},-2}$ on $\mathcal T'_s$ such that the set $\{\|\cdot\|_{\mathbb{T},-2}\leq B\}$ is relatively weak-* compact for each $B>0$. 
The norm $\|\cdot\|_{\mathbb{T},-2}$ involves a certain family of functions $\psi_m\in\mathcal T$ ($m\in\N^2$), 
and reads:
\begin{equation}
\|\Pi\|^2_{\mathbb{T},-2} 
:= 
\sum_{m\in\N^2}\frac{1}{(1+|m|^2)^2}\, |\Pi(\psi_m)|^2,
\qquad \Pi\in \mathcal T'_s
.
\end{equation}
Above, $|m|^2 = m_1^2+m_2^2$ for $m=(m_1,m_2)\in\N^2$. 
It is therefore enough to prove the existence of $c(\cdot)>0$ with $\lim_{B\rightarrow\infty}c(B)=\infty$ such that:
\begin{equation}\label{eq_to_prove_exp_tightness}
\sup_{N\in\N^*}\Prob^{\nu^N_{g_0}}\Big(\Big\|\frac{1}{T}\int_0^T\Pi^N_t \, dt\Big\|_{\mathbb{T},-2}\geq A\Big) 
\leq 
e^{-c(A)T}
.
\end{equation}
To prove~\eqref{eq_to_prove_exp_tightness},
let $\epsilon\in(0,1)$ and define $c_\epsilon := \sum_m (1+|m|^2)^{1+\epsilon}<\infty$. 
A union bound gives:
\begin{align}
&\Prob^{\nu^N_{g_0}}\Big(\Big\|\frac{1}{T}\int_0^T\Pi^N_t \, dt\Big\|_{\mathbb{T},-2}\geq A\Big) 
\nonumber\\
&\hspace{3cm}\leq 
\sum_{m\in\N^2} 
\Prob^{\nu^N_{g_0}}\Big(\Big|\frac{1}{T}\int_0^T\Pi^N_t(\psi_m)\, dt\Big|\geq \frac{A}{c_\epsilon^{1/2}} (1+|m|^2)^{\frac{1-\epsilon}{2}}\Big)
. 
\end{align}
The $\psi_m$ are just restrictions to $\rhd$ of the eigenvectors of the torus Laplacian on $[-2,2)^2$ (see Appendix~\ref{app_sobolev_spaces}). 
In particular, $\|\psi_m\|_\infty \leq C$ where $C$ does not depend on $N,m$. 
Each probability in the above sum is thus estimated by Corollary~\ref{coro_Boltzmann_gibbs_sec3} according to:
\begin{align}
\Prob^{\nu^N_{g_0}}\Big(\Big\|\frac{1}{T}\int_0^T\Pi^N_t \, dt\Big\|_{\mathbb{T},-2}\geq A\Big) 
&\leq
\sum_{m\in\N^2} 
\exp\Big[- \frac{\zeta_2 A}{c_\epsilon^{1/2}} \frac{(1+|m|^2)^{\frac{1-\epsilon}{2}}}{\|\psi_m\|_{\infty}} T + C(\rho_\pm)T\Big]
\nonumber\\
&\leq 
\sum_{m\in\N^2} 
\exp\Big[- \frac{\zeta_2 A}{c_\epsilon^{1/2} C} \frac{(1+|m|^2)^{\frac{1-\epsilon}{2}}}{\|\psi_m\|_{\infty}} T + C(\rho_\pm)T\Big]
,
\end{align}
where the parameter $\zeta_2$ is a universal constant defined in item 3 of Lemma~\ref{lemm_type}. 
The expression in the right-hand side is summable for any $\epsilon<1$. 
It is moreover bounded by $e^{-c(A)T}$, 
with $c(A)>0$ for $A$ large enough and $\lim_{B\rightarrow\infty}c(B)=+\infty$.  
This completes the proof.
\end{proof}

\section{Lower bound for smooth trajectories}\label{sec_lower_bound}
In this section, 
we give a lower bound on $\mathbb Q^{\pi^N_{inv}}_T(\mathcal O)$ 
(defined in~\eqref{eq_def_Q_T}) 
when $\mathcal O$ is an open subset of $(\mathcal T'_s,*)$, in terms of the kernels $k_h,h\in\s(\epsilon_B)$, 
with $\epsilon_B$ the quantity appearing in Theorems~\ref{theo_large_devs}--\ref{theo_entropic_problem}. 
As for standard large deviations (see Chapter 10 in~\cite{Kipnis1999}), 
we consider the tilted dynamics $\Prob_h,h\in\s(\epsilon_B)$ such that $\frac{1}{4}\big<k_h,\cdot\big>\in \mathcal O$, and obtain a lower bound by proving that the measure $\Prob_h$ concentrates on $\mathcal O$. 
In the following, let $\Q_{T,h}$ denote the law of $\frac{1}{T}\int_0^T\Pi^N_tdt$ under $\Prob_h$.\\

We first change the initial measure to $\nu^N_{g_h}$, with $g_h$ chosen according to~\eqref{eq_main_equation}, 
and tilt the dynamics by $h$. 
Using Jensen inequality to obtain the last line below, one finds, 
with ${\rm D}_h = d\Prob_h/d\Prob$ the Radon-Nikodym derivative on $[0,T]$ (see~\eqref{eq_def_der_radon_nikodym}):
\begin{align}
\log \Q_T^{\pi^N_{inv}}(\mathcal O) 
&= 
\log\E_h^{\nu^N_{g_h}}\Big[ {\bf 1}_{\mathcal O} \frac{\pi^N_{inv}(\eta_0)}{\nu^N_{g_h}(\eta_0)} ({\rm D}_h)^{-1}\Big] 
\nonumber\\
&= 
\log\E^{\nu^N_{g_h}}_{T,h,\mathcal O}\Big[\frac{\pi^N_{inv}(\eta_0)}{\nu^N_{g_h}(\eta_0)}({\rm D}_h)^{-1}\Big] + \log \Q_{T,h}^{\nu^N_{g_h}}(\mathcal O)
\nonumber\\
&\geq 
\E^{\nu^N_{g_h}}_{T,h,\mathcal O}\big[-\log {\rm D}_h\big] + \E^{\nu^N_{g_h}}_{T,h,\mathcal O}\Big[\log\Big(\frac{\pi^N_{inv}(\eta_0)}{\nu^N_{g_h}(\eta_0)}\Big)\Big] +\log \Q_{T,h}^{\nu^N_{g_h}}(\mathcal O)
\nonumber\\
&\geq \E^{\nu^N_{g_h}}_{T,h,\mathcal O}\big[-\log {\rm D}_h\big] - C(\rho_\pm,h)N +\log \Q_{T,h}^{\nu^N_{g_h}}(\mathcal O)
.
\label{eq_entropic_term_borne_inf_0}
\end{align}
Above, $C(g,g_h)>0$, 
while $\Prob^{\nu^N_{g_h}}_{T,h,\mathcal O}$ is the probability $\Prob^{\nu^N_{g_h}}_h$ conditional to $\big\{\frac{1}{T}\int_0^T\Pi^N_t \, dt\in \mathcal O\big\}$:
\begin{equation}
\Prob^{\nu^N_{g_h}}_{T,h,\mathcal O}(\cdot)  
= 
\Big(\Q^{\nu^N_{g_h}}_{T,h}(\mathcal O)\Big)^{-1}\Prob^{\nu^N_{g_h}}_h\Big(\cdot \cap  \Big\{ \frac{1}{T}\int_0^T\Pi^N_t \, dt\in \mathcal O\Big\}\Big)
,
\label{eq_lower_bound_0}
\end{equation}
with:
\begin{equation}
\E^{\nu^N_{g_h}}_{T,h,\mathcal O}[\,\cdot\,] 
= 
\int \cdot\ d\Prob^{\nu^N_{g_h}}_{T,h,\mathcal O}
.
\end{equation}
The terms appearing in~\eqref{eq_entropic_term_borne_inf_0} are of three types: the change of initial condition corresponding to the constant $-C(\rho_\pm,g_h)N$, 
which will vanish upon division by $T$ when $T$ is large; 
the dynamical part with $\log {\rm D}_h$, 
and the term $\Q^{\nu^N_{g_h}}_{T,h}(\mathcal O)$. 
The latter is well controlled only if $h$ is such that, 
under $\Q_{T,h}^{\nu^N_{g_h}}$, 
correlations are typically in $\mathcal O$ when $N,T$ are large. 
For such an $h$, upon dividing by $T$ and taking the large $T$ limit, 
only the dynamical part will contribute. 
The limit of $\Q^{\nu^N_{g_h}}_{T,h}(\mathcal O)$ is worked out in the next section, 
and the dynamical part is studied in Section~\ref{sec_dynamical_part_lower_bound}.
\subsection{Law of large numbers and Poisson equation}\label{sec_LLN}
\begin{proposition}\label{coro_limite_hydro}
Let $h\in\s(\epsilon_B)$, and let $k_h$ be the large $N$ limit of the correlations under $\nu^N_{g_h}$, where $g_h$ solves the main equation~\eqref{eq_main_equation}. If $\mathcal O\subset (\mathcal T'_s,*)$ is an open set containing $\frac{1}{4}\big<k_h,\cdot\big>$, then:
\begin{equation}
\liminf_{N\rightarrow\infty}\liminf_{T\rightarrow\infty}\Q_{T,h}^{\nu^N_{g_h}}(\mathcal O) = 1.\label{eq_LLN_Q(O)}
\end{equation}
\end{proposition}
An open set in $(\mathcal T'_s,*)$ is a (possibly uncountable) union of finite intersections of sets of the form $\Big\{\Big|\frac{1}{T}\int_0^T\Pi^N_t(\phi)\, dt - \frac{1}{4}\big<k_h,\phi\big>\Big| \in U\Big\}$, for an open set $U\subset \R$ and $\phi\in\mathcal T$. 
It is therefore enough to prove~\eqref{eq_LLN_Q(O)} for those sets, with $U = (-\epsilon,\epsilon)$ for $\epsilon>0$.

Let us first check that it suffices to prove~\eqref{eq_LLN_Q(O)} for symmetric $\phi\in C^2(\bar\rhd)$ with $\phi_{|\partial\square}=0$.  
Recall the notation $\|\phi\|^2_{2,N} = N^{-2}\sum_{i,j}\phi_{i,j}^2$. 
For each $a>0$ and $\phi,\psi\in\mathcal T$, 
the entropy inequality and Markov inequality together with the second item of Lemma~\ref{lemm_size_controllable} give, 
for some $C(h,\rho_\pm)>0$:
\begin{align}
\Prob^{\nu^N_{g_h}}_{h}\bigg(\Big|\frac{1}{T}\int_0^T\Pi^N_t(\phi-\psi)\, dt\Big|
\geq 
a\bigg)
&\leq 
\frac{C(h,\rho_\pm)\|\phi-\psi\|_{2,N}}{a} 
\nonumber\\
&= 
\frac{C(h,\rho_\pm)\|\phi-\psi\|_{2}}{a} + o_N(1)
.
\end{align}
In addition, $\big|\big<k_h,\phi-\psi\big>\big|
\leq 
\|k_h\|_2\|\phi-\psi\|_2$. 
Since one can approximate $\phi\in\mathcal T$ in $\mathbb L^2(\squaredash)$ with arbitrary precision by some $\psi\in \mathcal T\cap C^2(\bar\rhd)$ with $\psi_{\partial\square}=0$,  
it is enough to focus on such $\psi$ as claimed. 
This is done in the next proposition, by means of a Poisson problem associated with the large $N$ limit of the generator $N^2L_h$.
\begin{proposition}\label{prop_poisson_equation}
Let $\phi\in \mathcal T\cap C^2(\bar\rhd)$ be a symmetric function with $\phi_{\partial\square}=0$, 
where we recall $\mathcal T=\mathbb H^2(\squaredash)$. 
Let $h\in\s(\epsilon_B)$. 
Then, for any $\epsilon\in(0,1)$ and any $T>0$, 
there are positive constants $C(h,\phi),C'(h,\phi)>0$ (independent of $T$) such that:
\begin{equation}
\Prob^{\nu^N_{g_h}}_h\bigg(\Big|\frac{1}{T}\int_{0}^ T\Pi^N_t(\phi)\, dt
-\frac{1}{4}\big<k_h,\phi\big>\Big|
\geq \epsilon\|\phi\|_2\bigg) 
\leq 
\frac{C(h,\phi)}{T}\big(\epsilon^{-2} +N\epsilon^{-1}\big) + \frac{C'(h,\phi)}{\epsilon^2N^{1/2}}
.
\label{eq_quantitative_LLN}
\end{equation}
\end{proposition}
\begin{proof}
Fix a symmetric $\phi\in \mathcal T\cap C^2(\bar\rhd)$. To prove Proposition~\ref{prop_poisson_equation}, we express the difference appearing in the probability in~\eqref{eq_quantitative_LLN} as a time integral involving the generator $N^2L_h$, plus a martingale term. The martingale term is then proven to fluctuate like $\sqrt{T}$ when $N,T$ are large. 
It thus vanishes in the large $T$ limit upon dividing by $T$. Recall that $\mathcal M(u,v)(x,y) = \int_{(-1,1)}u(z,x)\bar\sigma(z)v(z,y)dz$ for any $u,v\in\mathbb L^2(\squaredash)$. The key ingredient is the following Poisson equation:
\begin{equation}
\begin{cases}
\frac{1}{2}\Delta f(x,y)
+ \frac{1}{2}\mathcal M(\partial_1f,\partial_1h) + \frac{1}{2}\mathcal M(\partial_1h,\partial_1 f) = \displaystyle{\frac{\phi(x,y)}{\|\phi\|_2}}\quad &\text{for }(x,y)\in\squaredash\, , \\
f = 0 \quad&\text{on }\partial \square,\\
(\partial_1-\partial_2)f(x_\pm,x) = 0\quad &\text{for }x\in (-1,1).\label{eq_Poisson_pour_partie_PiN_du_generateur}
\end{cases}
\end{equation}
For $\phi\in \mathcal T\cap C^2(\bar\rhd)$ satisfying $\phi_{\partial\square}=0$, 
we prove in Appendix~\ref{app_Poisson} that~\eqref{eq_Poisson_pour_partie_PiN_du_generateur} has a unique solution $f_\phi\in \mathcal T \cap C^3(\bar\rhd)$, 
a symmetric function on $\squaredash$. 
It satisfies $\|f_\phi\|_\infty\leq C$ for a constant independent of $\phi$, 
see Proposition~\ref{prop_solving_P_triangle}. 
The semi-martingale decomposition for $\Pi^N(f_\phi)$ reads:
\begin{equation}
\forall T\geq 0,\qquad 
\Pi^N_T(f_\phi) 
= 
\Pi^N_0(f_\phi) + \int_0^TN^2L_h\Pi^N_t(f_\phi)\, dt + M^{N,f_\phi}_t
.
\label{eq_decomp_martingale_borne_inf}
\end{equation}
Let us first use the Poisson equation~\eqref{eq_Poisson_pour_partie_PiN_du_generateur} to express $N^2 L_h\Pi^N(f_\phi)$ in terms of $\Pi^N(\phi)$. 
\begin{lemma}\label{lemm_lien_generateur_Poisson}
\begin{equation}
N^2L_h\Pi^N(f_\phi) = \frac{1}{\|\phi\|_2}\Big(\Pi^N( \phi) - \frac{1}{4}\big<k_h,\phi\big>\Big) +\theta^N(f_\phi),\label{eq_dans_lemm_moments_4_3}
\end{equation}
where $\theta^N(f_\phi)$ is an error term of size $N^{-1/2}$ of large type (recall Definitions~\ref{def_controllability}--\ref{def_LS_type} of error terms and of large type). 
\end{lemma}
Assuming Lemma~\ref{lemm_lien_generateur_Poisson} for the moment, 
let us prove Proposition~\ref{prop_poisson_equation}. 
For each $T>0$ and $\epsilon>0$, integrate~\eqref{eq_dans_lemm_moments_4_3} between $0$ and $T$ and use the martingale decomposition~\eqref{eq_decomp_martingale_borne_inf} to find:
\begin{align}
&\Prob^{\nu^N_{g_h}}_h\bigg(\Big|\frac{1}{T}\int_{0}^ T\Pi^N_t(\phi)\, dt-\frac{1}{4}\big<k_h,\phi\big>\Big|\geq \epsilon\|\phi\|_2\bigg)
\nonumber\\
 &\qquad 
 = 
 \Prob^{\nu^N_{g_h}}_h\bigg(\frac{1}{T}\Big|\Pi^N_T(f_\phi)-\Pi^N_0(f_\phi)-M_T^{N,f_\phi}-\int_0^T\theta^N_t(f_\phi)\, dt\Big|>\epsilon\bigg)
 . 
\end{align}
Let us estimate each of the terms appearing in the last probability. 
Proposition~\ref{prop_Boltzmann_gibbs_sec3} takes care of $\theta^N(f_\phi)$: 
there are $\gamma,C>0$ depending on $\rho_\pm,h$, but not $\phi$, 
such that for all $N\in\N^*$ and $T>0$:
\begin{equation}
\E^{\nu^N_{g_h}}_h\bigg[\Big|\frac{1}{T}\int_0^T\theta^N_t(f_\phi)\, dt\Big|\bigg] 
\leq 
\frac{1}{\gamma}\log \E^{\nu^N_{g_h}}_h\bigg[\exp\Big|\frac{\gamma}{T}\int_0^T\theta^N_t(f_\phi)\, dt\Big|\bigg] 
\leq
CN^{-1/2}
.
\end{equation}
By Markov- and Chebychev inequalities:
\begin{align}
&\Prob^{\nu^N_{g_h}}_h\bigg(\Big|\frac{1}{T}\int_{0}^ T\Pi^N_t(\phi)\, dt-\frac{1}{4}\big<k_h,\phi\big>\Big|\geq \epsilon\|\phi\|_2\bigg)
\nonumber\\
&\qquad\leq 
\frac{3C}{\epsilon N^{1/2}}+ \frac{3}{T\epsilon}\E^{\nu^N_{g_h}}_{h}\Big[\big|\Pi^N_T(f_\phi)-\Pi^N_0(f_\phi)\big|\Big] + \frac{9}{\epsilon^2 T^2}\E^{\nu^N_{g_h}}_h\Big[\big<M^{N,f_\phi}\big>_T\Big]
.
\label{eq_estimate_erreor_term_borne_inf}
\end{align}
The middle expectation involving the terms $\Pi^N_{T}(f_\phi)$, $\Pi^N_0(f_\phi)$ is bounded by $\|f_\phi\|_\infty N$ 
(which does not depend on $\phi$).  
Its contribution to~\eqref{eq_estimate_erreor_term_borne_inf} thus vanishes when $T$ is large. 
Let us prove that the quadratic variation of $M^{N,f_\phi}$ has average bounded linearly in time for $N$ large, 	
which will be enough to conclude the proof of Proposition~\ref{prop_poisson_equation}. 
The quadratic variation is given for each $t\geq0$ by:
\begin{equation}\label{eq_quadratic_variation_M_f_t}
\big<M^{N,f_\phi}\big>_t 
= 
\frac{1}{2}\int_0^t\sum_{i<N-1}c_h\big(\eta(s),i,i+1\big)\Big[\frac{\bar\rho'}{4N}(f_\phi)_{i,i+1}+\frac{1}{2N}\sum_{j\notin\{i,i+1\}}\bar\eta_j(s)\partial^N_1 (f_\phi)_{i,j}\Big]^2ds
. 
\end{equation}
Recall from the definition~\eqref{eq_def_jump_rates_H} of $c_h$ that $\sup_{N,\Omega_N}c_h \leq c(h)$ for some $c(h)>0$. 
Using the inequality $(a+b)^2 \leq 2a^2+2b^2$ for $a,b\in\R$ and the fact that $f_\phi$ is bounded, 
one has:
\begin{align}
\E^{\nu^N_{g_h}}_h\Big[\big<M^{N,f_\phi}\big>_T\Big]
&\leq 
c(h)\E^{\nu^N_{g_h}}_h\bigg[\int_0^T\sum_{i<N-1}\Big(\frac{1}{2N}\sum_{j\notin\{i,i+1\}}\bar\eta_j(t)\partial^N_1(f_\phi)_{i,j}\Big)^2dt\bigg]
\nonumber\\
&\qquad +\frac{Tc(h) \|f_\phi\|^2_\infty}{N}
.
\end{align}
The integrand at each time $t\leq T$ is of the form $N^{-1}X_{2,\{0\}}^{u}$ in the notations of Lemma~\ref{lemm_size_controllable}, 
with $u=(u_{i,j})_{i,j\in\Lambda_N}$ given by:
\begin{align}
u_{i,j} 
&= 
\frac{1}{N}\sum_{k\notin\{i-1,i,j-1,j,N-1\}}\partial^N_1 (f_{\phi})_{k,i}\partial^N_1 (f_{\phi})_{k,j} 
\nonumber\\
&
= \int_{(-1,1)}\partial_1 f_\phi(z,i/N)\partial_1 f_\phi(z,j/N)dz + \frac{v^N_{i,j}}{N}
,
\end{align}
where $v^N$ is a discretisation error bounded by $C(\phi)>0$. 
$N^{-1}X_{2,\{0\}}^{u}$ is thus controllable with size $1$, 
and Proposition~\ref{prop_Boltzmann_gibbs_sec3} yields the existence of $C(h,\phi)>0$, $C'(h,\phi)>0$ such that:
\begin{equation}\label{eq_bound_quadratic_variation}
\frac{9}{\epsilon^2 T^2}\E^{\nu^N_{g_h}}_h\Big[\big<M^{N,f_\phi}\big>_T\Big]\leq \frac{1}{\epsilon^2}\Big(\frac{C(h,\phi)}{ T} + \frac{C(h,\phi)}{N}\Big)
.
\end{equation}
This estimate and~\eqref{eq_estimate_erreor_term_borne_inf} conclude the proof of Proposition~\ref{prop_poisson_equation} assuming Lemma~\ref{lemm_lien_generateur_Poisson}, 
proven below.
\end{proof}
\begin{proof}[Proof of Lemma~\ref{lemm_lien_generateur_Poisson}.] 
The starting point is the expression of $N^2L_h\Pi^N(f_\phi)$ worked out in Corollary~\ref{coro_comparaison_D_H_martingale}: 
there is a $\Gamma$-controllable error term $\tilde\epsilon^N(h,f_\phi)$ with size $N^{-1/2}$, 
of large type, such that:
\begin{align}
N^2L_h \Pi^N(f_\phi) 
&= 
\frac{1}{2}\Pi^N\Big(\Delta f_\phi + \mathcal M(\partial_1f_{\phi},\partial_1h) + \mathcal M(\partial_1 h,\partial_1 f_{\phi})\Big) -\frac{(\bar\rho')^2}{4}\int_{(-1,1)}f_\phi(x,x)dx \nonumber\\
&\quad  + \frac{1}{4}\int_{\squaredash}\bar\sigma(x)\bar\sigma(y)\big[\partial_1 f_\phi(x,y)\partial_1h(x,y)\big]dxdy + \tilde\epsilon^N(h,f_\phi).\label{eq_Poisson_problem_0}
\end{align}
Note the absence of the diagonal term $\sum_{i<N-1}\bar\eta_i \bar\eta_{i+1}(\partial_1 (f_\phi)_{i_+,i}-\partial_1 (f_\phi)_{i_-,i})$: 
recall from Remark~\ref{remark_normal_derivative} that this term corresponds to the derivative of $f_\phi$ in the normal direction to the diagonal, 
which vanishes according to~\eqref{eq_Poisson_pour_partie_PiN_du_generateur}. 
By Proposition~\ref{prop_Boltzmann_gibbs_sec3}, 
there are $\gamma,C>0$ depending on $h,\phi$ with:
\begin{equation}
\E^{\nu^N_{g_h}}_h\bigg[\Big|\int_0^T\tilde \epsilon^N_t(h,f_\phi)\, dt\Big|\bigg] 
\leq 
\frac{1}{\gamma}\log \E^{\nu^N_{g_h}}_h\bigg[\exp\Big|\gamma\int_0^T\tilde \epsilon^N_t(h,f_\phi)\, dt\Big|\bigg] 
\leq 
C(h,\phi)N^{-1/2}.\label{eq_size_epsilon_alpha_preuve_limite_hydro}
\end{equation}
As $f_\phi$ solves~\eqref{eq_Poisson_pour_partie_PiN_du_generateur},~\eqref{eq_Poisson_problem_0} can be written as:
\begin{align}
N^2L_h \Pi^N(f_\phi) &= \frac{1}{\|\phi\|_2}\Pi^N(\phi)  -\frac{(\bar\rho')^2}{4}\int_{(-1,1)}f_\phi(x,x)\, dx\nonumber\\
&\quad    + \frac{1}{4}\int_{\squaredash}\bar\sigma(x)\bar\sigma(y)\big[\partial_1 f_\phi(x,y)\partial_1h(x,y)\big]\, dxdy+  \tilde\epsilon^N(h,f_\phi).\label{eq_Poisson_problem_0_bis}
\end{align}
Equation~\eqref{eq_Poisson_problem_0_bis} will correspond to~\eqref{eq_dans_lemm_moments_4_3} with $\theta^N(f_\phi) = \tilde \epsilon^N(h,f_\phi)$, if we can prove:
\begin{align}
-\frac{1}{4\|\phi\|_2}\big<k_h,\phi\big> 
&= 
-\frac{(\bar\rho')^2}{4}\int_{(-1,1)}f_\phi(x,x)dx 
\nonumber\\
&\quad + 
\frac{1}{4}\int_{\squaredash}\bar\sigma(x)\bar\sigma(y)\big[\partial_1 f_\phi(x,y)\partial_1h(x,y)\big]\, dx\, dy
.
\label{eq_k_h_is_constant_terms_Poisson}
\end{align}
By definition, $k_h +\bar\sigma =(\bar\sigma^{-1}-g_h)^{-1}$ with $g_h$ solving the main equation~\eqref{eq_main_equation}. 
We prove in Appendix~\ref{app_Poisson} that this is equivalent to saying that $k_h$ solves the Euler-Lagrange equation~\eqref{eq_Euler-Lagrange_appendix}. 
The identity~\eqref{eq_k_h_is_constant_terms_Poisson} is then straightforwardly obtained by taking $f_\phi$ as a test function in the weak formulation~\eqref{eq_Euler-Lagrange_weak} of the Euler-Lagrange equation, 
then integrating by parts.
\end{proof}
\subsection{Estimation of the dynamical part and conclusion of the lower bound}\label{sec_dynamical_part_lower_bound}
Let $h\in\s(\epsilon_B)$. 
In this section, we estimate the term $\E^{\nu^N_{g_h}}_{T,h,\mathcal O}\big[-\log {\rm D}_h\big]$ arising in~\eqref{eq_entropic_term_borne_inf_0}, 
and prove:
\begin{equation}
\liminf_{N\rightarrow\infty}\liminf_{T\rightarrow\infty} 
\frac{1}{T}
\E^{\nu^N_{g_h}}_{T,h,\mathcal O}\big[-\log {\rm D}_h\big] 
=
-J_h(k_h)
=
- \mathcal I_{\epsilon_B}(k_h)
=
- \mathcal I_{\infty}(k_h)
,
\label{eq_ccl_lower_bound}
\end{equation}
where $\mathcal I_\infty$ is defined in~\eqref{eq_def_rate_function} and $\mathcal I_{\epsilon_B}(k_h)$ is short for $\mathcal I_{\epsilon_B}\big(\frac{1}{4}\big<k_h,\cdot\big>\big)$, 
idem for $\mathcal I_{\infty}(k_h)$ and $J_h(k_h)$. 
This concludes the proof of the lower bound for regular kernels close to that of the steady state.\\

Recall from Corollary~\ref{coro_comparaison_D_H_martingale} the definition of the martingale $M^{N,h}_T$. 
The average $\E^{\nu^N_{g_h}}_{T,h,\mathcal O}\big[-\log {\rm D}_h\big]$ reads:
\begin{align}
\E^{\nu^N_{g_h}}_{T,h,\mathcal O}\big[-\log {\rm D}_h\big]
&= 
 \Big(\Q_{T,h}^{\nu^N_{g_h}}(\mathcal O)\Big)^{-1}\E^{\nu^N_{g_h}}_{h}\bigg[{\bf 1}_{\mathcal O}\Big(-M^{N,h}_T - \frac{1}{2}\int_0^T\Pi^N_t\big( \mathcal M(\partial_1h,\partial_1h)\big)\, dt\Big)\bigg]
\nonumber\\
&\hspace{-0.5cm}
- \frac{T}{8}\int_{\squaredash}\bar\sigma(x)\bar\sigma(y)\big[\partial_1 h(x,y)\big]^2\, dx\, dy 
+ \frac{\E^{\nu^N_{g_h}}_{h}\Big[{\bf 1}_{\mathcal O}\int_0^T\hat \epsilon^N(h)(\eta_t)dt\Big]}{\Q_{T,h}^{\nu^N_{g_h}}(\mathcal O)}
,
\label{eq_to_estimate_dynamical_part}
\end{align}
with $\hat \epsilon^N(h)$ the error term defined in Corollary~\ref{coro_comparaison_D_H_martingale}. 
Proposition~\ref{coro_limite_hydro} establishes the convergence of $\Q^{\nu^N_{g_h}}_{T,h}(\mathcal O)$ to $1$ as $T$, then $N$ become large.
The error term $\hat\epsilon^N(h)$ is controlled through
Proposition~\ref{prop_Boltzmann_gibbs_sec3}:
\begin{equation}
\limsup_{N\rightarrow\infty}\sup_{T>0}\E^{\nu^N_{g_h}}_{h}\Big[\Big|\frac{1}{T}\int_0^T\hat \epsilon^N(h)(\eta_t)dt\Big|\Big] = 0
.
\end{equation}
In particular,
\begin{equation}
\limsup_{N\rightarrow\infty}\limsup_{T\rightarrow\infty}
\Bigg|\frac{\E^{\nu^N_{g_h}}_{h}\Big[{\bf 1}_{\mathcal O}\int_0^T\hat \epsilon^N(h)(\eta_t)dt\Big]}{\Q_{T,h}^{\nu^N_{g_h}}(\mathcal O)}\Bigg|
=0
.
\label{eq_error_term_lower_bound}
\end{equation}
Consider now the expectation in the first line of~\eqref{eq_to_estimate_dynamical_part}. 
By Cauchy-Schwarz inequality and the fact that $\big<M^{N,h}\big>_T$ grows at most linearly in $T$ up to a small error vanishing with $N$ (see~\eqref{eq_bound_quadratic_variation}), 
the contribution of the martingale term vanishes:
\begin{equation}
\limsup_{N\rightarrow\infty}\limsup_{T\rightarrow\infty}\frac{\E^{\nu^N_{g_h}}_{h}\big[{\bf 1}_{\mathcal O}\big|M^{N,h}_T\big|\big]}{T\Q_{T,h}^{\nu^N_{g_h}}(\mathcal O)} = 0
.
\end{equation}
Consider now the last remaining term in~\eqref{eq_to_estimate_dynamical_part}, 
the term involving $\Pi^N_\cdot(\mathcal M(\partial_1h,\partial_1h))$. 
To estimate it, we state the following moment bound, 
proven afterwards.
\begin{lemma}\label{lemm_moment_3_halves}
For any continuous and bounded $F:\squaredash\rightarrow\R$, 
there is $C(h,\|F\|_\infty)>0$ such that:
\begin{equation}
\sup_{N\in\N^*}\sup_{t\geq 0}\E^{\nu^N_{g_h}}_h\big[|\Pi^N_t(F)|^{3/2}\big] 
\leq 
C(h,\|F\|_{\infty})
.
\label{eq_moment_bound_3_2}
\end{equation}
\end{lemma}
Lemma~\ref{lemm_moment_3_halves} implies that, 
to prove convergence of the term involving $\Pi^N_\cdot(\mathcal M(\partial_1h,\partial_1h))$ in~\eqref{eq_to_estimate_dynamical_part}, 
it is enough to prove that it converges in probability. 
Proposition~\ref{coro_limite_hydro} shows convergence of ${\bf 1}_{\mathcal O}$ to $1$ in probability, 
while Proposition~\ref{prop_poisson_equation} implies convergence of the correlation field. 
We deduce:
\begin{align}
\limsup_{N\rightarrow\infty}&\limsup_{T\rightarrow\infty} 
\bigg|\E^{\nu^N_{g_h}}_h\Big[{\bf 1}_{\mathcal O}\frac{1}{T}\int_0^T\Pi^N_t\big(\mathcal M(\partial_1h,\partial_1h)\big)\, dt\Big] - \frac{1}{4}\big<k_h,\mathcal M(\partial_1 h,\partial_1h)\big>
\bigg| 
= 0
.
\label{eq_first_line_dynamical_part_computed}
\end{align}
From equations~\eqref{eq_error_term_lower_bound} to~\eqref{eq_first_line_dynamical_part_computed}, 
we obtain:
\begin{align}
&\liminf_{N\rightarrow\infty}\liminf_{T\rightarrow\infty}
\frac{1}{T}\E^{\nu^N_{g_h}}_{T,h,\mathcal O}\big[-\log {\rm D}_h\big] 
\nonumber\\
&\hspace{2cm}=  
-\frac{1}{8}\int_{\squaredash}dx dy\ \bar\sigma(x)\bar\sigma(y)\big[\partial_1 h(x,y)\big]^2 - \frac{1}{8}\big<k_h, \mathcal M(\partial_1 h,\partial_1 h)\big>
\nonumber\\
&\hspace{2cm}= 
- \frac{1}{8}\int_{(-1,1)}dz\, \bar\sigma(z) \big<\partial_1 h(z,\cdot),(\bar\sigma+k_h)\partial_1 h(z,\cdot)\big>
.
\label{eq_final_estimate_dynamical_part_lower_bound}
\end{align}
Let us check that the right-hand side in~\eqref{eq_final_estimate_dynamical_part_lower_bound} is indeed equal to $-\mathcal I_{\epsilon_B}(k_h)$, 
and that $\mathcal I_{\infty}(k_h) = \mathcal I_{\epsilon_B}(k_h)$.  
The weak formulation~\eqref{eq_EL} of the Euler-Lagrange equation on $k_h$ with test function $h$ gives:
\begin{equation}
J_h(k_h)
=
\frac{1}{8}\int_{(-1,1)}dz\, \bar\sigma(z) \big<\partial_1  h(z,\cdot),(\bar\sigma+k_h)\partial_1 h(z,\cdot)\big>
.
\label{eq_J_h_k_h}
\end{equation}
Again using~\eqref{eq_EL} for test functions $f\in\s(\epsilon_B)$, 
the rate function $\mathcal I_{\epsilon'}(k_h)
:=
\sup_{f\in\s(\epsilon')}J_{f}(k_h)
$ reads, for any $\epsilon'\geq \epsilon_B$:
\begin{align}
\mathcal I_{\epsilon'}(k_h)
&=
\sup_{f\in\s(\epsilon')}\bigg\{ \frac{1}{4}\int_{(-1,1)}dz\, \bar\sigma(z) \big<\partial_1 f(z,\cdot),(\bar\sigma+k_h)\partial_1 h(z,\cdot)\big> 
\nonumber\\
&\hspace{4cm}- \frac{1}{8}\int_{(-1,1)}dz\, \bar\sigma(z) \big<\partial_1 f (z,\cdot),(\bar\sigma+k_h)\partial_1 f(z,\cdot)\big>\bigg\}
\label{eq_formule_representation_rate_function_0}\\
&=
J_h(k_h)
-
\sup_{f\in\s(\epsilon')}\bigg\{ \frac{1}{8}\int_{(-1,1)}dz\, \bar\sigma(z) \big<\partial_1 (f-h)(z,\cdot),(\bar\sigma+k_h)\partial_1 (f-h)(z,\cdot)\big> \bigg\}
\nonumber\\
&=
J_h(k_h)
.
\nonumber
\end{align}
This in particular yields $J_h(k_h) = \mathcal I_{\epsilon_B}(k_h)=\mathcal I_\infty(k_h)$ and, 
together with~\eqref{eq_final_estimate_dynamical_part_lower_bound}--\eqref{eq_J_h_k_h}, proves~\eqref{eq_ccl_lower_bound}.
\begin{proof}[Proof of Lemma~\ref{lemm_moment_3_halves}]
Let $F:\squaredash\rightarrow\R$ be bounded. 
Let $t\geq 0$, and let $\epsilon>0$ that will eventually be chosen as $\epsilon=\frac{1}{2}$. 
The moment bound is obtained by a careful application of the entropy inequality, 
putting to good use the $O(N^{-1/2})$ estimate on the size of the relative entropy of Theorem~\ref{theo_entropic_problem}. 
Note that $|\Pi^N_t(F)|\leq \|F\|_\infty N$. As a result, fixing $c>0$ to be chosen later and applying the entropy inequality to $c\lambda {\bf 1}_{|\Pi^N(F)|>\lambda}$ for each $\lambda>1$ in the second line below:
\begin{align}
&\E^{\nu^N_{g_h}}_h\big[\big|\Pi^N_t(F)\big|^{1+\epsilon}\big] 
\leq 
1+ (1+\epsilon)\int_{1}^{\|F\|_\infty N} \lambda^{\epsilon}\Prob^{\nu^N_{g_h}}\big( |\Pi^N_t(F)|>\lambda\big)d\lambda\\
&\leq 
1+ (1+\epsilon)\int_{1}^{\|F\|_\infty N} c^{-1}\lambda^{-1+\epsilon} \Big[H(f_t\nu^N_{g_h}|\nu^N_{g_h}) + \log \Big(1+(e^{c\lambda}-1)\nu^N_{g_h}\big(|\Pi^N(F)|>\lambda\big)\Big)\Big]d\lambda
.
\nonumber
\end{align}
By Theorem~\ref{theo_entropic_problem}, 
$H(f_t\nu^N_{g_h}|\nu^N_{g_h}) \leq CN^{-1/2}$ for some $C=C(h,\rho_\pm)>0$. 
Moreover, 
by Corollary~\ref{coro_estimate_W_section_Ising}, 
the probability involving $\nu^N_{g_h}$ above is bounded by $C(h,\|F\|_\infty)e^{-c(h,\|F\|_\infty)\lambda}$ for some $C(h,\|F\|_\infty),c(h,\|F\|_\infty)>0$. 
Choosing $c=c(h,\|F\|_\infty)/2$, one obtains the existence of $C'(h,\|F\|_\infty)>0$ such that:
\begin{align}
\E^{\nu^N_{g_h}}_h\big[\big|\Pi^N_t(F)\big|^{1+\epsilon}\big] 
&\leq  
1+\frac{(1+\epsilon)C'(h,\|F\|_\infty)}{\epsilon N^{1/2-\epsilon}}  \nonumber\\
&\qquad+ 
\frac{2(1+\epsilon)}{c(h,\|F\|_\infty)}\int_1^{\|F\|_\infty N} \lambda^{-1+\epsilon}e^{-c(h,\|F\|_\infty)\lambda/2}\, d\lambda. 
\end{align}
The integral is bounded with $N$ whatever the choice of $\epsilon>0$. 
Overall, the right-hand side above is therefore bounded with $N$ as soon as $\epsilon\in(0,1/2]$, 
which yields~\eqref{eq_moment_bound_3_2}. 
This concludes the proof of the lemma. 
\end{proof}
\subsection{Towards a lower bound for non regular correlations}\label{sec_extension_non_regular}
In this section, we discuss how to extend the lower bound to non regular correlations close to $k_0$, 
i.e. to kernels $k\in\mathbb H^1(\squaredash)$ which solve the Euler-Lagrange equation~\eqref{eq_EL} associated with a possibly non regular bias $h$, 
and for which $\|\nabla (k-k_0)\|_{2}$ is small enough. 
As we shall see below, 
this extension is much simpler here than in the case of large deviations of the density for the open SSEP~\cite{Bertini2003}. 
This is due to the fact that the Euler-Lagrange equation~\eqref{eq_EL_section_perspectives} is linear as an equation with fixed $h$ and unknown $k$, 
and to the smallness assumption on admissible biases $h$. \\

Let us sketch an informal proof of this extension.  
Take a symmetric $k\in\mathbb H^1(\squaredash)$ with $k_{|\partial\square}=0$ and recall that the Euler-Lagrange equation~\eqref{eq_EL}, 
viewed as an equation with unknown $h$, 
reads for all test functions $\phi\in \mathbb H^1(\squaredash)$ with $\phi_{|\partial \square} =0$:
\begin{equation}
\frac{1}{2}\int_{\squaredash}\nabla (k-k_0)(x,y)\cdot \nabla\phi(x,y)\,dx\,dy 
+ \frac{1}{2}\int_{(-1,1)}\bar\sigma(z)\big<\partial_1 h(z,\cdot), C_k \partial_1\phi(z,\cdot)\big>\, dz 
=
0
.
\label{eq_EL_section_perspectives}
\end{equation}
For $k$ close to $k_0$, 
one has $\bar\sigma + k = (\bar\sigma + k_0) + (k-k_0)\geq \alpha\,\text{id}$ for some $\alpha>0$, 
and the second term is a coercive bilinear form for the norm $\|\nabla\cdot \|_2$. 
There is thus a unique symmetric $h\in\mathbb H^1(\squaredash)$ with $h_{|\partial\square}=0$ associated with $k$. 
Using~\eqref{eq_EL_section_perspectives} with test function $h$, 
this $h$ has small norm (recall $c(\rho_\pm)\leq \bar\sigma \leq 1/4$ and $\|\nabla h\|_2^2 = 2\|\partial_1 h\|_2^2$ as $h$ is symmetric):
\begin{equation}
\frac{\alpha c(\rho_\pm)}{2} \|\nabla h\|^2_2
\leq
\alpha \int_{\squaredash} \bar\sigma(z)\partial_1h(x,z)^2 \,dx\,dz
\leq 
\|\nabla h\|_2\|\nabla(k-k_0)\|_2
.
\label{eq_bound_h_sec_perspectives}
\end{equation}
One would expect the estimate~\eqref{eq_bound_h_sec_perspectives} to be enough to approximate $h$ in $\|\nabla\cdot\|_2$ norm by more regular biases. 
However, to avoid technical difficulties, 
our set $\s(\epsilon_B)$ of regular biases was defined with an assumption on the sup norm of derivatives rather than on $\|\nabla\cdot \|_2$. 
This assumption is technical and could be lifted at the cost of lengthening some arguments. 
Assuming therefore that one can take $(h_n)\in\s(\epsilon_B)^{\N}$ with $\lim_n\|\nabla(h_n-h)\|_2=0$, 
let $k_n$ denote the kernel associated with $h_n$ ($n\in\N$) through~\eqref{eq_EL_section_perspectives}. 
By linearity of Equation~\eqref{eq_EL_section_perspectives}, 
written for $k$ and $k_n$ with test function $k-k_n$, 
the bounds $C_k=\bar\sigma+k \leq c\, \text{id}$, $\bar\sigma\leq 1/4$ and $\|\nabla h_n\|^2_2=2\|\partial_1 h_n\|_2^2$ give:
\begin{equation}
\|\nabla(k-k_n)\|_2^2 
\leq
\frac{c}{8}\|\nabla(h-h_n)\|_2\|\nabla (k-k_n)\|_2
+\frac{1}{8}\|k-k_n\|_2\|\nabla (k-k_n)\|_2\|\nabla h_n\|_2
.
\end{equation}
By assumption $k$ is close to $k_0$. 
This means that $\|\nabla h\|_2$ is small by~\eqref{eq_bound_h_sec_perspectives}. 
Since $(\|\nabla h_n\|_2)_n$ converges to $\|\nabla h\|_2$, 
$\|\nabla h_n\|_2$ is also small for large enough $n$. 
Using Poincaré inequality on $\square$ with 0 Dirichlet boundary conditions on $\partial\square$ to bound $\|k-k_n\|_2$ by $(2\pi)^{-1}\|\nabla (k-k_n)\|_2$, 
we thus find:
\begin{equation}
\|\nabla(k-k_n)\|_2\Big(1-\frac{\|\nabla h_n\|_2}{16\pi}\Big)
\leq 
\frac{c}{8}\|\nabla(h-h_n)\|_2
.
\end{equation}
Assuming that $k$ is close enough to $k_0$ for the parenthesis on the left-hand side to be, say, larger than $1/2$ when $n$ is large enough; 
we conclude that $\lim_n\|\nabla(k_{n}-k)\|_2 =0$. 
This yields the alleged general lower bound:
\begin{equation}
\mathcal I_{\epsilon_B}(k_{n})
=
\frac{1}{8}\int_{(-1,1)}dz\, \bar\sigma(z)\big<\partial_1 h_{n}(z,\cdot),\big(\bar\sigma +k_n\big)\partial_1 h_n(z,\cdot)\big>
\underset{n\rightarrow\infty}{\longrightarrow}\mathcal I_{\epsilon_B}(k)
.
\end{equation}

\begin{appendices}
\renewcommand{\thesection}{\Alph{section}}

\section{Correlations and concentration under discrete Gaussian measures}\label{app_discrete_gaussian_measures}
In this section, we investigate the measures $\nu^N_g$, defined for $g:\square\rightarrow\R$ as follows:
\begin{equation}
\bar\nu^N 
=
\bigotimes_{i\in\Lambda_N}\text{Ber}(\bar\rho_i),
\qquad \forall\eta\in\Omega_N,\quad 
\nu^N_g(\eta) 
= 
\frac{1}{\mathcal Z^N_g} \exp\bigg[\frac{1}{2N}\sum_{i\neq j \in \Lambda_N}g_{i,j}\bar\eta_i\bar\eta_j\bigg]\bar\nu^N(\eta)
,
\label{eq_def_bar_nu_G}
\end{equation}
where the partition function $\mathcal Z^N_g$ is a normalising constant. 
For simplicity, we assume throughout that $g\in\mathbb L^2(\square)$ is a negative kernel, i.e.:
\begin{align}
\forall u\in\mathbb L^2((-1,1)),\qquad 
\int_{\squaredash} u(x) g(x,y)u(y)dxdy 
\leq 0
.
\label{eq_g_negative_kernel_appendix}
\end{align}
None of the result in Appendix~\ref{app_discrete_gaussian_measures} would be modified if we instead assumed that $g= g_-+ g_+$ where $\pm g_\pm$ is a positive kernel, 
and $\|g_+\|_{2}\leq c$ for a small enough constant $c$
.
\subsection{Bound on the partition function and correlations}\label{sec_correl_bar_nu_G}
\begin{lemma}[Bound on the partition function]\label{lemm_bound_Z_g}
Let $g:\square\rightarrow\R$ be a continuous, bounded and symmetric function, and assume that $g$ is a negative kernel (see~\eqref{eq_g_negative_kernel_appendix}). Then:
\begin{equation}
\sup_{N\geq 1}\mathcal Z^N_g <\infty,\qquad \mathcal Z^N_g := \bar\nu^N\Big[\exp\Big[\frac{1}{2N}\sum_{i\neq j\in\Lambda_N}\bar\eta_i\bar\eta_jg_{i,j}\Big]\Big].\label{eq_bound_partition_function}
\end{equation}
\end{lemma}
\begin{proof}
As $g$ is continuous and a negative kernel, the matrix $(g_{i,j})_{(i,j)\in\Lambda_N^2}$ has negative eigenvalues, and one has:
\begin{equation}
\forall \eta\in\Omega_N,\qquad 2\Pi^N(g) + \frac{1}{2N}\sum_{i\in\Lambda_N}(\bar\eta_i)^2 g_{i,i} = \frac{1}{2N}\sum_{(i,j)\in\Lambda_N^2}\bar\eta_i\bar\eta_j g_{i,j} \leq 0.
\end{equation}
The continuity of $g$ on the diagonal and the bound $|\bar\eta_\cdot|\leq 1$ imply:
\begin{align}
\mathcal Z^N_g  
&= 
\bar\nu^N\Big[\exp\Big[\frac{1}{2N}\sum_{(i,j)\in\Lambda_N^2}\bar\eta_i\bar\eta_j g_{i,j} - \frac{1}{2N}\sum_{i\in\Lambda_N}(\bar\eta_i)^2g_{i,i}\Big]\Big]\nonumber\\
&\leq e^{\|g\|_\infty}.
\end{align}
\end{proof}
\begin{lemma}\label{lemm_bound_correlations_bar_nu_G}
Let $g:\square\rightarrow\R$ satisfy the hypotheses of Lemma~\ref{lemm_bound_Z_g}. Then:
\begin{equation}
\forall n\in\N^*,\qquad \sup_{I\subset\Lambda_N : |I|=n}\Big|\nu^N_g\Big(\prod_{a\in I}\bar\eta_a\Big) \Big|= 
O(N^{-n/2})
.
\label{eq_bound_correl_dans_lemm_bound_correl}
\end{equation}
Moreover, if $g\in C^1(\bar\rhd)\cap C^1(\bar\lhd)$ with $g(\pm 1,\cdot) = 0$, then, for each $n\in\N^*$ and each $\epsilon\in\{+,-\}$:
\begin{equation}
\sup_{\substack{I\subset\Lambda_N\\ \epsilon(N-1)\in I\text{ and }|I|=n}}\big|\nu^N_g\Big(\bar\eta_{\epsilon(N-1)}\prod_{a\in I\setminus\{\epsilon(N-1)\}}\bar\eta_a\Big)\big| = 
O(N^{-n/2-1}).\label{eq_correlations_bord_dans_lemm_bound_correl}
\end{equation}
\end{lemma}
\begin{proof}
Let $I_n\subset \Lambda_N$ with $|I|=n\in\N^*$. The proof relies on a development of the exponential defining $\nu^N_g$ in~\eqref{eq_def_bar_nu_G}, and the observation that the measure $\bar\nu$ is product, so that, for each $p\in\N^*$,
\begin{equation}
\forall(i_1,...,i_p)\in\Lambda_N^p,
\qquad 
|\{i_1,...,i_p\}|>p/2
\quad \Rightarrow\quad 
\bar\nu^N\Big(\prod_{j=1}^p\bar\eta_{i_j}\Big) = 0
.
\label{eq_moyenne_sous_bar_nu_is_zero}
\end{equation}
In other words, in a product $\bar\eta_{i_1}...\bar\eta_{i_n}$, each $\bar\eta$ must be paired with at least another $\bar\eta$ with the same index $i_1,...,i_n$. Since $\bar\eta_\cdot$ is bounded, a direct consequence is that, if $1\leq p\leq n$ and $w^p = (w^p_{i})_{i\in\Lambda_N}$ satisfies $\sup_{\Lambda_N}|w^p|=O_N(1)$, then:
\begin{equation}
\bar\nu^N\Big[\prod_{p=1}^n \Big(\frac{1}{N}\sum_{j\in\Lambda_N}w^p_j\bar\eta_j\Big)\Big] = O(N^{-\lceil n/2\rceil}).\label{eq_moment_fluct_sous_bar_nu}
\end{equation}
Write $I_p = \{i_1,...,i_p\}$ for each $1\leq p\leq n$ (and by extension write $I_{i_p} := I_p$). 
Define, for $J\subset\Lambda_N$:
\begin{equation}\label{eq_def_G_minus_I}
G_{J^c}(\eta) := \frac{1}{4N}\sum_{i\neq j \in \Lambda_N\setminus\{ J\}}\bar\eta_i\bar\eta_jg_{i,j},
\qquad \eta\in\Omega_N
.
\end{equation}
Then, for each $\eta\in\Omega_N$:
\begin{align}
2\Pi^N(g) &= 2G_{I_1^c}(\eta) + \bar\eta_{i_1}\sum_{j\neq i_1}\bar\eta_j g_{i_1,j}\nonumber\\
&= 2G_{I_p^c}(\eta) + \frac{1}{N}\sum_{q=1}^p \bar\eta_{i_q}\sum_{j\in\Lambda_N\setminus I_q}\bar\eta_j g_{i_q,j},\qquad 1\leq p \leq n.\label{eq_splitting_Pi_de_g_en_G_minus_a_et_reste}
\end{align}
As a result, one can write:
\begin{equation}\label{eq_def_m_one_point_0}
\mathcal Z^N_g \nu^N_g\Big(\prod_{p=1}^n\bar\eta_{i_p}\Big) 
= 
\bar\nu^N\bigg[\Big(\prod_{p=1}^n\bar\eta_{i_p}\Big)e^{2G_{I_n^c}(\eta)}\exp\Big[ \frac{1}{N}\sum_{p=1}^n\bar\eta_{i_p}\sum_{j\in\Lambda_N\setminus I_p}\bar\eta_j g_{i_p,j}\Big]\bigg]
,
\end{equation}
and $G_{I_n^c}$ does not involve any $\bar\eta_{i_p}$ ($1\leq p \leq n$). 
Note also that $\sup_N\mathcal Z^N_g<\infty$  by Lemma \ref{lemm_bound_Z_g}, 
thus establishing Lemma~\ref{lemm_bound_correlations_bar_nu_G} only boils down to estimating the right-hand side above. 
To do so, 
recall first the following identity:
\begin{align}
\forall x\in\R,\qquad e^{x} 
= 
\sum_{p=0}^{n-1} \frac{x^m}{m!} + \int_0^1 \frac{x^{m}}{(m-1)!} (1-t)^{m}e^{tx}dt
.
\end{align}
Expanding the second exponential in~\eqref{eq_def_m_one_point_0} to order $n$ thus yields, writing $n_0 := \lfloor n/2\rfloor$:
\begin{align}
\mathcal Z^N_g \nu^N_g&\Big(\prod_{p=1}^n\bar\eta_{i_p}\Big) 
= 
\bar\nu^N\bigg[\Big(\prod_{p=1}^n\bar\eta_{i_p}\Big)e^{2G_{I_n^c}(\eta)} \Big\{T_1(\eta) + T_2(\eta) + T_3(\eta)\Big\}\bigg]
,
\label{eq_splitting_average_T_1T_2_T3}
\end{align}
where:
\begin{align}
T_1(\eta) &:= 
\sum_{m=0}^{n_0-1}\frac{1}{m!}\Big(\frac{1}{N}\sum_{q=1}^n\bar\eta_{i_q}\sum_{j\in\Lambda_N\setminus I_q}\bar\eta_j g_{i_q,j}\Big)^m \label{eq_first_line_bound_correl}
,\\
T_2(\eta) 
&:= 
\sum_{m=n_0}^{n-1}\frac{1}{m!}\Big(\frac{1}{N}\sum_{q=1}^n\bar\eta_{i_q}\sum_{j\in\Lambda_N\setminus I_q}\bar\eta_j g_{i_q,j}\Big)^m \label{eq_second_line_bound_correl}
,\\
T_3(\eta) 
&:=
\int_0^1 \frac{(1-t)^{n-1}}{(n-1)!}\Big(\frac{1}{N}\sum_{q=1}^{n}\bar\eta_{i_q}\sum_{j\in\Lambda_N\setminus I_q}\bar\eta_j g_{i_q,j}\Big)^{n} \exp\Big[\frac{t}{N}\sum_{q=1}^n\bar\eta_{i_q}\sum_{j\in\Lambda_N\setminus I_q}\bar\eta_j g_{i_q,j}\Big]dt.\label{eq_third_line_bound_correl}
\end{align}
Let us estimate the average of $T_1,T_2, T_3$ separately. 
Consider first $T_1$ in~\eqref{eq_first_line_bound_correl}. 
For each $0\leq m\leq n_0-1$, developing the term elevated to the power $m$ yields products of at most $2m = 2n_0-2<n$ $\bar\eta$'s. 
In view of Equation~\eqref{eq_moyenne_sous_bar_nu_is_zero}, this is not enough to achieve a pairing of all indices in $\prod_{p=1}^n \bar\eta_{i_p}$, 
thus the average of $T_1$ vanishes identically. \\
Let us count the factors of $\bar\eta$'s appearing in $T_2$ in~\eqref{eq_second_line_bound_correl}. 
For each $n_0\leq m <n$, the term to the power $m$ in~\eqref{eq_second_line_bound_correl} is a sum of terms of the form:
\begin{align}
\prod_{i\in K}\Big(\frac{\bar\eta_{i}}{N}\sum_{j\in\Lambda_N\setminus I_i}\bar\eta_j g_{i,j}\Big)^{n_{i}},
\qquad 
K\subset I_n,(n_{i})_{i\in K}\in\N^{|K|}\text{ with }1\leq |K|\leq m,\ \sum_{i\in K}n_{i} = m.\label{eq_forme_terme_developpe_second_line_bound_correl}
\end{align}
Again by~\eqref{eq_moyenne_sous_bar_nu_is_zero}, 
each term appearing in~\eqref{eq_forme_terme_developpe_second_line_bound_correl}
gives a non-vanishing contribution to the average~\eqref{eq_splitting_average_T_1T_2_T3} provided each element $i_1,...,i_n$ of $I_n$ appears at least once in the expansion. 
Since the term~\eqref{eq_forme_terme_developpe_second_line_bound_correl} contains the product $\prod_{i\in K}\bar\eta_{i}$ with $|K|<n$,
$n-|K|$ indices $j$ of the sums in~\eqref{eq_forme_terme_developpe_second_line_bound_correl} must be fixed to an element of $I_n\setminus K$. 
A contribution $1/N$ arises each time an index $j$ is singled out in a sum. 
It follows that the contribution of $T_2$ to the average~\eqref{eq_splitting_average_T_1T_2_T3} is bounded by a sum of terms of the form:
\begin{equation}\label{eq_to_estimate_second_line_bound_correl_0}
\sum_{m=n_0}^{n-1}\frac{C(m)}{ N^{n-|K|}}\Big(\prod_{i\in I_n\setminus K}\|g_{i,\cdot}\|_\infty\Big)\bar\nu^N\bigg[e^{2G_{I_n^c}(\eta)}\prod_{i\in K}\Big|\frac{1}{N}\sum_{j\in\Lambda_N\setminus I_i}\bar\eta_jg_{i,j}\Big|^{n'_{i}}\bigg]
,
\end{equation}
where $C(m)>0$, $K\subset I_n$, and the $n'_{i}$ satisfy:
\begin{equation}
0\leq n'_{i}\leq n_{i} \quad \text{for } i\in K,
\quad \text{and}\quad 
\sum_{i\in K}n'_{i} = (m+|K|-n){\bf 1}_{m+|K|-n\geq 0}.
\end{equation}
Let us now estimate the average in~\eqref{eq_to_estimate_second_line_bound_correl_0}. 
Recall that $g$ is continuous, and a negative kernel 
(as defined in~\eqref{eq_g_negative_kernel_appendix}). 
As a result, $(g_{i,j})_{i,j\in \Lambda_N\setminus I_n}$ has negative eigenvalues, and for each $\eta\in\Omega_N$:
\begin{align}
2G_{I^c_n}(\eta) &= \frac{1}{2N}\sum_{i,j\in \Lambda_N\setminus I_n}\bar\eta_i\bar\eta_j g_{i,j} - \frac{1}{2N}\sum_{i\in\Lambda_N\setminus I_n}(\bar\eta_i)^2g_{i,i}\nonumber\\
&\leq - \frac{1}{2N}\sum_{i\in\Lambda_N\setminus I_n}(\bar\eta_i)^2g_{i,i}\leq \|g\|_\infty.\label{eq_g_restreint_hors_I_forme_quadra_negative}
\end{align}
Using then Cauchy-Schwarz inequality followed by the estimate~\eqref{eq_moment_fluct_sous_bar_nu} 
for moments of correlations under the product measure $\bar\nu^N$, 
the average in~\eqref{eq_to_estimate_second_line_bound_correl_0} is  bounded by:
\begin{align}\label{eq_to_estimate_second_line_bound_correl_1}
\bar\nu^N\bigg[e^{2G_{I^c_n}(\eta)}\prod_{i\in K}\Big|\frac{1}{N}\sum_{j\in\Lambda_N\setminus K}\bar\eta_jg_{i,j}\Big|^{n'_i}\bigg] 
&\leq 
e^{\|g\|_\infty}\bar\nu^N\bigg[\prod_{i\in K}\Big(\frac{1}{N}\sum_{j\in\Lambda_N\setminus K}\bar\eta_jg_{i,j}\Big)^{2n'_i}\bigg]^{1/2} 
\nonumber\\
&= 
{\bf 1}_{m+|K|-n\geq 0}\, O(N^{-(m+|K|-n)/2})
.
\end{align}
Putting together~\eqref{eq_to_estimate_second_line_bound_correl_0}-\eqref{eq_to_estimate_second_line_bound_correl_1} and summing on all possible choices for $K,(n_q)_{q\in K},(n'_q)_{q\in K}$ yields:
\begin{align}
\bigg|\bar\nu^N\bigg[e^{2G_{I^c_n}(\eta)}\Big(\prod_{q\in I}\bar\eta_{i_q}\Big)\cdot T_2(\eta)\bigg]\bigg| 
&= 
\sum_{m=n_0}^{n-1}\sum_{K:1\leq |K|\leq m}\sum_{(n_i),(n'_i)} O(N^{-(n-|K|)})O(N^{-(m+|K|-n)/2}) \nonumber\\
&= 
O(N^{-n/2}).\label{eq_estimate_second_line_bound_correl}
\end{align}
It remains to estimate the contribution of the average of $T_3$, 
defined in~\eqref{eq_third_line_bound_correl}. 
Since there are already $n$ sums involving $\bar\eta$, 
it is enough to use the following bound, 
valid for each $\eta\in\Omega_N$:
\begin{equation}
\bigg|\Big(\prod_{p=1}^n\bar\eta_{i_p}\Big) T_3(\eta)\bigg| 
\leq 
\sum_{q=1}^n\Big|\frac{1}{N}\sum_{j\in\Lambda_n\setminus I_q}\bar\eta_j g_{i_q,j}\Big|^n 
\int_0^1 \frac{(1-t)^{n-1}}{(n-1)!} \exp\Big[\frac{t}{N}\sum_{q=1}^n\bar\eta_{i_q}\sum_{j\in\Lambda_N\setminus I_q}\bar\eta_j g_{i_q,j}\Big]dt.
\end{equation}
The integral on $t$ is taken care of with the following identity: for each $t\geq 0$ and $\eta\in\Omega_N$,
\begin{equation}
2G_{I^c_n}(\eta) + \frac{t}{N}\sum_{q=1}^n\bar\eta_{i_q}\sum_{j\in\Lambda_N\setminus I_q}\bar\eta_j g_{i_q,j} = 2\Pi^N(g) + \frac{(t-1)}{N}\sum_{q=1}^n\bar\eta_{i_q}\sum_{j\in\Lambda_N\setminus I_q}\bar\eta_j g_{i_q,j}
.
\label{eq_completion_G_I_c_into_G}
\end{equation}
Indeed, since $g$ is a negative kernel, it follows that:
\begin{align}
\bigg|\bar\nu^N\bigg[e^{2G_{I_n^c}(\eta)}\Big(\prod_{q=1}^n\bar\eta_{i_q}\Big) \, T_3\bigg]\bigg| 
&\leq 
\frac{e^{(n+1)\|g\|_\infty}}{n!}\bar\nu^N\bigg[\Big|\frac{1}{N}\sum_{q=1}^{n}\bar\eta_{i_q}\sum_{j\in\Lambda_N\setminus I_q}\bar\eta_j g_{i_q,j}\Big|^n\bigg]\nonumber\\
&= O(N^{-n/2}).\label{eq_estimate_third_line_bound_correl}
\end{align}
Putting~\eqref{eq_estimate_second_line_bound_correl}-\eqref{eq_estimate_third_line_bound_correl} together yields the first part of Lemma~\ref{lemm_bound_correlations_bar_nu_G}, 
i.e.~\eqref{eq_bound_correl_dans_lemm_bound_correl}. \\

Let us now prove the improved estimate~\eqref{eq_correlations_bord_dans_lemm_bound_correl} if correlations include points at the boundaries. 
If $g$ is assumed to be $C^1$ with $g(\pm 1,\cdot) = 0$, then $\|g_{\pm(N-1),\cdot}\|_\infty = O(N^{-1})$. 
By~\eqref{eq_moyenne_sous_bar_nu_is_zero}, the index $\pm(N-1)$ must arise in the terms~\eqref{eq_first_line_bound_correl}-\eqref{eq_third_line_bound_correl} to compensate the $\bar\eta_{\pm(N-1)}$ present in the product $\prod_{q\in I}\bar\eta_{i_q}$. 
It follows that the $O(N^{-n/2})$ bound in~\eqref{eq_estimate_second_line_bound_correl}-\eqref{eq_estimate_third_line_bound_correl} is improved to $O(N^{-n/2-1})$, which concludes the proof of Lemma~\ref{lemm_bound_correlations_bar_nu_G}.
\end{proof}
\subsection{Exponential moments of higher order correlations}\label{sec_Ising_measure}
Let $g\in g_0+\s(\infty)$ be a negative kernel, 
where this set is defined in~\eqref{eq_def_s_infty}. 
In this section, we give bounds on the size of exponential moments, 
under $\nu^N_g$, of random variables involving $n$-point correlations, $n\geq 1$. 
These are useful when applying the entropy inequality.\\
Such concentration results are established in the literature by means of a logarithmic Sobolev inequality, 
see~\cite{Goetze2019}-\cite{Sambale2020}. 
That $g$ be a negative kernel implies that, 
for any $F: \Omega_N\rightarrow\R$:
\begin{equation}
\nu^N_g\big(\exp[F]\big) 
\leq 
e^{\|g\|_\infty}\big(\mathcal Z^N_g\big)^{-1}\bar\nu^N\big(\exp[F]\big)
,
\label{eq_concentration_under_product_measure_is_enough}
\end{equation}
and it is enough to estimate exponential moments under the product measure $\bar\nu^N$. 
To do so, let us fix some notations. 
For $d \in\N$, let $A:\Lambda_N^d\rightarrow\R$ be a tensor. 
Define its Hilbert-Schmidt norm by:
\begin{equation}
\|A\|_{HS} 
= 
\Big[\sum_{(i_0,...,i_{d-1})\in\Lambda_N^d} \big(A(i_0,...,i_{d-1})\big)^2\Big]^{1/2}.\label{eq_def_Hilbert_Schmidt_norm}
\end{equation}
For $J\subset \Z$ containing $0$, let $X^A_{d,J}(\eta)$ be defined for $\eta\in\Omega_N$ as:
\begin{equation} 
X^A_{d,J}(\eta) 
= 
\sum_{\substack{(i_0,...,i_{d-1}) \in\Lambda_N^d \\ i_0+J\subset\Lambda_N}} A(i_0,...,i_{d-1})\bar\eta_{i_0+J}\prod_{p=1}^{d-1}\bar\eta_{i_p},
\qquad 
\bar\eta_{i_0+J} := \prod_{j\in J}\bar\eta_{i_0+j}
.
\label{eq_def_X_A_d_J_appendice}
\end{equation}
The next theorem gives concentration estimates of $X^A_{d,J}$ under $\bar\nu^N$ for $J\subset \Z$ and $d\in\N^*$. 
The case $J = \{0\}$, $d\leq 4$ corresponds to Theorem 1.4 in~\cite{Goetze2019}, 
but their proof extends to any $d\in\N^*$. 
However, we also use in the article the case $J=\{0,1\}$, 
for which a proof is needed.
\begin{theorem}\label{theo_concentration_Ising}
Let $J\subset\Z$ contain $0$ and $d\in\N^*$. 
Assume that $A$ is such that $A(i_0,...,i_{d-1})$ vanishes whenever the same site appears twice in $\bar\eta_{i_0+J}\prod_{p=1}^{d-1}\bar\eta_{i_p}$ for $\eta\in\Omega_N$, 
i.e. assume:
\begin{equation}
\forall (i_0,...,i_{d-1})\in\Lambda_N^{d},\quad 
\Big(\exists j\in J, |\{i_0+j,i_1,...,i_{d-1}\}|<d\Big)
\ \Rightarrow\ 
A(i_0,...,i_{d-1}) = 0
.
\label{eq_assumption_diagonal_tensor}
\end{equation}
There are then constants $c_d>0$ depending only on $d$ such that, 
for any $c\in(0,c_d)$ and any $N$ with $J\subset\Lambda_N$:
\begin{equation}
\bar\nu^N \bigg(\exp\bigg[ \frac{c|X^A_{d,J}|^{2/d}}{\|A\|_{HS}^{2/d}}\bigg]\bigg)
\leq 
2.\label{bound_exp_moment_power_2_over_d}
\end{equation}
\end{theorem} 
\begin{proof}
The proof for general $J$ and the $J=\{0\}$ case in~\cite{Goetze2019} are very similar, so we only give a sketch. Without loss of generality, $A$ can be assumed to be invariant under permutation of its last $d-1$ indices. The idea is to proceed by recursion on $d$, noticing that, for each $\ell\in\Lambda_N$ and $\eta\in\Omega_N$:
\begin{align}
\nabla_\ell X^A_{d,J}(\eta) :&= X^A_{d,J}(\eta^\ell) - X^A_{d,J}(\eta) 
\nonumber\\
&= 
(1-2\eta_\ell)
\Big[\sum_{i_0: i_0+J\ni \ell}\sum_{i_1,...,i_{d-1}}A(i_0,...,i_{d-1})\bar\eta_{(i_0 + J)\setminus\{\ell\}}\prod_{a=1}^{d-1}\bar\eta_{i_a}  
\nonumber\\
&\hspace{3cm}
+ (d-1)\sum_{\substack{i_0,...,i_{d-2} \\ i_0+J\subset\Lambda_N}}A(i_0,...,i_{d-2},\ell)\bar\eta_{i_0+J}\prod_{a=1}^{d-2}\bar\eta_{i_a}\Big] 
.
\label{eq_nabla_X_d}
\end{align}
Fix $d\in\N$, $J\subset\Z$ with $0\in J$, and $N$ with $J\subset\Lambda_N$. For brevity, simply write $X_d$ for $X^A_{d,J}$.\\

\noindent\textbf{Step 1: reduction to moment bound.} To prove \eqref{bound_exp_moment_power_2_over_d}, it is enough to prove the existence of $C_d>0$ such that:
\begin{equation}
\forall p\in\N^*,\qquad M_p(X_d) 
:= 
\bar\nu^N\big[|X_d|^p\big]^{1/p} \leq C_d \|A\|_{HS}\  p^{\frac{d}{2}}
.\label{eq_bound_moment_desired_for_bound_X_d}
\end{equation}
Indeed, assuming such a bound, one has, for each $c>0$, using Jensen inequality when $d\geq 2$ for the convex function $f(x) = x^{d/2}$, $x\in\R$:
\begin{align}
\bar\nu^N\big[\exp \big(c|X_d|^{2/d}/\|A\|_{HS}^{2/d}\big)\big] 
&= 
1 +\sum_{p=1}^\infty \frac{c^p}{p!}\bar\nu^N\big(|X_d|^{2p/d}/\|A\|_{HS}^{2p/d}\big) 
\nonumber\\
&\leq 
1+\sum_{p=1}^\infty \frac{\big(c C_d^{2/d} p\big)^p}{p!}
.
\end{align}
As $p^p\leq p!e^p$ for each $p\geq 1$, taking $c\leq c_d := (2C_d^{2/d}e )^{-1}$ yields \eqref{bound_exp_moment_power_2_over_d}.\\

\noindent\textbf{Step 2: moment estimate.} 
It is enough to prove \eqref{eq_bound_moment_desired_for_bound_X_d} for $p\geq 2$, since the first moment can be estimated by Cauchy-Schwarz inequality. 
We will restrict to $p\geq 2$ at some point in the computation. 
For now, we treat $p$ as a continuous variable in $\R_+^*$ and differentiate $M_\cdot(X_d)$. 
For each $p>0$, one has:
\begin{align}
\frac{dM_p(X_d)}{dp} 
&= 
\frac{d}{dp}\Big(\exp\Big[\frac{1}{p}\log\bar\nu^N\big(|X_d|^p\big)\Big]\Big) \nonumber\\
&= 
-\frac{\log \bar\nu^N\big(|X_d|^p\big)}{p^2}M_p(X_d)+\frac{1}{p}\frac{\bar\nu^N\big(|X_d|^p\log|X_d|\big)}{\bar\nu^N\big(|X_d|^p\big)}M_p(X_d)
,
\end{align}
which can be written as:
\begin{equation}
\forall p> 0,\qquad \frac{dM_p(X_d)}{dp} 
= 
\frac{M_p(X_d)^{1-p}}{p^2}\text{Ent}(|X_d|^p)
,\label{eq_derivee_moment_preuve_concentration}
\end{equation}
with $\text{Ent}(F^2)$ the entropy of $F^2$ against $\bar\nu^N$, given by:
\begin{equation}
\forall F:\Omega_N\rightarrow\R,\qquad \text{Ent}(F^2) 
= 
\bar\nu^N(F^2\log F^2) - \bar\nu^N(F^2)\log\bar\nu^N( F^2)
.
\end{equation}
The entropy on the right-hand side of \eqref{eq_derivee_moment_preuve_concentration} is estimated by means of a logarithmic Sobolev inequality, 
satisfied by $\bar\nu^N$ for the Glauber dynamics on $\Omega_N$ (see e.g. Theorem A.1. in~\cite{Diaconis1996}): 
there is $C_{LS}>0$, independent of $N$, such that:
\begin{equation}
\forall F:\Omega_N\rightarrow\R_+,\qquad 
\text{Ent}(F^2)
\leq 
C_{LS} \bar\nu^N \Big(\sum_{i\in\Lambda_N} \big[\nabla_i F\big]^2\Big)
,
\label{eq_log_sob_nu_bar}
\end{equation}
where, for $i\in\Lambda_N$ and $\eta\in\Omega_N$, 
$\nabla_i F(\eta) = F(\eta^i) - F(\eta)$ for each $F:\Omega_N\rightarrow\R$. 
Similarly, a Poincare inequality holds with constant $C_{LS}/2$:
\begin{equation}
\forall F:\Omega_N\rightarrow\R_+,\qquad 
\bar\nu^N(F^2)-\bar\nu^N(F)^2 
\leq 
\frac{C_{LS}}{2} \bar\nu^N \Big(\sum_{i\in\Lambda_N} \big[\nabla_i F\big]^2\Big)
.
\label{eq_Poincare_nu_bar}
\end{equation}
Injecting \eqref{eq_log_sob_nu_bar} in \eqref{eq_derivee_moment_preuve_concentration} and proceeding as in~\cite{Goetze2019}, 
one successively obtains, restricting to $p>2$:
\begin{align}
\forall p> 2,\qquad \frac{dM_p(X_d)^2}{dp} 
&\leq 
\frac{2C_{LS} M_p(X_d)^{2-p}}{p^2}\bar\nu^N\Big(\sum_{i\in\Lambda_N} \big[\nabla_i (|X_d|^{p/2})\big]^2\Big)\nonumber\\
&\leq 
C_{LS}M_p(X_d)^{2-p}\bar\nu^N\Big(|X_d|^{p-2}\sum_{i\in\Lambda_N} \big[\nabla_i |X_d|\big]^2\Big)
.
\end{align}
Applying Hölder inequality with exponents $(p/(p-2),p/2)$ then yields, 
for each $p>2$:
\begin{equation}
\frac{dM_p(X_d)^2}{dp}
\leq 
C_{LS}M_{p/2}\Big(\sum_{i\in\Lambda_N} \big[\nabla_i |X_d|\big]^2\Big) 
\leq 
C_{LS} M_{p/2}\Big(\sum_{i\in\Lambda_N} \big[\nabla_i X_d\big]^2\Big)
.
\end{equation}
The function $M_\cdot(X_d)$ is increasing for $p>0$ by \eqref{eq_derivee_moment_preuve_concentration}. 
As a result, 
integrating between $2$ and $p$ and using the Poincare inequality \eqref{eq_Poincare_nu_bar} to estimate $M_2(X_d)$ yields:
\begin{align}
\forall p\geq 2,\qquad 
M_p(X_d)^2
&\leq 
M_2(X_d)^2 + C_{LS}(p-2) M_{p/2}\Big(\sum_{i\in\Lambda_N} \big[\nabla_i X_d\big]^2\Big)\nonumber\\
&\leq 
C_{LS}\ p M_{p/2}\Big(\sum_{i\in\Lambda_N} \big[\nabla_i X_d\big]^2\Big)
.
\label{eq_recursion_moments_X_d}
\end{align}
\noindent\textbf{Step 3: recursion on $d$.} Let $p\geq 2$. 
For $j\in\Lambda_N$ and $0\leq a\leq d-1$, 
define $A^{(i_a=j)}$ as the $d-1$-tensor $\big(A(i_0,...,i_{a-1},j,i_{a+1},...,i_{d-1})\big)_{(i_q)_{q\neq a}}$, 
and note that $A^{(i_a=j)}$ also satisfies the assumption \eqref{eq_assumption_diagonal_tensor}. 
Recall that $X_d$ was short for $X^A_{d,J}$. 
Let us prove by recursion on $d$:
\begin{equation}
\forall p\geq 2,\qquad 
M_p\big(X^A_{d,J}\big)^2 
\leq 
2^{d-1}(C_{LS}\ p)^d (d!)^2|J|^2\|A\|_{HS}^2
,
\label{eq_bound_moments_X_d}
\end{equation}
%
Proving such a result would conclude the proof. 
In the $d=1$ case, 
\eqref{eq_recursion_moments_X_d} yields, 
for each $p\geq 2$:
\begin{equation}
M_p(X_{1,J}^A)^2 
\leq 
C_{LS}\ p M_{p/2}\Big(\sum_{i\in\Lambda_N} \big[\nabla_i X_{1,J}^A\big]^2\Big)
\leq 
C_{LS}\ p \sum_{i\in\Lambda_N} M_{p}\big(\nabla_\ell X_{1,J}^A\big)^2
.
\end{equation}
By \eqref{eq_nabla_X_d}, 
bounding $\bar\eta_\cdot$ by $1$, 
the result for $d=1$ is proven:
\begin{equation}
\forall p\geq 2,\qquad 
M_p(X_{1,J}^A)^2 
\leq 
C_{LS}\ p \sum_{\ell\in\Lambda_N}\Big[\sum_{i_0+J\ni\ell}|A(i_0)|\Big]^2 \leq C_{LS}\ p |J|^2\|A\|_{HS}^2
.
\end{equation}
For $d\geq 2$, \eqref{eq_recursion_moments_X_d} similarly gives, 
bounding $\bar\eta_\cdot$ by $1$:
\begin{align}
M_{p/2}\Big(\sum_{\ell\in\Lambda_N} \big[\nabla_\ell X^A_{d,J}\big]^2\Big) 
&\leq 
\sum_{\ell\in\Lambda_N} M_p\big(\nabla_\ell X^A_{d,J}\big)^2 \nonumber\\
&\leq 
\sum_{\ell\in\Lambda_N}\Big[(d-1) M_p\big(X^{A^{(i_{d-1}=\ell)}}_{d-1,J}\big) + \sum_{i : i+J\ni \ell} M_p\big(X^{A^{(i_0=i)}}_{d-1,\{0\}}\big)\Big]^2
.
\end{align}
For each $\ell\in\Lambda_N^d$, 
the recursion hypothesis at rank $d-1$ applied to $X^{A^{(i_{d-1}=\ell)}}_{d-1,J}$, 
and to $X^{A^{(i_0=i)}}_{d-1,\{0\}}$ for each $i + J\ni\ell$, yields:
\begin{align}
M_{p/2}&\Big(\sum_{\ell\in\Lambda_N} \big[\nabla_\ell X^A_{d,J}\big]^2\Big)\nonumber\\
& \leq 
2^{d-2}(C_{LS}\ p)^{d-1} ((d-1)!)^2\sum_{\ell\in\Lambda_N}\Big[|J|(d-1)\|A^{(i_{d-1}=\ell)}\|_{HS} + \sum_{i:i +J\ni\ell} \|A^{(i_0=i)}\|_{HS}\Big]^2\nonumber\\
& \leq 
2^{d-1}(C_{LS}\ p)^{d-1} ((d-1)!)^2d^2|J|^2\|A\|_{HS}^2 = 2^{d-1}(C_{LS}p)^{d-1} (d!)^2|J|^2\|A\|_{HS}^2
,
\end{align}
where we used $(a+b)^2\leq 2a^2 + 2b^2$ for $a,b\in\R$ and bounded $(d-1)^2+1$ by $d^2$ to obtain the second inequality. Injecting this bound in \eqref{eq_recursion_moments_X_d} yields \eqref{eq_bound_moments_X_d} at rank $d$, concluding the proof.
\end{proof}
In the next corollary, 
we use Theorem~\ref{theo_concentration_Ising} to establish the controllability results of Lemma~\ref{lemm_size_controllable} on the variables $X^{A}_{d,J}$. 
\begin{corollary}\label{coro_estimate_W_section_Ising}
Let $J\subset\Z$ contain $0$ and, 
for $N\in\N^*$ and $d\in\N^*$, 
let $A = A(d,N)$ be a $d$-tensor. 
We do not assume that $A$ satisfy~\eqref{eq_assumption_diagonal_tensor}, 
but instead that $\sup_{N\in\N^*}\|A\|_\infty<\infty$. 
There are constants $\gamma_d,C_d(g)>0$ that are independent of $A$, 
with $\gamma_1=+\infty$, such that, 
for each $\gamma\leq \gamma_d$:
\begin{align}
&\text{for each }N\in\N^*\text{ with }J\subset\Lambda_N,\nonumber\\
&\hspace{2.5cm}\quad 
\begin{cases}
\displaystyle{\log\nu^N_g\Big(\exp\Big[\frac{\gamma}{\|A\|_\infty} N^{-1/2}X^{A}_{0,J}\Big]\Big) \leq C_0(g)\gamma^2}
\quad &\text{if }d=1,\\ \\
\displaystyle{\log\nu^N_g\Big(\exp\Big[\frac{\gamma}{\|A\|_\infty}N^{-(d-1)}X^{A}_{d,J}\Big]\Big) \leq \frac{C_d (g)\gamma }{N^{\frac{d-2}{2}}} }\quad &\text{if }d\geq 2.
\end{cases}
\label{eq_moments_exponentiels_appendice}
\end{align}
If instead $A$ is assumed to satisfy~\eqref{eq_assumption_diagonal_tensor}, 
then~\eqref{eq_moments_exponentiels_appendice} is valid with $\|A\|_\infty$ replaced by $\|A\|_{2,N} := \|N^{-d/2}A\|_{HS}\leq \|A\|_\infty$ everywhere.
\end{corollary}
\begin{proof}
Fix $J\subset \Z$ with $0\in J$, and $N$ such that $J\subset\Lambda_N$. 
We will use Theorem~\ref{theo_concentration_Ising} to obtain~\eqref{eq_moments_exponentiels_appendice}. 
Let us first explain why it is not necessary to assume the condition~\eqref{eq_assumption_diagonal_tensor} on the tensor $A$ if it is bounded. 
Assume that Corollary~\ref{coro_estimate_W_section_Ising} is proven for tensors satisfying the condition~\eqref{eq_assumption_diagonal_tensor}, 
and assume that $A$ is bounded but does not satisfy the condition. 
Then we claim that $N^{-(d-1)}X^A_{d,J}$ can be written as $N^{-(d-1)}X^{\tilde A}_{d,J}$ with $\tilde A$ satisfying condition~\eqref{eq_assumption_diagonal_tensor}, 
plus a sum of terms of the form $N^{-(d-1)}X^{B}_{d',J'}$ for bounded $d'$-tensors $B$, 
$1\leq d'\leq d-1$, and $J'\subset J$. 
A recursion then gives Corollary~\ref{coro_estimate_W_section_Ising} also in the case of bounded tensors that do not satisfy~\eqref{eq_assumption_diagonal_tensor} provided the corollary holds under~\eqref{eq_assumption_diagonal_tensor}. 
The existence of the decomposition of $X^A_{d,J}$ is not difficult to see: 
it is enough to sum separately on the indices $i_0,...,i_{d-1}$ in $X^A_{d,J}$ for which the same site appears twice in the collection $I_J = (i_0+J,i_1,...,i_{d-1})$, 
and use the formula $\bar\eta_i^2 = \bar\sigma_i + \bar\eta_i\sigma'(\bar\rho_i)$ ($i\in\Lambda_N$).\\

We thus prove Corollary~\ref{coro_estimate_W_section_Ising} for tensors $A$ that satisfy~\eqref{eq_assumption_diagonal_tensor}.\\ 
We first obtain estimates of tail probabilities of $X_d := X^A_{d,J}$ under $\nu^N_g$. 
These are then used to obtain~\eqref{eq_moments_exponentiels_appendice}. 
By Theorem~\ref{theo_concentration_Ising}, 
there is $\zeta_d>0$ such that, 
for any $\zeta\in(0,\zeta_d)$:
\begin{equation}
\nu^N_g\bigg(
\exp\bigg[\frac{\zeta|X_d|^{2/d}}{\|A\|_{HS}^{2/d}}\bigg]
\bigg)
\leq 
\frac{e^{\|g\|_\infty}}{\mathcal Z_g^N}\bar\nu^N\bigg(
\exp\bigg[\frac{\zeta|X_d|^{2/d}}{\|A\|_{HS}^{2/d}}\bigg]\bigg)
\leq 
2e^{\|g\|_\infty}\big(\mathcal Z^N_g\big)^{-1}
,
\label{eq_bound_exp_mom_0}
\end{equation}
Define:
\begin{equation}
\|A\|_{2,N} 
:=
\Big(\frac{1}{N^d}\sum_{i_0,...,i_{d-1}}A(i_0,...,i_{d-1})^2\Big)^{1/2}
=
\|N^{-d/2}A\|_{HS}
.
\end{equation}
By assumption on $A$,
\begin{equation}
\sup_{N\in\N^*}\|A\|_{2,N}
<
\infty
.
\end{equation}
For $\zeta$ as in~\eqref{eq_bound_exp_mom_0}, 
one has then, for each $t\geq 0$:
\begin{align}
\nu^N_g\Big(\big|X_d\big|>t\Big) 
&\leq 
2\big(\mathcal Z^N_g\big)^{-1} e^{\|g\|_\infty} \exp\Big[-\frac{\zeta t^{2/d}}{\|A\|_{HS}^{2/d}}\Big]
\nonumber\\
&\leq 
2\big(\mathcal Z^N_g\big)^{-1} e^{\|g\|_\infty}\exp\Big[-\frac{\zeta t^{2/d}}{N\|A\|_{2,N}^{2/d}}\Big]
.
\label{eq_def_zeta_size_moments_exp}
\end{align}
At this point the proof is the same for each $d\geq 1$. 
We focus on the $d\geq 2$ case. 
Let $\gamma>0$ and write:
\begin{equation}
\nu^N_g \Big[\exp\Big(\frac{\gamma}{\|A\|_{2,N}} N^{-(d-1)}|X_d|\Big)\Big]
=  
1+\int_{0}^\infty e^t \nu^N_g\Big(\frac{\gamma}{\|A\|_{2,N}} N^{-(d-1)}|X_d|>t\Big)\, dt
.
\label{eq_esperance_f(X)_bounded_by}
\end{equation}
Note that $N^{-(d-1)}X_d$ is bounded by $C\|A\|_{2,N} N$ for some numerical constant $C>0$ and each $N\in\N^*$. 
Thus, with $\zeta$ given by~\eqref{eq_def_zeta_size_moments_exp}:
\begin{align}
\nu^N_g\Big(\exp \Big[\frac{\gamma}{\|A\|_{2,N}} N^{-(d-1)}|X_d|\Big]\Big) 
&\leq 
1 + 2\big(\mathcal Z^N_g\big)^{-1} 
e^{\|g\|_\infty}\int_0^{\gamma CN} 
\exp\Big[
t-\frac{\zeta t^{2/d}N^{(d-2)/d}}{\gamma^{2/d} }\Big]
\, dt \nonumber\\
&=: 
1+2\big(\mathcal Z^N_g\big)^{-1} e^{\|g\|_\infty}
\int_0^{\gamma C N} q_\gamma(t)\, dt.
\end{align}
If $\gamma$ is small enough, 
we claim that the negative part of the exponential is dominant. 
Indeed, one has:
\begin{align}
\forall t\leq \gamma C N,\qquad 
&q_\gamma(t)
\leq 
\exp\Big[
-\frac{\zeta t^{2/d}N^{(d-2)/d}}{2\gamma^{2/d}}
\Big]\nonumber\\
&\Leftrightarrow\qquad  
\gamma 
\leq 
\gamma_d,
\quad \gamma_d = \gamma_d(\zeta) := \frac{C^{(2-d)/d}\zeta}{2}.
\end{align}
For any $\gamma<\gamma_d$, 
one has then:
\begin{align}
&\nu^N_g\Big(\exp \Big[\frac{\gamma}{\|A\|_{2,N}} N^{-(d-1)}|X_d|\Big]\Big)
\nonumber\\
&\hspace{2cm}\leq 
1 + 2\big(\mathcal Z^N_g\big)^{-1} e^{\|g\|_\infty}
\int_0^{\infty} 
\exp\Big[
-\frac{\zeta t^{2/d}N^{(d-2)/d}}{2\gamma^{2/d}}\Big]
dt.
\end{align}
The change of variable $u = t \zeta^{d/2}\gamma^{-1} N^{(d-2)/2}$ and the boundedness of $(\mathcal Z^N_g)_N$ 
(see Lemma~\ref{lemm_bound_Z_g}) conclude the proof of Corollary~\ref{coro_estimate_W_section_Ising} for $d\geq 2$: 
for each $\gamma<\gamma_d$,
\begin{align}
&\nu^N_g\Big(\exp \Big[\frac{\gamma}{\|A\|_{2,N}} N^{-(d-1)}|X_d|\Big]\Big) 
\nonumber\\
&\hspace{2cm}\leq
1 + \frac{2\gamma \zeta^{-d/2}}{N^{(d-2)/2}}\big(\mathcal Z^N_g\big)^{-1} e^{\|g\|_\infty}
\int_0^{\infty} \exp\Big[-\frac{u^{2/d} }{2}\Big]\, du
.
\end{align}
\end{proof}
\subsection{Log-Sobolev inequality under $\nu^N_g$}\label{app_LSI}
The following proposition extends to the measures $\nu^N_g$ a similar log-Sobolev inequality derived in~\cite{Goncalves2022} for product measures. 
\begin{proposition}\label{prop_LSI}
Let $0\leq \rho_-\leq \rho_+<1$ and $h,g:\squaredash\rightarrow\R$ be bounded with $\|g\|_\infty<1/4$.  
There is $C_{LS}(\rho_\pm,\|g\|_\infty,\|h\|_\infty)>0$ such that, 
for each density $f$ for $\nu^N_g$,
\begin{equation}
H(f\nu^N_g|\nu^N_g) 
\leq
C_{LS}N^2 \nu^N_g\big(\Gamma_h(\sqrt{f})\big).
\end{equation}
If in particular $h\in\s(\epsilon)$ and $g\in g_0+\s(\epsilon')$ with $\epsilon,\epsilon'\in(0,1/4)$, 
then $C_{LS}$ can be taken to depend only on $\rho_\pm$.
\end{proposition}
\begin{proof}
We shall see that $C_{LS}(\rho_\pm,\|g\|_\infty,\|h\|_{\infty})$ can be taken to be an increasing function of $\|g\|_\infty,\|h\|_{\infty}$, 
which proves the second claim.

As $\| h\|_\infty<\infty$, 
the jump rates $c_h$ satisfy:
\begin{align}
c(\eta,i) 
&\leq 
e^{\|h\|_\infty} c(\eta,i),\qquad i\in\{\pm(N-1)\},\nonumber\\
c(\eta,j,j+1)
&\leq 
e^{\|h\|_\infty}c_h(\eta,j,j+1),\qquad j<N-1.
\end{align}
It is therefore enough to prove the proposition for $h=0$. 
By assumption $\|g\|_\infty<1/4$, 
thus $\|g\|_2\leq 2\|g\|_\infty$, 
thus the difference $\lambda^N_g$ between largest and smallest eigenvalues of the matrix $(N^{-1}g_{i,j})_{i\neq j}$ is strictly below $1$ for large enough $N$. 
By Theorem 1 in~\cite{Bauerschmidt2019}, 
the following log-Sobolev inequality for the Glauber dynamics associated with $\nu^N_g$ holds: 
for each density $f$ with respect to $\nu^N_g$,
\begin{equation}\label{eq_LSI_Glauber}
H(f\nu^N_g|\nu^N_g) 
\leq
\frac{1}{2}\Big(1 + \frac{2\lambda^N_g}{1- \lambda^N_g}\Big) \sum_{i\in\Lambda_N}\nu^N_g\Big(\big[ \nabla_i\sqrt{f(\eta)}\, \big]^2\Big)
.
\end{equation}
The claim of Proposition~\ref{prop_LSI} is then an immediate consequence of~\eqref{eq_LSI_Glauber} and the following bound: for some $c(\rho_\pm,\|g\|_\infty)>0$,
\begin{equation}\label{eq_bound_Glauber_by_Kawasaki}
\forall i\in\Lambda_N,\qquad
\nu^N_g\Big(\big[ \nabla_i\sqrt{f(\eta)}\big]^2\Big)
\leq
c(\rho_\pm,\|g\|_\infty)N \nu^N_g\big(\Gamma(\sqrt{f})\big)
.
\end{equation}
Let us now prove~\eqref{eq_bound_Glauber_by_Kawasaki}. 
It is the claim of Lemma 4.2 in~\cite{Goncalves2021} when $g\equiv 0$. 
When $g\not\equiv0$, 
the proof is very similar and we only explain what changes. 
To prove~\eqref{eq_bound_Glauber_by_Kawasaki}, 
the idea is that changing the occupation number at a site $i$ requires one to bring the particle or hole at $i$ all the way to a reservoir (say the one on the right), 
perform an exchange, 
then bring back to $i$ the new hole/particle. 
This is expressed rigorously by the following formula, 
which gives a recursive description of the above procedure: 
for any $F:\Omega_N\rightarrow\R$,
\begin{equation}
\nabla_{i}F(\eta) 
= 
\nabla_{i,i+1}F(\eta) 
+ \nabla_{i+1}F(\eta^{i,i+1})
+ \nabla_{i,i+1}F\big((\eta^{i,i+1})^{i+1}\big)
.
\end{equation}
Using the identity $(x+y+z)^2 \leq 2(1+\beta)(x^2+y^2) + (1+\beta^{-1})z^2$, 
valid for any $\beta>0$ and $x,y,z\in\R$, 
one finds, after changes of variables:
\begin{align}
\nu^N_g\Big(\big[\nabla_i F(\eta)\big]^2\Big)
&\leq 
2(1+\beta) \sum_{\eta\in\Omega_N}\nu^N_g(\eta) \big[\nabla_{i,i+1}F(\eta)\big]^2\bigg[ 1 + \frac{\nu^N_g\big((\eta^i)^{i,i+1}\big)}{\nu^N_g(\eta)}\bigg]
\nonumber\\
&\quad+ 
(1+\beta^{-1})\sum_{\eta\in\Omega_N}\nu^N_g(\eta) \big(\nabla_{i+1}F(\eta)\big)^2\frac{\nu^N_g\big(\eta^{i,i+1}\big)}{\nu^N_g(\eta)}
.
\end{align}
The computation of the above ratios is the only place where $g$ plays a role. 
Since $\|g\|_\infty<\infty$, 
one can check that, for each $\eta\in\Omega_N$ (recall that $\bar\nu^N$ is the product Bernoulli measure defined in~\eqref{eq_def_bar_nu_g_intro_bis}):
\begin{equation}
\frac{\nu^N_g\big((\eta^{i})^{i,i+1}\big)}{\nu^N_g(\eta)}
\leq
\frac{\bar\nu^N\big((\eta^{i})^{i,i+1}\big)}{\bar\nu^N(\eta)}e^{2\|g\|_\infty}
,\qquad
\frac{\nu^N_g\big(\eta^{i,i+1}\big)}{\nu^N_g(\eta)}
\leq 
\frac{\bar\nu^N\big(\eta^{i,i+1}\big)}{\bar\nu^N(\eta)}e^{\|g\|_\infty}
.
\end{equation}
The rest of the proof is then identical to that of Lemma~4.2 in~\cite{Goncalves2021}. 
In particular,
\begin{equation}
\frac{\bar\nu^N\big((\eta^{i})^{i,i+1}\big)}{\bar\nu^N(\eta)}
\leq
\frac{1}{\min\{1-\rho_+,\rho_-\}}-1
,
\qquad
\frac{\bar\nu^N\big(\eta^{i,i+1}\big)}{\bar\nu^N(\eta)}
\leq 
1+ \frac{\bar\rho'}{\min\{1-\rho_+,\rho_-\}^2 N}
.
\end{equation}
This yields the following bound:
\begin{align}
\nu^N_g\Big(\big[\nabla_i F(\eta)\big]^2\Big)
&\leq 
2(1+\beta)\Big[1+\Big(\frac{1}{\min\{1-\rho_+,\rho_-\}}-1\Big)e^{2\|g\|_\infty}\Big] \nu^N_g\Big(\big[\nabla_{i,i+1} F(\eta)\big]^2\Big)
\nonumber\\
&\quad+
\big(1+\beta^{-1}\big)\Big(1+ \frac{\bar\rho'}{\min\{1-\rho_+,\rho_-\}^2 N}e^{\|g\|_\infty}\Big) \nu^N_g\Big(\big[\nabla_{i+1} F(\eta)\big]^2\Big)
.
\end{align}
Iterating the bound with the choice $\beta = N$ concludes the proof of~\eqref{eq_bound_Glauber_by_Kawasaki}, 
thus of Proposition~\ref{prop_LSI}.
\end{proof}
\section{Integration by parts formulae}\label{sec_IPP}
Fix $h\in\s(\infty)$ and $g\in g_0+\s(\infty)$ 
(the set $\mathcal S(\infty)$ is defined in~\eqref{eq_def_s_infty}) throughout. 
In this section, we provide integration by parts formulas under the measure $\nu^N_g$, both in the bulk and close to the reservoirs. 
These formulas are key to the renormalisation procedure used to estimate error terms, in Section~\ref{sec_estimate_bulk_error_terms}. 
In particular, they are essential to proving the $\Gamma$-controllability of the variables $U_0^\pm,U_1^\pm, X^{\phi_2}_{2,\{0,1\}}$ encountered in Lemma~\ref{lemm_size_controllable}.
\subsection{Integration by parts in the bulk}
Before stating the result, let us give some notations and explain what we mean by an integration by parts formula. 
Fix a density $f:\Omega_N\rightarrow\R$ for $\nu^N_g$. 
For $i<N-1$, let $\Gamma_h^{i,i+1}$ be defined as:
\begin{equation}
\forall\eta\in\Omega_{N},\qquad \Gamma_h^{i,i+1}(\sqrt{f})(\eta) = \frac{1}{2}c_h(\eta,i,i+1)\big[\nabla_{i,i+1}\sqrt{f}(\eta)\big]^2,
\end{equation}
with, for any $u : \Omega_{N}\rightarrow\R$ and any $i<N-1$:
\begin{equation}
\forall\eta\in\Omega_{N},\qquad \nabla_{i,i+1} u(\eta) = u(\eta^{i,i+1}) - u(\eta).
\end{equation}
The jump rates $c_h$ are defined in~\eqref{eq_def_jump_rates_H}. \\
Consider a family $(\omega_i)_{i\in\Lambda_N}$ of functions on $\Omega_N$. 
To estimate certain error terms in the adjoint $L^*_h{\bf 1}$ in Section~\ref{sec_computation_L_star}, 
a renormalisation scheme is used in Section~\ref{sec_estimate_W_3} below. 
For each density $f$ for $\nu^N_g$, 
some $\epsilon>0$ and $i\leq N-\epsilon N$, 
this scheme consists in estimating the replacement:
\begin{align}
&\nu^N_g\Big(f \Big[\omega_i - \frac{1}{\epsilon N}\sum_{j= i}^{i+\epsilon N -1}\omega_j\Big]\Big)\nonumber\\
&\hspace{2cm}= 
\nu^N_g\Big(f \sum_{j=i}^{i+\epsilon N-2}\frac{\epsilon N - 1-(j-i)}{\epsilon N}\big[\omega_{j}-\omega_{j+1} \big]\Big) 
\end{align}
in terms of the entropy $H(f\nu^N_g|\nu^N_g)$ and the averaged carré du champ $\nu^N_g\big(\Gamma_h(f^{1/2})\big)$. 
The right-hand side in the last equation is obtained via a simple resummation. 
The key issue, then, 
is to understand how the space gradient $\omega_j-\omega_{j+1}$ ($j<N-1$) can be turned into difference $f(\eta')-f(\eta)$, 
for a transition $\eta\rightarrow\eta'$ allowed by the dynamics (i.e. one of the differences contained in $\Gamma_h(f^{1/2})$). 
Solving this issue is what we mean by finding an integration by parts formula, 
typically a formula of the form:
\begin{equation}
\sum_{\eta\in\Omega_N}(\omega_{i+1}-\omega_i) f(\eta)\nu^N_g(\eta) = \sum_{\eta\in\Omega_N}q(\eta) [f(\eta^{i,i+1}) - f(\eta)]\nu^N_g(\eta) + \nu^N_g (fX),\label{eq_IPP_formelle}
\end{equation}
where $q,X$ are explicit functions, 
and the average of $X$ can be estimated via the entropy inequality or another integration by parts formula.
The natural choice for $\omega_\cdot$ in our case is $\omega_\cdot = \bar\eta_\cdot$, however a simpler formula~\eqref{eq_IPP_formelle} is obtained, following~\cite{Jara2018}, through the choice:
\begin{equation}
\forall i\in\Lambda_N,\qquad \omega_i = \frac{\bar\eta_i}{\bar\sigma_i},\quad \bar\eta_i = \eta_i-\bar\rho_i,\quad \bar\sigma_i = \bar\rho_i(1-\bar\rho_i).\label{eq_def_omega_IPP}
\end{equation}
\begin{lemma}\label{lemm_IPP}
Let $f$ be a $\nu^N_g$-density. Fix $i<N-1$ and let $u:\Omega_{N} \rightarrow\R$ be such that $\nabla_{i,i+1}u=0$. Then:
\begin{align}
\nu^N_g\big[ u(\omega_{i+1}-\omega_{i})f\big] 
&= 
\nu^N_g\big(u q\nabla_{i,i+1} f \big)
\nonumber\\
&\quad- (\bar\rho_{i+1}-\bar\rho_i)\nu^N_g\Big[\omega_i\omega_{i+1}e^{-(\eta_{i+1}-\eta_i)C^g_i/N}uf \Big] 
\nonumber\\
&\quad 
+ \nu^N_g\Big[(\omega_{i+1}-\omega_i)\Big(1-e^{-(\eta_{i+1}-\eta_i)C^g_i/N}\Big)uf\Big]
,
\label{eq_IPP_omega_x}
\end{align}
where the function $q = q_{i,i+1}(\eta)$ is given by:
\begin{equation}
\forall\eta\in\Omega_N,\qquad 
q(\eta) 
= 
\frac{\eta_i(1-\eta_{i+1})}{\bar\rho_i(1-\bar\rho_{i+1})}e^{-(\eta_{i+1}-\eta_i) C^g_i/N}.\label{eq_def_q_i_i+1}
\end{equation}
Recall that, for each $i<N-1$, $N^{-1}(\eta_{i}-\eta_{i+1})C^g_i = \nabla_{i,i+1}\Pi^N(g)$ is defined in~\eqref{eq_def_B_x_D_x}, 
and satisfies $\max_i|C_i^g|\leq 2\|g\|_\infty$.
\end{lemma}
\begin{proof}
Let $i<N-1$ and $q:\Omega_N\rightarrow\R$. 
Notice that, by definition of $\nabla_{i,i+1}$:
\begin{equation}
\nu^N_g\big[uq \nabla_{i,i+1} f\big] 
=
\sum_{\Omega_N} uq\nabla_{i,i+1}f \nu^N_g 
=  
\sum_{\Omega_N} f \frac{\nabla_{i,i+1}(uq\nu^N_g)}{\nu^N_g}\nu^N_g 
= 
\nu^N_g\Big[u f \frac{\nabla_{i,i+1}(q\nu^N_g)}{\nu^N_g}\Big],\label{eq_IPP_0}
\end{equation}
where we used $\nabla_{i,i+1}u=0$ to obtain the second equality. 
The gradient in the right-hand side reads:
\begin{equation}
\forall \eta\in\Omega_N,\qquad \frac{\nabla_{i,i+1}(q\nu^N_g)(\eta)}{\nu^N_g(\eta)} 
= 
q(\eta^{i,i+1})\exp\Big[-\frac{(\eta_{i+1}-\eta_i)}{N}\big(2C^g_i + \partial^N\lambda_i\big)\Big]-q(\eta),
\end{equation}
where $\lambda_i = \log(\bar\rho_i/(1-\bar\rho_i))$. We need to choose a suitable $q$ in order to have a difference $\omega_{i+1}-\omega_i$ arise above. 
In the $g=0$ case, corresponding to~\cite{Jara2018}, one can take:
\begin{equation}
\tilde q(\eta) 
:= 
\frac{\eta_i(1-\eta_{i+1})}{\bar\rho_i(1-\bar\rho_{i+1})}.
\end{equation}
When $g\neq 0$, the exponential of $C^g_i$ does not change things much, 
and if $q$ is taken as in~\eqref{eq_def_q_i_i+1}, then:
\begin{align}
\frac{\nabla_{i,i+1}(q\nu^N_g)}{\nu^N_g}(\eta) 
&= 
\bigg[\frac{\eta_{i+1}(1-\eta_i)}{\bar\rho_{i+1}(1-\bar\rho_i)}- \frac{\eta_i(1-\eta_{i+1})}{\bar\rho_i(1-\bar\rho_{i+1})}\bigg]e^{-(\eta_{i+1}-\eta_i)C^g_i/N}.\label{eq_first_formule_q_i_iplusone}
\end{align}
The variables $\omega_\cdot$ 
(see~\eqref{eq_def_omega_IPP}) are tailored to give the above bracket a nice expression 
(see (A.3) in~\cite{Jara2018}): 
\begin{equation}
\bigg[\frac{\eta_{i+1}(1-\eta_i)}{\bar\rho_{i+1}(1-\bar\rho_i)}- \frac{\eta_i(1-\eta_{i+1})}{\bar\rho_i(1-\bar\rho_{i+1})}\bigg] 
= 
\omega_{i+1}-\omega_i +(\bar\rho_{i+1}-\bar\rho_i)\omega_i\omega_{i+1}.
\end{equation}
This formula can be checked by looking for the left-hand side as a polynomial in 
$\omega_i,\omega_{i+1}$, of the form $a+b\omega_i + c\omega_{i+1} + d\omega_i\omega_{i+1}$ for real numbers $a,b,c,d$. 
Equation~\eqref{eq_first_formule_q_i_iplusone} then becomes:
\begin{align}
\frac{\nabla_{i,i+1}(q\nu^N_g)}{\nu^N_g}(\eta) &= \omega_{i+1}-\omega_i + (\omega_{i+1}-\omega_i)\Big(e^{-(\eta_{i+1}-\eta_i)C^g_i/N}-1\Big) \nonumber\\
&\quad + (\bar\rho_{i+1}-\bar\rho_i)\omega_i\omega_{i+1}e^{-(\eta_{i+1}-\eta_i)C^g_i/N},
\end{align}
which proves the lemma when plugged into~\eqref{eq_IPP_0}.
\end{proof}
The next lemma is a rewriting of Lemma~\ref{lemm_IPP} in terms of the carré du champ operator.
\begin{lemma}\label{lemm_IPP_forme_dir}
Let $i<N-1$ and let $u:\Omega_{N}\rightarrow\R$ be such that $\nabla_{i,i+1} u = 0$. 
There is then a constant $C = C(h,g,\bar\rho)>0$ such that, for any $\delta>0$:
\begin{align}
\nu^N_g\big[ u (\omega_{i+1}-\omega_{i})f\big]&\leq \delta N^2\nu^N_g\big[\Gamma_h^{i,i+1}(f^{1/2})\big] + \frac{C}{\delta N^2} \nu^N_g\big[f u^2\big]\label{eq_IPP_forme_dir_bulk}\\
&\quad-(\bar\rho_{i+1}-\bar\rho_i)\nu^N_g\Big[\omega_i\omega_{i+1}e^{-(\eta_{i+1}-\eta_i)C^g_i/N} uf\Big] \nonumber\\
&\quad+ \nu^N_g\Big[(\omega_{i+1}-\omega_i)\Big(1-e^{-(\eta_{i+1}-\eta_i)C^g_i/N}\Big)uf\Big].\nonumber
\end{align}
\end{lemma}
\begin{proof}
Let $i<N-1$ and $\beta>0$. In~\eqref{eq_IPP_omega_x}, write, for each $\eta\in\Omega_{N}$:
\begin{align}
\nabla_{i,i+1}f(\eta)u(\eta)q(\eta) 
&=
\beta^{1/2} [f^{1/2}(\eta^{i,i+1})-f^{1/2}(\eta)]
\nonumber\\
&\qquad\cdot\beta^{-1/2}u(\eta)q(\eta)[f^{1/2}(\eta^{i,i+1})+f^{1/2}(\eta)]
.
\end{align}
Apply then Cauchy-Schwarz inequality to obtain:
\begin{align}
\nu^N_g\big[ u (\omega_{i+1}-\omega_{i})f\big]
&\leq 
\frac{\beta}{2} \nu^N_g\big[\big(\nabla_{i,i+1}f^{1/2}\big)^2\big] + \frac{1}{\beta} \nu^N_g\big[\big(f(\eta)+f(\eta^{i,i+1})\big)u^2 q^2\big]
\label{eq_IPP_forme_dir_0}
\\
&\qquad
-(\bar\rho_{i+1}-\bar\rho_i)\nu^N_g\Big[\omega_i\omega_{i+1}e^{-(\eta_{i+1}-\eta_i)C^g_i/N} uf\Big] 
\nonumber\\
&\qquad 
+ \nu^N_g\Big[(\omega_{i+1}-\omega_i)\Big(1-e^{-(\eta_{i+1}-\eta_i)C^g_i/N}\Big)uf\Big].\nonumber
\end{align}
Changing variables $\eta\leftarrow \eta^{i,i+1}$, since $\bar\rho_i\in[\rho_-,\rho_+]\subset (0,1)$, the second expectation in~\eqref{eq_IPP_forme_dir_0} reads:
\begin{align}
\frac{1}{\beta}\nu^N_g\big[\big(f(\eta)+f(\eta^{i,i+1})\big)u^2 q^2\big] &= \frac{1}{\beta}\sum_{\eta\in\Omega_N}f(\eta)u^2(\eta)\Big[q^2(\eta)+q^2(\eta^{i,i+1})\frac{\nu^N_g(\eta^{i,i+1})}{\nu^N_g(\eta)}\Big]\nu^N_g(\eta)\nonumber\\
&\leq \frac{e^{4\|g\|_\infty}M}{\beta}\nu^N_g\big[fu^2\big],
\end{align}
where $M = \max \{(\bar\rho_j(1-\bar\rho_{j+1}))^{-2},(\bar\sigma_j\bar\sigma_{j+1})^{-1} : j<N-1\}$ depends only on $\rho_\pm$.\\
Consider now the first term in the right-hand side of~\eqref{eq_IPP_forme_dir_0}. Since $c_h(\eta,i,i+1)\geq c(\eta,i,i+1)e^{-2\|h\|_\infty}$ for each $\eta\in\Omega_N$, it reads:
\begin{equation}
\frac{\beta}{2} \int \big[\nabla_{i,i+1}f^{1/2}\big]^2 d\nu^N_g \leq \beta e^{2\|h\|_\infty}\nu^N_g\big(\Gamma_h^{i,i+1}(\sqrt{f})\big).
\end{equation}
Taking $\beta = \delta N^2 e^{2\|h\|_\infty}$ concludes the proof.
\end{proof}
\subsection{Integration by parts at the boundary and boundary correlations}\label{sec_boundary_correl_dirichlet}
Here, we estimate dynamical correlations involving sites close the reservoirs, i.e. correlations of the form $\bar\eta_{\pm(N-1)}X_N$ for a function $X_N:\Omega_N\rightarrow\R$. Recall that $h\in\s(\infty)$, 
$g\in g_0+\s(\infty)$ are fixed, and the definition~\eqref{eq_def_L_pm} of the jump rates at the boundary. Define, for $f:\Omega_N\rightarrow\R$:
\begin{equation}
\forall\eta\in\Omega_N,\qquad \Gamma_h^{\pm}(f) = \frac{c_h(\eta,\pm (N-1))}{2}\big[f(\eta^{\pm (N-1)})-f(\eta)\big]^2,\label{eq_def_Gamma_bord}
\end{equation}
and observe:
\begin{equation}
c_h(\eta,\pm(N-1))\geq e^{-2\|h\|_\infty}\min\{\rho_{\epsilon_1},(1-\rho_{\epsilon_2}):\epsilon_1,\epsilon_2\in\{-,+\}\}.\label{eq_lower_bound_jump_rate_boundary}
\end{equation}
\begin{lemma}\label{lemm_small_boundary_correl_sous_bar_nu_G}
For $n,N\in\N^*$, let $\phi_n:\Lambda_N^n\rightarrow\R$ satisfy $\sup_{N\in\N^*} \|\phi_n\|_{\infty} <\infty$. 
For $\eta\in\Omega_N$ and $\epsilon\in\{+,-\}$, 
define $U^\epsilon_0(\eta) = \bar\eta_{\epsilon(N-1)}$ ($\eta\in\Omega_N$) 
and, for $n\geq 1$:
\begin{equation}
U_n^\epsilon(\eta) 
= 
\bar\eta_{\epsilon(N-1)}V_n^\epsilon(\eta)
,\qquad 
V_n^\epsilon(\eta) 
= 
\frac{1}{N^{n-1}}\sum_{i_1,...,i_n\neq \epsilon(N-1)}\bar\eta_{i_1}...\bar\eta_{i_n}\phi_n(i_1,...,i_n)
.
\end{equation}
Then, for each $n\in\N$, $U_n^\epsilon$ is $\Gamma$-controllable with size $N^{-1}$ and of vanishing type in the sense of Lemmas~\ref{lemm_size_controllable}--\ref{lemm_type}. 
More precisely, there is $C>0$ such that, for any $\nu^N_g$-density $f$ and any $\delta>0$:
\begin{align}
\nu^N_g\big(fN^{1/2} U_n^\epsilon\big) 
&\leq 
\delta N^2\nu^N_g\big(\Gamma_h^{\epsilon}(\sqrt{f})\big) 
\nonumber\\
&\quad 
+ \nu^N_g\bigg(f \frac{C}{N}\Big[\frac{(V^\epsilon_n)^2}{\delta} + \Big|\frac{V^\epsilon_n}{N^{1/2}}\sum_{j\neq \epsilon(N-1)}\bar\eta_j(Ng_{\epsilon(N-1),j})\Big|\Big]\bigg) 
+ \alpha^N_{n}
,
\label{eq_explicit_bound_boundaryy}
\end{align}
setting $V_0^\epsilon :=1$ by convention and with $\alpha^N_n = O(N^{-1/2})$.
\end{lemma} 
\begin{remark}
The estimate on the size of $U^\epsilon_n$ is optimal only if $n\leq 1$. 
$U^\epsilon_{n}$ with $n\leq 1$ are used in the computation of the adjoint in Section~\ref{sec_computation_L_star}, 
while $U^{\epsilon}_2$ is useful in Section~\ref{sec_estimate_bulk_error_terms}.
\demo
\end{remark}
\begin{proof}
Let $n\in\N^*$. 
Using the notations of Corollary~\ref{coro_estimate_W_section_Ising}, 
the term in the expectation in the right-hand side of~\eqref{eq_explicit_bound_boundaryy} is of the form $N^{-(2n-1)}X^{\psi_{2n}}_{2n,\{0\}}+N^{-n-\frac{1}{2}}X_{n+1,\{0\}}^{\psi_{n+1}}$ with $\sup_{N\in\N^*}\|\psi_\ell\|_{\infty}<\infty$, 
$\ell\in\{n,2n\}$. 
It is thus controllable with size (at most) $\delta^{-1}$. 
The estimate~\eqref{eq_explicit_bound_boundaryy} then implies that 
$U_n^\epsilon$ is $\Gamma$-controllable with size $N^{-1}$ by taking $\delta=N^{1/2}$, 
but also that $N^{1/2}U_n^\epsilon$ is $\Gamma$-controllable with size $1$ and thus that $U^n_\epsilon$ is of vanishing type. \\

Let us therefore prove~\eqref{eq_explicit_bound_boundaryy}. 
We do so for the left boundary $\epsilon=-$, the proof for the right boundary being identical. 
The goal is to create a gradient of $f$ of the form~\eqref{eq_def_Gamma_bord}. We use the shorthand notation $b= -(N-1)$. Notice that $V_n^-(\eta^b) = V_n^-(\eta)$ for any $\eta$. 
 The mapping $\eta\mapsto \eta^{b}$ is bijective on $\Omega_N$, thus:
\begin{align}
\nu^N_g(f\bar\eta_{b}V^-_n) = \frac{1}{2}\sum_{\eta\in\Omega_N} \nu^N_g(\eta)V^-_n(\eta)\Big[f(\eta)(\eta_{b}-\bar\rho_{b}) + f(\eta^{b})(1-\eta_{b}-\bar\rho_{b})\frac{\nu^N_g(\eta^{b})}{\nu^N_g(\eta)}\Big].\label{eq_moyenne_nu_f_boundary}
\end{align}
For $\eta\in\Omega_N$, let us compute the ratio $\nu^N_g(\eta^{b})/\nu^N_g(\eta)$:
\begin{align}
\frac{\nu^N_g(\eta^{b})}{\nu^N_g(\eta)}& = \frac{(1-\eta_{b})\bar\rho_{b}+ \eta_{b}(1-\bar\rho_{b})}{\eta_{b}\bar\rho_{b} + (1-\eta_{b})(1-\bar\rho_{b})} e^{2\Pi^N(g)(\eta^{b}) - 2\Pi^N(g)(\eta)}\nonumber\\
&= \Big(\frac{\bar\rho_{b}}{1-\bar\rho_{b}}\Big)^{1-2\eta_{b}}\exp\Big[\frac{(1-2\eta_{b})}{N}\sum_{j\neq b}\bar\eta_j g_{b,j}\Big].\label{eq_gradient_mesure_au_bord}
\end{align}
For future reference, notice that~\eqref{eq_gradient_mesure_au_bord} is bounded by $C(\rho_\pm)e^{2\|g\|_\infty}$ for some $C(\rho_\pm)>0$. Forgetting $g$ for a second in~\eqref{eq_moyenne_nu_f_boundary}, notice also that, for each $\eta\in\Omega_N$:
\begin{align}
f(\eta)&(\eta_b-\bar\rho_b)+f(\eta^b)(1-\eta_b-\bar\rho_b)\Big(\frac{\bar\rho_{b}}{1-\bar\rho_{b}}\Big)^{1-2\eta_{b}}\nonumber \\
&\qquad = \eta_b(1-\bar\rho_b)\big[f(\eta)-f(\eta^b)\big] - (1-\eta_b)\bar\rho_b\big[f(\eta) - f(\eta^b)\big] \nonumber\\
&\qquad = \bar\eta_b \big(f(\eta)-f(\eta^b)\big),\label{eq_ordre_0_gradient_mesure_bord}
\end{align}
which involves a gradient of $f$ as desired. Coming back to~\eqref{eq_gradient_mesure_au_bord}, note that, since $g(\pm 1,\cdot) = 0$, the argument of the exponential in~\eqref{eq_gradient_mesure_au_bord} is bounded by $O(N^{-1})$. Equation~\eqref{eq_ordre_0_gradient_mesure_bord} and the existence of $C(g)>0$ such that $|e^{x}-1-x|\leq C(g)x^2$ holds for $x\leq 2\|g\|_\infty$ therefore yield the bound:
\begin{align}
&\Bigg|
N^{1/2}\nu^N_g(f\bar\eta_{b}V^-_n) - \frac{N^{1/2}}{2}\sum_{\eta\in\Omega_N} \nu^N_g(\eta)V^-_n(\eta)\bar\eta_b \big[f(\eta)-f(\eta^b)\big] 
\nonumber\\
&\hspace{1cm}- 
\frac{N^{1/2}}{2N^2}\sum_{\eta\in\Omega_N}\nu^N_g(\eta)V^-_n(\eta)f(\eta^b)
(1-\eta_b-\bar\rho_b)\Big(\frac{\bar\rho_{b}}{1-\bar\rho_{b}}\Big)^{1-2\eta_{b}}(1-2\eta_b)\sum_{j\neq b}\bar\eta_j (Ng_{b,j})
\Bigg| 
\nonumber\\
&\hspace{3cm}\leq 
\frac{C(g)\|V_n^-\|_{\infty}}{N^{3/2}}=:\alpha^n_N
.
\label{eq_dvplt_exponentiel_IPP_bord}
\end{align}
Since $\|V_n^-\|_\infty = O(N)$, 
$\alpha^n_N=O(N^{-1/2})$ as claimed. 
It remains to estimate the second and third terms in the left-hand side. 
Consider first the third term. Using the bijection $\eta\mapsto \eta^b$ to turn $f(\eta^b)$ into $f(\eta)$, 
recalling that $V^-_n(\eta^b) = V^-_n(\eta)$ for each $\eta\in\Omega_N$, 
and bounding the ratio~\eqref{eq_gradient_mesure_au_bord} by $C(\rho_\pm)e^{2\|g\|_\infty}$, one finds:
\begin{align}
\Big|&\frac{N^{1/2}}{2N^2}\sum_{\eta\in\Omega_N}\nu^N_g(\eta)V^-_n(\eta)f(\eta^b)(1-\eta_b-\bar\rho_b)\Big(\frac{\bar\rho_{b}}{1-\bar\rho_{b}}\Big)^{1-2\eta_{b}}\sum_{j\neq b}\bar\eta_j (Ng_{b,j})\Big)\Big| \nonumber\\
&\hspace{4cm}\leq 
\frac{C(\rho_\pm)e^{2\|g\|_\infty}}{N^{3/2}}\nu^N_g\Big(f\Big|V^-_n\sum_{j\neq b}\bar\eta_j(Ng_{b,j})\Big|\Big),\label{eq_un_des_termes_RHS_IPP_boundary}
\end{align}
which is one of the terms appearing in~\eqref{eq_explicit_bound_boundaryy}. 
Consider now the second term in the left-hand side of~\eqref{eq_dvplt_exponentiel_IPP_bord}. 
For $\beta>0$ and $\eta\in\Omega_N$, 
split as follows:
\begin{equation}
N^{1/2}[f(\eta)-f(\eta^b)]
=
\beta N^{1/2}[f^{1/2}(\eta)-f^{1/2}(\eta^b)]\beta^{-1}[f^{1/2}(\eta)+f^{1/2}(\eta^b)]
.
\end{equation}
Apply then Cauchy-Schwarz inequality twice to find, bounding $(\bar\eta_b)^2$ by $1$:
\begin{align}
\frac{N^{1/2}}{2}&\sum_{\eta\in\Omega_N} V^-_n(\eta)\bar\eta_b \big[f(\eta)-f(\eta^b)\big] 
\\
&\qquad\leq 
\frac{\beta N}{4}\sum_{\eta\in\Omega_N} \big[f^{1/2}(\eta)-f^{1/2}(\eta^b)\big]^2 + \frac{1}{2\beta}\sum_{\eta\in\Omega_N}\nu^N_g(\eta)\big[f(\eta)+f(\eta^{b})\big](V_n^-)^2(\eta)
\nonumber.
\end{align}
As for~\eqref{eq_un_des_termes_RHS_IPP_boundary}, the last expectation is bounded by $\beta^{-1}C(g)\nu^N_g(f(V^-_n)^2)$. 
To conclude the proof, recall from~\eqref{eq_def_Gamma_bord} the expression of $\Gamma_h^-$ and from~\eqref{eq_lower_bound_jump_rate_boundary} the lower bound $C(\rho_\pm)e^{-\|h\|_\infty}$ on the jump rates. 
Choose $\beta = 2C(\rho_\pm)e^{-2\|h\|_\infty}\delta N$ for $\delta>0$. Then:
\begin{align}
\frac{1}{2}\sum_{\eta\in\Omega_N} V^-_n(\eta)\bar\eta_b \big[f(\eta)-f(\eta^b)\big] 
&\leq 
\delta N^2 \nu^N_g\big(\Gamma_h^-(\sqrt{f})\big) 
\nonumber\\
&\quad+ 
\frac{C(g)C(\rho_\pm)^{-1}e^{2\|h\|_\infty}}{2\delta N} \nu^N_g\big(f(V^-_n)^2\big)
.
\end{align}
The last equation together with~\eqref{eq_dvplt_exponentiel_IPP_bord} and~\eqref{eq_un_des_termes_RHS_IPP_boundary} is precisely the right-hand side of~\eqref{eq_explicit_bound_boundaryy}.
\end{proof}
\section{Control of the error terms}\label{sec_estimate_bulk_error_terms}
Fix $h\in\s(\infty)$ 
(this set is defined in~\eqref{eq_def_s_infty}) and an associated $g_h\in g_0+ \s(\infty)$, 
solution of the main equation~\eqref{eq_main_equation}. 
In this section we estimate, 
for each density $f$ for $\nu^N_{g_h}$, the average against $f\nu^N_{g_h}$ of the function $X^\phi_{2,\{0,1\}}$, 
$\phi:\Lambda_N^2\rightarrow\R$, 
defined below in~\eqref{eq_def_W_3_appendix}. 
This proves the last item of Lemma~\ref{lemm_size_controllable}. 
We also estimate the expectation of the time average of any error term encountered in the text, 
proving Proposition~\ref{prop_Boltzmann_gibbs_sec3} and Corollary~\ref{coro_Boltzmann_gibbs_sec3}.
\begin{proposition}\label{prop_error_terms}
Let $\phi \in \Lambda_N^{2}\rightarrow\R$ satisfy $\sup_{N\in\N^*}\|\phi_2\|_\infty<\infty$. Recall that $X^{\phi}_{2,\{0,1\}}$, abbreviated as $X_2$, was defined in~\eqref{eq_def_X_A_d_J_appendice} by:
\begin{equation}
\forall \eta\in\Omega_N,\qquad X_{2}(\eta) = \sum_{i<N-1}\sum_{j\notin\{i,i+1\}}\bar\eta_{i}\bar\eta_{i+1}\bar\eta_{j} \phi(i,j).\label{eq_def_W_3_appendix}
\end{equation}
Then $N^{-1}X_{2}$ is $\Gamma$-controllable with size $N^{-1/2}$ (recall Definition~\ref{def_controllability} of controllability).
\end{proposition}
The next proposition was stated as Proposition~\ref{prop_Boltzmann_gibbs_sec3}, 
and is proven in Section~\ref{sec_proof_boltzmann_gibbs}.
\begin{proposition}\label{prop_Boltzmann_Gibbs}
Let $0<\rho_-\leq \rho_+<1$ and $\epsilon_B(\rho_-)$ be chosen as in Lemma~\ref{lemm_L_star_as_e_plus_carre_du_champ}, 
and assume $\bar\rho'\leq \epsilon_B$, $h\in\s(\epsilon_B)$. 
Let $E^N:\Omega_N\rightarrow\R$ be an error term with size $a_N =o_N(1)$, 
and let $F^N$ be controllable with size $1$. 
There are then $\gamma,C$ and $\gamma',C'>0$ depending on $h,\rho_\pm$ such that, 
for each $T>0$ and $N\in\N^*$:
\begin{align}\label{eq_bound_time_expectation_W_3_W_4}
\frac{1}{T}\log\E^{\nu^{N}_{g_h}}\Big[\exp\Big|\gamma\int_0^T  E^N(\eta_t)\, dt\Big|\Big] 
&\leq 
C a_N,\qquad\nonumber\\ 
\frac{1}{T}\log \E^{\nu^{N}_{g_h}}\Big[\exp\Big|\gamma'\int_0^T  F^N(\eta_t)\, dt\Big|\Big] 
&\leq 
C'
.
\end{align}
\end{proposition}
Propositions~\ref{prop_error_terms}-\ref{prop_Boltzmann_Gibbs} are proven in the next two sections.
\subsection{Estimate of $X_2$}\label{sec_estimate_W_3}
In this section, we prove Proposition~\ref{prop_error_terms}. Fix $\phi : \Lambda_N^2\rightarrow\R, N\in\N^*$ as in the proposition. 
Fix also a density $f$ for $\nu^N_{g_h}$ once and for all. 
The proof of Proposition~\ref{prop_error_terms} being very technical, we first present its general structure. \\
The idea is to smoothen the product $\bar\eta_i\bar\eta_{i+1}$ into a quantity that depends on all $\bar\eta$'s in a box of size $\ell$ with $\ell$ sufficiently large, then use the entropic inequality to estimate the resulting term. The cost of this replacement will be estimated by an integration by parts formula, see Section~\ref{sec_IPP}. \\
We need room between the indices $i,i+1$ and $j$ in the definition~\eqref{eq_def_W_3_appendix} of $X_2$ to take averages in a box. Let $I_\ell$ be the segment $\{0,...,\ell-1\}$ and split the sum on $j$ in~\eqref{eq_def_W_3_appendix} as follows:
\begin{align}
\forall\eta\in\Omega_N,\qquad 
\frac{1}{N}X_2(\eta) 
= 
\overrightarrow{X}_2^{\ell} + \overleftarrow{X}_2^{\ell},\qquad \overrightarrow{X}_2^{\ell} 
&= 
\frac{1}{N}\sum_{i<N-1}\sum_{\substack{j\in\Lambda_N\setminus\{i\} \\ j\notin i+1+I_\ell}}\bar\eta_i\bar\eta_{i+1}\bar\eta_j\phi(i,j)
,
\\
\overleftarrow{X}_2^{\ell} 
&= 
\frac{1}{N}\sum_{i<N-1}\sum_{\substack{j\in\Lambda_N\setminus\{i,i+1\}\\ j\in i+1+I_\ell }}\bar\eta_i\bar\eta_{i+1}\bar\eta_j \phi(i,j) 
.
\label{eq_def_X_left_X_right_W_3}
\end{align}
The direction of the arrow indicates the direction in which the replacement of $\bar\eta_i$ 
($\leftarrow$) or $\bar\eta_{i+1}$ ($\rightarrow$) by averages on sites to the left of $i$ ($\leftarrow$) or to the right of $i+1$ ($\rightarrow$) is going to be performed. 
Estimates for $\overleftarrow{X}_2 ^{\ell}$ and $\overrightarrow{X}_2 ^{\ell}$ are identical, so we only estimate the latter. 
In practice, the replacement is made thanks to the integration by parts Lemma~\ref{lemm_IPP_forme_dir}, 
which uses $\omega_\cdot = \bar\eta_\cdot/\bar\sigma_\cdot$ as main variable. 
Write: 
\begin{equation}
A(i,j) := \bar\sigma_{i+1} \phi(i,j),\qquad i<N-1,j\in\Lambda_N.
\end{equation}
Then:
\begin{equation}
\forall\eta\in\Omega_{N},\qquad 
\overrightarrow{X}_2 ^{\ell}(\eta) 
= 
\frac{1}{N}\sum_{i<N-1}\sum_{\substack{j\in\Lambda_N\setminus\{i,i+1\} \\ j\notin i+1+I_\ell }}\bar\eta_i\omega_{i+1}\bar\eta_jA(i,j),
\end{equation}
and we replace $\omega_{i+1}$ by $\frac{1}{\ell}\sum_{a\in i+1+I_\ell}\omega_a$. If $i+1$ is too close to the reservoirs, i.e. if $i+\ell>N-1$, 
then this replacement does not make sense. 
In this case, 
we spread the unit mass at $i+1$ to $1/\ell$ at each site in $i+1+I_\ell\cap \Lambda_N$, 
and leave the remaining $\frac{N-1-(i+1)}{\ell}$ mass at the boundary. 
This is summarised in the following definition of the replacement $\overrightarrow\omega^\ell_{i+1}$ of $\omega_{i+1}$:
\begin{equation}
\forall i<N-1,\qquad 
\overrightarrow\omega^\ell_{i+1} 
= 
\frac{1}{\ell}\sum_{a=i+1}^{\min\{i+\ell,N-1\}}\omega_a + {\bf 1}_{i+\ell>N-1}\Big(1-\frac{N-1-i}{\ell}\Big)\omega_{N-1}
.
\label{eq_def_rightarrow_omega}
\end{equation}
\noindent\textbf{Choice of $\ell$.} Let $\overrightarrow{Y}_2^{\ell}$ denote the averaged version of $\overrightarrow{X}_2 ^{\ell}$:
\begin{align}
\forall\eta\in\Omega_{N},\qquad 
&\overrightarrow{Y}_2 ^{\ell}(\eta) 
= 
\overrightarrow{Z}_2 ^{\ell}(\eta) \label{eq_first_line_Y_fleche_N} \\
&\quad + \frac{1}{N}\sum_{i<N-1}\sum_{\substack{j\in\Lambda_N\setminus\{i,i+1\} \\ j\notin i+1+I_\ell }}\bar\eta_{i}\bar\eta_j{\bf 1}_{i+\ell>N-1}\Big(1-\frac{N-1-i}{\ell}\Big)\omega_{N-1}A(i,j)
,
\nonumber
\end{align}
with:
\begin{equation}
\forall \eta\in\Omega_N,\qquad \overrightarrow{Z}_2 ^{\ell}(\eta) 
= 
\frac{1}{N}\sum_{i<N-1}\sum_{\substack{j\in\Lambda_N\setminus\{i,i+1\}\\ j\notin i+1+I_\ell }}\bar\eta_{i}\bar\eta_j\Big(\frac{1}{\ell}\sum_{a=i+1}^{\min\{i+1+\ell,N-1\}}\omega_a\Big)A(i,j)
.
\end{equation}
The last term in~\eqref{eq_first_line_Y_fleche_N} is $\Gamma$-controllable with size $N^{-1}$ by Lemma~\ref{lemm_small_boundary_correl_sous_bar_nu_G}. For the replacement of $\overrightarrow{X}_2 ^{\ell}$ by $\overrightarrow{Y}_2 ^{\ell}$ to be useful, $\overrightarrow{Z}_2 ^{\ell}$ should be controllable with size $o_N(1)$. This requirement will fix the choice of $\ell$. Looking at Corollary~\ref{coro_estimate_W_section_Ising}, we see that any $\ell$ such that $\ell = o(N)$ fails, so we take:
\begin{equation}
\ell := N.\label{eq_choice_of_ell_first_step_W_3}
\end{equation}
The entropy inequality is then effective on $\overrightarrow{Z}_2^{N}$. Indeed, it is of the form:
\begin{equation}
\forall\eta\in\Omega_{N},\qquad 
\overrightarrow{Z}_2^{N}(\eta) 
= 
\frac{1}{N^2}\sum_{\substack{(i,j,a)\in\Lambda_N^3 \\ |\{(i,j,a)\}|=3}}\bar\eta_i\bar\eta_j\bar\eta_a\tilde A(i,j,a)
,
\end{equation}
for some function $\tilde A$ satisfying $|\tilde A(i,j,a)|\leq |A(i,j)|$ for each $(i,j,a)\in\Lambda_N^3$. By Corollary~\ref{coro_estimate_W_section_Ising}, $\overrightarrow{Z}_2^{N}$ is therefore controllable with size $N^{-1/2}$: there are $\gamma,C>0$ such that:
\begin{equation}
\nu^N_{g_h}\big(f\overrightarrow{Z}_2^{N}\big) \leq \frac{H(f\nu^N_g|\nu^N_g)}{\gamma} +\frac{1}{\gamma}\log \nu^N_g\Big[\exp \big(\gamma \overrightarrow{Z}_2^{N}\big)\Big] \leq \frac{H(f\nu^N_g|\nu^N_g)}{\gamma} +\frac{C}{ N^{1/2}}.
\end{equation}
\noindent \textbf{Cost of the replacement.} Let us estimate the cost of replacing $\overrightarrow{X}_2^{N}$ by $\overrightarrow{Y}_2^{N}$, defined in~\eqref{eq_first_line_Y_fleche_N}. To do so, we use of the following integration by parts identity, which explicitly describes how to spread the unit mass at $i+1$, $i<N-1$, to $1/N$ on every site up to the boundary, where the remaining mass is then left. One has:
\begin{equation}
\omega_{i+1} - \overrightarrow{\omega}_{i+1}^N = \sum_{a=i+1}^{\min\{i+N,N-1\}-1} \phi_N(a-(i+1))(\omega_{a}-\omega_{a+1}),
\end{equation}
with:
\begin{equation}
\phi_N(b) 
= 
\frac{N-1-b}{N}{\bf 1}_{0\leq b<N},
\qquad b\in\Z
.
\end{equation}
For brevity, for $a\in\Lambda_{N}$, let $u_a$ denote the quantity:
\begin{equation}
\forall \eta\in\Omega_{N},\qquad 
u_a(\eta) 
= 
\frac{1}{N}\sum_{\substack{i<N-1 \\ a-N<i<a}}\phi_{N}(a-(i+1))\bar\eta_i\sum_{\substack{j\in\Lambda_N\setminus\{i,i+1\} \\ j\notin i+1+I_N}}\bar\eta_jA(i,j)
.
\label{eq_def_u_y}
\end{equation}
Then, for each $\eta\in\Omega_{N}$:
\begin{equation}
\overrightarrow{X}^{N}_2-\overrightarrow{Y}^{N}_2 
= 
\sum_{a<N-1}(\omega_{a}-\omega_{a+1})u_a(\eta)
.
\end{equation}
To estimate the expectation of the right-hand side above under $f\nu^N_g$, apply, for each $a<N-1$, 
the integration by parts formula of Lemma~\ref{lemm_IPP_forme_dir}, 
with $u=-u_a$. 
There is thus a constant $C>0$ such that, for each $\delta>0$:
\begin{align}
&\sum_{a<N-1}\nu^N_{g_h}\big(f(\omega_a-\omega_{a+1})u_a\big) \leq 
\delta N^2 \nu^N_{g_h}\big(\Gamma_h(\sqrt{f})\big) 
\nonumber\\
&\qquad
+ \frac{C}{\delta N^2} \sum_{a<N-1}\nu^N_{g_h}\big(u_a^2f\big)
\quad 
(:= \nu^N_{g_h}(fR_1))
\label{eq_IPP_pour_W_3_u_square_term}\\
&\qquad
+\sum_{a<N-1}(\rho_{a+1}-\rho_a)\nu^N_{g_h}\big(\omega_a\omega_{a+1}e^{-(\eta_{a+1}-\eta_a)C^{g_h}_a/N }u_a f\big)
\quad 
(:= \nu^N_{g_h}(fR_2))
\label{eq_IPP_pour_W_3_omega_y_omega_y+1_term}\\
&\qquad
- \sum_{a<N-1}\nu^N_{g_h}\Big((\omega_{a+1}-\omega_{a})\Big(1-e^{-(\eta_{a+1}-\eta_a)C^{g_h}_a/N }\Big)u_a f\Big)
\quad 
(:= \nu^N_{g_h}(fR_3))
.
\label{eq_IPP_pour_W_3_C_g_term}
\end{align}
Let us estimate one by one each of~\eqref{eq_IPP_pour_W_3_u_square_term}-\eqref{eq_IPP_pour_W_3_omega_y_omega_y+1_term}-\eqref{eq_IPP_pour_W_3_C_g_term}. 
Consider first $\nu^N_{g_h}(fR_3)$, and note that:
\begin{equation}
\sup_N\sup_{a\in\Lambda_N}|(\omega_{a+1}-\omega_a)(\eta_{a+1}-\eta_a)|
\leq 
C(\rho_\pm)
.
\end{equation}
As a result, 
using the identity $e^x =1 + \int_0^1 xe^{tx}dt$ for $x\in\R$ and the fact that $|C^{g_h}_\cdot|\leq 2\|g_h\|_\infty$, $\nu^N_{g_h}(fR_3)$ can be bounded as follows:
\begin{align}
\big|\nu^N_{g_h}(fR_3)\big|
\leq 
C(\rho_\pm)e^{2\|g\|_\infty} \nu^N_{g_h}\bigg(\frac{1}{N}\sum_{a<N-1}f\big|C^{g_h}_a u_a\big|\bigg)
.
\label{eq_bound_R_3}
\end{align}
By definition of $C^{g_h}_\cdot$ (see e.g. Lemma~\ref{lemm_IPP}) and of $u_\cdot$ in~\eqref{eq_def_u_y}, the product $C^{g_h}_\cdot u_\cdot$ is of the form:
\begin{equation}
\forall a\in\Lambda_N,\qquad 
C^{g_h}_a u_a 
= 
\frac{1}{N^2}\sum_{(i,j,b)\in\Lambda_N^3}\bar\eta_i\bar\eta_j\bar\eta_b \psi^a_{i,j,b}
+
D^{g_h}_a u_a,
\quad  
\sup D^{g_h}_\cdot = O(N^{-1})
,
\end{equation}
where the functions $\psi^a :(-1,1)^3\rightarrow\R$ are bounded uniformly in $a$. 
It follows by Corollary~\ref{coro_estimate_W_section_Ising} that $R_3$ is controllable with size $N^{-1/2}$.\\
Consider now $R_2$, defined in~\eqref{eq_IPP_pour_W_3_omega_y_omega_y+1_term}. 
Again using $e^x = 1+\int_0^1xe^{tx}dt$ for $x\in\R$, we can bound $\nu^N_{g_h}(f R_2)$ as follows:
\begin{align}
\nu^N_{g_h}(f R_2) 
&\leq 
\frac{\bar\rho'}{N}\sum_{a<N-1}\nu^N_{g_h}\big(f\omega_a \omega_{a+1}u_a\big)+\frac{\bar\rho'e^{2\|g_h\|_\infty}}{N^2}\sum_{a<N-1}\frac{1}{\bar\sigma_a\bar\sigma_{a+1}}\nu^N_{g_h}\Big(f\big|C^{g_h}_au_a\big|\Big).
\label{eq_second_line_R_2}
\end{align}
Recalling the definition of $u_\cdot$ from~\eqref{eq_def_u_y}, the first term in~\eqref{eq_second_line_R_2} is of the form $N^{-2}X^{B}_{3,\{0,1\}}$ in the notations of Theorem~\ref{theo_concentration_Ising}, i.e. of the form $N^{-2}\sum_{i,j,b}\bar\eta_i\bar\eta_{i+1}\bar\eta_j\bar\eta_b B(i,j,b)$, with $B$ bounded. Corollary~\ref{coro_estimate_W_section_Ising} tells us that this function does not behave worse than a sum of three-point correlations, and is therefore controllable with size $N^{-1/2}$. In addition, the second term in~\eqref{eq_second_line_R_2} has the same structure as $N^{-1} R_3$, and is therefore controllable with size $N^{-3/2}$.\\
Consider finally $R_1$ in~\eqref{eq_IPP_pour_W_3_u_square_term}. It reads:
\begin{align}
\nu^N_{g_h}(fR_1) 
= 
\frac{C}{\delta N^4}\sum_{a<N-1} \sum_{\substack{i,j<N-1 \\ a-N< i,j<a}}\bigg[\phi_N(a-(i+1))&\phi_N(a-(j+1)) \bar\eta_i\bar\eta_j
\label{eq_bound_R_1}\\
&\times\sum_{\substack{b,c\in \Lambda_N\setminus\{i,i+1\} \\ b\notin i+1+I_N \\ c\notin j+1+I_N}}\bar\eta_b\bar\eta_c A(i,b)A(j,c)\bigg]
.
\nonumber
\end{align}
In particular, it is of the form $N^{-3}X^{v_4}_{4,\{0\}}$ for a $v_4:\Lambda_N^4\rightarrow\R$ satisfying $\sup_{N}\sup_{\Lambda_N^4}|v_4|<\infty$, and therefore controllable with size $N^{-1}$ by Corollary~\ref{coro_estimate_W_section_Ising}. For each $\delta>0$, we have proven the existence of a controllable function $R_\delta:\Omega_N\rightarrow\R$ such that:
\begin{equation}
\Big|\nu^N_{g_h}\big(\overrightarrow X^{N}_2\Big)\Big| 
\leq 
2\delta N^2\nu^N_{g_h}\big(\Gamma_h(\sqrt{f})\big)
+ 
\nu^N_{g_h}\big(f( \overrightarrow Y^{N}_2 + R_\delta )\big),
\qquad 
R_\delta 
:= 
R_1^\delta +  R_2+ R_3
.
\end{equation}
The arguments above do not depend on the sign of $A$ in the definition~\eqref{eq_def_X_left_X_right_W_3} of $\overrightarrow X^{N}_2$. This implies that $\overrightarrow X^{N}_2$ is $\Gamma$ controllable with size $N^{-1/2}$ in the sense of Definition~\ref{def_controllability}. Since the same arguments also apply to $\overleftarrow X^{N}_2$, Proposition~\ref{prop_error_terms} is proven.
\subsection{Proof of Proposition~\ref{prop_Boltzmann_gibbs_sec3} and Corollary~\ref{coro_Boltzmann_gibbs_sec3}}\label{sec_proof_boltzmann_gibbs}
\begin{proof}
Let $T>0$. 
Corollary~\ref{coro_Boltzmann_gibbs_sec3} is obtained as a side product of the proof of Proposition~\ref{prop_Boltzmann_Gibbs},  
which we focus on.
Write $G^N$ for either the error term $E^N$ with size $a_N$, 
or the $(\Gamma$-$)$controllable function $F^N$ with size $1$.  
For each $\gamma>0$, 
Feynman-Kac inequality~\eqref{eq_FK_sous_mu_general} and the bound~\eqref{eq_bound_adjoint_dansl_lemma_L_star} on the adjoint $L^*_h{\bf 1}$ in $\mathbb{L}^2(\nu^N_{g_h})$ imply:
\begin{align}
&\frac{1}{T}\log \E^{\nu^N_{g_h}}_h\Big[\exp\Big[\gamma \int_0^T G^N(\eta_t)dt\Big]\Big]
\nonumber\\
&\hspace{3cm}\leq 
\sup_{f\geq 0: \nu^N_{g_h}(f)=1}\Big\{ 
\nu^N_{g_h}\Big(f\big(\gamma G^N + \e/2\big)\Big) - \frac{N^2}{4}\nu^N_{g_h}\big(\Gamma_h(\sqrt{f})\big)
\Big\}
.
\end{align}
In addition, 
the function $\e$ satisfies, 
for some constant $C = C(\rho_\pm,h)>0$:
\begin{equation}
\nu^N_{g_h}\big(f\e\big) 
\leq 
\frac{H(f\nu^N_{g_h}|\nu^N_{g_h})}{4C_{LS}} 
+ \frac{2C}{N^{1/2}}
\leq
\frac{N^2}{4}\nu^N_{g_h}\big(\Gamma_h(\sqrt{f})\big)
+ \frac{2C}{N^{1/2}}
, 
\end{equation}
where we used the entropy inequality first, 
then the log-Sobolev inequality of Proposition~\ref{prop_LSI} to get the right-hand side. 
Thus:
\begin{align}
&\frac{1}{T}\log \E^{\nu^N_{g_h}}_h\Big[\exp\Big[\gamma \int_0^T G^N(\eta_t)dt\Big]\Big]
\nonumber\\
&\hspace{3cm}\leq 
\sup_{f\geq 0: \nu^N_{g_h}(f)=1}\Big\{ 
\nu^N_{g_h}\big(f\gamma G^N \big) - \frac{N^2}{8}\nu^N_{g_h}\big(\Gamma_h(\sqrt{f})\big) + \frac{C}{N^{1/2}}
\Big\}
.
\label{eq_preuve_BG_0}
\end{align}
Suppose first that $G^N$ is controllable with size $s_N$ (including both $s_N = O_N(1)$ and $s_N = o_N(1)$ cases): 
for some $\gamma_0>0$,
\begin{equation}
\forall N\in\N^*,\qquad
\log \nu^N_{g_h}\big[e^{\gamma_0 G^N}\big] 
\leq s_N
.
\end{equation}
By the entropy- and log-Sobolev inequalities, 
the quantity inside the supremum in
Equation~\eqref{eq_preuve_BG_0} is bounded above,
for each density $f$ for $\nu^N_{g_h}$, 
by:
\begin{equation}
\Big(\frac{\gamma C_{LS}}{\gamma_0}- \frac{1}{8}\Big)N^2\nu^N_{g_h}\big(\Gamma_h(\sqrt{f})\big) + \frac{\gamma s_N}{\gamma_0} + \frac{C}{N^{1/2}}
.
\end{equation}
Taking any $\gamma<\frac{\gamma_0}{8C_{LS}}$ ensures that the first term is negative, 
and concludes the proof in the controllable case:
\begin{equation}
\forall N\in\N^*,\qquad
\frac{1}{T}\log \E^{\nu^N_{g_h}}_h\Big[\exp\Big[\gamma \int_0^T G^N(\eta_t)dt\Big]\Big]
\leq 
\frac{\gamma s_N}{\gamma_0} + {C}{N^{1/2}}.
\end{equation}
If $G^N$ is not only controllable but also of LS type (recall Definition~\ref{def_LS_type}), 
corresponding to the controllable case in Corollary~\ref{coro_Boltzmann_gibbs_sec3}, 
then $\gamma_0>2^{10}C_{LS}>8 C_{LS}$ by assumption, 
thus one can take $\gamma >1$.\\

If now $G^N$ is $\Gamma$-controllable with size $s_N$, 
the idea is the same, 
except that one first bounds $\nu^N_{g_h}(fG^N)$ from above using Definition~\ref{def_controllability} of $\Gamma$-controllability:
\begin{equation}
\forall \delta>0,\qquad
\nu^N_{g_h}(fG^N) 
\leq 
\delta N^2\nu^N_{g_h}\big(\Gamma_h(\sqrt{f})\big) 
+ \frac{1}{\delta}\nu^N_g(Y_{G^N})
,
\end{equation}
with $Y_{G^N}$ controllable with size $s_N$. 
Choosing $\delta = \frac{1}{16}$, 
the proof becomes identical to the controllable case. 
\end{proof}
\begin{remark}
If one is interested only in estimating the expectation of the time integral of $G^N$ rather than its exponential moment, 
then the log-Sobolev inequality is not necessary. 
One can instead directly rely on Theorem~\ref{theo_entropic_problem} that bounds the relative entropy along the dynamics.
\demo
\end{remark}
\section{The Neumann condition on the diagonal}\label{sec_Neumann}
Let $h\in\s(\infty)$ 
(this set is defined in~\eqref{eq_def_s_infty}) 
Assume that $\bar\rho'\leq \epsilon_B$, $h\in\s(\epsilon_B)$ so that the solution $g_h\in g_0+\s(\infty)$ of the main equation~\eqref{eq_main_equation} exists and the conclusions of Lemma~\ref{lemm_L_star_as_e_plus_carre_du_champ} hold. 
In this section, we rewrite the term:
\begin{equation}
\frac{1}{4}\sum_{i<N-1}\bar\eta_i\bar\eta_{i+1}\big(\partial_1 h_{i_+,i} - \partial_1 h_{i_-,i}\big)
=
\frac{1}{4}\sum_{i<N-1}\omega_i\omega_{i+1}\bar\sigma_i\bar\sigma_{i+1}\big(\partial_1 h_{i_+,i} - \partial_1 h_{i_-,i}\big),\qquad \omega_{\cdot} = \frac{\bar\eta_{\cdot}}{\bar\sigma_{\cdot}},
\label{eq_terme_Neumann_appendix}
\end{equation}
as a function of the two-point correlations field $\Pi^N$, defined in~\eqref{eq_def_Pi}. 
This is necessary in the proof of upper-bound large deviations, in order to obtain a closed expression of the Radon-Nikodym derivative in terms of the field $\Pi^N$. 
It is done through the integration by parts Lemma~\ref{lemm_IPP},
replacing $\omega_i\omega_{i+1}$ by local averages of $\omega$'s. 
As $h\in\s(\infty)$, 
the function $\delta_h(x) = \partial_1 h(x_+,x)-\partial_1 h(x_-,x)$, 
$x\in(-1,1)$ can be extended into an element of $C^2([-1,1])$, 
still denoted by $\delta_h$, 
which satisfies $\delta_h(\pm 1) = 0$. 
Let $\epsilon\in(0,1)\in \N^*$ and $I_{\epsilon N} := \{0,...,N\epsilon-1\}$, 
writing $\epsilon N$ for $\lfloor \epsilon N\rfloor$. 
In the large $N$ limit, the correlation fields $\Pi$ we consider act on $\mathcal T$, 
which only contains functions with a certain regularity. 
We cannot simply replace $\omega_{i+1}$ by a uniform average of $\omega_\cdot$ on $i+1+I_{\epsilon N}$ and obtain an element of $\mathcal T$ 
(the indicator function ${\bf 1}_{[0,\epsilon)}$ is not regular enough). 
Consider instead a smooth function $\chi^\epsilon\in C^\infty(\bar\square)$ with $\chi^\epsilon = 0$ on $\partial\square$, 
$0\leq \chi^\epsilon \leq 2/\epsilon$, 
and such that $\chi^\epsilon(x,\cdot)$ approximates $\epsilon^{-1}{\bf 1}_{(x,x+\epsilon)\cap(-1,1)}$ in the following sense: 
$\chi^\epsilon(x,\cdot)$ is supported on $(x,x+\epsilon)\cap(-1,1)$ for each $x\in(-1,1)$, and:
\begin{equation}
\int_{\square}\big|\chi^\epsilon(x,y) - \epsilon^{-1}{\bf 1}_{(x,x+\epsilon)\cap(-1,1)}(y)\big|^2\, dx\, dy
\leq
\epsilon
.
\label{eq_L2_convergence_Neumann}
\end{equation}
Define then, 
recalling that $\bar\sigma(x) = \bar\rho(x)(1-\bar\rho(x))$ for $x\in[-1,1]$:
\begin{equation}
\forall (x,y)\in\squaredash,\qquad 
\mathcal N^{\delta_h}_\epsilon(x,y) 
= 
\frac{\bar\sigma(x)}{\bar\sigma(y)}\delta_h(x)\chi^\epsilon(x,y)
.
\end{equation}
Note that $\mathcal N^{\delta_h}_\epsilon$ belongs to $\mathcal T$, 
defined in~\eqref{eq_def_T}, 
thus $\Pi(\mathcal N^{\delta_h}_\epsilon)$ is now a well defined object for each $\Pi\in\mathcal T'_s$.
At the microscopic level, $\Pi^N$ also acts on functions with less regularity and we have the following more general result.
\begin{proposition}\label{prop_neumann_condition_averaging}
Let $q\in C^0([-1,1])$. 
Define, for $\eta\in\Omega_{N}$:
\begin{equation}
W_q(\eta) 
= 
\frac{1}{4}\sum_{i<N-1}\bar\eta_i\bar\eta_{i+1}q_i
.
\end{equation}
For each $\epsilon\in(0,1)$ smaller than some $\epsilon_0(\rho_\pm,h,q)>0$,  
there are constants $C_1(\rho_\pm,h,q)>0$ and $C_2(\rho_\pm,h,q,\epsilon)>0$ such that, 
for each $N$ larger than some constant depending on $\epsilon$ and each $T>0$:
\begin{align}
&\E^{\nu^N_{g_h}}_h\bigg[\exp\Big[\int_0^T \big[W_q(\eta_t) - \Pi^N_t\big(\mathcal N^q_\epsilon\big)\big]dt\Big] \bigg] 
\nonumber\\
&\hspace{3cm}\leq 
\exp\Big[ C_1(\rho_{\pm},h,q)\epsilon^{1/2} T + \frac{C_2(\rho_{\pm},h,q,\epsilon)T}{N^{1/2}}\Big]
.
\label{eq_estimate_replacement_neumann}
\end{align}
\end{proposition}
\begin{proof}
Let $T>0$, $\epsilon\in(0,1)$, 
and write:
\begin{equation}
W_q - \Pi^N(\mathcal N^q_\epsilon)
=
W_q - \overrightarrow{W}^{\epsilon N}_q+  \overrightarrow{W}^{\epsilon N}_q -\Pi^N(\mathcal N^q_\epsilon)
,
\label{eq_splitting_Neumann}
\end{equation}
where $\overrightarrow{W}^{\epsilon N}_q$ corresponds to $W_q$ in which the unit mass at each $i+1<N$ has been replaced by a mass $(\epsilon N)^{-1}$ at each site in $\{i+1,...,i+\epsilon N\wedge N-1\}$ 
(recall that $\epsilon N = \lfloor \epsilon N\rfloor$):
\begin{equation}
\overrightarrow{W}^{\epsilon N}_q 
:=
\Pi^N(A_\epsilon),
\qquad
A_\epsilon(x,y) := \frac{\bar\sigma(x)}{\bar\sigma(y)}q(x)\, \epsilon^{-1}{\bf 1}_{(x,x+\epsilon)\cap(-1,1)}(y)
,\quad (x,y)\in\square
.
\end{equation}
Up to applying Cauchy-Schwarz inequality in the exponential moment in~\eqref{eq_estimate_replacement_neumann}, 
it is enough to separately estimate the contribution of each difference in~\eqref{eq_splitting_Neumann}.\\

Consider first the contribution of $\overrightarrow{W}^{\epsilon N}_q -\Pi^N(\mathcal N^q_\epsilon)$. 
According to Corollary~\ref{coro_estimate_W_section_Ising},  for some $\gamma,C>0$ independent of $\epsilon,N,q$:
\begin{equation}
\frac{1}{T}\log \E^{\nu^N_{g_h}}_h\bigg[\exp\Big[\gamma \int_0^T \|A_\epsilon -\mathcal N^q_\epsilon\|_{2,N}^{-1}\big|\overrightarrow{W}^{\epsilon N}_q(\eta_t) -\Pi^N_t(\mathcal N^q_\epsilon)\big|\, dt\Big] \bigg] 
\leq
C
.
\label{eq_controllability_replacement_Neumann}
\end{equation}
Hölder inequality applied to~\eqref{eq_controllability_replacement_Neumann} 
then yields:
\begin{equation}
\frac{1}{2T}\log \E^{\nu^N_{g_h}}_h\bigg[\exp\Big[2\int_0^T \big|\overrightarrow{W}^{\epsilon N}_q(\eta_t) -\Pi^N_t(\mathcal N^q_\epsilon)\big|\, dt\Big] \bigg] 
\leq
\frac{1}{\gamma}\|A_\epsilon -\mathcal N^q_\epsilon\|_{2,N}
.
\label{eq_first_exp_mom_Neumann}
\end{equation}
Let us estimate $\|A_\epsilon -\mathcal N^q_\epsilon\|_{2,N}$. 
Since $\chi^\epsilon$ is smooth, 
\begin{equation}
\max_{(i,j)\in\Lambda_N^2}\sup_{(r,s)\in[0,1)^2}\Big|\chi^\epsilon_{i,j} -  \chi^\epsilon\Big(\frac{i+r}{N},\frac{j+s}{N}\Big)\Big|
\leq 
\frac{C(\epsilon)}{N}
.
\label{eq_discrete_L2_norm_chi_epsilon}
\end{equation}
It is also not difficult to check that, for some different $C(\epsilon)>0$:
\begin{align}
&\frac{1}{N^2}\sum_{(i,j)\in\Lambda_N^2} \int_{[0,1]^2} \epsilon^{-2}\bigg[{\bf 1}\Big\{\frac{j}{N}\in\Big(\frac{i}{N},\frac{i}{N}+\epsilon\Big)\Big\} - {\bf 1}\Big\{\frac{j}{N}\in\Big(\frac{i}{N},\frac{i}{N}+\epsilon\Big)\Big\}\bigg]^2\, du\, dv
\nonumber\\
&\hspace{9cm}\leq
\frac{3}{\epsilon^2 N}
.
\label{eq_discrete_L2_norm_indicator}
\end{align}
From~\eqref{eq_discrete_L2_norm_chi_epsilon}--\eqref{eq_discrete_L2_norm_indicator} and property~\eqref{eq_L2_convergence_Neumann} of $\chi^\epsilon$, 
using also the elementary bound $(a+b)^2\leq 2a^2+2b^2$ for real $a,b$; 
the difference $A_\epsilon -\mathcal N^q_\epsilon$ therefore satisfies:
\begin{align}
\|A_\epsilon -\mathcal N^q_\epsilon\|_{2,N}^2
:=
\frac{1}{N^2}\sum_{(i,j)\in\Lambda_N^2}\big(A_\epsilon -\mathcal N^q_\epsilon\big)_{i,j}^2
&\leq 
2\|A_\epsilon-\mathcal N^q_\epsilon\|_2^2 + \frac{C(\rho_\pm,q,\epsilon)}{N}
\nonumber\\
&\leq 
2\epsilon + \frac{C(\rho_\pm,q,\epsilon)}{N}
.
\end{align}
This yields a bound on~\eqref{eq_first_exp_mom_Neumann} of the same form as the right-hand side in Proposition~\ref{prop_neumann_condition_averaging} for any $\epsilon>0$ and any $N$ large enough depending on $\epsilon,\rho_\pm,q$.

Consider now the contribution of $W_q - \overrightarrow{W}^{\epsilon N}_q$ to~\eqref{eq_estimate_replacement_neumann}. 
The idea is the same as in the proof of Proposition~\ref{prop_Boltzmann_Gibbs}: 
we express $W_q - \overrightarrow{W}^{\epsilon N}_q$ in terms of the carré du champ operator and explicit controllable functions with size vanishing when $\epsilon$ is small. 
We start from the bound~\eqref{eq_preuve_BG_0} on exponential moments:
\begin{align}\label{eq_starting_point_Neumann}
&\frac{1}{T}\log \E^{\nu^N_{g_h}}_h\bigg[\exp\Big[\int_0^T \big[W_q(\eta_t) - \overrightarrow{W}^{\epsilon N}_q(\eta_t)\big)\big]dt\Big] \bigg]\\
&\qquad \leq 
\sup_{f\geq 0 : \nu^N_{g_h}(f) = 1}\bigg\{ \nu^N_{g_h}\Big(f\big[W_q - \overrightarrow{W}^{\epsilon N}_q\big]\Big) - \frac{N^2}{8}\nu^N_{g_h}\big(\Gamma_h(\sqrt{f})\big) + \frac{C(\rho_\pm,h)}{N^{1/2}}\bigg\}\nonumber.
\end{align}
Recall Definition~\ref{def_controllability} of controllability. 
To obtain the claim~\eqref{eq_controllability_replacement_Neumann}, 
it is enough to prove:
\begin{equation}
\epsilon^{-1/2}\big[W_q - \overrightarrow{W}^{\epsilon N}_q\big]
\text{ is }(\Gamma\text{-})\text{controllable with size 1}
.
\label{eq_sufficient_claim_Neumann}
\end{equation}
Indeed, if so, $W_q - \overrightarrow{W}^{\epsilon N}_q$ will be $(\Gamma$-)controllable with size $\epsilon^{1/2}$, 
and of LS type (see Definition~\ref{def_LS_type}) for $\epsilon$ small enough depending on $q,\rho_\pm,h$. \\

To prove~\eqref{eq_sufficient_claim_Neumann}, 
we use the integration by parts Lemma~\ref{lemm_IPP}. 
It is formulated with the variables $\omega_i = \bar\eta_i/\bar\sigma_i$, $i\in\Lambda_N$, 
for which $W_q$ becomes:
\begin{equation}
\forall \eta\in\Omega_N,\qquad 
W_q(\eta) 
= 
\frac{1}{4}\sum_{i<N-1}\omega_{i}\omega_{i+1}\bar\sigma_i\bar\sigma_{i+1}q_i = \frac{1}{4}\sum_{i<N-1}\omega_{i}\omega_{i+1}(\bar\sigma_i)^2q_i + \theta^{N,0}(\eta)
,
\end{equation}
where $\theta^{N,0}$ is the error term:
\begin{equation}
\forall\eta\in\Omega_N,\qquad 
\theta^{N,0}(\eta) 
= 
\frac{1}{4N}\sum_{i<N-1}\omega_{i}\omega_{i+1}N\bar\sigma_i\big[\bar\sigma_{i+1}-\bar\sigma_i\big]q_i
.
\end{equation}
It is of the form $N^{-1}X^{v}_{1,\{0,1\}}$ with the notations of Theorem~\ref{theo_concentration_Ising}, 
thus controllable with size $N^{-1}$ and of vanishing type 
(recall the terminology of Definitions~\ref{def_controllability}--\ref{def_LS_type}).\\ 
Recall from~\eqref{eq_def_rightarrow_omega} the definition of the quantity $\overrightarrow{\omega}^{\epsilon N}_{i+1}$:
\begin{equation}
\overrightarrow{\omega}^{\epsilon N}_{i+1}
=
\frac{1}{\epsilon N} \sum_{j=i+1}^{\min\{i+\epsilon N,N-1\}}\omega_j 
+ {\bf 1}_{i+\epsilon N>N-1}\Big(1-\frac{N-1-i}{\epsilon N}\Big)\omega_{N-1}
,\qquad
i<N-1
.
\end{equation}
For each $i<N-1$, one can write as before:
\begin{equation}
\omega_{i+1} - \overrightarrow{\omega}_{i+1}^{\epsilon N} 
= 
\sum_{a=i+1}^{\min\{i+N,N-1\}-1} \phi_{\epsilon N}(a-(i+1))(\omega_{a}-\omega_{a+1}),
\end{equation}
with:
\begin{equation}
\phi_{\epsilon N}(b) 
= 
\frac{\epsilon N-1-b}{\epsilon N}{\bf 1}_{0\leq b<\epsilon N},\quad b\in\Z
.
\end{equation}
Define then $u_j:\Omega_N\rightarrow\R$ for $j>-(N-1)$ similarly to~\eqref{eq_def_u_y}:
\begin{equation}
\forall \eta\in\Omega_N,\qquad 
u_j(\eta) 
= 
\sum_{i<j} \bar\sigma_i^2q_i\omega_i \phi_{\epsilon N}(j-(i+1))
.
\label{eq_def_u_j_Neumann}
\end{equation}
With this definition, the quantity $W_q - \overrightarrow{W}^{\epsilon N}_q$ reads, 
for each $\eta\in\Omega_N$:
\begin{align}
W_q(\eta) - \overrightarrow{W}^{\epsilon N}_q(\eta)
&= 
\sum_{j>-(N-1)}\big(\omega_j - \omega_{j+1}\big) u_j(\eta) \nonumber\\
&\quad
+ \omega_{N-1}\sum_{i<N-1}{\bf 1}_{i+\epsilon N>N-1}\Big(1-\frac{N-1-i}{\epsilon N}\Big)\omega_i(\bar\sigma_i)^2
+ \theta^{N,0}(\eta)
.
\label{eq_IPP_Neumann}
\end{align}
Fix a density $f$ for $\nu^N_{g_h}$. 
The first term in the second line involves boundary correlations. 
According to Lemma~\ref{lemm_small_boundary_correl_sous_bar_nu_G},  
it is $\Gamma$-controllable with size $N^{-1}$ and of vanishing type: 
for each $\delta>0$, 
there is a function $\mathcal D_{\delta}$, 
controllable with size $N^{-1}$ and of vanishing type, 
such that:
\begin{align}
&\nu^N_{g_h}\bigg(f\omega_{N-1}\sum_{i<N-1}{\bf 1}_{i+\epsilon N>N-1}\Big(1-\frac{N-1-i}{\epsilon N}\Big)(\bar\sigma_i)^2q_i\omega_i\bigg) 
\nonumber\\
&\hspace{3cm}\leq 
\delta N^2 \nu^N_{g_h}\big(\Gamma_h(\sqrt{f})\big) + \nu^N_{g_h}\big(f\mathcal D_{\delta}\big)
.
\label{eq_def_D_lambda}
\end{align}
It therefore remains to estimate the other term in the right-hand side of~\eqref{eq_IPP_Neumann}. 
By the integration by parts Lemma~\ref{lemm_IPP_forme_dir} applied to $-u_j$ for each $j<N-1$, 
there is a constant $C>0$ such that, 
for any $\delta>0$:
\begin{align}
\nu^N_{g_h}\Big(&f\sum_{j<N-1}\big(\omega_i - \omega_{i+1}\big) u_j\Big) 
\leq 
\delta N^2 \nu^N_{g_h}\big(\Gamma_h(\sqrt{f})\big) 
\nonumber\\
&\quad+ 
\frac{C}{\delta N^2} \sum_{j<N-1}\int f|u_j|^2d\nu^N_{g_h} 
\quad(:= \frac{1}{\delta}\nu^N_{g_h}(f\mathcal N_1))\nonumber\\
&\quad+
\sum_{j<N-1}(\bar\rho_{j+1}-\bar\rho_j)\int \omega_j\omega_{j+1} e^{-(\eta_{j+1}-\eta_j)C^{g_h}_j/N }f u_jd\nu^N_{g_h} 
\quad(:= \nu^N_{g_h}(f\mathcal N_2))\nonumber\\
&\quad- 
\sum_{j<N-1}\int \big(\omega_{j+1}-\omega_{j}\big)\Big(1-e^{-(\eta_{j+1}-\eta_j)C^{g_h}_j/N}\Big) f u_jd\nu^N_{g_h} 
\quad(:= \nu^N_{g_h}(f\mathcal N_3))
.
\label{eq_def_mathcal_N_3}
\end{align}
The functions $\mathcal N_i$ ($1\leq i \leq 3$) are then estimated as in~\eqref{eq_bound_R_1}--\eqref{eq_second_line_R_2}--\eqref{eq_bound_R_3} respectively, 
but not all of them are error terms and one has to be careful to get a bound that vanishes with $\epsilon$. 
Let us check that each of them indeed satisfies~\eqref{eq_sufficient_claim_Neumann}. 
For $\mathcal N_1$,
~\eqref{eq_def_u_j_Neumann} and the bound $\|u_j\|_\infty\leq \|q\|_\infty\epsilon N$ imply that it is of the form $N^{-1}X^{v_2}_{2,\{0\}}$ with $v_2$ given by:
\begin{equation}
v_2(a,b)
=
\frac{C}{\delta N}\sum_{j<N-1}{\bf 1}_{j>\max\{a,b\}} q_aq_b\bar\sigma_a^2\bar\sigma_b^2\phi_{\epsilon N}(j-(a+1))\phi_{\epsilon N}(j-(b+1))
\leq
\frac{C\|q\|_\infty^2\epsilon}{\delta }
.
\end{equation}
It follows that $\epsilon^{-1} \mathcal N_1$ satisfies~\eqref{eq_sufficient_claim_Neumann}, 
thus $\epsilon^{-1/2}\mathcal N_1$ as well.\\

Consider now $\mathcal N_2$. 
From~\eqref{eq_second_line_R_2}, we get:
\begin{equation}
\nu^N_{g_h}\big(f\mathcal N_2\big) 
\leq 
\frac{\bar\rho'}{N}\sum_{j<N-1}\nu^N_{g_h}\big(\omega_j\omega_{j+1}u_j \big)
+ 
\frac{\bar\rho'e^{2\|g_h\|_\infty}}{N}\sum_{j<N-1}\frac{1}{\bar\sigma_j\bar\sigma_{j+1}}\nu^N_{g_h}\big(f|C_j^{g_h}u_j|\big)
.
\label{eq_estimate_N_2}
\end{equation}
Here we do not even need the $\Gamma$-controllability with size $N^{-1/2}$ of the first term established in Proposition~\ref{prop_Boltzmann_Gibbs}. 
Instead, 
recall that $C^{g_h}_j = B^{g_h}_j + D^{g_h}_j$ ($j<N-1$) with:
\begin{equation}
B^{g_h}_j(\eta) 
= 
\frac{1}{2N}\sum_{\ell\notin\{j,j+1\}}\bar\eta_\ell\partial^N_1 (g_h)_{j,\ell},
\qquad 
D^{g_h}_j(\eta) 
= 
\frac{\partial^N\bar\rho_{j} \, (g_h)_{j,j+1}}{2N} 
=
O(N^{-1})
.
\end{equation}
It is then enough, 
recalling that $u_j$ is given by~\eqref{eq_def_u_j_Neumann}, 
to notice that $\nu^N_{g_h}(f\mathcal N_2)$ satisfies:
\begin{align}
\nu^N_{g_h}\big(f\mathcal N_2\big)
&\leq 
\nu^N_{g_h}\big(f \big[N^{-1}X^{w_2}_{2,\{0,1\}}\big]\big) 
\nonumber\\
&\qquad+ 
\frac{\bar\rho' e^{2\|g_h\|_\infty}}{2N}\sum_{j<N-1}\frac{1}{\bar\sigma_j\bar\sigma_{j+1}}\nu^N_{g_h}\Big(f\big[ |N^{-1}X^{w^j_2}_{2,\{0\}}| + Y_2\big]\Big)
,
\end{align}
where:
\begin{align}
w_2(i,j) 
&=
{\bf 1}_{i<j}\frac{\bar\rho'\bar\sigma_i}{\bar\sigma_j\bar\sigma_{j+1}}q_i\phi_{\epsilon N}(j-(i+1))
,
\nonumber\\
w^j_2(a,b) 
&=
\frac{1}{2}{\bf 1}_{a\notin\{j,j+1\}}{\bf 1}_{b<j} \bar\sigma_bq_b\, \phi_{\epsilon N}(j-(b+1)) (\partial^N_1 g_h)_{j,a}
,
\\
Y_2 &= \sum_{j<N-1}|u_j D^{g_h}_j| 
\leq 
C(\rho_\pm,q,g_h) \epsilon
\nonumber
.
\end{align}
Due to the fact that $\phi_{\epsilon N}$ is non-zero only for $\epsilon N$ different integers, 
$w_2$ and each $w^j_2$ have $\|\cdot\|_{2,N}$-norm bounded by $C\epsilon^{1/2}$ for some $C=C(\rho_\pm,q,g_h)>0$ that does not depend on $j<N-1$. 
It follows from Corollary~\ref{coro_estimate_W_section_Ising} that $\epsilon^{-1/2}\mathcal N_2$ is controllable with size $1$ and of large type, 
i.e. satisfies~\eqref{eq_sufficient_claim_Neumann}. \\

Consider finally $\mathcal N_3$, 
defined in~\eqref{eq_def_mathcal_N_3}. 
The bound~\eqref{eq_bound_R_3} shows that $\nu^N_{g_h}(f\mathcal N_3)$ is bounded by a constant times the second term in~\eqref{eq_estimate_N_2}. 
It follows that $\epsilon^{-1/2}\mathcal N_3$ also satisfies~\eqref{eq_sufficient_claim_Neumann}, 
which concludes the proof.
\end{proof}
\section{Sobolev spaces}\label{app_sobolev_spaces}
\begin{definition}
Let $U\subset \R^2$ be a bounded open set with Lipschitz boundary. 
For $n\in\N$ and $p\geq 1$, let $\mathbb W^{n,p}(U)$ be the following space. 
If $n=0$, it is simply $\mathbb L^p(U)$. 
If $n\geq 1$, $\mathbb W^{n,p}(U)$ is the set of functions $f\in\mathbb L^p(U)$ such that, 
for any $(n_1,n_2)\in\N^2$ with $n_1+n_2\leq n$, there is $f^{n_1,n_2}\in\mathbb L^p(U)$ satisfying:
\begin{align}
&\forall u\in C^\infty_c(U),\qquad 
\int_{U}f(x,y) \partial_1^{n_1}\partial_2^{n_2} u(x,y)\, dx\, dy 
\nonumber\\
&\hspace{5cm}
= 
(-1)^{n_1+n_2}\int_{U}f^{n_1,n_2}(x,y)u(x,y)\, dx\, dy
.
\end{align}
$\mathbb W^{n,p}(U)$ is a separable Banach space for the norm:
\begin{equation}
\forall f\in \mathbb W^{n,p}(U),\qquad \|f\|_{\mathbb W^{n,p}(U)} = \Big[\sum_{\substack{(n_1,n_2) \in \N^2 \\ n_1+n_2 \leq n}}\|f^{n_1,n_2}\|^2_{\mathbb L^p(U)}\Big]^{1/2}.
\end{equation}
Moreover, the set $C^\infty(\bar U)$ of restrictions of elements of $C^\infty(\R^2)$ to $\bar U$ is dense in $\mathbb W^{n,p}(U)$ for $\|\cdot\|_{\mathbb W^{n,p}(U)}$. 
In the special case $p=2$, define $\mathbb H^{n}(U):= \mathbb W^{n,p}(U)$. This is a Hilbert space.
\end{definition}
Along the text, we make use of the following Sobolev embedding results (see Theorem 4.12 in~\cite{Adams2003} and Theorem 1.4.4.1 in~\cite{Grisvard2011}).
\begin{proposition}\label{prop_sobolev_embeddings}
Let $U\subset \R^2$ be a bounded set with Lipschitz boundary. The following embeddings hold.
\begin{itemize}
	\item Let $p>2$ and $n\in\N^*$, then $\mathbb W^{n,p}(U)\subset C^{n}(\bar U)$.
	\item Let $p\geq 2$ and $n\in\N^*$, then $\mathbb W^{n,p}(U)\subset \mathbb W^{\ell,q}(U)$ for any $\ell\leq n-1$ and any $q\geq 1$. 
\end{itemize}
\end{proposition}
In our case, $U= \squaredash = \lhd\cup \rhd$, where we recall that $\square = (-1,1)^2$, $\squaredash = \square\setminus D$ and $\rhd = \{(x,y)\in\squaredash:x<y\} = \squaredash\setminus\{\lhd\}$. 
We are interested in the subset $\mathcal T'_s$, defined in~\eqref{eq_def_T_prime}, 
of the topological dual $\mathcal T'$ of $\mathcal T$. 
\begin{definition}\label{def_weak_star_convergence}
If $(X,\|\cdot\|_X)$ is a Banach space, let $X'$ be its topological dual, equipped with the norm:
\begin{equation}
\forall L\in X',\qquad N_{X}(L) = \sup_{ \phi\in X\setminus\{0\}}\frac{|L(\phi)|}{\|\phi\|_X}.\label{eq_def_norm_duale}
\end{equation}
If $\phi :\squaredash\rightarrow \R$, let $\phi_s(x,y) = [\phi(x,y)+\phi(y,x)]/2$ denote its symmetric part, and let $\mathcal T'_s\subset \mathcal T$ be the subset of elements $\Pi$ satisfying $\Pi(\phi) = \Pi(\phi_s)$ for any $\phi\in\mathcal T$. Then:
\begin{align}
\forall \Pi\in\mathcal T'_s,
\qquad \mathcal N_{\mathcal T}(\Pi) 
\, &:= 
\sup_{ \phi\in \mathcal T\setminus\{0\}}\frac{|\Pi(\phi)|}{\|\phi\|_{\mathbb H^2(\squaredash)}} 
= 
\sup_{ \phi\in \mathcal T\setminus\{0\}}\frac{|\Pi(\phi_s)|}{\|\phi_s\|_{\mathbb H^2(\squaredash)}} 
\nonumber\\
&= 
\sup_{ \phi\in \mathcal T\setminus\{0\}}\frac{|\Pi(\phi_{|\rhd})|}{\|\phi_{|\rhd}\|_{\mathbb H^2(\rhd)}}
,
\label{eq_norm_on_T_prime_s}
\end{align}
where $\phi_{|\rhd}$ is the restriction of $f$ to $\rhd$.
$\mathcal T'_s$ is closed for the norm~\eqref{eq_norm_on_T_prime_s}. \\
The weak$^*$ topology on $\mathcal T'$ is the topology of simple convergence: a sequence $(\Pi_n) \in (\mathcal T')^{\N}$ weak$^*$ converges to $\Pi\in \mathcal T'$ if and only if:
\begin{equation}
\forall \phi\in \mathcal T,\qquad 
\lim_{n\rightarrow\infty}\Pi_n(\phi) 
= 
\Pi(\phi)
\quad \Leftrightarrow\quad 
\forall \phi\in\mathcal T_{\rhd},\qquad 
\lim_{n\rightarrow\infty}\Pi_n(\phi_{|\rhd}) 
= 
\Pi(\phi_{|\rhd})
.
\end{equation}
The set $\mathcal T'_s$ is also closed for the weak$^*$ topology. 
We write $\big(\mathcal T'_s,* \big)$ when explicitly referring to this topology.
\end{definition}
\subsection{Compact sets}
In this section, 
we give a sufficient condition for compactness in $\big(\mathcal T'_s,* \big)$. 
Recall that $\mathcal T=\mathbb H^2(\squaredash)$. 
Banach-Alaoglu's theorem characterises compactness in $(\mathcal T'_s,*)$:
\begin{proposition}[Banach-Alaoglu]\label{prop_Banach_Alaoglu}
Let $\mathcal K\subset (\mathcal T'_s,*)$ be such that $\sup_{\mathcal K} \mathcal N_{\mathcal T}<\infty$. 
Then $\mathcal K$ is relatively weak$^*$ compact.
\end{proposition}
The norm $\mathcal N_{\mathcal T}$ is defined through a supremum, 
which is difficult to work with. 
Instead, we formulate a sufficient condition for compactness which involves a sum. 
Such a characterisation is known to hold when the underlying space is periodic, 
e.g. on the torus $\mathbb{T}^2_{-2,2}=[-2,2)^2$: 
a linear form $\pi\in(\mathbb H^{2}(\mathbb{T}^2_{-2,2}))'$ is bounded if and only if
\begin{equation}
\|\pi\|_{\mathbb{T},-2} 
:=
\sum_{m\in\N^2}(1+|m|^2)^{-2}\big|\pi(\phi_m)\big|^2<\infty,\label{eq_def_norme_moins_n}
\end{equation}
with $|m|^2 = m_1^2+m_2^2$ for $m=(m_1,m_2)\in\N^2$, 
and where $(\phi_m)_{m\in\N^2}$ is an orthonormal basis of $\mathbb L^2(\mathbb T^2_{-2,2})$ made of real eigenvalues of the Laplacian: 
for $(x,y)\in\mathbb T^2_{-2,2}$, writing $\N^*:=\N\setminus\{0\}$:
\begin{equation}
\phi_m(x,y) 
= 
\varphi_{m_1}(x)\varphi_{m_2}(y),\quad \varphi_{m_1}(x) = \begin{cases}
1/2\quad &\text{if }m_1 = 0,\\
2^{-1/2}\cos\Big(\frac{m'\pi x}{2}\Big) \quad &\text{if }m_1=2m'-1\in\N^*,\\
2^{-1/2}\sin\Big(\frac{m'\pi x}{2}\Big)\quad &\text{if }m_1=2m'\in\N^*
.
\end{cases}\label{eq_def_eigenvalues_Laplace_sobolev}
\end{equation}	 
The equivalence between~\eqref{eq_def_norme_moins_n} and $\mathcal N_{\mathcal T}$ in the periodic setting comes from the fact that a function has the same regularity as its Fourier transform. 
In our case, however, this property does not hold because of the boundaries in $\squaredash$, 
and $\mathcal N_{\mathcal T}$ is not equivalent to the norm $\|\cdot\|_{\mathbb{T},-2}$.

We look for a sufficient condition for compactness that can nonetheless be stated in terms of the norm $\|\cdot\|_{\mathbb{T},-2}$ defined in~\eqref{eq_def_norme_moins_n}. 
To do so, note first that, by~\eqref{eq_norm_on_T_prime_s}, 
it is sufficient to work on the triangle $\rhd = \{(x,y)\in\squaredash:x<y\}$. 
The idea is then to extend elements $\Pi\in \mathcal T'_s$ to linear forms acting on the larger space $\mathbb H^{2}(\mathbb{T}^2_{-2,2})$ of test functions, 
then check that the norms of $\Pi$ and its extension are comparable.  
Define then:
\begin{equation}\label{eq_prolongement_Pi}
\Pi^{\text{ext}}(u) = \Pi\big(u_{|\rhd}\big), 
\qquad 
u\in \mathbb H^2(\mathbb{T}^2_{-2,2})
.
\end{equation}
Clearly, $\Pi^{\text{ext}}$ is a linear form on $\mathbb H^2(\mathbb{T}^2_{-2,2})$, 
although it may not be bounded any more. 
The sufficient condition for compactness can now be stated.
\begin{proposition}\label{prop_cara_compact_sets}
Let $A>0$, 
and let $\mathcal K_A := \big\{\Pi\in\mathcal T' : \|\Pi^{\text{ext}}\|_{\T,-2}\leq A\}$. 
Then $\mathcal K_A$ is weak$^*$ relatively compact in $\mathcal T'_s$.
\end{proposition}
\begin{proof}
The goal is to bound the original norm $\mathcal N_{\mathcal T}$ by the norm $\|\cdot\|_{\mathbb{T},-2}$ from~\eqref{eq_def_norme_moins_n}, 
then use the Banach-Alaoglu theorem 
(Proposition~\ref{prop_Banach_Alaoglu}) to conclude. 
In view of~\eqref{eq_norm_on_T_prime_s}, 
it is sufficient to work with test functions defined on the triangle. 
Consider:
\begin{equation}
\mathcal T_{\rhd} = \big\{\phi_{|\rhd} : \phi\in\mathcal T \big\}.
\end{equation}
We first explain how to embed $\mathcal T_\rhd$ into $\mathbb H^2(\mathbb{T}^2_{-2,2})$ and bound the norm $\mathcal N_{\mathcal T}$ by the strong dual norm $\|\cdot\|_{\mathbb H^{-2}(\mathbb{T}^2_{-2,2})}$ (recall~\eqref{eq_def_norm_duale}). 
Since this norm is equivalent to the norm~\eqref{eq_def_norme_moins_n}, 
this will be enough to conclude. \\

By Theorem 1.4.3.1 in~\cite{Grisvard2011}, there is a continuous linear extension $P$ from $\mathcal T_\rhd$ to $\mathbb H^2(\R^2)$, 
i.e. there is $C(\rhd)>0$ such that:
\begin{equation}
\forall u\in\mathcal T_\rhd,\qquad 
P u \in\mathbb H^2(\R^2),\quad (Pu)_{|\rhd} 
= u,
\quad \|Pu\|_{\mathbb H^2(\R^2)}\leq C(\rhd)\|u\|_{\mathbb H^2(\rhd)}
.
\end{equation}
Let $\chi\in C^\infty(\R^2)$ be equal to $1$ on $\rhd$, and be compactly supported in $(-2,2)^2$. 
By Theorem 1.4.4.2 in~\cite{Grisvard2011}, $\mathbb H^n(\R^2)$ is a Banach algebra as soon as $n\geq 2$. It follows that there is $C(\chi,\rhd)>0$ such that:
\begin{align}
\forall u\in\mathcal T_\rhd,\qquad
\|\chi Pu\|_{\mathbb H^2(\R^2)}
&\leq 
C(\rhd)\|\chi\|_{{\mathbb H^2(\R^2)}}\|Pu\|_{\mathbb H^2(\R^2)}
=
C(\chi,\rhd)\|u\|_{\mathbb H^2(\rhd)}
\nonumber\\
&\leq 
C(\chi,\rhd)\|\chi Pu\|_{\mathbb H^2(\R^2)}
,\label{eq_equiv_norme__T_et_extension}
\end{align}
where the last inequality comes from the inclusion $\rhd\subset(-2,2)^2$ and the fact that $(\chi Pu)_{|\rhd} = u$. 
The mapping $\chi P$ is an embedding from $\mathcal T_\rhd$ to $\mathcal T^{\text{ext}}_{\rhd} := \mathbb H^2_0((-2,2)^2)$, 
the closure of $C^\infty$, compactly supported functions on $(-2,2)^2$ for the norm of $\mathbb H^2((-2,2)^2)$.  
An element $u$ of $\mathcal T^{\text{ext}}_{\rhd}$ can be turned into a periodic function in $\mathbb H^2(\T^2_{-2,2})$ with the same norm, 
by setting $u(\cdot + (4a,4b)) = u(\cdot)$ for each $(a,b)\in \Z^2$. \\
Let us now compare the elements of $\mathcal T'_s$ and their extensions to $\mathbb H^2(\mathbb{T}^2_{-2,2})$. 
Take $\Pi\in\mathcal T'_s$, and extend it to a linear form $\Pi^{\text{ext}}$ on $\mathbb H^2(\mathbb{T}^2_{-2,2})$ (possibly unbounded) through~\eqref{eq_prolongement_Pi}. 
Then:
\begin{align}
\mathcal N_{\mathcal T}(\Pi) 
:\hspace{-0,1cm}&= 
\sup_{u\in\mathcal T_\rhd\setminus \{0\}}\frac{|\Pi(u)|}{\|u\|_{\mathbb H^2(\rhd)}} 
= 
\sup_{u\in \mathcal T_\rhd\setminus \{0\}}\frac{|\Pi^{\text{ext}}(\chi P u)|}{\|u\|_{\mathbb H^2(\rhd)}}\nonumber\\
&\leq 
C(\chi,\rhd)\sup_{v\in\mathcal T^{\text{ext}}_{\rhd}\setminus \{0\}}\frac{|\Pi^{\text{ext}}(v)|}{\|v\|_{\mathbb H^2(\R^2)}} \nonumber\\
&\leq 
C(\chi,\rhd)\sup_{w\in\mathbb H^2(\mathbb{T}^2_{-2,2})\setminus \{0\}}\frac{|\Pi^{\text{ext}}(w)|}{\|w\|_{\mathbb H^2(\mathbb{T}^2_{-2,2})}}
=: 
C(\chi,\rhd)\|\Pi^{\text{ext}}\|_{\mathbb H^{-2}(\mathbb{T}^2_{-2,2})}
.
\label{eq_borne_norme_Pi}
\end{align}
The fact that $\|\cdot \|_{\mathbb H^{-2}(\T^2_{-2,2})}\leq c\|\cdot \|_{\T,-2}$ for some $c>0$ concludes the proof.
\end{proof}
\section{Poisson equations}\label{app_Poisson}
In this section, we give conditions for the existence and uniqueness of solutions of the various Poisson equations 
-- among which the Euler-Lagrange equation~\eqref{eq_EL} and the main equation~\eqref{eq_main_equation} 
-- encountered along the text. 
To do so, we prove that finding a kernel $k$ that solves the (linear) Euler-Lagrange equation or a solution $g$ of the (non-linear) main equation is the same, 
and rewrite all linear equations in a common framework.\\

Throughout Appendix~\ref{app_Poisson}, 
a function $u\in\mathbb L^2(\squaredash)$ is identified with the kernel operator $u\phi(\cdot) = \int u(\cdot,y)\phi(y)\,dy$ for $\phi\in\mathbb L^2((-1,1))$. 
View also $\bar\sigma$ as a multiplication operator: 
if $\phi\in\mathbb L^2((-1,1))$,
\begin{align}
(\bar\sigma u)(x,y) &= \bar\sigma(x)u(x,y),
\quad
(u\bar\sigma )(x,y) = u(x,y)\bar\sigma(y),
\nonumber\\
(u\bar\sigma \phi)(x) 
&= 
\int_{(-1,1)}u(x,y)\bar\sigma(y)\phi(y)\, dy
.
\end{align}
For a symmetric function $f\in\mathcal T$, 
write $\delta_f$ for the operator:
\begin{align}
\delta_f(x) 
&=
(\partial_2-\partial_1)f(x,x_+)
=
\partial_2f(x,x_+)-\partial_2 f(x,x_-)
\nonumber\\
&=
-(\partial_2-\partial_1)f(x_+,x)
,\qquad
x\in(-1,1)
.
\end{align}
\subsection{Euler-Lagrange equation}\label{app_EL}
The next proposition is classical, 
and proves Proposition~\ref{prop_Euler-Lagrange equation}.
\begin{proposition}\label{prop_EL_appendix}
Let $\Pi\in\mathcal T'_s$ be associated with a kernel $k$ via $\Pi =\frac{1}{4}\big<k,\cdot\big>$, 
and write $C_k := \bar\sigma +k$.  
Write for short $\mathcal I_{\infty}(k)$ for $\mathcal I_{\infty}(\Pi)$, 
and idem for $J_h(k)$. 
Assume $\mathcal I_{\infty}(k)<\infty$. 
There is then a generalised bias $h\in\mathbb L^2(\squaredash)$, 
with $h$ a symmetric function admitting a weak derivative that satisfies:
\begin{equation}
\int_{(-1,1)}\bar\sigma(z)\big<\partial_1h(z,\cdot),C_k\partial_1 h(z,\cdot)\big>\, dz 
<
\infty
,
\end{equation}
such that $\mathcal I_\infty(k) = J_h(k)$. 
Moreover, $k$ and $h$ are related through the following Euler-Lagrange equation: 
for any test function $\phi\in\mathbb H^1(\squaredash)$ with $\phi_{|\partial\square}=0$, 
\begin{equation}
\frac{1}{2}\int_{\squaredash}\nabla (k-k_0)\cdot\nabla\phi 
- \int_{(-1,1)}\bar\sigma(z)\big<\partial_1h(z,\cdot),C_k\partial_1 \phi(z,\cdot)\big>\, dz 
=
0
.
\label{eq_Euler-Lagrange_weak}
\end{equation}
This is a weak formulation of:
\begin{align}
\begin{cases}
&\Delta k  (x,y) - \big[\bar\sigma(x)\delta_h(x) + \bar\sigma(y)\delta_h(y)\big]k(x,y) 
\\ 
&\hspace{1.6cm}
- \displaystyle{\int_{(-1,1)} \Big[\partial_x\big(\bar\sigma(x)\partial_1 h(x,z)\big)k(z,y) + k(x,z)\partial_y\big(\partial_2 h(z,y)\bar\sigma(y)\big)\Big]dz}
\\
&\hspace{2cm}
=   \partial_x  \big(  \bar\sigma(x) \bar\sigma(y)   \; \partial_1  h(x,y) \big)
+  \partial_y  \big(  \bar\sigma(x) \bar\sigma(y)   \; \partial_2  h(x,y) \big) \quad \text{for }(x,y)\in\squaredash\, ,
\\
&h_{|\partial\square} = k_{|\partial\square} =0,
\\
&\bar\sigma(x)^2(\partial_2-\partial_1)h(x_+,x) - (\partial_2-\partial_1)k(x_+,x) 
=
(\bar\rho')^2\quad \text{for }x\in(-1,1).
\end{cases}
\label{eq_Euler-Lagrange_appendix}
\end{align}
\end{proposition}
\begin{proof}
The existence (and uniqueness) of the generalised bias $h$ such that $J_h(k) = \mathcal I_\infty(k)$ is classical and follows from arguments similar to those of Section~\ref{sec_dynamical_part_lower_bound}, 
see Lemma 5.3 in Chapter 10 of~\cite{Kipnis1999}. 
Let us show that $k$ satisfies the Euler-Lagrange equation with bias $h$.  
Since $h$ is such that $J_h(k) = \inf_{\tilde h\in\s(\infty)}J_{\tilde h}(k)$, 
one has $\epsilon^{-1}(J_{h\pm \epsilon\phi}-J_h)\leq 0$ for any $\phi\in\mathbb H^1(\squaredash)$ with $\phi_{|\partial\square}=0$ and small enough $\epsilon>0$. 
Thus:
\begin{align}
&\frac{1}{8} \int_{\squaredash} \nabla k\cdot\nabla\phi  
+\frac{1}{4}\int_{(-1,1)} (\bar\rho')^2 \phi(x,x)\,dx
\nonumber\\
&\quad
-\frac{1}{4}\int_{(-1,1)}\bar\sigma(z)\big< \partial_1 h(z,\cdot),(\bar\sigma + k)(\partial_1 \phi(z,\cdot))\big>\, dz 
=
0
.
\end{align}
Applying this equation in the $h=0$ case corresponding to the steady-state kernel $k_0$ (or just recalling the expression~\eqref{eq_def_k_0}), one gets:
\begin{equation}
\frac{1}{8} \int_{\squaredash} \nabla k_0\cdot\nabla\phi   
= 
- \frac{1}{4}\int_{(-1,1)} (\bar\rho')^2 \phi(x,x)\,dx
.
\end{equation}
This yields~\eqref{eq_Euler-Lagrange_appendix}. 
Equation~\eqref{eq_EL} follows by careful integration by parts. 
For instance, 
the term $\bar\sigma(y)\delta_h(y)k(x,y)$ comes from the contribution of $k$ in the integral involving $C_k = \bar\sigma +k$:
\begin{align}
&\int_{(-1,1)^2}\bar\sigma(z)\partial_1 h (z,x) (k\partial_1\phi)(z,x)\, dx\,dz
=
\int_{(-1,1)^2}\bar\sigma(z)\partial_1 h (z,x) \partial_1(k\phi)(z,x)\, dx\,dz
\nonumber\\
&\hspace{4cm}=
\int_{(-1,1)} \bar\sigma(x)\Big[\partial_1 h(x_-,x) - \partial_1 h(x_+,x)\Big](k\phi)(x,x)\, dx 
\nonumber\\
&\hspace{4cm}\qquad - 
\int_{(-1,1)^2} \partial_1\big(\bar\sigma(z)\partial_1 h(z,x)\big)(k\phi)(z,x)\, dx \, dz
.
\end{align}
Note that $k\phi(x,x)$ is simply $\int_{(-1,1)}k(y,x)\phi(y,x)\, dy$ as $k$ is symmetric, 
thus:
\begin{align}
&\int_{(-1,1)} \bar\sigma(x)\Big[\partial_1 h(x_-,x) - \partial_1 h(x_+,x)\Big](k\phi)(x,x)\, dx 
\nonumber\\
&\hspace{3cm}=
\int_{(-1,1)^2}k(x,y)\phi(x,y) \bar\sigma(y)\Big[\partial_1 h(y_-,y) - \partial_1 h(y_+,y)\Big]\, dx \, dy
.
\end{align}
Using the symmetry of $h$ and the fact that $\partial_1h(x_-,x)=\partial_1h(x,x_+)$ 
(it is the same point in the same triangle), 
the bracket involving $\partial_1 h$ is $-\delta_h$ as claimed:
\begin{equation}
\partial_1 h(x_-,x)-\partial_1h(x_+,x) 
= 
(\partial_2-\partial_1)h(x_+,x) 
= 
-\delta_h(x),
\qquad 
x\in(-1,1)
.
\end{equation}
\end{proof}
\subsection{Equivalence of the Euler-Lagrange equation and the main equation}\label{app_equivalence_EL_main_equation}
\begin{proposition}\label{prop_linear_main_equation}
Let $h\in\s(\infty)$, where this set is defined in~\eqref{eq_def_s_infty}.  
Recall that $\square := (-1,1)^2$, $\squaredash = \square\setminus D$ and $\rhd = \{(x,y)\in\squaredash:x<y\}$, $\lhd=\squaredash\setminus\rhd$. 
\begin{itemize}
	\item Let $k\in C^3(\bar\rhd)\cap C^3(\bar\lhd)$ be symmetric and suppose $k$ solves the Euler-Lagrange equation~\eqref{eq_Euler-Lagrange_appendix}. 
Assume that 
$k$ satisfies:
\begin{equation}
\int_{\squaredash}\bar\sigma^{-1}(x)k(x,y)^2\bar\sigma^{-1}(y)\, dx\,dy 
<
1
.
\label{eq_size_condition_k_appendix_F}
\end{equation}
Then 
the correlation operator $C_k=\bar\sigma+k$ is invertible in $\mathbb L^2(\squaredash)$.  
Define $g\in\mathbb L^2(\squaredash)$ through $C_k = \bar\sigma +k = (\bar\sigma^{-1}-g)^{-1}$.  
Then $g\in C^3(\bar\rhd)\cap C^3(\bar\lhd)$ is symmetric and solves the main equation~\eqref{eq_main_equation}.
	\item Conversely, let $g\in C^3(\bar\rhd)\cap C^3(\bar\lhd)$ be symmetric and solve the main equation~\eqref{eq_main_equation}, and assume:
\begin{equation}
	\int \bar\sigma(x)g(x,y)^2\bar\sigma(y)\, dx\, dy
	< 
	1
	.
	\label{eq_size_condition_g_appendix_F}
	\end{equation}	
	Then $\bar\sigma^{-1}-g$ is invertible, and 
	$k\in C^3(\bar\rhd)\cap C^3(\bar\lhd)$ defined through $\bar\sigma+k = (\bar\sigma^{-1}-g)^{-1}$ is symmetric and solves the Euler-Lagrange equation~\eqref{eq_Euler-Lagrange_appendix}.
\end{itemize}
\end{proposition}
\begin{remark}
Conditions~\eqref{eq_size_condition_k_appendix_F}--\eqref{eq_size_condition_g_appendix_F} hold if $\bar\rho'\leq \epsilon, h\in\s(\epsilon)$ for small enough $\epsilon$.
\demo
\end{remark}
\begin{proof}
We only show that $k\in C^3(\bar\rhd)\cap C^3(\bar\lhd)$ symmetric solving the Euler-Lagrange equation implies that $g$ solves the main equation, 
is symmetric and is in $C^3(\bar\rhd)\cap C^3(\bar\lhd)$, 
the other implication being similar.

By assumption~\eqref{eq_size_condition_k_appendix_F}, 
$g$ admits a series expansion:
\begin{align}
g 
&= 
-\bar\sigma^{-1/2}\big(1+\bar\sigma^{-1/2}k\bar\sigma^{-1/2}\big)^{-1}\bar\sigma^{-1/2}\, +\, \bar\sigma^{-1}
\nonumber\\
&=
-\bar\sigma^{-1/2}\circ\sum_{n\geq 1}(-1)^n\big(\bar\sigma^{-1/2}k\bar\sigma^{-1/2}\big)^{\circ n} \circ\bar\sigma^{-1/2}
\nonumber\\
&=
-\bar\sigma^{-1}\circ\sum_{n\geq 1} (k\circ\bar\sigma^{-1})^{\circ n}
,
\label{eq_series_expansion_g}
\end{align}
where $\circ$ denotes composition and $\circ n$ $n$-times composition, $n\geq 1$. 
In particular, $g$ is symmetric, satisfies~\eqref{eq_size_condition_g_appendix_F}, 
and it inherits the regularity of $k$: $g\in C^3(\bar\rhd)\cap C^3(\bar\lhd)$. 
Moreover,~\eqref{eq_series_expansion_g} already shows that $g_{|\partial \square}=0$ if $k_{|\partial \square}=0$.

To check the Neumann condition on the diagonal and the fact that $g$ satisfies the main equation~\eqref{eq_main_equation}, 
let us write derivatives of $k$ in terms of $g$. 
We henceforth drop the symbol $\circ$. 
Using the inverse operator $C_k^{-1} = (\bar\sigma^{-1}-g)$, 
it holds that:
\begin{equation}
(\bar\sigma +k) (\bar\sigma^{-1}-g) = \text{id} 
\quad \Rightarrow \quad 
k \bar\sigma^{-1} = C_k g
\quad \Rightarrow \quad 
k  = C_k g \bar\sigma.
\end{equation}
In the same way $k  =  \bar\sigma g  C_k$.
Differentiating $k$ with respect to the second variable
(or alternatively integrating against the derivative $\phi'$ of a test function $\phi$ on $(-1,1)$ and integrating by parts) yields:
\begin{equation}
\partial_2 k 
= 
\partial_2 (C_k g \bar\sigma)
=
C_k\partial_2 (g\bar\sigma)
.
\label{eq_partial_2k}
\end{equation}
Note that since $k$ is symmetric, 
$\partial_2 k(x,y)=\partial_1 k(y,x)$ for $(x,y)\in\squaredash$ means that $\partial_1 k$ is the operator adjoint to $\partial_2 k$ in $\mathbb L^2((-1,1))$, 
thus one has also:
\begin{equation}
\partial_1k 
= 
\partial_1(\bar\sigma g) C_k
.
\label{eq_partial_1k}
\end{equation}
In particular, this gives the Neumann condition for $g$ on the diagonal:
\begin{equation}
(\partial_1 -\partial_2)k(x_+,x)
= 
\bar\sigma(x)^2(\partial_1 -\partial_2)g(x_+,x)
,\qquad 
x\in(-1,1)
.
\label{eq_bc_Neumann_k_g_link}
\end{equation}
Let us now prove that $g$ satisfies the main equation~\eqref{eq_main_equation}. 
Recall the convention: 
if $\phi:\squaredash\rightarrow\R$ and $q:(-1,1)\rightarrow\R$, 
then $\phi q$ is the function $(x,y)\mapsto \phi(x,y)q(y)$, 
while $(q\phi)(x,y) := q(x)\phi(y)$.\\
The Euler-Lagrange equation~\eqref{eq_Euler-Lagrange_appendix} then reads:
\begin{align}
\Delta k    
-  \Big[\partial_1 \big(\bar\sigma \partial_1 h \big) C_k  +
C_k \partial_2\big(\partial_2 h \bar\sigma \big)\Big] \,- \, \big[\bar\sigma \delta_h k + k  \bar\sigma \delta_h\big]
= 0 
,
\label{eq_Euler-Lagrange_appendix simplifie}
\end{align}
where the above is an equality between functions on $\squaredash$. 
In~\eqref{eq_Euler-Lagrange_appendix simplifie}, 
the right-hand side in~\eqref{eq_Euler-Lagrange_appendix} has been included in the term in the bracket by using the operator $C_k = \bar \sigma + k$. 
The partial derivatives $\partial_2 k,\partial_1k$ were obtained in~\eqref{eq_partial_2k}--\eqref{eq_partial_1k}. 
As $k$ is not $C^2$ across the diagonal, 
taking the second derivative is more subtle, 
and we do so against a test function $\phi\in C^\infty(\square)$ vanishing on $\partial\square$:
\begin{equation}
- \partial_2k\partial_2\phi 
= 
\partial^2_2 k\phi + \delta_k\phi
=
- C_k \; \partial_2 (g \bar\sigma)\partial_2\phi
=
C_k \partial^2_2(g\bar\sigma) + C_k\bar\sigma\delta_g\phi
.
\end{equation}
Above, recall that, 
for a symmetric function $\psi:\squaredash\rightarrow\R$, 
$\delta_\psi$ is the function operator:
\begin{align}
\delta_\psi(x) 
&=
(\partial_2-\partial_1)\psi(x,x_+)
=
(\partial_1-\partial_2)\psi(x_+,x)
\nonumber\\
&=
\partial_1\psi(x_+,x)-\partial_1 \psi(x_-,x)
=
\partial_2\psi(x,x_+)-\partial_2 \psi(x,x_-)
,\qquad
x\in(-1,1)
.
\label{eq_def_delta_h_appendix}
\end{align}
Since $\delta_k = \bar\sigma^2\delta_g$, we get:
\begin{equation}
\partial^2_2 k 
=
C_k \partial^2_2(g\bar\sigma) 
+k\bar\sigma\delta_g
.
\end{equation}
By symmetry, one has also 
\begin{equation}
\partial_1 k  
=  
\partial_1 ( \sigma g ) \; C_k 
\quad \Rightarrow \quad 
\partial^2_1 k  
=  
\partial^2_1 ( \sigma g ) \; C_k    + \delta_g\bar\sigma k 
.
\end{equation}
Thus:
\begin{equation}
\Delta k  
= 
( \bar\sigma \, \partial^2_1 g  + 2 \bar\sigma'   \, \partial_1 g  + \bar\sigma''  \, g) C_k 
+ C_k   (  \partial^2_2 g \, \bar\sigma   + 2 \partial_2 g  \, \bar\sigma'  +  g\bar\sigma''  ) 
+
\delta_g\bar\sigma k + k\bar\sigma\delta_g
.
\end{equation}
Compose~\eqref{eq_Euler-Lagrange_appendix simplifie} by $C_k^{-1}$ on each side.  
Using $C_k^{-1} = \bar\sigma^{-1}-g$, 
this yields for $\Delta k$, 
still as an identity between functions on $\squaredash$:
\begin{align}
C_k^{-1} \, \Delta k  \, C_k^{-1} 
& = 
\bar\sigma^{-1} \partial^2_1(\bar\sigma\, g) - g \partial^2_1 ( \sigma g ) 
+    \partial^2_2( g \, \bar\sigma ) \bar\sigma^{-1}
-   \partial^2_2 (g \bar\sigma) g 
\nonumber\\
&\qquad + C_k^{-1}(\delta_g\bar\sigma k + k\bar\sigma\delta_g)C_k^{-1}
.
\end{align}
From $kC_k^{-1}=\bar\sigma g$ and the following identities obtained by integration by parts 
\begin{align}
\label{eq: IPP}
- g \partial^2_1 ( \sigma g ) = \partial_2  g \bar\sigma \, \partial_1 g      +   g \bar\sigma  \delta_g
 + \partial_2g\bar\sigma' g
,
\end{align}
we obtain:
\begin{align}
C_k^{-1}\Delta kC_k^{-1}
&=   
\Delta g  + 2 \frac{\bar\sigma'}{ \bar\sigma}  \, \partial_1 g  +   2 \partial_2 g  \, \frac{\bar\sigma'}{ \bar\sigma} 
+ 2 \partial_2  g \bar\sigma \, \partial_1 g      
+  g \bar\sigma  \, \delta_g  + \bar\sigma  \, \delta_g \, g
 +    g  \, \frac{\bar\sigma''}{ \bar\sigma} 
 +     \, \frac{\bar\sigma''}{ \bar\sigma} g 
 \nonumber\\
 &\quad 
  + \partial_2g\bar\sigma' g + g\bar\sigma'\partial_1g +\delta_g\bar\sigma g + g\delta_g\bar\sigma -2g\delta_g\bar\sigma^2g
.
\label{eq: appendice-relation A1}
\end{align}
We next show that:
\begin{align}
C_k^{-1} \,  \Big[  \partial_1 \big(\bar\sigma \partial_1 h \big) C_k  +
C_k \partial_2\big(\partial_2 h \bar\sigma \big) \Big]  \, C_k^{-1} 
& = \Delta h   +  \frac{\bar\sigma'}{\bar\sigma} \partial_1  h 
+ \partial_2  h \frac{\bar\sigma'}{\bar\sigma} \nonumber \\
& \quad +  \partial_2 g \bar\sigma \partial_1  h 
+  \partial_2 h \bar\sigma \partial_1  g
+ g \bar \sigma \delta_h +   \delta_h \bar \sigma g.
\label{eq: appendice-relation A2}
\end{align}
This follows by noticing that 
\begin{align}
C_k^{-1} \,  \partial_1 \big(\bar\sigma \partial_1 h \big) 
&= 
\bar\sigma^{-1} \, \big( \bar\sigma' \partial_1  h 
+ \bar\sigma \partial^2_1  h  \big) - g  \,  \partial_1 \big(\bar\sigma \partial_1 h \big)
\nonumber\\ 
&= 
\frac{\bar\sigma'}{\bar\sigma} \partial_1  h +  \partial^2_1  h 
+ \partial_1 g \bar\sigma \partial_1  h   + g \bar \sigma \delta_h
,
\end{align}
where in the last step, we used an integration by parts as in \eqref{eq: IPP}. 
Part of the boundary terms cancels thanks to the final identity:
\begin{align}
& C_k^{-1} \, [\bar\sigma \delta_h k + k \bar\sigma \delta_h] \, C_k^{-1}
 = 
 (\bar\sigma^{-1}-g) \, \bar\sigma \delta_h + \bar\sigma \delta_h (\bar\sigma^{-1}-g)
- 2 (\bar\sigma^{-1}-g)  \bar\sigma \delta_h \bar \sigma (\bar\sigma^{-1}-g) \nonumber \\
& \qquad = 2  \delta_h - g  \bar\sigma \delta_h -  \bar\sigma \delta_h g
- 2 (1-g \bar\sigma)  \delta_h  (1- \bar\sigma g)
= 
g  \bar\sigma \delta_h +  \bar\sigma \delta_h g -   2 g  \bar\sigma^2  \delta_h  g
.
\label{eq: appendice-relation A3}
\end{align}
Plugging the relations~\eqref{eq: appendice-relation A1},~\eqref{eq: appendice-relation A2},~\eqref{eq: appendice-relation A3} 
in the Euler-Lagrange equation~\eqref{eq_Euler-Lagrange_appendix simplifie}, 
we get:
\begin{align}
& 0 =   \Delta (g-h)  +  \frac{\bar\sigma'}{ \bar\sigma}  \, \partial_1  ( 2g-h)  
+   \partial_2 (2 g -h) \, \frac{\bar\sigma'}{ \bar\sigma} 
+  \partial_2  (g-h) \bar\sigma \, \partial_1 g       
 +  \partial_2 g \bar\sigma \partial_1 (g- h) 
 \nonumber\\
& +   2 g \bar\sigma^2  \delta_{h-g}  g  +    g \bar\sigma \, \left(  \frac{\bar\sigma''}{ \bar\sigma^2 } 
+2 \delta_g -2 \delta_h \right)
+   \bar\sigma   \left(  \frac{\bar\sigma''}{ \bar\sigma^2 }    +2\, \delta_g -2  \delta_h
\right) g  + \partial_2g\bar\sigma' g + g\bar\sigma'\partial_1g
.
\label{eq_EL_main_eq_0}
\end{align}
In addition, the boundary conditions of \eqref{eq_Euler-Lagrange_appendix} imply that 
$\delta_{g-h} = -  \frac{\sigma''}{ 2 \sigma^2}$ 
(recall $\bar\sigma'' = -2(\bar\rho')^2$). 
Finally, an integration by parts gives:
\begin{equation}
\partial_2g\bar\sigma' g 
= 
-g\bar\sigma''g - g\bar\sigma'\partial_1g
.
\end{equation}
The second line in~\eqref{eq_EL_main_eq_0} thus vanishes, 
and the first line is precisely the main equation~\eqref{eq_main_equation}.
\end{proof}
In the next sections, we focus on establishing existence, uniqueness and regularity for solutions of the Euler-Lagrange equation~\eqref{eq_Euler-Lagrange_appendix}. 
In view of Proposition~\ref{prop_EL_appendix}, 
this will prove Proposition~\ref{prop_main_equation}.
\subsection{Existence and uniqueness}
We now focus on solving the Euler-Lagrange equation~\eqref{eq_Euler-Lagrange_appendix} and the Poisson equation~\eqref{eq_Poisson_pour_partie_PiN_du_generateur}. 
To do so, we rewrite them in a common framework. 
Both equations are formulated as equations on $\squaredash$ involving symmetric functions. 
To solve them, it is therefore enough to look at the equation in a single triangle, 
say $\rhd$. 
To do so, let us introduce some notations.

Recall that the function $\mathcal M$ acts on $\phi,\psi\in\mathbb L^2(\squaredash)$ according to:
\begin{equation}
\mathcal M(\phi,\psi)(x,y) 
=
\int_{-1}^1 \phi(x,z)\bar\sigma(z)\psi(y,z)\, dz
,\qquad
(x,y)\in\squaredash\, 
.
\label{eq_def_mathcalW_appendix}
\end{equation}
Define by extension $\mathcal M$ on the triangle as follows: if $(\phi,\psi)\in \mathbb L^2(\rhd)^2$,
\begin{equation}
\mathcal M(\phi,\psi)
:=
\mathcal M(\phi_s,\psi_s),
\qquad 
\phi_s(x,y) 
= 
\begin{cases}
\phi(x,y)\quad \text{if } x<y,\\
\phi(y,x)\quad \text{if } x>y.
\end{cases}
\end{equation}
For $\phi\in \mathbb L^2(\rhd)$ and $q:(-1,1)\rightarrow \R$, 
we often write:
\begin{equation}
(qf)(x,y) 
:=
q(x) f(x,y),
\quad 
(fq)(x,y) 
:= 
f(x,y)q(y),
\qquad
(x,y)\in\rhd
.
\end{equation}
Given symmetric $\phi,\psi,\xi:\squaredash\rightarrow\R$ and 
$d:(-1,1)\rightarrow\R$, 
we say that $f:\squaredash\rightarrow\R$ solves the Poisson problem $(P) = (P_{\phi,\psi,\xi,d})$ if $f$ is symmetric, 
and $f_{|\rhd}$ is a classical solution of $(P_\rhd)$, where:
\begin{align}
(P_{\rhd}) : 
\begin{cases}
&\displaystyle{\Delta f(x,y) + [d(x)+d(y)]f(x,y) - \mathcal M(f,\xi)} 
\hspace{3cm}\text{ for }(x,y)\in\rhd,\\ 
&\ 
\qquad\displaystyle{  - \mathcal M(\xi,f) - \mathcal M(\partial_1 f,\partial_1\psi) - \mathcal M(\partial_1\psi,\partial_1 f) = \phi(x,y)}
\\
&\partial_{\nu_\rhd}f = 0\quad \text{on }D,
\\
&f = 0\quad \text{on }(\partial \rhd)\setminus D.
\end{cases}
\label{eq_probleme_entropique_appendice}
\end{align}
where $\partial_{\nu_\rhd}$ stands for the normal derivative on the diagonal.
\begin{remark}\label{rmk_lien_P_triangle_main_equation}
\begin{itemize}
	\item[(i)] If $\psi=h\in\s(\infty)$ (recall~\eqref{eq_def_test_functions_s_B}), 
	$d=0$, $\xi=0$ and $\phi\in \mathcal T\cap C^2(\bar\rhd)$ has norm $2$, 
then $(P)$ corresponds to the Poisson problem encountered in the proof of large deviations, 
in~\eqref{eq_Poisson_pour_partie_PiN_du_generateur}.
	\item[(ii)] Let $h\in\s(\infty)$, 
	$\psi =0$, 
	$d = \bar\sigma \delta_h$ with $\delta_h$ given in~\eqref{eq_def_delta_h_appendix}. 
	Take $\xi$ as follows:
	\begin{equation}
	\xi(x,y) 
	= 
		\partial_1\big(\bar\sigma\partial_1 h\big)(x,y)
		=		
	\partial_x\big(\bar\sigma(x)\partial_1 h(x,y)\big)
	,
	\qquad 
	(x,y)\in\squaredash\,
	.
	\label{eq_xi_for_k}
	\end{equation}
	Define then $\phi(x,y)$ for $(x,y)\in\squaredash$ as:
	\begin{align}
&\phi(x,y)
=
-\partial_1\big(\bar\sigma h\big)(x,y)-\partial_2\big(h\bar\sigma'\big)(x,y)
+\mathcal M(k_0 + \bar\sigma h \bar\sigma, \xi)(x,y)  
\nonumber\\
&\quad
+\mathcal M(k_0 + \bar\sigma h \bar\sigma, \xi)(x,y)
-\big[k_0(x,y) + (\bar\sigma h\bar\sigma)(x,y)\big]\big[\bar\sigma(x)\delta_h(x) + \bar\sigma(y)\delta_h(y)\big]
.
\label{eq_phi_for_k}
	\end{align}
	Then $(P)$ is the Euler-Lagrange equation~\eqref{eq_Euler-Lagrange_appendix} written with unknown $f = k-k_0 -\bar\sigma h \bar\sigma$. 
	Note that $f$ is chosen so that the boundary conditions in~\eqref{eq_probleme_entropique_appendice} imply the Neumann condition $\partial_{\nu_\rhd}f =0$.
	\end{itemize}
	\demo
\end{remark}
In the remainder of the section, we study existence, uniqueness and regularity of solutions of $(P)$. 
For $\star\in\{\rhd,\lhd\}$, 
we write $\big<\cdot,\cdot\big>_\star$ for the usual scalar product on $\mathbb L^2(\star)$, 
and simply $\big<\cdot,\cdot\big>$ for the scalar product on $\mathbb L^2(\square) = \mathbb L^2(\squaredash)$. 
The norm on $\mathbb L^2(\squaredash)$ is denoted by $\|\cdot\|_2$. 
Let also $\text{tr}$ denote the trace operator on the boundary of $\squaredash$. 
When interested only in a portion $\Gamma$ of the boundary, we may write $\text{tr}_\Gamma$. \\
We will use the fact that the Laplacian with $0$ Dirichlet boundary condition on $(\partial\rhd)\setminus D$ and $0$ Neumann boundary condition on the diagonal $D$ has a gap $\alpha > \pi^2/4>0$, 
see e.g. Equation 5 in Section 3.3. of~\cite{Siudeja2016}. 
This means that, 
for any symmetric $f\in \mathcal T$ satisfying the boundary conditions of $(P_{\rhd})$, one has:
\begin{align}
\|f_{|\rhd}\|^2_{\rhd} 
\leq 
\alpha^{-1} \|\nabla f_{|\rhd}\|^2_{\rhd} 
\quad \Rightarrow\quad 
\|f\|^2_{2}\leq \alpha^{-1}\|\nabla f\|^2_{2}
.
\label{eq_Poincare_inequality}
\end{align}
We first obtain existence and uniqueness of solutions of $(P)$ in the set $\mathcal T_{(P)}\subset \mathcal T=\mathbb H^2(\squaredash)$ of functions satisfying the boundary conditions of $(P)$ by a fixed point argument. 
The set $\mathcal T_{(P)}$ and its counterpart $\mathcal T_{(P_\rhd)}$ for functions on $\rhd$ are defined as follows:
\begin{align}
\mathcal T_{(P)} 
&= 
\mathbb H^2(\squaredash) \cap \big\{ f : \text{tr}_D(\partial_{\nu_\rhd}f) = 0, \text{tr}_{\partial\square}(f)=0\big\},
\nonumber\\
\mathcal T_{(P_\rhd)} 
&= 
\mathbb H^2(\rhd)\cap\big\{f : \text{tr}(f) = 0\text{ on }\partial(\rhd)\setminus \bar D, \text{tr}_D(\partial_{\nu_\rhd}f) = 0\big\} = \big\{f_{|\rhd} : f\in \mathcal T_{(P)}\big\}
.
\label{eq_def_sets_H_H_star}
\end{align}
\begin{proposition}[Solving $(P)$]\label{prop_solving_P_triangle}
Let $\phi,\xi\in\mathbb L^2(\squaredash)$, 
$\psi\in\mathcal T$ be symmetric functions. 
Let also $d:(-1,1)\rightarrow\infty$ be bounded.  
For $f\in \mathcal T_{(P)}$, 
define $Sf$ as the symmetric function such that, for $(x,y)\in \rhd$,
\begin{align}
Sf(x,y) 
&= 
(-\Delta_{\rhd})^{-1}\bigg[-\phi + df + fd -\mathcal M(f,\xi) -\mathcal M(\xi,f)
\nonumber\\
&\hspace{3cm}
- \mathcal M(\partial_1f,\partial_1\psi) - \mathcal M(\partial_1\psi,\partial_1f)\bigg](x,y)
. 
\label{eq_def_S_P_triangle}
\end{align}
Above, $\Delta_{\rhd}^{-1}$ is the inverse of the Laplacian on $\rhd$ with $0$ Dirichlet condition on $(\partial\rhd)\setminus D$, and $0$ Neumann conditions on the diagonal $D$. \\
Then $S f\in \mathcal T_{(P)}$. 
Moreover, if $d=0$, $\xi=0$ and $\|\nabla \psi\|_{2}\leq 1$, 
$\|\phi\|_{2}\leq 2$, 
then $S$ has a unique fixed point $f_{\phi,\psi}\in \mathcal T_{(P)}$ 
with $\|f_{\phi,\psi}\|_{\mathbb H^1(\squaredash)}\leq C$ and $\|f_{\phi,\psi}\|_\infty \leq C'$, 
for constants $C,C'>0$ independent of $\phi,\psi$.\\
If instead $\psi=0$, 
there is a fixed point $f_{d,\phi,\xi}\in \mathcal T_{(P)}$ provided $\|\xi\|_2, \|d\|_\infty$ are small enough,  with $\|f_{d,\phi,\xi}\|_{\mathbb H^1(\squaredash)}\leq \delta(\|d\|_\infty,\|\xi\|_2,\|\phi\|_2)$ and $\delta$ vanishes when $(\|d\|_\infty,\|\xi\|_2,\|\phi\|_2)$ vanishes.
\end{proposition}
\begin{proof}
We prove that $S$ is a contraction on $\mathcal T_{(P)}$ for the norm $\|\nabla \cdot\|_{2}$. 
Let us start by showing that $S$ is well defined. 
The inverse operator $\Delta_{\rhd}^{-1}$ exists by Lemma 4.4.3.1 in~\cite{Grisvard2011} and, 
by Theorem 4.4.3.7 in~\cite{Grisvard2011}, 
maps $\mathbb L^2(\rhd)$ onto $\mathcal T_{(P_\rhd)}$. 
It follows that $S(\mathcal T_{(P)})\subset \mathcal T_{(P)}$. \\
We now prove that $S$ is a contraction. For $f\in \mathcal T_{(P)}$, one has:
\begin{align}
\|\nabla S f \|_{2}^2 
&= 
\big<(\nabla Sf)_{|\rhd}, (\nabla Sf)_{|\rhd}\big>_\rhd +\big<( \nabla Sf)_{|\lhd}, (\nabla Sf)_{|\lhd}\big>_{\lhd} 
\nonumber\\
&= 
-\big<(Sf)_{|\rhd},(\Delta Sf)_{|\rhd}\big>_{\rhd} -\big<(Sf)_{|\lhd},(\Delta Sf)_{|\lhd}\big>_{\lhd} = -\big<Sf,\Delta Sf\big>
\label{eq_norme_gradient_Sf_square}
.
\end{align}
The integration by parts is legitimate by Theorem 1.5.3.1 in~\cite{Grisvard2011}. 
Let us compute the right-hand side. 
One has, using $\bar\sigma\leq 1/4$:
\begin{align}
\forall (x,y)\in \squaredash,\qquad 
\mathcal M(f,\xi)(x,y) 
&= 
\int_{(-1,1)}f(z,x)\bar\sigma(z)\xi(z,y)\, dz
\nonumber\\
&\leq 
\frac{1}{4}\Big(\int_{(-1,1)}f(z,x)^2\, dz\Big)^{1/2}\Big(\int_{(-1,1)}\psi(z,y)^2\, dz\Big)^{1/2}
,
\end{align}
and the same holds for $\mathcal M(\partial_1f,\partial_1\psi)$
As a result, 
by Cauchy-Schwarz inequality and using $\|\partial_1 f \|_{2} = 2^{-1/2}\|\nabla f\|_{2}$ as implied by the symmetry of $f$:
\begin{align}
\Big|\big<S f, \mathcal M(f,\xi)\big>\Big|
&\leq 
\frac{1}{4}\|Sf\|_{2}\|f\|_{2}\|\xi\|_{2},
\nonumber\\
\Big|\big<S f, \mathcal M(\partial_1f,\partial_1\psi)\big>\Big|
&\leq 
\frac{1}{8}\|Sf\|_{2}\|\nabla f\|_{2}\|\nabla \psi\|_{2}
.
\end{align}
Recalling the expression of $S f$ from~\eqref{eq_def_S_P_triangle},
~\eqref{eq_norme_gradient_Sf_square} is therefore bounded as follows:
\begin{align}
\|\nabla S f \|_{2}^2 
&\leq 
\|S f\|_{2}\Big(\|\phi\|_{2} + 2\|d\|_\infty \|f\|_2  +  \frac{\|f\|_{2}\|\xi\|_2}{4}+ \frac{\|\nabla f\|_{2}\|\nabla\psi\|_2}{8}\Big)
.
\end{align}
Since $f,Sf \in \mathcal T_{(P)}$, 
the Poincaré inequality~\ref{eq_Poincare_inequality} can be applied and yields:
\begin{align}
\|\nabla Sf \|_{2} 
&\leq 
\alpha^{-1/2}\Big(\|\phi\|_{2} + 2\alpha^{-1/2}\|d\|_\infty \|\nabla f\|_2  
\nonumber\\
&\hspace{3cm}+  \frac{\alpha^{-1/2}\|\nabla f\|_{2}\|\xi\|_2}{4}+ \frac{\|\nabla f\|_{2}\|\nabla\psi\|_2}{8}\Big)
.
\label{eq_general_bound_gradient_S_f}
\end{align}
By similar computations, if $f_1,f_2\in \mathcal T_{(P)}$, one obtains:
\begin{align}
\|\nabla (Sf_1-Sf_2) \|_{2} 
&\leq 
\alpha^{-1/2}\Big(2\alpha^{-1/2}\|d\|_\infty +  \frac{\alpha^{-1/2}\|\xi\|_2}{4}+ \frac{\|\nabla\psi\|_2}{8}\Big)\|\nabla (f_1-f_2)\|_2
.
\label{eq_general_bound_gradient_Sf1_minus_Sf2}
\end{align}
In particular, $S$ is a contraction as soon as:
\begin{equation}
\alpha^{-1/2}\Big(2\alpha^{-1/2}\|d\|_\infty +  \frac{\alpha^{-1/2}\|\xi\|_2}{4}+ \frac{\|\nabla\psi\|_2}{8}\Big)
<
1
.
\label{eq_sufficient_condition_contraction}
\end{equation}
Recall that $\alpha>\pi^2/4$. 
If $d=0$, $\xi=0$, 
$\|\nabla \psi\|_2\leq 1$ and $\|\phi\|_2\leq 2$ (which includes case (i) of Remark~\ref{rmk_lien_P_triangle_main_equation}), 
then:
\begin{equation}
\|\nabla S(f_1-f_2) \|_{2} 
\leq
\frac{1}{4\pi}\|\nabla(f_1-f_2)\|_2,
\qquad 
\|\nabla Sf\|_{2} 
\leq
\frac{1}{4\pi}\|\nabla f \|_2 + \frac{4}{\pi}
.
\end{equation}
There is thus a unique fixed point $f_{\phi,\psi}\in \mathcal T_{(P)}$, 
and since $\|\phi\|_2\leq 2$, 
it belongs to the ball $B(0,c) = \{ u\in \mathbb H^1(\squaredash):\text{tr}(u) = 0\text{ on }\partial\square, \|\nabla u\|_{2}\leq c\}$ with $c = 8/(2\pi-1)\leq 2$. 
Poincaré inequality~\eqref{eq_Poincare_inequality} yields $\|f_{\phi,\psi}\|_{\mathbb H^1(\squaredash)}\leq (4+16/\pi^2)^{1/2}\leq 2$, 
and Theorem 4.3.1.4 in~\cite{Grisvard2011} yields $\|f_{\phi,\psi}\|_{\mathbb H^2(\squaredash)}\leq C$ for some universal $C>0$. 
The Sobolev embedding $\mathbb H^2(\squaredash)\subset C^0(\bar\rhd)\cap C^0(\bar\lhd)$ then implies $\|f_{\phi,\psi}\|_{\infty}\leq C'$ for a universal $C'>0$ as claimed.\\

Consider now the case $\psi=0$ which includes item (ii) of Remark~\ref{rmk_lien_P_triangle_main_equation}. 
Then:
\begin{align}
\|\nabla S(f_1-f_2) \|_{2} 
&\leq
\frac{4}{\pi^2}\Big(2\|d\|_\infty + \frac{\|\xi\|_2}{4}\Big)\|\nabla(f_1-f_2)\|_2,
\nonumber\\
\|\nabla Sf\|_{2} 
&\leq
\frac{4}{\pi^2}\Big(2\|d\|_\infty + \frac{\|\xi\|_2}{4}\Big)\|\nabla f \|_2 + \frac{2\|\phi\|_2}{\pi}
.
\end{align}
If $\|d\|_\infty,\|\xi\|_2$ are small enough, $S$ is a contraction and the norm of its fixed point vanishes when $\|d\|_\infty,\|\xi\|_2$ and $\|\phi\|_2$ vanish. 
This concludes the proof.
\end{proof}
\subsection{Regularity estimates}
In Proposition~\ref{prop_solving_P_triangle}, the solution of $(P)$ has been shown to be in $\mathbb H^2(\squaredash)$. 
In this section, 
we use results of~\cite{Grisvard2011} to argue that the solution of $(P)$ is more regular if the data $\phi,\psi$ are regular. 
This concludes the proof of Proposition~\ref{prop_main_equation}. 
The study of regularity is made very complicated by the presence of corners.
\begin{proposition}[Theorem 5.1.3.1. in~\cite{Grisvard2011}]\label{prop_regularite_P_triangle_general}
Let $b\in\N$, $p>2$ and let $\zeta\in \mathbb W^{b,p}(\rhd)$. Let $S_1,S_2,S_3$ denote the corners of $\rhd$ numbered in a counter-clockwise fashion, 
with $S_1$ the upper left corner. 
Consider on $\rhd$ the equation $\Delta f = \zeta$, 
with the boundary conditions of $(P_\rhd)$. 
If $b= 0$, then $f\in \mathbb W^{2,p}(\rhd)$. 
If $b\leq 3$, $f\in \mathbb W^{b+2,p}(\rhd)$ provided $\zeta$ vanishes at the corners, i.e. provided:
\begin{equation}
\forall j\in\{1,2,3\},\qquad 
\zeta\big(S_j\big) 
= 
0
.
\label{eq_conditionb_is_1_grisvard}
\end{equation}
\end{proposition}
\begin{remark}
Though the statement of Proposition~\ref{prop_regularite_P_triangle_general} makes no mention of them, we recall notations from~\cite{Grisvard2011} so that the reader may check that Theorem 5.1.3.1 applies to our case.\\
Label by $j\in\{1,...,3\}$ the line segments composing $\partial\rhd$ in a counter clockwise fashion, with the convention that $j=1$ for the $y=-1$ segment. $S_j$ is then the point joining segments $j,j+1$ in $\partial \rhd$. Let $\omega_j$ be the counter-clockwise measure of the inwards angle at $S_j$:
\begin{equation}
\omega_j = \begin{cases}
\pi/2\quad &\text{if }j=1,\\
\pi/4\quad &\text{if }j\in\{2,3\}.
\end{cases}
\end{equation}
Let $\nu_j=\mu_j$ denote the unit outwards normal and $\tau_j$ be the (counter clockwise) unit tangent vector on the line segment $j$. Define also $\Phi_j = \pi/2$ if $j\in\{1,2\}$, $\Phi_j = 0$ if $j=3$ and $\Phi_{3+1} := \Phi_1$. Finally, for $m\in\Z$ and each $j$, define:
\begin{equation}
\lambda_{j,m} = \frac{\Phi_j- \Phi_{j+1}+m\pi}{\omega_j} = \begin{cases}
2m\quad &\text{if }j=1,\\
2+4m\quad &\text{if }j=2,\\
-2+4m\quad &\text{if }j=3.
\end{cases}\label{eq_def_lambda_m_j}
\end{equation}
\demo
\end{remark}
In our context, Proposition~\ref{prop_regularite_P_triangle_general} translates to the following result.
\begin{proposition}[Regularity of solutions of $(P_\rhd)$]\label{prop_regularity_P_triangle}
Let $\phi,\psi,d,\xi$ be such that the solution $f$ given by Proposition~\ref{prop_solving_P_triangle} exists.
\begin{enumerate}
	\item[(i)] Assume $d=0$, $\xi=0$, let $\phi\in C^2(\bar\rhd)\cap C^2(\bar\lhd)$ be symmetric and let $\psi\in \s(\infty)$, 
	corresponding to item (i) of Remark~\ref{rmk_lien_P_triangle_main_equation}. 
	If $\phi_{|\partial\square} = 0$, 
	then $f \in \mathbb W^{4,p}(\squaredash)$ for any $p>2$.
	\item[(ii)] Recall the definition~\eqref{eq_def_test_functions_s_B} of $\s(\epsilon)$ for $\epsilon>0$. 
	Take $\psi=0$, and $d,\phi,\xi$ defined in terms of $h\in\s(\epsilon)$
	as in item (ii) of Remark~\ref{rmk_lien_P_triangle_main_equation}. 
	Then $f \in \mathbb W^{4,p}(\squaredash)$ for any $p>2$, 
	and $k = f+ k_0 + \bar\sigma h \bar\sigma \in k_0 + \s(\epsilon')$ for some $\epsilon'>0$ depending only on $\rho_\pm$ and $\epsilon$.
\end{enumerate}
\end{proposition}
\begin{proof}
Since $f$ is symmetric, we work only on $\rhd$. 
Let us first assume $f$ has the alleged regularity and treat all claims of item $2$ that do not have to do with the regularity of $f$.\\
Define:
\begin{equation}
\|k-k_0\|_{C^1}
:=
\max\Big\{\|k-k_0\|_\infty,\|\partial_1(k-k_0)\|_{\infty}\Big\}
,
\end{equation}
Notice first that $k-k_0 = f+\bar\sigma h \bar\sigma$ with $h\in\s(\epsilon)$ implies that $k-k_0\in\mathbb W^{4,p}(\squaredash)$ for all $p>2$, 
thus $k$ also as $k_0$ is regular. \\
Let us now prove that $\|k-k_0\|_{C^1}$ vanishes when $d,\phi,\xi$ vanish. 
By Sobolev embedding (see Proposition~\ref{prop_sobolev_embeddings}), 
it is enough to prove the same for $\|k-k_0\|_{\mathbb W^{2,p}(\squaredash)}$. 
Since $k-k_0 = f+\bar\sigma h \bar\sigma$, 
it is enough to bound $\|f\|_{\mathbb W^{2,p}(\squaredash)}$. 
Theorem 4.3.2.4 in~\cite{Grisvard2011} implies that, for a universal constant $C>0$:
\begin{equation}
\|f\|_{\mathbb W^{2,p}} 
\leq 
C \Big(\|\Delta f\|_{\mathbb W^{0,p}} + \|f\|_{\mathbb W^{1,p}} \Big)
.
\end{equation}
Recalling the expression~\eqref{eq_def_S_P_triangle} of the mapping $S$, 
it holds that:
\begin{align}
\|\Delta f\|_{\mathbb W^{0,p}(\squaredash)}
&\leq 
\big(\|\phi\|_p + 2\|d\|_\infty \|f\|_p + \|\xi\|_\infty\big) \|f\|_{\mathbb W^{1,p}(\squaredash)}
\nonumber\\
&\leq
C'\big(\|\phi\|_p + 2\|d\|_\infty \|f\|_p + \|\xi\|_\infty\big) \|f\|_{\mathbb H^{2}(\squaredash)}
,
\end{align}
where the second inequality is again a Sobolev embedding. 
The fact that $\|f\|_{\mathbb H^2(\squaredash)}$ vanishes with $d,\phi,\xi$ now follows from Proposition~\ref{prop_solving_P_triangle}.

We now prove that $k-k_0$ satisfies the boundary conditions of elements of $\s(\infty)$. 
Since $f,k_0,h$ vanish on $\partial\square$, so does $k$. 
Moreover, $\partial_{\rhd} h=0$ at the extremities $S_2,S_3$ of the diagonal $D$ and
$\partial_{\rhd} f=0$ on $D$ means $\partial_{\rhd} (k-k_0)(S_2)=0 = \partial_{\rhd} (k-k_0)(S_3)$. 
We have thus shown that if $f$ has the alleged regularity, then $k-k_0\in\s(\epsilon')$ for some $\epsilon'>0$ depending on $d,\phi,\xi$, 
i.e. on $\rho_\pm,\epsilon$ recalling~\eqref{eq_phi_for_k}--\eqref{eq_xi_for_k} and $h\in\s(\epsilon)$. \\

Let us now prove the regularity of $f$. 
By definition, $\Delta f = \Delta Sf$, with $S$ defined in~\eqref{eq_def_S_P_triangle}. 
Let $p>2$. 
The idea is classical: if $f\in \mathbb W^{2+n,p}(\squaredash)$, $n\in\N$, 
we want to prove that $\Delta S f\in \mathbb W^{1+n,p}(\squaredash)$, 
from which $f\in \mathbb W^{3+n,p}(\squaredash)$ by Proposition~\ref{prop_regularite_P_triangle_general} provided $\Delta Sf$ satisfies suitable boundary conditions. 
To implement this recursion scheme, we first prove that $\Delta Sf\in\mathbb L^p(\squaredash)$. 
By assumption on $\phi,\psi$ in case of item (i) 
(using the Sobolev embedding $C^2(\bar\rhd)\subset \mathbb W^{2,s}(\rhd)$ for any $s>2$), 
and from~\eqref{eq_xi_for_k}--\eqref{eq_phi_for_k} for item (ii), 
we see that it is the regularity of $f$ only that limits the regularity of $\Delta Sf$. 
The fact that $\Delta Sf\in\mathbb L^p(\squaredash)$ then follows from the Sobolev embedding $\mathbb H^1(\squaredash)\subset\mathbb L^{s}(\squaredash)$, 
valid for any $s>2$, 
see Proposition~\ref{prop_sobolev_embeddings}. 
It follows that $f\in\mathbb W^{2,p}(\squaredash)$ by Proposition~\ref{prop_regularite_P_triangle_general}.

To obtain further regularity on $f$, 
let us check that $\Delta Sf$ is in $\mathbb W^{b-1,p}(\squaredash)$ whenever $f\in \mathbb W^{b,p}(\squaredash)$ for $b\in\N^*$, 
and also that $\Delta Sf$ satisfies the condition~\eqref{eq_conditionb_is_1_grisvard} in Proposition~\ref{prop_regularite_P_triangle_general}. \\
The regularity of $\Delta Sf$ boils down to proving that $\mathcal W(u,v)$, 
defined in~\eqref{eq_def_mathcalW_appendix}, 
is in $\mathbb W^{b,p}(\rhd)$ ($b\in\N^*$) 
whenever $u,v\in \mathbb W^{b,p}(\rhd)$. 
This is the claim of the following lemma.
\begin{lemma}\label{lemm_terme_integral_W_k_si_data_W_k}
Let $p\geq 2$ and $1\leq b\leq 4$. 
Let $u,v \in \mathbb W^{b,p}(\rhd)$, 
and recall from~\eqref{eq_def_mathcalW_appendix} the definition of $\mathcal M $. 
Then $\mathcal M(u,v)\in \mathbb W^{b,p}(\rhd)$.
\end{lemma}
Lemma~\ref{lemm_terme_integral_W_k_si_data_W_k} is easily proven by approximating $u,v$ in $\mathbb W^{b,p}(\rhd)$ by sequences in $C^\infty(\bar\rhd)$, 
and integrating by parts.

It remains to prove that $\Delta Sf$ satisfies the condition of Proposition~\ref{prop_regularite_P_triangle_general}, 
i.e. that $\Delta Sf(S_j) = 0$ for $j\in\{1,2,3\}$. 
By assumption in the case of item (i), 
and from the expression~\eqref{eq_phi_for_k} and the definition~\eqref{eq_def_test_functions_s_B} of $\s(\epsilon)$ for item (ii), 
we know that $\phi_{|\partial\square} =0$.
The fact that $f_{|\partial\square}=0$ gives $\mathcal M(f,\xi)(S_j) = 0 = \mathcal M(\xi,f)(S_j)$ for each $j$. 
Integrating by parts and since $\psi\in\s(\infty)$, 
one has also $\mathcal M(\partial_1 f,\partial_1\psi)(S_j)=0=\mathcal M(\partial_1 \psi,\partial_1f)(S_j)$.
It follows that $\Delta Sf(S_j) = 0$ for $j\in\{1,2,3\}$, 
thus $f\in \mathbb W^{4,p}(\squaredash)$. 
Since $p>2$ was arbitrary, this concludes the proof.
\end{proof}
\subsection{Bounds on the solution and definition of $\epsilon_B$}\label{app_choice_epsilonB}
Let $h\in\s(\infty)$ and let $g_h = \bar\sigma^{-1} - (\bar\sigma + k_h)^{-1}$ be the solution of the main equation~\eqref{eq_main_equation} obtained from the solution $k_h$ of the Euler-Lagrange equation~\eqref{eq_Euler-Lagrange_appendix}.  
In this section, 
we define the $\epsilon_B$ arising in Theorems~\ref{theo_large_devs}--\ref{theo_entropic_problem}, 
and show that the $C^1$ norm of $g_h-g_0$ can be controlled by $\|h\|_{C^1}$ and $\bar\rho'$ as claimed in Theorem~\ref{theo_entropic_problem}.  
Writing $C_{k_0} = \bar\sigma + k_0 = (\bar\sigma^{-1}-g_0)^{-1}$, 
notice first that, 
as soon as $\|k_h-k_0\|_2$ is sufficiently small:
\begin{align}
g_h-g_0 
&= 
C_{k_0}^{-1}- (\bar\sigma + k_h)^{-1} 
=
C_{k_0}^{-1} - (C_{k_0} + (k_h-k_0))^{-1}
\nonumber\\
&=
C_{k_0}^{-1} \sum_{n\geq 1} (-1)^n\big((k_h-k_0)C_{k_0}^{-1}\big)^{\circ n}
.
\end{align}
For $k_h-k_0\in\s(\epsilon')$ for some suitably small $\epsilon'>0$, 
one can take derivatives inside the sum, which directly yields:
\begin{equation}
\|g_h-g_0\|_{C^1}
\leq 
\tilde \delta(\rho_\pm, \|k_h-k_0\|_{C^1}),
\qquad
\lim_{(x,y)\rightarrow 0}\tilde \delta(x,y) 
= 
0
.
\label{eq_norm_g_as_norm_k}
\end{equation}
Since $f=k_h-k_0-\bar\sigma h \bar\sigma$ is the solution of $(P)$, $\|k_h-k\|_{C^1}$ in turn only depends on $\rho_\pm,\|h\|_{C^1}$, 
thus $\|g_h-g_0\|_{C^1}$ only depends on $\bar\rho',\|h\|_{C^1}$, and vanishes when they both vanish as claimed in Theorem~\ref{theo_entropic_problem}.

We now define $\epsilon_B$. 
\begin{definition}[Definition of $\epsilon_B$]\label{def_epsilon_B_and_G}
Let $h\in\s(\infty)$.
Let $\rho_-\in(0,1)$, and choose $\epsilon_B=\epsilon_B(\rho_-)>0$ and $\rho_+\in[\rho_-,1)$ such that, 
if $\bar\rho':=\frac{\rho_+-\rho_-}{2}\leq \epsilon_B$ and $h\in\s(\epsilon_B)$, then:
\begin{itemize}
	\item The contraction $S$ has a fixed point $f$. 
	\item $g_h$ is a negative kernel (as defined in~\eqref{eq_g_negative_kernel_appendix}).
	\item To ensure good concentration properties in Section~\ref{sec_computations}, 
	$\|g_h-h\|_{C^1}\|g_h\|_{C^1}\leq (2^{10}C_{LS})^{-1}$, 
	$\bar\rho' \|g_h-h\|_{C^1}\leq (2^{10}C_{LS})^{-1}$ and $\bar\rho'\epsilon_B \leq (2^{10}C_{LS})^{-1}$,  
	where $C_{LS}$ is the log-Sobolev constant appearing in  Lemma~\ref{lemm_LSI_sec3}.
\end{itemize}
\end{definition}

\end{appendices}


\providecommand{\bysame}{\leavevmode\hbox to3em{\hrulefill}\thinspace}
\providecommand{\MR}{\relax\ifhmode\unskip\space\fi MR }
\providecommand{\MRhref}[2]{%
  \href{http://www.ams.org/mathscinet-getitem?mr=#1}{#2}
}
\providecommand{\href}[2]{#2}

\end{document}